%% file: Ainf-coho.tex
\documentclass[11pt, reqno]{amsart}

\usepackage[OT2, T1]{fontenc}

\usepackage[usenames,dvipsnames]{color}
\newcommand{\revise}[1]{{\color{Bittersweet}   [#1]}}   
\newcommand{\ready}[1]{{\color{Brown}   [#1]}}	    

\renewcommand{\revise}[1]{#1}   
\renewcommand{\ready}[1]{#1}	    

\usepackage{url}
\usepackage{amsmath}
\usepackage{array}
\usepackage{graphicx}
\usepackage{amssymb}
\usepackage{amsthm}
\definecolor{link}{RGB}{11,0,128}
\usepackage[colorlinks=true,citecolor=NavyBlue,linkcolor=Bittersweet,urlcolor=link]{hyperref} 
\usepackage{paralist}
\usepackage{colonequals}		
\usepackage{color}
\usepackage{amscd}
\usepackage[all,cmtip]{xy}	
\usepackage{sseq}			
\usepackage{verbatim}
\usepackage{parskip}
\usepackage[in]{fullpage}
\usepackage{graphicx}
\usepackage{mathrsfs} 			
\usepackage{cleveref}
\usepackage{enumitem}
\usepackage{moreenum}			
\usepackage[alphabetic,lite]{amsrefs} 	
\usepackage{stmaryrd}	
\usepackage[foot]{amsaddr}
\usepackage{subfiles}
\usepackage{ragged2e}

\DeclareSymbolFont{cyrletters}{OT2}{wncyr}{m}{n}
\DeclareMathSymbol{\Sha}{\mathalpha}{cyrletters}{"58}

\newcommand{\gL}{\lambda}

\newcommand{\bA}{\mathbb{A}}
\newcommand{\bB}{\mathbb{B}}

\newcommand{\bF}{\mathbb{F}}
\newcommand{\bG}{\mathbb{G}}

\newcommand{\bL}{\mathbb{L}}
\newcommand{\bN}{\mathbb{N}}

\newcommand{\bQ}{\mathbb{Q}}

\newcommand{\bT}{\mathbb{T}}

\newcommand{\bZ}{\mathbb{Z}}

\newcommand{\cF}{\mathcal{F}}

\newcommand{\cH}{\mathcal{H}}

\newcommand{\cN}{\mathcal{N}}
\newcommand{\cO}{\mathcal{O}}

\newcommand{\cU}{\mathcal{U}}

\newcommand{\cX}{\mathcal{X}}
\newcommand{\cY}{\mathcal{Y}}

\newcommand{\fm}{\mathfrak{m}}

\newcommand{\fo}{\mathfrak{o}}

\newcommand{\fM}{\mathfrak{M}}

\newcommand{\fO}{\mathfrak{O}}

\newcommand{\fU}{\mathfrak{U}}

\newcommand{\fX}{\mathfrak{X}}

\newcommand{\sT}{\mathscr{T}}

\newcommand{\ra}{\rightarrow}

\newcommand{\xra}{\xrightarrow}

\newcommand{\hra}{\hookrightarrow}

\newcommand{\Infty}{\infty}
\newcommand{\wt}{\widetilde}
\newcommand{\wh}{\widehat}
\newcommand{\eps}{\epsilon}
\newcommand{\del}{\partial}
\newcommand{\pr}{^{\prime}}

\newcommand{\ce}{\colonequals}
\newcommand{\ov}{\overline}

\newcommand{\sm}{\mathrm{sm}}

\renewcommand{\b}{\textbf}
\newcommand{\surjects}{\twoheadrightarrow}

\newcommand{\tensor}{\otimes} 		
\newcommand{\isomto}{\overset{\sim}{\longrightarrow}}

\newcommand{\cont}{{\mathrm{cont}}}		
\newcommand{\cris}{{\mathrm{cris}}}		
\newcommand{\dR}{{\mathrm{dR}}}		
\newcommand{\st}{{\mathrm{st}}}		
\newcommand{\et}{\mathrm{\acute{e}t}}	
\newcommand{\Zar}{\mathrm{Zar}}		
\newcommand{\llb}{\llbracket}		
\newcommand{\rrb}{\rrbracket}		
\newcommand{\tors}{\mathrm{tors}}		
\newcommand{\psh}{\mathrm{psh}}		

\renewcommand{\i}{^{-1}}
\renewcommand{\th}{^{\mathrm{th}}}

\newcommand{\PD}{{\mathrm{PD}}}	

\newcommand{\leftexp}[2]{{\vphantom{#2}}^{#1}{#2}}

\providecommand{\In}[1]{\left\langle#1\right\rangle}
\providecommand{\p}[1]{\left(#1\right)}
\providecommand{\up}[1]{{\upshape(}#1{\upshape)}}
\providecommand{\uref}[1]{{\upshape\ref{#1}}}
\providecommand{\uS}{{\upshape\S}}
\providecommand{\ucolon}{{\upshape:} }
\providecommand{\uscolon}{{\upshape;} }
\providecommand{\sqp}[1]{\left[#1\right]}
\providecommand{\f}[2]{\frac{#1}{#2}}
\DeclareMathOperator{\Ker}{Ker}			
\DeclareMathOperator{\Coker}{Coker}		
\DeclareMathOperator{\im}{Im}			
\DeclareMathOperator{\Spec}{Spec}		
\DeclareMathOperator{\Spf}{Spf}		
\DeclareMathOperator{\Spa}{Spa}		
\DeclareMathOperator{\Fitt}{Fitt}			
\DeclareMathOperator{\Frac}{Frac}		
\DeclareMathOperator{\id}{id}			
\DeclareMathOperator{\Cone}{Cone}		
\DeclareMathOperator{\Tor}{Tor}			
\DeclareMathOperator{\Gal}{Gal}	
\DeclareMathOperator{\ord}{ord}	
\DeclareMathOperator{\length}{length}		
\DeclareMathOperator{\Fil}{Fil}			
\newcommand{\ba}{\begin{aligned}}
\newcommand{\ea}{\end{aligned}}
\newcommand{\be}{\begin{equation}}
\newcommand{\ee}{\end{equation}}
\newcommand{\pf}{\begin{proof}}
\newcommand{\bpf}{\begin{proof}}
\newcommand{\epf}{\end{proof}}
\newcommand{\bthm}{\begin{thm}}
\newcommand{\ethm}{\end{thm}}

\newcommand{\bprop}{\begin{prop}}
\newcommand{\eprop}{\end{prop}}
\newcommand{\bcor}{\begin{cor}}
\newcommand{\ecor}{\end{cor}}

\newcommand{\brem}{\begin{rem}}
\newcommand{\erem}{\end{rem}}

\newcommand{\brems}{\begin{rems} \hfill \begin{enumerate}[label=\b{\thesubsection.},ref=\thesubsection]}
\newcommand{\bremstweak}{\begin{rems-tweak} \hfill \begin{enumerate}[label=\b{\thesubsection.},ref=\thesubsection]}
\newcommand{\bremst}{\begin{rems-tweak} \hfill \begin{enumerate}[label=\b{\thesubsection.},ref=\thesubsection]}
\newcommand{\remi}{\addtocounter{subsection}{1} \item}

\newcommand{\erems}{\end{enumerate} \end{rems}}
\newcommand{\eremstweak}{\end{enumerate} \end{rems-tweak}}
\newcommand{\eremst}{\end{enumerate} \end{rems-tweak}}
\newcommand{\blem}{\begin{lemma}}
\newcommand{\elem}{\end{lemma}}

\newcommand{\bconj}{\begin{conj}}
\newcommand{\econj}{\end{conj}}
\newcommand{\bprob}{\begin{Problem}}
\newcommand{\eprob}{\end{Problem}}

\newcommand{\bq}{\begin{Q}}
\newcommand{\eq}{\end{Q}}
\newcommand{\benum}{\begin{enumerate}[label={{\upshape(\alph*)}}]}
\newcommand{\benuma}{\begin{enumerate}[label={{\upshape(\arabic*)}}]}
\newcommand{\benumr}{\begin{enumerate}[label={{\upshape(\roman*)}}]}
\newcommand{\eenum}{\end{enumerate}}
\newcommand{\bc}{}
\newcommand{\bd}{\begin{defn}}
\newcommand{\ed}{\end{defn}}

\newcommand{\beg}{\begin{eg}}
\newcommand{\eeg}{\end{eg}}

\newcommand{\bcl}{\begin{claim}}
\newcommand{\ecl}{\end{claim}}
\newcommand{\lab}{\label}

\newcommand{\q}{\quad}
\newcommand{\qq}{\quad\quad}
\newcommand{\qqq}{\quad\quad\quad}

\newcommand{\tst}{\textstyle}

\newcommand{\val}{\mathrm{val}}
\newcommand{\Inf}{\mathrm{inf}}
\newcommand{\ad}{\mathrm{ad}}			
\newcommand{\proet}{\mathrm{pro\acute{e}t}}			

\newcommand*{\QED}{\hfill\ensuremath{\qed}}

\theoremstyle{plain}
\newtheorem{thm}[subsection]{Theorem}
\Crefname{thm}{Theorem}{Theorems}

\Crefname{rethm}{Theorem}{Theorem}
\newtheorem{prop}[subsection]{Proposition}
\Crefname{prop}{Proposition}{Propositions} 
\newtheorem{Q}[subsection]{Question}
\Crefname{Q}{Question}{Questions}
\newtheorem{Problem}[subsection]{Problem}
\Crefname{Problem}{Problem}{Problems}
\newtheorem{conj}[subsection]{Conjecture}
\Crefname{conj}{Conjecture}{Conjectures}
\newtheorem{cor}[subsection]{Corollary}
\Crefname{cor}{Corollary}{Corollaries}
\newtheorem{lemma}[subsection]{Lemma}

\Crefname{subprop}{Proposition}{Propositions}

\Crefname{subcor}{Corollary}{Corollaries}
\newtheorem{sublem}[equation]{Lemma}
\Crefname{sublem}{Lemma}{Lemmas}
\newtheorem{claim}[equation]{Claim}

\theoremstyle{remark}
\Crefname{claim}{Claim}{Claims}

\Crefname{subrem}{Remark}{Remarks}

\theoremstyle{definition}
\newtheorem{defn}[subsection]{Definition}
\Crefname{defn}{Definition}{Definitions}

\Crefname{conv}{Convention}{Conventions}
\newtheorem{eg}[subsection]{Example}
\Crefname{eg}{Example}{Examples}
\newtheorem{rem}[subsection]{Remark}
\Crefname{rem}{Remark}{Remarks}
\newtheorem*{rems}{Remarks}

\theoremstyle{plain}
\newtheorem{thm-tweak}[subsection]{Theorem}
\Crefname{thm-tweak}{Theorem}{Theorems}
\newtheorem{lemma-tweak}[subsection]{Lemma}
\Crefname{lemma-tweak}{Lemma}{Lemmas}
\newtheorem{cor-tweak}[subsection]{Corollary}
\Crefname{cor-tweak}{Corollary}{Corollaries}
\newtheorem{prop-tweak}[subsection]{Proposition}
\Crefname{prop-tweak}{Proposition}{Propositions} 
\newtheorem{conj-tweak}[subsection]{Conjecture}
\Crefname{conj-tweak}{Conjecture}{Conjectures} 

\theoremstyle{definition}
\newtheorem{defn-tweak}[subsection]{Definition}
\Crefname{defn-tweak}{Definition}{Definitions}
\newtheorem{eg-tweak}[subsection]{Example}
\Crefname{eg-tweak}{Example}{Examples}
\newtheorem*{rems-tweak}{Remarks}
\newtheorem{rem-tweak}[subsection]{Remark}
\Crefname{rem-tweak}{Remark}{Remarks}

\newtheoremstyle{subsection-tweak}
   {11pt}
   {3pt}%
   {}
   {}%
   {\bfseries}
   {}%
   {.5em}
   {\thmnumber{\@{#1}{}\@{#2}.}%
    \thmnote{~{\bfseries#3.}}}    
    
\theoremstyle{subsection-tweak}
\newtheorem{pp}[subsection]{}
\newcommand{\bpp}{\begin{pp}}
\newcommand{\epp}{\end{pp}}

\theoremstyle{subsection-tweak}
\newtheorem{pp-tweak}[subsection]{}

\numberwithin{equation}{subsection}

\makeatletter
\def\@tocline#1#2#3#4#5#6#7{
    \begingroup 
    \@ifempty{#4}{%
    }{%
    }%

    \parindent\z@ \leftskip#3\relax \advance\leftskip\@tempdima\relax
    #5\hskip-\@tempdima
      \ifcase #1
       \or\or \hskip 2em \or \hskip 1em \else \hskip 3em \fi%
      #6\nobreak\relax
    \dotfill\hbox to\@pnumwidth{\@tocpagenum{#7}}\par
    \nobreak
    \endgroup
  }
 \def\l@section{\@tocline{1}{0pt}{1pc}{}{}}

\renewcommand{\tocsection}[3]{%
  \indentlabel{\@ifnotempty{#2}{\makebox[1.3em][l]{%
    \ignorespaces#1 \bfseries{#2}.\hfill}}}\bfseries{#3}
    \vspace{1.5pt}}

\renewcommand{\tocsubsection}[3]{%
  \indentlabel{\@ifnotempty{#2}{\hspace*{-0.5em}\makebox[2.1em][l]{%
    \ignorespaces#1#2.\hfill}}}#3
    \vspace{1.5pt}}

\makeatother

\begin{document}

\title{The $A_\Inf$-cohomology in the semistable case }

\author{K\k{e}stutis \v{C}esnavi\v{c}ius and Teruhisa Koshikawa}
\address{Laboratoire de Math\'{e}matiques d'Orsay, Univ.~Paris-Sud, CNRS, Universit\'{e} Paris-Saclay, 91405 Orsay, France}
\address{Research Institute for Mathematical Sciences, Kyoto University, Kyoto 606-8502, Japan}
\email{kestutis@math.u-psud.fr}
\email{teruhisa@kurims.kyoto-u.ac.jp}


\date{\today}


\begin{abstract} \revise{For a proper, smooth scheme $X$ over a $p$-adic field $K$, we show that any proper, flat, semistable $\cO_K$-model $\cX$ of $X$ whose logarithmic de Rham cohomology is torsion free determines the same $\cO_K$-lattice inside $H^i_\dR(X/K)$ and, moreover, that this lattice is functorial in $X$. For this, we extend the results of Bhatt--Morrow--Scholze on the construction and the analysis of an $A_\Inf$-valued cohomology theory of $p$-adic formal, proper, smooth $\cO_{\ov{K}}$-schemes $\fX$ to the semistable case.  The relation of the $A_\Inf$-cohomology to the $p$-adic \'{e}tale and the logarithmic crystalline cohomologies allows us to reprove the semistable conjecture of Fontaine--Jannsen.} \end{abstract}


\maketitle

\hypersetup{
    linktoc=page,     
}

\renewcommand*\contentsname{}
\q\\
\tableofcontents


\section{Introduction}

\revise{
\ready{
\bpp[Integral relations between $p$-adic cohomology theories] \lab{integral-rels}
For a proper, smooth scheme $X$ over a complete, discretely valued extension $K$ of $\bQ_p$ with a perfect residue field $k$, comparison isomorphisms of $p$-adic Hodge theory relate the $p$-adic \'{e}tale, de Rham, and, in the case of semistable reduction, also crystalline cohomologies of $X$. For instance, they show that for $i \in \bZ$, the $\Gal(\ov{K}/K)$-representation $H^i_\et(X_{\ov{K}}, \bQ_p)$ functorially determines the filtered $K$-vector space $H^i_\dR(X/K)$. Even though the integral analogues of these isomorphisms are known to fail in general, one may still consider their hypothetical consequences, for instance, one may ask the following.
\begin{itemize}
\item
For proper, flat, semistable $\cO_K$-models $\cX$ and $\cX'$ of $X$ endowed with their standard log structures, do the images of $H^i_{\log\dR}(\cX/\cO_K)$ and $H^i_{\log\dR}(\cX'/\cO_K)$ in $H^i_\dR(X/K)$ agree?
\end{itemize}
One of the goals of the present paper is to show that the answer is positive if the logarithmic de Rham cohomology of the models $\cX$ and $\cX'$ is torsion free (see \eqref{log-dR-X-Xpr} and \Cref{lat-dR-id}): in this case, both $H^i_{\log\dR}(\cX/\cO_K)$ and $H^i_{\log\dR}(\cX'/\cO_K)$ agree with the $\cO_K$-lattice in $H^i_\dR(X/K)$ that is functorially determined by $H^i_\et(X_{\ov{K}}, \bZ_p)$. The good reduction case of this result may be derived from the work of Bhatt--Morrow--Scholze \cite{BMS16} on integral $p$-adic Hodge theory, and our approach, as well as the bulk of this paper, is concerned with extending the framework of \emph{op.~cit.}~to the semistable case. 
\epp

\bpp[The $A_\Inf$-cohomology in the semistable case] \lab{Ainf-sst-case}
To approach the question above, we set $C \ce \wh{\ov{K}}$, let $A_\Inf \ce W(\cO_{C}^\flat)$ be the basic period ring of Fontaine, and, for a semistable $\cO_K$-model $\cX$ of $X$, similarly to the smooth case treated in \cite{BMS16}, construct the $A_\Inf$-cohomology~object
\[
R\Gamma_{A_\Inf}(\cX) \in D^{[0,\, 2\dim(X)]}(A_\Inf)
\]
that is quasi-isomorphic to a bounded complex of finite free $A_\Inf$-modules and has finitely presented cohomology $H^i_{A_\Inf}(\cX)$. We show that base changes of $R\Gamma_{A_\Inf}(\cX)$ recover other cohomology theories:
\be \lab{Ainf-coho-bc-ann}
\ba
\tst R\Gamma_{A_\Inf}(\cX) \tensor_{A_\Inf}^\bL  W(C^\flat) &\cong \tst R\Gamma_\et(X_{\ov{K}}, \bZ_p) \tensor_{\bZ_p}^\bL W(C^\flat), \\
 R\Gamma_{A_\Inf}(\cX) \tensor^\bL_{A_\Inf,\, \theta} \cO_C &\cong R\Gamma_{\log\dR}(\cX/\cO_K) \tensor^\bL_{\cO_K} \cO_C, \\
 R\Gamma_{A_\Inf}(\cX) \tensor^\bL_{A_\Inf} W(\ov{k}) &\cong R\Gamma_{\log\cris}(\cX_k/W(k)) \tensor^\bL_{W(k)} W(\ov{k}),
 \ea
\ee
see \S\ref{Ainf-coho}; here $R\Gamma_{\log\cris}$ denotes the logarithmic crystalline (that is, Hyodo--Kato) cohomology, $W(k)$ (resp.,~$\cO_K$) carries the log structure associated to $\bN_{\ge 0} \xra{1\,\mapsto\, 0,\, 0\, \mapsto\, 1} W(k)$ (resp.,~$\cO_K \setminus\{ 0\} \hra \cO_K$), and $\cX_k$ is endowed with the base change of the standard log structure $\cO_{\cX,\, \et} \cap (\cO_{\cX,\, \et}[\f{1}{p}])^\times$ of $\cX$. 

If the cohomology of $R\Gamma_{\log\dR}(\cX/\cO_K)$ is torsion free, then each $H^i_{A_\Inf}(\cX)$ is $A_\Inf$-free and the base changes \eqref{Ainf-coho-bc-ann} hold in each individual cohomological degree (see \S\ref{indiv-BC} and \Cref{free-free}). In this case, the Fargues equivalence for Breuil--Kisin--Fargues $\Gal(\ov{K}/K)$-modules  allows us to prove that
\[
\text{the $\Gal(\ov{K}/K)$-representation} \q H^i_\et(X_{\ov{K}}, \bZ_p) \q \text{determines} \q H^i_{A_\Inf}(\cX)
\]
(see \Cref{lat-dR-id}). Then $H^i_\et(X_{\ov{K}}, \bZ_p)$ also determines\footnote{\ready{The implicit functor is nonexact, as it must be: there exists a nonexact sequence of abelian schemes over $\bZ_2$ that is short exact over $\bQ_2$ (see \cite{BLR90}*{7.5/8}), so there is no \emph{exact} functor $F$ with $F(H^1_\et((-)_{\ov{\bQ}_2}, \bZ_2)) = H^1_\dR(-/\bZ_2)$.}} $H^i_{\log\dR}(\cX/\cO_K)$ (and $H^i_{\log\cris}(\cX_k/W(k))$) and, since the same reasoning applies to another model $\cX'$, the result claimed in \S\ref{integral-rels} follows.

The base changes \eqref{Ainf-coho-bc-ann} also allow us to extend the cohomology specialization results obtained in the good reduction case in \cite{BMS16}. Qualitatively, in \Cref{free-free} we show that $H^i_{\log \dR}(\cX/\cO_K)$ is torsion free if and only if so is $H^i_{\log \cris}(\cX_k/W(k))$, in which case $H^i_\et(X_{\ov{K}}, \bZ_p)$ is torsion free. Quantitatively, in \Cref{et-cris-tors,et-dR-tors} we show that for every $n \ge 0$, 
\[
\ba
\length_{\bZ_p}((H^i_\et(X_{\ov{K}}, \bZ_p)_\tors)/p^n) &\le \length_{W(k)}((H^i_{\log \cris}(\cX_k/W(k))_\tors)/p^n), \\
\length_{\bZ_p}((H^i_\et(X_{\ov{K}}, \bZ_p)_\tors)/p^n) &\le \f{1}{\length_{\cO_K}(\cO_K/p)} \cdot \length_{\cO_K}((H^i_{\log \dR}(\cX/\cO_K)_\tors)/p^n).
\ea
\]
\epp

\bpp[The semistable comparison isomorphism]
The analysis of $R\Gamma_{A_\Inf}(\cX)$, specifically, its relation to the $p$-adic \'{e}tale and the logarithmic crystalline cohomologies, permits us to reprove in \Cref{sst-comp} the semistable conjecture of Fontaine--Jansen \cite{Kat94a}*{Conj.~1.1}:
\be \lab{sst-comp-ann}
R\Gamma_\et (X_{\ov{K}}, \bZ_p)\otimes^\bL_{\bZ_p} B_\st \cong R\Gamma_{\log\cris}(\cX_{k} /W(k)) \otimes^\bL_{W(k)} B_\st.
\ee
Other proofs of this conjecture have been given in \cite{Tsu99}, \cite{Fal02}, \cite{Niz08}, 
\cite{Bha12}, \cite{Bei13}, and \cite{CN17}, whereas \cite{BMS16} used $R\Gamma_{A_\Inf}(\cX)$ to reprove the crystalline conjecture. Similarly to \cite{CN17}, we establish \eqref{sst-comp-ann} for $p$-adic \emph{formal} $\cO_K$-schemes $\cX$ that are proper, flat, and ``semistable.'' 

A key result that leads to \eqref{sst-comp-ann} is the absolute crystalline comparison isomorphism
\be \lab{abs-crys-ann}
R\Gamma_{A_\Inf}(\cX) \tensor^\bL_{A_\Inf} A_\cris \cong R\Gamma_{\log\cris}(\cX_{\cO_{\ov{K}}/p}/A_\cris)
\ee
of \Cref{Acris-BC}, whose construction in \S\ref{section-Acris} forms the technical core of this paper. This construction is based on an ``all possible coordinates'' technique that is a variant of its analogue used to establish \eqref{abs-crys-ann} in the smooth case in \cite{BMS16}*{\S12}. The presence of singularities and log structures creates additional complications that do not appear in the smooth case and are overviewed in \S\ref{section-Acris}. 

Using the absolute crystalline comparison isomorphism, in \Cref{BdR-comp} we compare the $A_\Inf$-cohomology of $\cX$ with the $B_\dR^+$-cohomology of $X$ defined by Bhatt--Morrow--Scholze in \cite{BMS16}*{\S13}:
\be \lab{BdR-comp-ann}
R\Gamma_{A_\Inf}(\cX) \tensor_{A_\Inf}^\bL B_\dR^+ \cong R\Gamma_\cris (X^\ad_C /B_{\dR}^{+}).
\ee
The identification \eqref{BdR-comp-ann} is important for ensuring that the semistable comparison \eqref{sst-comp-ann} is compatible with the de Rham comparison proved in \cite{Sch13}, and hence that it respects filtrations.

As for the question posed in \S\ref{integral-rels}, even though it only involves the \'{e}tale and the de Rham cohomologies, the resolution of its torsion free case outlined in \S\ref{Ainf-sst-case} uses both \eqref{abs-crys-ann} and \eqref{BdR-comp-ann} (so also the bulk of the material of this paper). This is because  we need to ensure that the determination of $H^i_\dR(X/K)$ by $H^i_\et(X_{\ov{K}}, \bQ_p)$ via the de Rham comparison of $p$-adic Hodge theory is compatible with the determination of $H^i_{\log \dR}(\cX/\cO_K)$ and $H^i_{\log \dR}(\cX'/\cO_K)$ by $H^i_\et(X_{\ov{K}}, \bZ_p)$ via $A_\Inf$-cohomology and Breuil--Kisin--Fargues modules. In fact, even for showing that the cohomology modules of $R\Gamma_{A_\Inf}(\cX)$ are Breuil--Kisin--Fargues, we already use the absolute crystalline comparison \eqref{abs-crys-ann}.
\epp

\bpp[The object $A\Omega_\fX$ and its base changes]
Even though above we have focused on schemes, the construction and the analysis of $R\Gamma_{A_\Inf}(-)$ works for any $p$-adic formal $\cO_C$-scheme $\fX$ that is semistable in the sense described in \S\ref{fir-setup} (see \eqref{sst-coord}) and that, whenever needed, is assumed to be proper. Specifically, for such an $\fX$, in \S\ref{AOX-def} we use the (variant for the \'{e}tale topology of the) definition of Bhatt--Morrow--Scholze from \cite{BMS16} to build an object
\[
A\Omega_\fX \in D^{\ge 0}(\fX_\et, A_\Inf) \qq \text{and to set} \qq R\Gamma_{A_\Inf}(\fX) \ce R\Gamma(\fX_\et, A\Omega_\fX).
\]
As in the smooth case of \cite{BMS16}, the relation of $R\Gamma_{A_\Inf}(\fX)$ to the $p$-adic \'{e}tale cohomology of the adic generic fiber $\fX_C^\ad$ of $\fX$ follows from the results of \cite{Sch13} (see  \S\ref{etale-section}). In turn, the relations to the logarithmic de Rham and crystalline cohomologies are the subjects of \S\ref{AOmegaX} and \S\ref{section-Acris}, respectively, and rest on the following identifications established in \Cref{dR-spec,abs-crys-comp}:
\be \lab{AOX-id-ann}
A\Omega_\fX \tensor^\bL_{A_\Inf,\, \theta} \cO_C \cong \Omega^\bullet_{\fX/\cO_C,\, \log} \qq \text{and} \qq A\Omega_{\fX}\wh{\otimes}^{\bL}_{A_\Inf}A_{\cris} \cong Ru_{*} (\cO_{\fX_{\cO_C/p}/A_{\cris}}),
\ee
where $u \colon (\fX_{\cO_C/p}/ A_{\cris})_{\log\cris} \ra  \fX_{\et}$ is the forgetful map of topoi. The arguments for  \eqref{AOX-id-ann} build on the same general skeleton as in \cite{BMS16} but differ, among other aspects, in how they handle the interaction of the Deligne--Berthelot--Ogus d\'{e}calage functor $L\eta$ used in the definition of $A\Omega_\fX$ with the intervening base changes and with the almost isomorphisms supplied by the almost purity theorem. Namely, for this, the nonflatness over the singular points of $\fX$ of the explicit perfectoid pro\'{e}tale covers that we construct makes it difficult to directly adapt the arguments from \emph{op.~cit.} Instead, we take advantage of several general results about $L\eta$ from \cite{Bha16}. Verifying their assumptions in our case amounts to the analysis in \S\ref{local-analysis} of continuous group cohomology modules built using the aforementioned perfectoid cover. The typical conclusion of this analysis is that these modules have no nonzero ``almost torsion'' and that the element $\mu \in A_\Inf$ kills  their ``nonintegral parts.''

\ready{Further and more specific overviews of our arguments are given in the beginning parts of the sections that follow. In the rest of this introduction, we fix the precise notational setup for the rest of the paper (see \S\ref{fir-setup}), discuss the logarithmic structure on $\fX$ that we later use without notational explication (see \S\ref{log-str-def}), and review the relevant general notational conventions (see~\S\ref{conv}).}
\epp

}

\bpp[The setup] \lab{fir-setup}
\ready{In what follows, we fix the following notational setup.
\begin{itemize}
\item
We fix an algebraically closed field $k$ of characteristic $p > 0$, let $C$ be the completed algebraic closure of $W(k)[\f{1}{p}]$,
and let $\fm \subset \cO_C$ be the maximal ideal in the valuation ring of $C$.

\item
For convenience, we fix an embedding $p^{\bQ} \subset C$, that is, for every prime $\ell$, we fix a system of compatible $\ell^n$-power roots $p^{1/\ell^\infty} \ce ( p^{1/\ell^n})_{n \ge 0}$ of $p$ in $\cO_C$.

\item
We fix a $p$-adic formal scheme $\fX$ over $\cO_C$ that in the \'{e}tale topology may be covered by open affines $\fU$ which admit an \'{e}tale $\cO_C$-morphism
\be  \lab{sst-coord} 
\qqq \fU = \Spf (R) \ra \Spf (R^{\square}) \q \text{with} \q R^\square \ce \cO_C\{t_0, \dotsc, t_r, t_{r + 1}^{\pm 1}, \dotsc,  t_d^{\pm 1}\}/(t_0\cdots t_r - p^q)
\ee
for some $d \ge 0$, some $0 \le r \le d$, and some $q \in \bQ_{> 0}$ (where $d$, $r$, and $q$ may depend on $\fU$).
\end{itemize}
For example, $C$ could be the completed algebraic closure of any discretely valued field $K$ of mixed characteristic $(0, p)$ with a perfect residue field. In addition, no generality is gained by allowing $p^q$ in \eqref{sst-coord} to be \emph{any} nonunit $\pi \in \cO_C \setminus \{ 0 \}$. The role of the embedding $p^{\bQ} \subset C$ is to simplify arguments with explicit charts for the log structure on $\fX$ (defined in \S\ref{log-str-def}); this is particularly useful in \S\ref{section-Acris}, especially in \S\S\ref{chart-target}--\ref{conv-chart-sm}. Our $C$ is less general than in \cite{BMS16}, where any complete algebraically closed nonarchimedean extension of $\bQ_p$ is typically allowed. One of the main reasons for this is that we want to be able to apply, especially in \S\ref{section-Acris}, certain auxiliary results from \cite{Bei13a} (besides, relations $t_0\cdots t_r - \pi$ in which $\pi$ has a nonrational valuation go beyond ``semistable reduction'').

The existence of \'{e}tale local \emph{semistable coordinates} \eqref{sst-coord} implies that each $\fX_{\cO_C/p^n}$ is flat and locally of finite presentation over $\cO_C/p^n$ and $\fX_{\cO_C/p^n}^\sm$ is dense in $\fX_{\cO_C/p^n}$. By \cite{SP}*{\href{http://stacks.math.columbia.edu/tag/04D1}{04D1}} and limit arguments, \eqref{sst-coord} is the formal $p$-adic completion of the $\ov{W(k)}$-base change of an \'{e}tale $\cO$-morphism 
\be \lab{sst-approx}
\tst U \ra \Spec\p{\cO[t_0, \dotsc, t_r, t_{r + 1}^{\pm 1}, \dotsc,  t_d^{\pm 1}]/(t_0\cdots t_r - p^q) }
\ee
for some discrete valuation subring $\cO \subset \ov{W(k)}$ that contains $p^q$. \emph{Loc.~cit.}~and \cite{GR03}*{7.1.6 (i)} also imply that $R$ is $R^\square$-flat. In addition, if $R \tensor_{\cO_C} k$ is not $k$-smooth, then $R$ determines $q$.\footnote{\ready{The following argument justifies this. Choose an $n \in \bZ_{> q}$ and let $A$ be the local ring of $\Spec(R/p^n)$ at some singular point. Without loss of generality, all the $t_i$ with $0 \le i \le r$ are noninvertible in $A$, so, in particular, $r \ge 1$. The $d$-th Fitting ideal $\Fitt_{d }(\Omega^1_{(R^\square/p^n)/(\cO_C/p^n)}) \subset R^\square/p^n$ is generated by the $r$-fold partial products $t_{0}\cdots \wh{t_{i}} \cdots t_{r}$ with $0 \le i \le r$, so the same holds for $\Fitt_{d}(\Omega^1_{A/(\cO_C/p^n)}) \subset A$ (see \cite{SGA7I}*{VI,~5.1~(a)}). Consequently, the quotient $(R^\square/p^n)/(\Fitt_{d }(\Omega^1_{(R^\square/p^n)/(\cO_C/p^n)}))$ is faithfully flat over $\cO_C/(p^q)$, and hence so is $A/(\Fitt_{d }(\Omega^1_{A/(\cO_C/p^n)}))$. It follows that $(p^q) \subset \cO_C$ is the preimage of $\Fitt_{d }(\Omega^1_{A/(\cO_C/p^n)}) \subset A$, to the effect that $R$ determines  $q$.
}}

Any smooth $p$-adic formal $\cO_C$-scheme $\fX$ meets the requirements above: indeed, then the cover $\{ \fU\}$ exists already for the Zariski topology with $r = 0$ and $q = 1$ for all $\fU$, see \cite{FK14}*{I.5.3.18}. Another key example is 
\be \lab{fX-sst}
\fX = \wh{\cX_{\cO_C}}
\ee
for some discrete valuation subring $\cO \subset \cO_C$ with a perfect residue field and a uniformizer $\pi \in \cO$ and  a locally of finite type, flat $\cO$-scheme $\cX$ that is \emph{semistable} in the sense that $\cX_{\cO/\pi}$ is a normal crossings divisor  in $\cX$ (as defined in \cite{SP}*{\href{http://stacks.math.columbia.edu/tag/0BSF}{0BSF}}), so that, in particular, $\cX$ is regular at every point of $\cX_{\cO/\pi}$.\footnote{To justify that any $\fX$ as in \eqref{fX-sst} meets the requirements, we first note that \'{e}tale locally on $\cX$ there exists a regular sequence such that the product its $r + 1$ first terms cuts out $\cX_{\cO/\pi}$. Thus, since any finite extension of $\cO/\pi$ is separable, the miracle flatness theorem \cite{EGAIV2}*{6.1.5} ensures that every $x \in \cX_{\cO/\pi}$ has an \'{e}tale neighborhood $U \ra \cX$ that admits an \'{e}tale $\cO$-morphism $U \ra \Spec(\cO[t_0, \dotsc, t_d]/(t_0\cdots t_r - \pi))$ or, equivalently, an \'{e}tale morphism
\be \lab{sst-def-map}
U \ra \Spec(\cO[t_0, \dotsc, t_r, t_{r + 1}^{\pm 1}, \dotsc, t_d^{\pm 1}]/(t_0\cdots t_r - \pi))).
\ee
\vspace{-11pt}}
Moreover, if $\cX$ is even \emph{strictly semistable} in the sense that $\cX_{\cO/\pi}$ is even a strict normal crossings divisor in $\cX$ (as defined in \cite{SP}*{\href{http://stacks.math.columbia.edu/tag/0BI9}{0BI9}}), then the \'{e}tale maps \eqref{sst-def-map} exist even Zariski locally on $\cX$, and so also the cover $\{\fU \}$ exists already for the Zariski topology of $\fX$.

\begin{itemize}
\item
We let $\fX_C^\ad$ denote the adic generic fiber of $\fX$. By \eqref{sst-coord} and \cite{Hub96}*{3.5.1}, the adic space $\fX_C^\ad$ is smooth over $C$; by \cite{Hub96}*{1.3.18~ii)}, if $\fX$ is $\cO_C$-proper, then $\fX_C^\ad$ is $C$-proper. 

\item
We let $(\fX_C^\ad)_\proet$ denote the pro\'{e}tale site of $\fX_C^\ad$ (reviewed in \cite{BMS16}*{\S5.1} and defined in \cite{Sch13}*{3.9} and \cite{Sch13e}*{(1)}) and let
\be\lab{nu-def}
\nu \colon (\fX_C^\ad)_\proet \ra \fX_\et
\ee
be the morphism to the \'{e}tale site of $\fX$ that sends any \'{e}tale $\fU \ra \fX$ to the constant pro-system associated to its adic generic fiber. 
By \cite{SP}*{\href{http://stacks.math.columbia.edu/tag/00X6}{00X6}}, this functor indeed defines a morphism of sites: by \cite{Hub96}*{3.5.1}, it preserves coverings, commutes with fiber products, and respects final objects. Thus, $\nu$ induces a morphism of topoi $(\nu\i, \nu_*)$ (see \cite{SP}*{\href{http://stacks.math.columbia.edu/tag/00XC}{00XC}}).
\end{itemize}
}
\epp

\ready{
\bpp[The logarithmic structure on $\fX$] \lab{log-str-def}
Unless noted otherwise, we always equip 
\benuma
\item \lab{OC-log-str}
the ring $\cO_C$ (resp.,~$\cO_C/p^n$ or $k$) with the log structure $\cO_C\setminus \{ 0 \} \hra \cO_C$ (resp.,~its pullback);

\item \lab{X-log-str}
the formal scheme $\fX$ (resp.,~$\fX_{\cO_C/p^n}$ or $\fX_k$) with the log structure given by the subsheaf associated to the subpresheaf\footnote{\ready{The subpresheaf and its associated subsheaf necessarily agree on every quasi-compact object $\fU$ of $\fX_\et$.}} $\cO_{\fX,\, \et} \cap (\cO_{\fX,\, \et}[\f{1}{p}])^\times \hra \cO_{\fX,\, \et}$ (resp.,~its pullback log structure).
\eenum
Both \ref{OC-log-str} and \ref{X-log-str} determine the same log structure on $\Spf(\cO_C)$, so the map $\fX \ra \Spf(\cO_C)$ is that of log formal schemes. Moreover, \'{e}tale locally on $\fX$, the log structure may be made explicit: in the presence of a coordinate morphism \eqref{sst-coord}, \Cref{log-claim,log-comp} below give an explicit chart for the log structure of $\fU$, namely, the chart \eqref{log-cl-prelog} in which we replace $\cO$ by $\cO_C$, replace $U$ by $\fU$, and set $\pi \ce p^q$.
\ready{This chart shows, in particular, that $\fU$ and $\cO_C$ may be endowed with fine log structures whose base changes along a ``change of log structure'' self-map of $\cO_C$ recover the log structures described in \ref{OC-log-str}--\ref{X-log-str} (for example, the fine log structure on $\cO_C$ could be the one determined by the chart $\bN_{\ge 0} \xra{a\, \mapsto (p^q)^a} \cO_C$, in which case the ``change of log structure'' self-map of $\cO_C$ is the identity on the underlying scheme $\Spec(\cO_C)$ and is determined on the log structures by the map of charts $\bN_{\ge 0} \xra{a\, \mapsto (p^q)^a} \cO_C \setminus \{ 0\}$). Since many common properties of maps of log schemes are stable under base change, in practice this means that we may often deal with the log structures in \ref{OC-log-str}--\ref{X-log-str} as if they were fine and, in particular, we may cite \cite{Kat89} for certain purposes. 

By the preceding discussion, all the log structures above are quasi-coherent and integral. Moreover, by \cite{Kat89}*{3.7~(2)}, each $\fX_{\cO_C/p^n}$ is log smooth over $\cO_C/p^n$, so that, by \cite{Kat89}*{3.10}, the $\cO_\fX$-module $\Omega^1_{\fX/\cO_C,\, \log}$ of logarithmic differentials is finite locally free.} We set 
\[
\tst \Omega^i_{\fX/\cO_C,\, \log} \ce \bigwedge^i \Omega^1_{\fX/\cO_C,\, \log},
\]
let $\Omega^\bullet_{\fX/\cO_C,\, \log}$ denote the logarithmic de Rham complex, and set
\[
R\Gamma_{\log\dR}(\fX/\cO_C) \ce R\Gamma(\fX_\et, \Omega^\bullet_{\fX/\cO_C,\, \log}).
\]
\epp

\bcl \lab{log-claim}
For a valuation subring $\cO \subset \ov{W(k)}$ and an $\cO$-scheme $U$ that has an \'{e}tale morphism
\[
U \ra \Spec\p{\cO[t_0, \dotsc, t_r, t_{r + 1}^{\pm 1}, \dotsc,  t_d^{\pm 1}]/(t_0\cdots t_r - \pi) }  \qq \text{for some nonunit} \qq \pi \in \cO \setminus \{ 0\},
\]
the log structure on $U$ associated to $\cO_{U,\, \et} \cap (\cO_{U,\, \et}[\f{1}{p}])^\times$ has the chart \be \lab{log-cl-prelog}
\bN_{\ge 0}^{r + 1} \sqcup_{\bN_{\ge 0}} (\cO \setminus \{ 0 \}) \ra \Gamma(U, \cO_U)
\ee
given by $(a_i)_{0 \le i \le r} \mapsto \prod_{0 \le i \le r} t_i^{a_i}$ on $\bN_{\ge 0}^{r + 1}$, the diagonal $\bN_{\ge 0} \ra \bN_{\ge 0}^{r + 1}$ and $\bN_{\ge 0} \xra{a\, \mapsto \pi^a} (\cO \setminus \{ 0\})$ on $\bN_{\ge 0}$, and the structure map $(\cO \setminus \{ 0 \}) \ra \Gamma(U, \cO_U)$ on $\cO \setminus \{ 0 \}$.
\ecl

\bpf
Without loss of generality, $U$ is affine, so, by a limit argument, 
we may assume that $\cO$ is discretely valued. Then $U$, endowed with the log structure associated to \eqref{log-cl-prelog}, is logarithmically regular in the sense of \cite{Kat94b}*{2.1} (compare with \cite{Bei12}*{\S4.1, proof of Lemma}). Therefore, since the locus of triviality of this log structure is $U[\f{1}{p}]$, the claim follows from \cite{Kat94b}*{11.6}.
\epf

\bcl \lab{log-comp}
For $\cO$ as in Claim \uref{log-claim}, a flat $\cO$-scheme $U$ \up{resp.,~and its formal $p$-adic completion $\fU$} endowed with the log structure associated to $\cO_{U,\, \et} \cap (\cO_{U,\, \et}[\f{1}{p}])^\times$ \up{resp.,~$\cO_{\fU,\, \et} \cap (\cO_{\fU,\, \et}[\f{1}{p}])^\times$}, 
\be \lab{formal-pc}
\text{the formal $p$-adic completion morphism } \q j \colon \fU \ra U \q \text{of log ringed \'{e}tale sites is strict.}
\ee
\ecl

\bpf
For a geometric point $\ov{u}$ of $\fU$, due to \cite{SP}*{\href{http://stacks.math.columbia.edu/tag/04D1}{04D1}}, the stalk map $\cO_{U,\, \ov{u}} \cong  j\i(\cO_{U,\, \ov{u}}) \ra \cO_{\fU,\, \ov{u}}$ induces an isomorphism $\cO_{U,\, \ov{u}}/p^n \cong \cO_{\fU,\, \ov{u}}/p^n$ for every $n > 0$. We consider the stalk map 
\be \lab{stalk-map}
\tst \cO_{U,\, \ov{u}} \cap (\cO_{U,\, \ov{u}}[\f{1}{p}])^\times \cong j\i(\cO_{U,\, \ov{u}} \cap (\cO_{U,\,\ov{u}}[\f{1}{p}])^\times) \ra \cO_{\fU,\, \ov{u}} \cap (\cO_{\fU,\, \ov{u}}[\f{1}{p}])^\times.
\ee
Every element $x$ of the target of \eqref{stalk-map} satisfies the equation $xy = p^n$ for some $n > 0$. We choose an $\wt{x} \in \cO_{U,\, \ov{u}}$ congruent to $x$ modulo $p^{n + 1}$, so that $\wt{x}\wt{y} = p^n + p^{n + 1}\wt{z}$ for some $\wt{y}, \wt{z} \in \cO_{U,\, \ov{u}}$. Since $1 + p\wt{z} \in \cO_{U,\, \ov{u}}^\times$, we adjust $\wt{y}$ to get $\wt{x}\wt{y} = p^n$, which shows that $\wt{x} \in \cO_{U,\, \ov{u}} \cap (\cO_{U,\, \ov{u}}[\f{1}{p}])^\times$ and $(p^n) \subset (\wt{x})$. Thus, the image of $\wt{x}$ in $\cO_{\fU,\, \ov{u}}$ and $x$ generate the same ideal, and hence are unit multiples of each other. Conversely, if $\wt{x}_1, \wt{x}_2 \in \cO_{U,\, \ov{u}} \cap (\cO_{U,\, \ov{u}}[\f{1}{p}])^\times$ are unit multiples of each other in $\cO_{\fU,\, \ov{u}}$, then, by reducing modulo $p^n$ for a large enough $n$, we see that they generate the same ideal in $\cO_{U,\, \ov{u}}$, so are unit multiples of each other already in $\cO_{U,\, \ov{u}}$. In conclusion, the map \eqref{stalk-map} induces an isomorphism
\[
\tst (\cO_{U,\, \ov{u}} \cap (\cO_{U,\, \ov{u}}[\f{1}{p}])^\times)/\cO_{U,\, \ov{u}}^\times \isomto (\cO_{\fU,\, \ov{u}} \cap (\cO_{\fU,\, \ov{u}}[\f{1}{p}])^\times)/\cO_{\fU,\, \ov{u}}^\times,
\]
to the effect that the map \eqref{formal-pc} is indeed strict, as claimed. \epf

\bpp[Conventions and additional notation] \lab{conv}
\ready{
For a field $K$, we let $\ov{K}$ be its algebraic closure (taken inside $C$ if $K$ is given as a subfield of $C$). If $K$ has a valuation, we let $\cO_K$ be its valuation subring and write $\ov{\cO}_K$ for the integral closure of $\cO_K$ in $\ov{K}$. In mixed characteristic, we normalize  the valuations by requiring that $v(p) = 1$. We let $(-)^\sm$ denote the smooth locus of a (formal) scheme over an implicitly understood base.  For power series rings, we use $\{ -\}$ to indicate decaying coefficients. For a topological ring $R$, we let $R^\circ$ denote the subset of powerbounded elements.

We let $W(-)$ (resp.~$W_n(-)$) denote $p$-typical Witt vectors (resp.,~their length $n$ truncation), and let $[-]$ denote Teichm\"{u}ller representatives. We let $\bZ_{(p)}$ be the localization of $\bZ$ at $p$, let $\mu_{p^n}$ be the group scheme of $p^n$-th roots of unity, and let $\zeta_{p^n}$ denote a primitive $p^n$-th root of unity. For brevity, we set $\bZ_p(1) \ce \varprojlim \p{\mu_{p^n}(C)}$. We let $\wh{M}$ denote the (by default, $p$-adic) completion of a module $M$ and, similarly, let $\wh{\bigoplus}$ denote the completion of a direct sum. Unless specified otherwise, we endow a $p$-adically complete module with the inverse limit of the discrete topologies.

We use the definition of a perfectoid ring given in \cite{BMS16}*{3.5} (the compatibility with prior definitions is discussed in \cite{BMS16}*{3.20}). Explicitly, by \cite{BMS16}*{3.9 and 3.10}, a $p$-torsion free ring $S$ is \emph{perfectoid} if and only if $S$ is $p$-adically complete and the divisor $(p) \subset S$ has a $p$-power root in the sense that there is a $\varpi \in S$ with $(\varpi^p) = (p)$ and $\xymatrix@C=30pt{S/\varpi S \ar[r]_{x \mapsto x^p}^{\sim} &S/p S.}$ In particular, for such an $S$, any $p$-adically formally \'{e}tale $S$-algebra $S'$ that is $p$-adically complete is also perfectoid.

For a ring object $R$ of a topos $\sT$, we write $D(\sT, R)$, or simply $D(R)$, for the derived category of $R$-modules.  For an object $M$ of a derived category, we denote its derived $p$-adic completion  by
\be \lab{der-comp-def}
\tst \wh{M} \ce R\lim_n (M \otimes^{\bL}_{\bZ} \bZ/p^n\bZ), \qq \text{and also set} \qq * \wh{\tensor}^\bL_\cdot - \ce R\lim_n( (* \tensor^\bL_\cdot - ) \tensor^\bL_\bZ \bZ/p^n\bZ )
\ee
(see \cite{SP}*{\href{http://stacks.math.columbia.edu/tag/0940}{0940}} for the definition of $R\lim$). 
For a morphism $f$ of ringed topoi, we use the commutativity of the functor $Rf_*$ with derived limits and derived completions, see \cite{SP}*{\href{http://stacks.math.columbia.edu/tag/0A07}{0A07} and \href{http://stacks.math.columbia.edu/tag/0944}{0944}}.

For a profinite group $H$ and a continuous $H$-module $M$, we write $R\Gamma_\cont(H, M)$ for the continuous cochain complex. Whenever convenient, we also view $R\Gamma_\cont(H, -)$ as the derived global sections functor of the site of profinite $H$-sets (see \cite{Sch13}*{3.7~(iii)} and \cite{Sch13e}*{(1)}). 

For commuting endomorphisms $f_1, \dotsc, f_n$ of an abelian group $A$, we recall the \emph{Koszul complex}:
\be \lab{Koszul-def}
\tst K_A(f_1, \dotsc, f_n) \ce A \tensor_{\bZ[x_1, \dotsc, x_n]} \bigotimes_{i = 1}^n \p{\bZ[x_1, \ldots, x_n] \xra{x_i} \bZ[x_1, \ldots, x_n]},
\ee
where $A$ is regarded as a $\bZ[x_1, \ldots, x_n]$-module by letting $x_j$ act as $f_{j}$, the tensor products are over $\bZ[x_1, \dotsc, x_n]$, and the factor complexes are concentrated in degrees $0$ and $1$.

For an ideal $I$ of a ring $R$ and an $R$-module complex $(M^\bullet, d^\bullet)$ with $M^j \cong 0$ for $j < 0$, the subcomplex
\be \lab{eta-def}
\eta_I(M^\bullet) \subset M^\bullet \qq \text{is defined by} \qq (\eta_I(M^\bullet))^j \ce \{ m \in I^jM^j \, \vert \, d^j(m) \in I^{j + 1}M^{j + 1}\}.
\ee
We will mostly (but \emph{a priori} not always, see \Cref{comp-map-0}) use $\eta_I(M^\bullet)$ as in \cite{BMS16}*{6.2}, namely, when $I$ is generated by a nonzerodivisor and the $M^j$ have no nonzero $I$-torsion.

A \emph{logarithmic divided power thickening} (or, for brevity, a \emph{log PD thickening}) is an exact closed immersion of logarithmic (often abbreviated to log) schemes equipped with a divided power structure on the quasi-coherent sheaf of ideals that defines the underlying closed immersion of schemes.}
\epp

\subsection*{Acknowledgements}
We thank Bhargav Bhatt and Matthew Morrow for writing the surveys \cite{Bha16} and \cite{Mor16}, which have been useful for preparing this paper. We thank the referee for helpful comments and suggestions. We thank Bhargav Bhatt, Pierre Colmez, Ravi Fernando, Luc Illusie, Arthur-C\'{e}sar Le Bras, Matthew Morrow, Wies\l awa Nizio\l, Arthur Ogus, Peter Scholze, Joseph Stahl, Jakob Stix, Jan Vonk, and Olivier Wittenberg for helpful conversations or correspondence. We thank the Kyoto Top Global University program for providing the framework in which this collaboration started. We thank the Miller Institute at the University of California Berkeley, the Research Institute for Mathematical Sciences at Kyoto University, and the University of Bonn for their support during the preparation of this article. 
}}


\subfile{etale}


\subfile{local-analysis}




\subfile{AOmega}


\subfile{Acris}


\subfile{BdR}


\subfile{RG}


\subfile{lattice}


\subfile{semistable}

\end{document}

%% file: etale.tex

\section{The object $A\Omega_\fX$ and the $p$-adic \'{e}tale cohomology of $\fX$} \lab{etale-section}

\ready{
\ready{As in the case when $\fX$ is smooth treated in \cite{BMS16}, the eventual construction of the $A_\Inf$-cohomology modules of $\fX$ rests on the object $A\Omega_\fX$ that lives in a derived category of $A_\Inf$-module sheaves on $\fX$. In this short section, we review the definition of $A\Omega_\fX$ in \S\ref{AOX-def} and then, in the case when $\fX$ is proper, review the connection between $A\Omega_\fX$ and the integral $p$-adic \'{e}tale cohomology of $\fX_C^\ad$ in \Cref{RG-et-id}. We begin by fixing the basic notation that concerns the ring $A_\Inf$ of integral $p$-adic Hodge theory. The setup of \S\S\ref{Ainf-not}--\ref{AOX-def} will be used freely in the rest of the paper.}

\bpp[The ring $A_\Inf$] \lab{Ainf-not}
\ready{We denote the tilt of $\cO_C$ by
\[
\tst \cO_C^\flat \ce \varprojlim_{y \mapsto y^p} \p{\cO_C/p}, \q \text{so that, by reduction mod $p$,} \q \varprojlim_{y \mapsto y^p} \cO_C \isomto \varprojlim_{y \mapsto y^p}\p{ \cO_C/p} = \cO_C^\flat
\]
as multiplicative monoids (see \cite{Sch12}*{3.4 (i)}). We regard $p^{1/p^\infty}$ fixed in \S\ref{fir-setup} as an element of $\cO_C^\flat$. Due to the fixed embedding $p^{\bQ_{\ge 0}} \subset \cO_C$, this element comes equipped with well-defined powers $(p^{1/p^\infty})^q \in \cO_C^\flat$ for $q \in \bQ_{\ge 0}$. For each $x \in \cO_C^\flat$, we let $(\ldots, x^{(1)}, x^{(0)})$ denote its preimage in $\varprojlim_{y\mapsto y^p} \cO_C$. The map $x \mapsto \val_{\cO_C}(x^{(0)})$ makes $\cO_C^\flat$ a complete valuation ring of height $1$ whose fraction field $C^\flat \ce \Frac(\cO_C^\flat)$ is algebraically closed (see \cite{Sch12}*{3.4~(iii), 3.7~(ii)}). We let $\fm^\flat$ denote the maximal ideal of $\cO_C^\flat$.

The basic period ring $A_\Inf$ of Fontaine is defined by
\[
\tst  A_\Inf \ce W(\cO_C^\flat) \qq \text{and comes equipped with the Witt vector Frobenius} \qq \varphi \colon A_\Inf \isomto A_\Inf.
\]
We equip the local domain $A_\Inf$ with the product of the valuation topologies via the Witt coordinate bijection $W(\cO_C^\flat) \cong \prod_{n = 1}^\infty \cO_C^\flat$. Then $A_\Inf$ is complete and its topology agrees with the $(p, [x])$-adic topology for any nonzero nonunit $x \in \cO_C^\flat$. We fix (once and for all) a compatible system $\eps = (\dotsc, \zeta_{p^2}, \zeta_p, 1)$ of $p$-power roots of unity in $\cO_C$, so that $\eps \in \cO_C^\flat$, and set 
\be \lab{mu-intro}
\mu \ce [\eps] - 1 \in A_\Inf.
\ee
 Since $(p, \mu) = (p, [\eps - 1])$, the topology of $A_\Inf$ is $(p, \mu)$-adic. By forming the limit of the sequences
\be \lab{mu-SES}
0 \ra W_n(\cO_C^\flat) \xra{\mu} W_n(\cO_C^\flat) \ra W_n(\cO_C^\flat)/\mu \ra 0,
\ee
we see that $A_\Inf/\mu$ is $p$-adically complete and that the ideal $(\mu) \subset A_\Inf$ does not depend on the choice of $\eps$ (use the fact that the valuation of $\zeta_p - 1$ does not depend on $\zeta_p$).

The assignment $[x] \mapsto x^{(0)}$ extends uniquely to a ring homomorphism
\be\lab{theta-intro}
\theta \colon A_\Inf \surjects \cO_C, \qq \text{the \emph{de Rham specialization} map of $A_\Inf$,}
\ee
which is surjective, as indicated, and intertwines the Frobenius $\varphi$ of $A_\Inf$ with the absolute Frobenius of $\cO_C/p$. Its kernel $\Ker(\theta) \subset A_\Inf$ is principal and  generated by the element 
\be \lab{xi-intro}
\tst \xi \ce \sum_{i = 0}^{p - 1} [\eps^{i/p}]
\ee
(see \cite{BMS16}*{3.16}). Analogues of the sequences \eqref{mu-SES} show that each $A_\Inf/\xi^n$ is $p$-adically complete. In fact, the map $\theta$ identifies $A_\Inf/\xi^{n}$ with the initial $p$-adically complete infinitesimal thickening of $\cO_C$ of order $n - 1$, see \cite{SZ17}*{3.13}. The composition 
\[ 
\theta \circ \varphi\i\colon A_\Inf \surjects \cO_C \qq \text{is the \emph{Hodge--Tate specialization} map of $A_\Inf$,}
\]
and its kernel is generated by the element $\varphi(\xi) = \sum_{i = 0}^{p - 1} [\eps^{i}]$. 

Due to the nature of our $C$ (see \S\ref{fir-setup}), the ring $\cO_C/p$ is a $k$-algebra, so $A_\Inf$ is a $W(k)$-algebra.}
\epp


\ready{
\bpp[The object $A\Omega_\fX$] \lab{AOX-def}
The operations that define $\cO_C^\flat$ and $A_\Inf$ make sense on the pro\'{e}tale site $(\fX_C^\ad)_\proet$: namely, as in \cite{Sch13}*{4.1, 5.10, and 6.1}, we have the integral completed structure sheaf 
\be \lab{O-hat-plus}
\tst \wh{\cO}_{\fX_C^\ad}^+ \ce \varprojlim_n (\cO_{\fX_C^\ad,\, \proet}^+/p^n), \qq \text{its tilt} \qq \wh{\cO}_{\fX_C^\ad}^{+,\, \flat} \ce \varprojlim_{y \mapsto y^p} (\cO_{\fX_C^\ad,\, \proet}^+/p),
\ee
and the basic period sheaf
\[
\bA_{\inf,\, \fX_C^\ad} \ce W(\wh{\cO}_{\fX_C^\ad}^{+,\, \flat}).
\]

For brevity, we often denote these sheaves simply by $\wh{\cO}^+$, $\wh{\cO}^{+,\, \flat}$, and $\bA_\Inf$. Affinoid perfectoids form a basis for $(\fX_C^\ad)_\proet$ (see \cite{Sch13}*{4.7}) and the construction of the map $\theta$ of \eqref{theta-intro} makes sense for any perfectoid $\cO_C$-algebra (see \cite{BMS16}*{\S3}). In particular,  $\bA_{\Inf,\, \fX_C^\ad}$ comes equipped with the map 
\be \lab{theta-sheaf-intro}
\tst \theta_{\fX_C^\ad}\colon \bA_{\Inf,\, \fX_C^\ad} \ra \wh{\cO}_{\fX_C^\ad}^+,
\ee
which, by construction, is compatible with the map $\theta\colon A_\Inf \ra \cO_C$, intertwines the Witt vector Frobenius $\varphi$ of $\bA_{\inf,\, \fX_C^\ad}$ with the absolute Frobenius of $\wh{\cO}_{\fX_C^\ad}^+/p$, and, by \cite{Sch13}*{6.3 and 6.5}, is surjective with $\Ker(\theta_{\fX_C^\ad}) = \xi \cdot \bA_{\Inf,\, \fX_C^\ad}$ (in addition, $\xi$ is not a zero divisor in $\bA_{\Inf,\, \fX_C^\ad}$).

The key object that we are going to study in this paper is
\be \lab{AOX-def-eq}
A\Omega_\fX \ce L\eta_{(\mu)}(R\nu_*( \bA_{\Inf,\, \fX_C^\ad})) \in D^{\ge 0}(\fX_{\et}, A_\Inf),
\ee
where the d\'{e}calage functor $L\eta$ of \cite{BMS16}*{\S6} is formed with respect to the ideal $(\mu)$ of the constant sheaf $A_\Inf$ of $\fX_\et$ (the definition of $L\eta_{(\mu)}$ builds on the formula \eqref{eta-def} for $\eta_{(\mu)}$). The formula \eqref{AOX-def-eq} may also be executed with the Zariski site $\fX_\Zar$ as the target of $\nu$, and it then defines the object
\be \lab{AO-Zar-def}
A\Omega_{\fX_\Zar} \in D^{\ge 0}(\fX_\Zar, A_\Inf),
\ee
which is  the $A\Omega_\fX$ that was used in \cite{BMS16}. We will only use $A\Omega_{\fX_\Zar}$ in \Cref{Zar-et-comp} (and in some results that lead to it) for comparison to $A\Omega_\fX$.

Since $\varphi(\mu) = \varphi(\xi)\mu$, by \cite{BMS16}*{6.11}, we have $L\eta_{(\varphi(\mu))} \cong L\eta_{(\varphi(\xi))} \circ L\eta_{(\mu)}$, so the Frobenius automorphism of $\bA_{\Inf,\, \fX_C^\ad}$ gives the Frobenius morphism
\be\lab{Frob-mor}
\tst A\Omega_\fX \tensor^\bL_{A_\Inf,\, \varphi} A_\Inf \cong L\eta_{(\varphi(\xi))}(A\Omega_\fX) \xra{\text{ \cite{BMS16}*{6.10 and 3.17~(ii)} }}   A\Omega_{\fX} \q \text{in} \q D^{\ge 0}(\fX_{\et}, A_\Inf),
\ee
which, by \cite{BMS16}*{6.14}, induces an isomorphism
\be \lab{Frob-iso}
\tst (A\Omega_\fX \tensor^\bL_{A_\Inf,\, \varphi} A_\Inf)[\f{1}{\varphi(\xi)}] \isomto (A\Omega_{\fX})[\f{1}{\varphi(\xi)}].
\ee
In addition, by \emph{loc.~cit.},~we also have 
\be \lab{inv-mu-no-L}
\tst A\Omega_\fX \tensor^\bL_{A_\Inf} A_\Inf[\f{1}{\mu}] \cong (R\nu_* (\bA_{\Inf,\, \fX_C^\ad})) \tensor_{A_\Inf}^\bL A_\Inf[\f{1}{\mu}],
\ee
so a result of Scholze \cite{BMS16}*{5.6} supplies the following relation to integral $p$-adic \'{e}tale cohomology:
\epp

\bthm \lab{RG-et-id}
If $\fX$ is proper over $\cO_C$, then there is an identification
\be \lab{RG-et-id-eq}
\tst R\Gamma(\fX_\et, A\Omega_\fX) \tensor^\bL_{A_\Inf} A_\Inf[\f{1}{\mu}] \cong  R\Gamma_\et(\fX_C^\ad, \bZ_p) \tensor_{\bZ_p}^\bL A_\Inf[\f{1}{\mu}]. \qedhere
\ee
\ethm

In broad strokes, the proof of \Cref{RG-et-id} given in \emph{loc.~cit.}~goes as follows: one considers the map
\be \lab{RG-under-map}
R\Gamma_\et(\fX_C^\ad, \bZ_p) \tensor^\bL_{\bZ_p} A_\Inf \cong R\Gamma_\proet(\fX_C^\ad, \bZ_p) \tensor^\bL_{\bZ_p} A_\Inf \ra R\Gamma_{\proet}(\fX_C^\ad, \bA_{\Inf,\, \fX_C^\ad})
\ee
induced by the inclusion $A_\Inf \hra \bA_{\Inf,\, \fX_C^\ad}$ and deduces from the almost purity theorem with, for instance, \Cref{pump-zero} below that the ideal
\be \lab{W-mflat-def}
W(\fm^\flat) \ce \Ker(W(\cO_C^\flat) \surjects W(k)) \q \text{of} \q A_\Inf
\ee
 kills the cohomology of its cone. Since $\mu$ lies in $W(\fm^\flat)$ and we have the identification \eqref{inv-mu-no-L}, it follows that the map \eqref{RG-under-map} induces the identification \eqref{RG-et-id-eq}.

\brem
In practice, $\fX$ often arises as the formal $p$-adic completion of a proper, finitely presented $\cO_C$-scheme $\cX$. In this situation, $\fX_C^\ad$ agrees with the adic space associated to $\cX_C$ (see \cite{Con99}*{5.3.1 4.},~\cite{Hub94}*{4.6~(i)}, and \cite{Hub96}*{1.9.2~ii)}) and, by \cite{Hub96}*{3.7.2}, we have
\[
R\Gamma_\et(\fX_C^\ad, \bZ_p) \cong R\Gamma_\et(\cX_{C}, \bZ_p).
\]
\erem
}}

%% file: local-analysis.tex

\section{The local analysis of $A\Omega_\fX$} \lab{local-analysis}

\ready{
Even though the definition of the object $A\Omega_\fX$ given in \eqref{AOX-def-eq} is global, the key computations that will eventually relate it to the logarithmic de Rham and crystalline cohomologies are local and are presented in this section. Under the assumption that $\fX$ has a coordinate morphism as in \eqref{sst-coord} (or \eqref{et-coord} below), their basic goal is to express the cohomology of the pro\'{e}tale sheaf $\bA_{\inf,\, \fX_C^\ad}$, at least after applying $L\eta_{(\mu)}$, in terms of continuous group cohomology formed using an explicit perfectoid pro\'{e}tale cover $\fX^\ad_{C,\,\infty}$ of $\fX_C^\ad$ (see \Cref{ANB-c}). The basic relation of this sort is supplied by the almost purity theorem, so the key point is to explicate the appearing group cohomology modules well enough in order to eliminate the ``almost'' ambiguities inherent in this theorem with the help of \Cref{Bha-lemma} below that comes from \cite{Bha16}. We first carry out this program for the simpler sheaf $\wh{\cO}_{\fX_C^\ad}^+$, and then build on this case to address $\bA_{\inf,\, \fX_C^\ad}$.

In comparison to the local analysis carried out in the smooth case in \cite{BMS16}, one complication is that the perfectoid cover of $\fX$ that gives rise to $\fX^\ad_{C,\,\infty}$ is not flat over the singular points of $\fX_k$. This makes it difficult to transfer various arguments with ``$q$-de Rham complexes'' across the coordinate morphism \eqref{et-coord}. In fact, we avoid $q$-de Rham complexes altogether and instead phrase the intermediate steps of the local analysis purely in terms of continuous group cohomology modules.

\bpp[The local setup] \lab{loc-loc-set-set}
We assume throughout \S\ref{local-analysis} that $\fX = \Spf( R)$ and for some $d \ge 0$, some $0 \le r \le d$, and some $q \in \bQ_{> 0}$, there is an \'{e}tale $\Spf(\cO_C)$-morphism as in \eqref{sst-coord}:
\be \lab{et-coord}
\fX = \Spf(R) \ra \Spf(R^{\square}) \equalscolon \fX^\square \q \text{with} \q R^\square \ce \cO_C\{t_0, \dotsc, t_r, t_{r + 1}^{\pm 1}, \dotsc,  t_d^{\pm 1}\}/(t_0\cdots t_r - p^q).
\ee
Due to our assumptions from \S\ref{fir-setup}, a general $\fX$ is of this form on a basis for its \'{e}tale topology.
\epp

\bpp[The perfectoid cover $\fX^\ad_{C,\, \infty}$] \lab{Rinfty}
\ready{For each $m \ge 0$, we consider the $R^\square$-algebra
\[
R_m^\square \ce \cO_C\{t_0^{1/p^m}, \dotsc, t_r^{1/p^m}, t_{r + 1}^{\pm 1/p^m}, \dotsc,  t_d^{\pm 1/p^m}\}/(t_0^{1/p^m}\cdots t_r^{1/p^m} - p^{q/p^m}), \q \text{and} \q R_\infty^\square \ce \p{\varinjlim R_m^\square}\wh{\ },
\]
where, as always unless mentioned otherwise (see \S\ref{conv}), the completion is $p$-adic.
Explicitly, we have the $p$-adically completed direct sum decomposition
\be \lab{Rinfty-dec}
R_\infty^\square \cong \wh{\bigoplus}_{\substack{ (a_0, \ldots, a_d) \in (\bZ[\f{1}{p}]_{\ge 0})^{\oplus (r + 1)} \oplus (\bZ[\f{1}{p}])^{\oplus(d - r)}, \\ \text{$a_{j} = 0$ for some $0 \le j \le r$}}}\ \  \cO_C \cdot t_0^{a_0}\cdots t_d^{a_d},
\ee
which shows that $R_\Infty^\square$ is perfectoid (see \S\ref{conv}) and that, for each $m \ge 0$, the map $R_m^\square \ra R_\infty^\square$  is an inclusion of an $R_m^\square$-module direct summand comprised of those summands $\cO_C \cdot t_0^{a_0}\cdots t_d^{a_d}$ of \eqref{Rinfty-dec} for which $p^ma_j \in \bZ$ for all $j$. 

The corresponding $R$-algebras are 
\[
R_m \ce  R \tensor_{R^\square} R_m^\square  \qq \text{and} \qq R_\infty \ce \p{\varinjlim R_m}\wh{\ }\, \cong (R \tensor_{R^{\square}} R_\infty^\square){\wh{\,\ }}.
\]
Each $R_m$ (resp.,~$R_\infty$) is a $p$-torsion free $p$-adically formally \'{e}tale $R_m^\square$-algebra (resp.,~$R_\infty^\square$-algebra). In particular, $R_\infty$ is perfectoid (see \S\ref{conv}). By \cite{GR03}*{7.1.6 (ii)}, each $R_m$ is $p$-adically complete. 

The summands in \eqref{Rinfty-dec} with $a_j \not\in \bZ$ for some $0 \le j \le d$ comprise an $R^\square$-submodule $M_\infty^\square$ of $R_\infty^\square$, and we set $M_\infty \ce R \wh{\tensor}_{R^\square} M_\infty^\square $. Thus, we have the $R^\square$-module (resp.,~$R$-module) decomposition
\be \lab{Rinfty-M-dec}
R_\infty^\square \cong R^\square \oplus M_\infty^\square \qq \text{(resp.,} \qq R_\infty \cong R \oplus M_\infty).
\ee

The profinite group
\[
\tst \Delta \ce \left\{ (\eps_0,  \dotsc, \eps_d) \in \p{\varprojlim_{m \ge 0} (\mu_{p^m}(\cO_C))}^{\oplus (d + 1)}\, \Big |\ \eps_0\cdots \eps_r = 1 \right\} \simeq \bZ_p^{\oplus d}
\]
acts $R^\square$-linearly on $R_m^\square$ by scaling each $t_j^{1/p^m}$ by the $\mu_{p^m}$-component of $\eps_j$. The induced actions of $\Delta$ on $R_\Infty^\square$ and $R_\infty$ are continuous, compatible, and preserve the decompositions \eqref{Rinfty-dec} and \eqref{Rinfty-M-dec}. 
In terms of the element $\eps$ fixed in \S\ref{Ainf-not}, $\Delta$ is topologically freely generated by the following $d$ elements:
\[
\ba
\  \delta_i &\ce (\eps\i, 1, \dotsc, 1, \eps, 1, \dotsc, 1) \q \text{for $i = 1, \ldots, r$,} \q \text{where the $0$-th and $i$-th entries are nonidentity;} \\
\  \delta_i &\ce (1, \dotsc, 1, \eps, 1, \dotsc, 1) \q \text{for $i = r + 1, \ldots, d$,} \q\, \text{where the $i$-th entry is nonidentity.} 
\ea
\]

After inverting $p$, for each $m \ge 0$, we have
\[
\tst R^\square_m[\f{1}{p}] \cong \bigoplus_{a_1, \dotsc, a_d \in \{0, \f{1}{p^m}, \dotsc, \f{p^m - 1}{p^m}\}} R^\square[\f{1}{p}] \cdot t_1^{a_1}\cdots t_d^{a_d}, 
\]
so $R_m^\square[\f{1}{p}]$ is the $R^\square[\f{1}{p}]$-algebra obtained by adjoining the $(p^m)\th$ roots of $t_1, \dotsc, t_d \in (R^\square[\f{1}{p}])^\times$, and hence is finite \'{e}tale over $R^\square[\f{1}{p}]$. Therefore, $\varinjlim_m\p{ R^\square_m[\f{1}{p}]}$ is a pro-(finite \'{e}tale) $\Delta$-cover of $R^\square[\f{1}{p}]$.  The explicit description \eqref{Rinfty-dec} implies that $R_m^\square = (R_m^\square[\f{1}{p}])^\circ$, so the pro-object 
\[
\tst (\fX^\square)^\ad_{C,\,\infty} \ce \varprojlim \Spa(R_m^\square[\f{1}{p}], R_m^\square), \q \text{which determines the perfectoid space} \q \Spa(R_\infty^\square[\f{1}{p}], R_\infty^\square),
\]
is an affinoid perfectoid pro-(finite \'{e}tale) $\Delta$-cover of the adic generic fiber $(\fX^\square)^\ad_C$ of $\Spf(R^\square)$; in particular, $(\fX^\square)^\ad_{C,\,\infty}$ is an affinoid perfectoid object of the pro\'{e}tale site $((\fX^\square)^\ad_C)_\proet$. Consequently, the $\fX_C^\ad$-base change of $(\fX^\square)^\ad_{C,\, \infty}$, namely, the tower 
\[
\tst \fX_{C,\,\infty}^{\ad} \ce \varprojlim \Spa(R_m[\f{1}{p}], R_m), \q \text{which determines the perfectoid space} \q \Spa(R_\infty[\f{1}{p}], R_\infty),
\]
is an affinoid perfectoid pro-(finite \'{e}tale) $\Delta$-cover of $\fX_C^{\ad}$, so, in particular, is an affinoid perfectoid object of $(\fX^\ad_C)_\proet$. 

By \cite{Sch13}*{4.10~(iii)}, the value on $\fX_{C,\,\infty}^{\ad}$ of the sheaf $\wh{\cO}_{\fX_C^\ad}^+$ reviewed in \eqref{O-hat-plus}  is the ring $R_\infty$.
}
\epp

\bpp[The cohomology of $\wh{\cO}^+$ and continuous group cohomology] \lab{alm-pur}
 By \cite{Sch13}*{3.5, 3.7 (iii) and its proof, 6.6} (see also \cite{Sch13e}), the \v{C}ech complex of the sheaf $\wh{\cO}^+_{\fX_C^\ad}$ with respect to the pro-(finite \'{e}tale) affinoid perfectoid cover 
 \[
\fX^\ad_{C,\, \infty} \surjects \fX_C^\ad 
\]
  is identified with the continuous cochain complex $R\Gamma_\cont(\Delta, R_\infty)$. 
In particular, by using \cite{SP}*{\href{http://stacks.math.columbia.edu/tag/01GY}{01GY}}, we obtain the edge map to the pro\'{e}tale cohomology of $\wh{\cO}^+_{\fX_C^\ad}$: 
\be \lab{e-ebox}
e \colon R\Gamma_\cont(\Delta, R_\infty) \ra R\Gamma_\proet(\fX_C^\ad, \wh{\cO}^+), 
\ee
which on the level of cohomology is described by the Cartan--Leray spectral sequence (see \emph{loc.~cit.}~or \cite{SGA4II}*{V.3.3}).
By the almost purity theorem \cite{Sch13}*{4.10 (v)}, the maximal ideal $\fm \subset \cO_C$ kills the cohomology groups of $\Cone(e)$. 
\epp

We will show in \Cref{e-iso} that $L\eta_{(\zeta_p - 1)}(e)$ is an isomorphism, so that $L\eta_{(\zeta_p - 1)}(R\Gamma_\proet(\fX_C^\ad, \wh{\cO}^+))$ is computed in terms of continuous group cohomology. For this, we will use the following lemma.

\blem[\cite{BMS16}*{8.11~(i)}] \lab{iso-on-L}
An $\cO_C$-module map $f\colon M \ra M'$ with $M[\fm] =\p{\f{M}{(\zeta_p - 1) M}}[\fm] = 0$ and both $\Ker f$ and $\Coker f$ killed by $\fm$  induces an isomorphism $\f{M}{M[\zeta_p - 1]} \isomto \f{M'}{M'[\zeta_p - 1]}$. \QED
\elem

In order to apply \Cref{iso-on-L}, we will check in \Cref{Rinfty-clean} that the cohomology modules $H^i_\cont(\Delta, R_\infty)$ have no nonzero $\fm$-torsion. This will use the following general lemmas.

\blem \lab{no-almost}
For an inclusion $\fo \subset \fO$ of a discrete valuation ring into a nondiscrete valuation ring of rank $1$, if $N$ is an $\fo$-module and $\fM \subset \fO$ denotes the maximal ideal, then $(N \tensor_\fo \fO)[\fM] = 0$.
\elem

\bpf
The $\fo$-flatness of $\fO$ reduces us to the case when $N$ is finitely generated, so it suffices to observe that $(\fO/(a))[\fM] = 0$ whenever $a \in \fO$.
\epf

\blem \lab{coho-big-sum}
Fix an $i \in \bZ_{\ge 0}$, let  $H$ be a profinite group, let $\{M_j\}_{j \in J}$ be $p$-adically complete, $p$-torsion free, continuous $H$-modules, and suppose that either
\benumr
\item \lab{CBS-i}
the group $H^i_\cont(H, M_j)$ is $p$-torsion free for every $j$, or

\item \lab{CBS-ii}
some $p^n$ kills $H^i_\cont(H, M_j)$ for every $j$.
\eenum
Then the following map is injective\ucolon
\[
\tst H^i_\cont(H, \wh{\bigoplus}_{j \in J} M_j) \hra \prod_{j \in J} H^i_\cont(H, M_j), \qq \text{where the completion is $p$-adic.}
\]
In particular, in the case \ref{CBS-i} \up{resp.,~\ref{CBS-ii}}, $H^i_\cont(H, \wh{\bigoplus}_{j \in J} M_j)$ is $p$-torsion free \up{resp.,~killed by $p^n$}.
\elem

\bpf
Let $c$ be a continuous $\p{\wh{\bigoplus}_{j \in J} M_j}$-valued $i$-cocycle that represents an element of the kernel. For each $j$, let $c_j$ be the ``$j$-th coordinate'' of $c$. We discard the $j$ with $c_j = 0$ and, for each remaining $j$, we choose the maximal $n_j \in \bZ_{\ge 0}$ such that $c_j$ is $(p^{n_j} M_j)$-valued, so that the function $j \mapsto n_j$ is finite-to-one.  Since each $M_j$ is $p$-torsion free, each $p^{-n_j} c_j $ is an $M_j$-valued continuous $i$-cocycle. 

In the case \ref{CBS-i}, the class of $p^{-n_j}c_j$ in $H^i_\cont(H, M_j)$ vanishes, so each $c_j$ is the coboundary of a $(p^{n_j}M_j)$-valued continuous $(i - 1)$-cochain $b_j$. In the case \ref{CBS-ii}, $p^n$ kills $H^i_\cont(H, M_j)$, so $c_j$ is the coboundary of a $(p^{n_j - n}M_j)$-valued continuous $(i - 1)$-cochain $b_j$ whenever $n_j \ge n$.

In both cases, the $b_j$'s exhibit $c$ as a continuous coboundary.
\epf

\blem[\cite{BMS16}*{7.3~(ii)}] \lab{Koszul}
Let $H$ be a profinite group isomorphic to $\bZ_p^{\oplus d}$ for some $d > 0$, and let $M \cong \varprojlim_{n \ge 1} M_n$ be a continuous $H$-module with each $M_n$ a discrete, $p^n$-torsion, continuous $H$-module. For any $\gamma_1, \dotsc, \gamma_d \in H$ that topologically freely generate $H$, there is a natural identification
\[
R\Gamma_\cont(H, M) \cong K_M(\gamma_1 - 1, \dotsc, \gamma_d - 1), \q \text{so also} \q H^j_\cont(H, M) \cong H^j(K_M(\gamma_1 - 1, \dotsc, \gamma_d - 1)),
\]
in the derived category {\upshape(}see \uS\uref{conv} for the notation{\upshape)}. \QED 
\elem


\bprop \lab{Rinfty-clean}
The element $\zeta_p - 1$ kills the $\cO_C$-modules $H^i_\cont(\Delta, M_\infty)$. Moreover, for each $b \in \cO_C$, the $\cO_C$-modules $R_\infty/b$ and $H^i_\cont(\Delta, R_\infty/b)$ have no nonzero $\fm$-torsion.
\eprop

\bpf
Let $S \ce \cO_C \cdot t_0^{a_0}\cdots t_d^{a_d}$ be a summand of \eqref{Rinfty-dec}. By \Cref{Koszul}, 
the $\cO_C$-module $H^i_\cont(\Delta, S)$ is the $i$-th cohomology of the $\cO_C$-tensor product of $d$ complexes of the form $\cO_C \xra{\zeta - 1} \cO_C$ for suitable $p$-power roots of unity $\zeta$. Moreover, since the $d$ complexes may be defined over some discrete valuation subring of $\cO_C$, \Cref{no-almost} ensures that 
\be \lab{RC-eq-0}
 H^i_\cont(\Delta, S) \qq \text{has no nonzero $\fm$-torsion.}
\ee
If $S$ contributes to $M_\infty$, that is, if $a_j \not\in \bZ$ for some $j$, then some $\zeta$ is not $1$, and the corresponding factor complex is quasi-isomorphic to $\cO_C/(\zeta - 1)$. Thus, in this case, 
\be \lab{RC-eq-00}
\zeta - 1, \q \text{and hence also} \q \zeta_p - 1, \q \text{kills} \q H^i_\cont(\Delta, S).
\ee
For $m > 0$, let $M_m^\square$ denote the $p$-adically completed direct sum of those summands $\cO_C \cdot t_0^{a_0}\cdots t_d^{a_d}$ of \eqref{Rinfty-dec} for which $m$ is the smallest \emph{nonnegative} integer with $p^m\cdot (a_0, \dotsc, a_d) \in \bZ^{\oplus(d + 1)}$. \Cref{coho-big-sum} and \eqref{RC-eq-0}--\eqref{RC-eq-00} imply that the $\cO_C$-module
\be \lab{RC-eq-1}
 H^i_\cont(\Delta, M_m^\square) \qq \text{has no nonzero $\fm$-torsion and is killed by $\zeta_p - 1$.}
\ee
Since $R$ is $R^\square$-flat and $R\tensor_{R^\square} M_m^\square $ is $p$-adically complete (see \S\ref{fir-setup} and \S\ref{Rinfty}), \Cref{Koszul} gives
\be \lab{RC-eq-2}
 H^i_\cont(\Delta, R \tensor_{R^\square} M_m^\square) \cong R \tensor_{R^\square} H^i_\cont(\Delta, M_m^\square).
\ee
Since $M_\infty \cong \wh{\bigoplus}_m (R \tensor_{R^\square} M_m^\square )$, \eqref{RC-eq-1}--\eqref{RC-eq-2} and \Cref{coho-big-sum} imply that $\zeta_p - 1$ kills $H^i_\cont(\Delta, M_\infty)$.

Since $R_\infty/b$ is $p$-adically complete and each of the summands of the decomposition
\[
\tst  R_\infty/(b, p^n) \cong R/(b, p^n) \oplus \bigoplus_{m > 0} (R \tensor_{R^\square}M_m^\square )/(b, p^n) \qq \text{for $n > 0$}
 \]
may be defined over a suitably large discrete valuation subring of $\cO_C$, \Cref{no-almost} ensures that $R_\infty/b$ has no nonzero $\fm$-torsion. In addition, the $\Delta$-action on each summand may be defined over a possibly larger such subring, so, by \Cref{no-almost,Koszul}, in the case $b \neq 0$ each
\[
 H^i_\cont(\Delta, (R \tensor_{R^\square} M_m^\square )/b), \q \text{so also} \q H^i_\cont(\Delta, M_\infty/b), \q  \text{has no nonzero $\fm$-torsion.}
\]
This conclusion extends to the case $b = 0$ because the $(\zeta_p - 1)$-annihilation of $H^i_\cont(\Delta, M_\infty)$ supplies the injection $H^i_\cont(\Delta, M_\infty) \hra H^i_\cont(\Delta, M_\infty/(\zeta_p - 1))$. It remains to observe that the $\cO_C$-module $H^i_\cont(\Delta, R/b)$ also has no nonzero $\fm$-torsion: $\Delta$ acts trivially on $R/b$, so \Cref{Koszul} ensures that $H^i_\cont(\Delta, R/b)$ is a direct sum of copies of $R/b$.
\epf

\bthm \lab{e-iso}
The edge map $e$ defined in \eqref{e-ebox} induces the isomorphism
\[
L\eta_{(\zeta_p - 1)}(e) \colon  L\eta_{(\zeta_p - 1)}(R\Gamma_\cont(\Delta, R_\infty)) \isomto L\eta_{(\zeta_p - 1)}(R\Gamma_\proet(\fX_C^\ad, \wh{\cO}^+)).
\]
\ethm

\bpf
\Cref{Rinfty-clean} ensures that the $\cO_C$-modules $H^i_\cont(\Delta, R_\infty)$ have no nonzero $\fm$-torsion and that $\f{H^i_\cont(\Delta,\, R_\infty)}{H^i_\cont(\Delta,\, R_\infty)[\zeta_p - 1]} \cong \f{H^i_\cont(\Delta,\, R)}{H^i_\cont(\Delta,\, R)[\zeta_p - 1]}$. Since $\Delta$ acts trivially on $R$, this last quotient is a finite direct sum of copies of $R$ (see \Cref{Koszul}), so, by \Cref{Rinfty-clean}, it has no nonzero $\fm$-torsion. Consequently, since $\fm$ kills the kernel and the cokernel of each map
\[
H^i(e) \colon H^i_\cont(\Delta, R_\infty) \ra H^i(\fX_C^\ad, \wh{\cO}^+)
\]
(see \S\ref{alm-pur}), \Cref{iso-on-L} applies to these maps and gives the desired conclusion. 
\epf

\brem \lab{e-iso-more}
\Cref{e-iso} extends as follows: for any profinite group $\Delta'$ equipped with a continuous surjection $\Delta' \surjects \Delta$ and any pro-(finite \'{e}tale) affinoid perfectoid $\Delta'$-cover
\[
\tst \Spa(R_\infty'[\f{1}{p}], R_\infty') \ra \Spa(R[\f{1}{p}], R) \cong \fX_C^\ad \qq \text{that refines the $\Delta$-cover} \qq \fX_{C,\, \infty}^\ad \ra \fX_C^\ad \qq \text{of \S\ref{Rinfty}}
\]
compatibly with the surjection $\Delta' \surjects \Delta$,
 the edge map $e'$ defined analogously to \eqref{e-ebox} induces the isomorphism
\[
L\eta_{(\zeta_p - 1)}(e') \colon  L\eta_{(\zeta_p - 1)}(R\Gamma_\cont(\Delta', R'_\infty)) \isomto L\eta_{(\zeta_p - 1)}(R\Gamma_\proet(\fX_C^\ad, \wh{\cO}^+)).
\]
Indeed, by the almost purity theorem \cite{Sch13}*{4.10~(v)}, the ideal $\fm$ kills the cohomology of $\Cone(e')$ (in addition to that of $\Cone(e)$), so the octahedral axiom (see \cite{BBD82}*{1.1.7.1}) ensures that it also kills the  cohomology of the cone of the map $R\Gamma_\cont(\Delta, R_\infty) \ra R\Gamma_\cont(\Delta', R'_\infty)$; \Cref{iso-on-L} then applies to this map and combines with \Cref{e-iso} to give the claim.
\erem

The main goal of this section is an analogue of \Cref{e-iso} for the sheaf $\bA_{\Inf,\, \fX_C^\ad}$ (see \Cref{ANB-c}). To prepare for it, in \S\ref{tilt} and \S\ref{Ainf} we describe the values of the sheaves $\wh{\cO}^{+,\, \flat}_{\fX_C^\ad}$ and $\bA_{\Inf,\, \fX_C^\ad}$ on~$\fX_{C,\, \infty}^\ad$.

\bpp[The tilt $R^{\flat}_\infty$] \lab{tilt}
Thanks to the explicit description \eqref{Rinfty-dec} of the perfectoid ring $R_\infty^\square$, its tilt $(R_\infty^{\square})^\flat \ce \varprojlim_{y \mapsto y^p} (R_\infty^\square/p)$ is described explicitly by the identification
\[
\ba
\tst (R_\infty^{\square})^\flat &\cong \tst \p{\varinjlim_m\, \p{ \cO_C^\flat [x_0^{1/p^m}, \dotsc, x_r^{1/p^m}, x_{r + 1}^{\pm 1/p^m}, \dotsc,  x_d^{\pm 1/p^m}]/(x_0^{1/p^m} \cdots x_r^{1/p^m} - (p^{1/p^\infty})^{q/p^m})}}^{\wh{\ \ \ }} \\
 & \cong \wh{\bigoplus}_{\substack{ (a_0, \ldots, a_d) \in (\bZ[\f{1}{p}]_{\ge 0})^{\oplus (r + 1)} \oplus (\bZ[\f{1}{p}])^{\oplus(d - r)}, \\ \text{$a_{j} = 0$ for some $0 \le j \le r$}}} \ \ \cO_C^\flat \cdot x_0^{a_0}\cdots x_d^{a_d},
\ea
\]
where $x_i^{1/p^m}$ corresponds to the $p$-power compatible sequence $(\ldots, t_i^{1/p^{m + 1}}, t_i^{1/p^m})$ of elements of $R_\infty^\square$, the completions are $p^{1/p^\infty}$-adic, and the decomposition is as $\cO_C^\flat$-modules. Thus,
\[
\tst \text{the tilt} \q R^{\flat}_\infty \ce \varprojlim_{y \mapsto y^p} (R_\infty/p) \q \text{of the perfectoid ring} \q R_\infty
\]
is identified with the $p^{1/p^\infty}$-adic completion of any lift of the \'{e}tale $R_\infty^\square/p$-algebra $R_\infty/p$ to an \'{e}tale $(R_\infty^\square)^\flat$-algebra (such a lift exists, see \cite{SP}*{\href{http://stacks.math.columbia.edu/tag/04D1}{04D1}}). By \cite{Sch13}*{5.11~(i)}, the value on $\fX_{C,\, \infty}^\ad$ of the sheaf $\wh{\cO}^{+,\, \flat}_{\fX_C^\ad}$ reviewed in \eqref{O-hat-plus} is the ring $R_\infty^\flat$.

By functoriality, the group $\Delta$ acts continuously and $\cO_C^\flat$-linearly on $(R_\infty^\square)^\flat$ and $R_\infty^\flat$. Explicitly, $\Delta$ respects the completed direct sum decomposition and an $(\eps_0, \dotsc, \eps_d) \in \Delta$ scales $x_j^{a_j}$ by $\eps_j^{a_j} \in \cO_C^\flat$.
\epp

Our analysis in \S\ref{Ainf} of the value on~$\fX_{C,\, \infty}^\ad$ of the sheaf $\bA_{\Inf,\, \fX_C^\ad}$  will hinge on the following lemmas.

\blem \lab{Rflat-clean}
Both $R_\infty^\flat/b$ and $H^i_\cont(\Delta, R_\infty^\flat/b)$ for each $b \in \cO_C^\flat\setminus \{ 0\}$ have no nonzero $\fm^\flat$-torsion.
\elem

\bpf
We may assume that $b \in \fm^\flat$, so, by using Frobenius, that $b\mid p^{1/p^\infty}$ in $\cO_C^\flat$. Then \Cref{Rinfty-clean} and the $\Delta$-isomorphism $R_\infty^\flat/b \cong R_\infty/b^\sharp$ for some $b^\sharp \in \cO_C$ gives the claim.
\epf

\blem \lab{A-Apr}
For any affinoid perfectoid $\Spa(R_\infty'[\f{1}{p}], R_\infty')$ over $\Spa(C, \cO_C)$, the ring
\[
\bA_\Inf(R_\infty') \ce W((R_\infty')^\flat) \qq \text{{\upshape(}resp.,} \qq \bA_{\Inf}(R'_{\infty})/\mu)
\]
is $(p, \mu)$-adically complete \up{resp.,~$p$-adically complete}. Moreover, for any $n, n' > 0$, the sequence $(p^n, \mu^{n'})$ is $\bA_\Inf(R_\infty')$-regular and the $A_\Inf/(p^n, \mu^{n'})$-algebra $\bA_{\Inf}(R'_{\infty})/(p^n, \mu^{n'})$ is flat.
\elem

\bpf
By its definition, the perfect $\cO_C^\flat$-algebra $(R_\infty')^\flat \ce \varprojlim_{y \mapsto y^p} (R_\infty'/p)$ has no nonzero $p^{1/p^\infty}$-torsion (that is, it is $\cO_C^\flat$-flat), so the regular sequence claim follows from \cite{SP}*{\href{http://stacks.math.columbia.edu/tag/07DV}{07DV}}. The formal criterion of flatness \cite{BouAC}*{Ch.~III, \S5.2, Thm.~1~(i)$\Leftrightarrow$(iv)} then implies the $A_\Inf/(p^n, \mu^{n'})$-flatness of $\bA_{\Inf}(R'_{\infty})/(p^n, \mu^{n'})$ (even with $n' = 0$). In addition, the short exact sequences \eqref{mu-SES} with $(R_\infty')^\flat$ in place of $\cO_C^\flat$ imply the $p$-adic completeness of $\bA_{\Inf}(R'_{\infty})/\mu$.

Analogously to the case of $A_\Inf$ discussed in \S\ref{Ainf-not}, we use the Witt coordinate bijection and the $\mu$-adic topology on $(R'_\infty)^{\flat}$ to topologize $\bA_\Inf(R'_\infty) \cong \prod_{n = 1}^\infty (R'_\infty)^{\flat}$ and we see that this topology agrees with the $(p, \mu)$-adic topology. Thus, $\bA_\Inf(R'_\infty)$ is $(p, \mu)$-adically complete. 
\epf

\bpp[The ring $\bA_\Inf(R_\infty)$] \lab{Ainf}
By \cite{Sch13}*{6.5~(i)}, the value on~$\fX_{C,\, \infty}^\ad$ of the sheaf $\bA_{\Inf,\, \fX_C^\ad}$ is the ring
\[
\bA_\Inf(R_\infty) \ce W(R^{\flat}_\infty).
\]
By \Cref{A-Apr} and the formal criterion of flatness, 
$\bA_\Inf(R_\infty)$ is $(p, \mu)$-adically formally flat as an $A_\Inf$-algebra and $(p, \mu)$-adically formally \'{e}tale as an $\bA_\Inf(R_\infty^\square)$-algebra. By using, in addition, \Cref{Rflat-clean}, we see that each quotient
\be \lab{Ainf-clean-1}
\bA_\Inf(R_\infty)/(p^n, \mu^{n'}), \qq \text{so also} \qq \bA_\Inf(R_\infty)/\mu, \qq \text{has no nonzero $W(\fm^\flat)$-torsion.}
\ee

In general, for a perfect $\bF_p$-algebra $A$, the Witt ring $W(A)$ is the unique $p$-adically complete $p$-torsion free $\bZ_p$-algebra $\wt{A}$ equipped with an isomorphism $\wt{A}/p\simeq A$ (see \cite{Bha16}*{2.5}). For an $a \in A$, the Teichm\"{u}ller $[a] \in \wt{A}$ is $\displaystyle \lim_{n \ra \infty} (\wt{a}_n^{p^n})$ where $\wt{a}_n \in \wt{A}$ is any lift of $a^{1/p^n}$ (see \cite{Bha16}*{2.4}). Therefore, 
\[\ba
\tst \bA_\Inf(R_\infty^\square) & \tst\cong \p{\varinjlim_m A_\Inf [X_0^{1/p^m}, \dotsc, X_r^{1/p^m}, X_{r + 1}^{\pm 1/p^m}, \dotsc,  X_d^{\pm 1/p^m}]/(\prod_{i = 0}^r X_i^{1/p^m} - [(p^{1/p^\infty})^{q/p^m}])}^{\wh{\ \ \ }} \\
&\cong \tst \wh{\bigoplus}_{\substack{ (a_0, \ldots, a_d) \in (\bZ[\f{1}{p}]_{\ge 0})^{\oplus (r + 1)} \oplus (\bZ[\f{1}{p}])^{\oplus(d - r)}, \\ \text{$a_{j} = 0$ for some $0 \le j \le r$}}} \ \ A_\Inf \cdot X_0^{a_0}\cdots X_d^{a_d},
\ea \]
where the completions are $(p, \mu)$-adic, the decomposition is as $A_\Inf$-modules, and, in terms of \S\ref{tilt}, we have $X_i^{1/p^m} = [x_i^{1/p^m}]$. 
The summands for which $a_i \in \bZ$ for all $i$ comprise a subring 
\be \lab{AR-box}
A(R^\square) \cong A_\Inf\{X_0, \dotsc, X_r, X_{r + 1}^{\pm 1}, \dotsc,  X_d^{\pm 1}\}/(X_0\cdots X_r - [(p^{1/p^\infty})^q]) \q \text{inside} \q \bA_\Inf(R_\infty^\square),
\ee
where the convergence is $(p, \mu)$-adic. The remaining summands, that is, those for which $a_i \not\in\bZ$ for some $i$, comprise an $A(R^\square)$-submodule $N_\infty^\square \subset \bA_\Inf(R_\infty^\square)$.

On sections over $\fX_{C,\, \infty}^\ad$, the map $\theta$ from \eqref{theta-sheaf-intro} is identified with the unique ring homomorphism  
\[
\theta\colon \bA_\Inf(R_\infty) \surjects R_\infty \qq \text{such that} \qq [x] \mapsto x^{(0)},
\] 
is surjective with the kernel generated by the regular element $\xi$ (see \cite{BMS16}*{3.10, 3.11}), and intertwines the Witt vector Frobenius of $\bA_\Inf(R_\infty)$ with the absolute Frobenius of $R_\infty/p$. Thus, 
\be \lab{AR-surj}
 \theta\colon A(R^\square) \surjects R^\square \qq \text{is described by} \qq X_i \mapsto t_i.
\ee
We use the surjection \eqref{AR-surj} to uniquely lift the \'{e}tale $R^\square/p$-algebra $R/p$ to a $(p, \mu)$-adically complete, formally \'{e}tale $A(R^\square)$-algebra $A(R)$.
By construction, we have the identification
\be \lab{A-etale-id}
\bA_\Inf(R_\infty) \cong \bA_\Inf(R^\square_\infty) \wh{\otimes}_{A(R^\square)} A(R),
\ee
where the completion is $(p, \mu)$-adic. Therefore, by setting $N_\infty \ce N_\infty^\square \wh{\otimes}_{A(R^\square)}A(R)$, we arrive at the decompositions of $\bA_{\Inf}(R_\infty^\square)$ and $\bA_{\Inf}(R_\infty)$ into ``integral'' and ``nonintegral'' parts:
\be \lab{dec-dec}
\bA_{\Inf}(R_\infty^\square) \cong A(R^\square)\oplus N_\infty^\square \qq \text{and} \qq 
\bA_{\Inf}(R_\infty) \cong A(R)\oplus N_\infty.
\ee
Modulo $\Ker \theta$ (that is, modulo $\xi$), these decompositions reduce to the decompositions \eqref{Rinfty-M-dec}.

The Witt vector Frobenius of $\bA_\Inf(R_\infty^\square)$ preserves $A(R^\square)$; explicitly: it is semilinear with respect to the Frobenius of $A_\Inf$ and raises each $X_i^{1/p^m}$ to the $p$-th power. By construction,  $A(R)$ inherits a Frobenius ring endomorphism from $A(R^\square)$, and the identification \eqref{A-etale-id} is Frobenius-equivariant.

The natural $\Delta$-action on $\bA_{\Inf}(R_\infty)$ is continuous and commutes with the Frobenius. Explicitly, $\Delta$ respects the completed direct sum decomposition and an $(\eps_0, \dotsc, \eps_d) \in \Delta$ scales $X_j^{a_j}$ by $[\eps_j^{a_j}] \in A_\Inf$. 
The $\Delta$-action on  $A(R^\square)$ lifts uniquely to a necessarily Frobenius-equivariant $\Delta$-action on  $A(R)$. In particular, $\Delta$ acts trivially on $A(R)/\mu$. The identifications \eqref{A-etale-id} and \eqref{dec-dec} are $\Delta$-equivariant. 
\epp

\bpp[The cohomology of $\bA_\Inf$ and continuous group cohomology] \lab{e-reedged}
Similarly to \S\ref{alm-pur}, 
the \v{C}ech complex of the sheaf $\bA_{\Inf,\, \fX_C^\ad}$ with respect to the pro-(finite \'{e}tale) affinoid perfectoid cover $\fX^\ad_{C,\,\infty} \ra \fX_C^\ad$ 
is identified with the continuous cochain complex $R\Gamma_{\cont}(\Delta, \bA_{\Inf}(R_\Infty))$. 
Thus, by using \cite{SP}*{\href{http://stacks.math.columbia.edu/tag/01GY}{01GY}}, we obtain the edge map to the pro\'{e}tale cohomology of $\bA_{\Inf,\, \fX_C^\ad}$:
\be \lab{another-e}
e \colon R\Gamma_{\cont}(\Delta, \bA_{\Inf}(R_\Infty)) \ra R\Gamma_\proet(\fX^\ad_{C}, \bA_{\Inf}). 
\ee
By the almost purity theorem, more precisely, by \cite{Sch13}*{6.5 (ii)}, the subset $[\fm^\flat] \subset A_\Inf$ that consists of the Teichm\"{u}ller lifts of the elements in the maximal ideal $\fm^\flat \subset \cO_C^\flat$ kills all the cohomology groups of $\Cone(e)$. 
Since $\mu \in W(\fm^\flat)$ \up{see \eqref{W-mflat-def}}, it will be useful to strengthen this annihilation as follows.
\epp

\blem \lab{kill-by-more}
The ideal $W(\fm^\flat) \subset A_\Inf$ defined in \eqref{W-mflat-def} kills each $H^i(\Cone(e))$.
\elem

\bpf
We argue similarly to \cite{BMS16}*{proof of Thm.~5.6}. 
Both the source and the target of $e$ 
are derived $p$-adically complete (see \S\ref{conv}), so, by \cite{BS15}*{3.4.4 and 3.4.14}, each $ H^i(\Cone(e))$ is also derived $p$-adically complete. Thus, the desired conclusion follows from the following lemma.
\epf

\blem \lab{pump-zero}
If $[\fm^\flat] A_\Inf$ kills a derived $p$-adically complete $A_\Inf$-module $H$, then so does $W(\fm^\flat)$.
\elem

\bpf
By the derived $p$-adic completeness, any free $A_\Inf$-module resolution $F^\bullet$ of $H$ satisfies
\[
\tst H\cong  \Coker\p{\varprojlim_n(F^{-1}/p^n) \ra \varprojlim_n(F^0/p^n)}.
\]
Moreover, for every $n \ge 1$ the ideals $[\fm^\flat]\cdot W_n(\cO_C^\flat)$ and $W_n(\fm^\flat) \ce \Ker(W_n(\cO_C^\flat) \ra W_n(k))$ of $W_n(\cO_C^\flat)$ agree. Thus, the $([\fm^\flat] A_\Inf)$-annihilation of $H$ implies that $W_n(\fm^\flat)$ kills both 
\[
H/p^n \cong H^0(F^\bullet \tensor_{A_\Inf} A_\Inf/p^n) \qq \text{and} \qq \Tor_1^{A_\Inf}(H, A_\Inf/p^n) \cong H^{-1}(F^\bullet \tensor_{A_\Inf} A_\Inf/p^n).
\]
Thus, since $[\fm^\flat]^2 = [\fm^\flat]$ and $F_0/p^n$ has no nonzero $m$-torsion for every nonzero $m \in [\fm^\flat]$, any element $x \in W_{n + 1}(\fm^\flat) \cdot (F_0/p^{n + 1})$ may be lifted to $W_{n + 1}(\fm^\flat) \cdot (F_{ - 1}/p^{n + 1})$, compatibly with a specified lift of its image $\ov{x} \in W_{n}(\fm^\flat) \cdot (F_0/p^{n})$ to $W_{n}(\fm^\flat) \cdot (F_{ -1}/p^{n})$. In particular, $W(\fm^\flat) \cdot (\varprojlim_n (F^0/p^n))$ lies in the image of $\varprojlim_n (F^{-1}/p^n)$, that is, $W(\fm^\flat)$ kills $H$, as desired.
\epf

We will show in \Cref{ANB-c} that $L\eta_{(\mu)}(e)$ is an isomorphism, so that continuous group cohomology computes $L\eta_{(\mu)}(R\Gamma_\proet(\fX_C^\ad, \bA_\Inf))$. For this, we will use the following lemma.

\blem \lab{Bha-lemma}
If $B \xra{b} B'$ is a morphism in $D(A_\Inf)$ such that each $H^i(B \tensor_{A_\Inf}^\bL A_\Inf/\mu)$ has no nonzero $W(\fm^\flat)$-torsion and $W(\fm^\flat)$ kills each $H^i(\Cone(b))$, then $L\eta_{(\mu)}(b)$ is an isomorphism.
\elem

\bpf
Since $L\eta$ is not a triangulated functor, the fact that $L\eta_{(\mu)}(\Cone(b)) \cong 0$ does not \emph{a priori} suffice. Instead, the ideal $(W(\fm^\flat))^2$ kills the cohomology of $\Cone(b) \tensor^\bL_{A_\Inf} A_\Inf/\mu$, so the sequences
\[
0 \ra H^i(B \tensor^\bL_{A_\Inf} A_\Inf/\mu) \ra H^i(B' \tensor^\bL_{A_\Inf} A_\Inf/\mu) \ra H^i(\Cone(b) \tensor^\bL_{A_\Inf} A_\Inf/\mu) \ra 0
\]
are short exact. By the Bockstein construction (see \cite{BMS16}*{6.12}), as $i$ varies, they comprise a short exact sequence whose terms are complexes that compute $L\eta_{(\mu)}(B) \tensor^\bL_{A_\Inf} A_\Inf/\mu$, etc. Thus, the vanishing of $L\eta_{(\mu)}(\Cone(b))$ implies that $(L\eta_{(\mu)}(b))\tensor^\bL_{A_\Inf} A_\Inf/\mu$ is an isomorphism. It follows that $\Cone(L\eta_{(\mu)}(b)) \tensor^\bL_{A_\Inf} A_\Inf/\mu \cong 0$, so $\mu$ acts invertibly on the cohomology of $\Cone(L\eta_{(\mu)}(b))$. But then, as we see after applying $- \tensor^\bL_{A_\Inf} A_\Inf[\f{1}{\mu}]$, this cohomology vanishes.
\epf

We now verify that the edge map $e$ defined in \eqref{another-e} also meets the first assumption of \Cref{Bha-lemma}.

\bprop \lab{coho-Ainf}
For each $i \in \bZ$, the $A_\Inf$-module $H^i_\cont(\Delta, \bA_{\Inf}(R_\Infty)/\mu)$ is $p$-torsion free and $p$-adically complete\uscolon moreover, the following natural maps are isomorphisms\ucolon
\be \lab{ANB-red}
 H^i_\cont(\Delta, \bA_\Inf(R_\Infty)/\mu) \tensor_{A_\Inf} A_\Inf/p^{n} \isomto H^i_\cont(\Delta, \bA_\Inf(R_\Infty)/(\mu, p^{n})) \q \text{for} \q n > 0
\ee
and 
\be \lab{coho-lim}
\tst H^i_\cont(\Delta, \bA_{\Inf}(R_\Infty)/\mu) \isomto \varprojlim_n \p{H^i_\cont(\Delta, \bA_{\Inf}(R_\Infty)/(\mu, p^n))}.
\ee
In addition, $H^i_\cont(\Delta, \bA_{\Inf}(R_\Infty)/(\mu, p^n))$ and $H^i_\cont(\Delta, \bA_{\Inf}(R_\Infty)/\mu)$ have no nonzero $W(\fm^\flat)$-torsion.
\eprop

\bpf
\ready{Since $A(R)/\mu$ is $p$-adically complete and has a trivial $\Delta$-action (see \Cref{A-Apr} and \S\ref{Ainf}), \Cref{Koszul} implies that $H^i_\cont(\Delta, A(R)/\mu)$ is a direct sum of copies of $A(R)/\mu$, and likewise for $H^i_\cont(\Delta, A(R)/(\mu, p^n))$. Consequently, since, by \eqref{Ainf-clean-1}, the rings $A(R)/(\mu, p^n)$ and $A(R)/\mu$ have no nonzero $W(\fm^\flat)$-torsion, the analogues of all the claims with $A(R)$ in place of $\bA_\Inf(R_\infty)$ follow. Thus, due to \eqref{dec-dec}, we only need to establish these analogues with $N_\infty$ in place of $\bA_\Inf(R_\infty)$.

To prepare for treating $N_\infty$, we start by building on the ideas of \cite{Bha16}*{proof of Lem.~4.6} to analyze a single summand $S \ce A_\Inf \cdot X_0^{a_0}\cdots X_d^{a_d}$ that, as in \S\ref{Ainf}, contributes to $N_\infty^\square$. We set  
\be \lab{bj-def}
\qq b_j \ce a_j - a_0 \q \text{for} \q 1 \le j \le r \qq \text{and} \qq b_j \ce a_j \q \text{for} \q r + 1 \le j \le d,
\ee
and let $m \in \bZ_{> 0}$ be the minimal such that $p^m b_j \in \bZ$ for all $j$. 
\Cref{Koszul} applied with the topological generators $\delta_1, \dotsc, \delta_d$ of $\Delta$ defined in \S\ref{Rinfty} gives an $A_\Inf$-isomorphism $H^i_\cont(\Delta, S/\mu) \simeq H^i(C^\bullet)$, where $C^\bullet$ is the $(A_\Inf/\mu)$-tensor product of the $d$ complexes 
\be \lab{factor-cx}
 [A_\Inf/\mu \xra{[\eps^{b_j}] - 1} A_\Inf/\mu] \cong A_\Inf/([\eps^{b_j}] - 1) \tensor^\bL_{A_\Inf} A_\Inf/\mu.
\ee
By reordering the $b_j$, we may assume that for all $j$ we have $b_j/b_1 \in \bZ_{(p)}$, so that $b_1 \not \in \bZ$ and both $[\eps^{b_1}] - 1 \mid [\eps^{b_j}] - 1$ and $[\eps^{b_1}] - 1 \mid \mu$. Then the object \eqref{factor-cx} with $j = 1$ is given by the complex $[A_\Inf/([\eps^{b_1}] - 1) \xra{0} A_\Inf/([\eps^{b_1}] - 1)]$ and, by using the left sides of \eqref{factor-cx} for the factors with $j \neq 1$, we see that $C^\bullet$ is quasi-isomorphic to a direct sum of shifts of $A_\Inf/([\eps^{b_1}] - 1) \cong A_\Inf/\varphi^{-m}(\mu)$. Thus, for $i \in \bZ$,
\be \lab{Smu-concl}
\tst H^i_\cont(\Delta, S/\mu) \simeq \bigoplus_I A_\Inf/\varphi^{-m}(\mu) \q \text{for some set $I$,} \q \text{and hence} \q H^i_\cont(\Delta, S/\mu)[p] = 0.
\ee
By \Cref{Koszul} and \cite{SP}*{\href{http://stacks.math.columbia.edu/tag/061Z}{061Z}, \href{http://stacks.math.columbia.edu/tag/0662}{0662}}, this implies that
\be \lab{S-bc}
H^i_\cont(\Delta, S/\mu) \tensor_{A_{\Inf}} A_{\Inf}/p^n \isomto H^i_\cont(\Delta, S/(\mu, p^n)).
\ee

We now analyze $N_\infty^\square$. Since $\bA_{\Inf}(R_\infty^\square)/\mu$ is $p$-adically complete, \S\ref{Ainf} gives the $\Delta$-decomposition
\[
\bA_{\Inf}(R_\infty^\square)/\mu \cong \wh{\bigoplus}_{\substack{ (a_0, \ldots, a_d) \in (\bZ[\f{1}{p}]_{\ge 0})^{\oplus (r + 1)} \oplus (\bZ[\f{1}{p}])^{\oplus(d - r)}, \\ \text{$a_{j} = 0$ for some $0 \le j \le r$}}} \ \ A_\Inf/\mu \cdot X_0^{a_0}\cdots X_d^{a_d}
\]
in which the completion is $p$-adic. \Cref{coho-big-sum}~\ref{CBS-i} then combines with \eqref{Smu-concl} to prove that 
\[
H^i_\cont(\Delta, N_\infty^\square/\mu)[p] = 0 \q \text{for each} \q i \in \bZ.
\]
Analogously to \eqref{S-bc}, this, in turn, implies that
\be \lab{Hi-box-bc}
H^i_\cont(\Delta, N_\Infty^\square/\mu) \tensor_{A_{\Inf}} A_{\Inf}/p^n \isomto H^i_\cont(\Delta, N_\Infty^\square/(\mu, p^n)).
\ee

Finally, we analyze $N_\infty$. The identification
\[
\tst N_\Infty/(\mu, p^{n}) \cong N_\Infty^\square/(\mu, p^{n}) \tensor_{A(R^\square)} A(R)
\]
is $\Delta$-equivariant and $A(R)/(\mu, p^n)$ is $(A(R^\square)/(\mu, p^n))$-flat, so \Cref{Koszul} gives the identifications
\be \lab{Hi-N-bc}
H^i_\cont(\Delta, N_\Infty/(\mu, p^{n})) \cong  H^i_\cont(\Delta, N_\Infty^\square/(\mu, p^{n})) \tensor_{A(R^\square)} A(R)  \q \text{for} \q n \ge 1,
\ee
which are compatible as $n$ varies. Consequently, for $n > 1$, the sequences
\be \lab{cont-n-hop}
0 \ra H^i_\cont(\Delta, N_\Infty/(\mu, p^{n}))[p] \ra H^i_\cont(\Delta, N_\Infty/(\mu, p^{n})) \ra H^i_\cont(\Delta, N_\Infty/(\mu, p^{n - 1})) \ra 0
\ee
are short exact because, by \eqref{Smu-concl} and \eqref{Hi-box-bc}, so are their analogues with $N_\infty^\square$ in place of $N_\infty$. By taking the inverse limit of these sequences for varying $n$ and using \cite{SP}*{\href{http://stacks.math.columbia.edu/tag/0D6K}{0D6K}}, 
we obtain
\be \lab{coho-N-lim}
\tst H^i_\cont(\Delta, N_\Infty/\mu) \isomto \varprojlim_n \p{H^i_\cont(\Delta, N_\Infty/(\mu, p^n))},
\ee
which is the sought analogue of \eqref{coho-lim}. The $p$-torsion freeness of $H^i_\cont(\Delta, N_\Infty/\mu)$ follows from \eqref{cont-n-hop}--\eqref{coho-N-lim} and, as in \eqref{S-bc}, it implies that
\[
H^i_\cont(\Delta, N_\Infty/\mu) \tensor_{A_{\Inf}} A_{\Inf}/p^n \isomto H^i_\cont(\Delta, N_\Infty/(\mu, p^n)).
\]
It remains to show that each $H^i_\cont(\Delta, N_\Infty/(\mu, p^n))$ has no nonzero $W(\fm^\flat)$-torsion.

The surjectivity aspect of the short exact sequences \eqref{cont-n-hop} implies that the sequences
\[
0 \ra N_\Infty/(\mu, p) \xra{p^{n - 1}} N_\Infty/(\mu, p^{n}) \ra N_\Infty/(\mu, p^{n - 1}) \ra 0 
\]
stay short exact after applying $H^i_\cont(\Delta, -)$. Thus, $H^i_\cont(\Delta, N_\Infty/(\mu, p^n))$ is a successive extension of copies of $H^i_\cont(\Delta, N_\Infty/(\mu, p))$. Consequently, it has no nonzero $W(\fm^\flat)$-torsion because, by \Cref{Rflat-clean}, neither does $H^i_\cont(\Delta, N_\infty/(\mu, p))$ (note that $N_\infty/(\mu, p)$ is a direct summand of $\bA_\Inf(R_\infty)/(\mu, p) \cong R_\infty^\flat/\mu$).} 
\epf

\bthm\lab{ANB-c}
The edge map $e$ defined in {\upshape\eqref{another-e}} induces the isomorphism
\[
 L\eta_{(\mu)}(e) \colon L\eta_{(\mu)}(R\Gamma_{\cont}(\Delta, \bA_{\Inf}(R_\Infty))) \isomto L\eta_{(\mu)}(R\Gamma_\proet(\fX^\ad_C, \bA_{\Inf,\, \fX^\ad_C})).
\]
\ethm

\bpf
By the projection formula \cite{SP}*{\href{http://stacks.math.columbia.edu/tag/0944}{0944}}, 
\be \lab{proj-proj}
R\Gamma_{\cont}(\Delta, \bA_{\Inf}(R_\Infty)) \tensor^\bL_{A_\Inf} A_\Inf/\mu \cong R\Gamma_{\cont}(\Delta, \bA_{\Inf}(R_\Infty)/\mu),
\ee
so \Cref{coho-Ainf} implies that the cohomology modules of $R\Gamma_{\cont}(\Delta, \bA_{\Inf}(R_\Infty)) \tensor^\bL_{A_\Inf} A_\Inf/\mu$ have no nonzero $W(\fm^\flat)$-torsion. Thus, the claim follows from Lemmas \ref{kill-by-more} and \ref{Bha-lemma}.
\epf

\brem \lab{ANB-c-more}
Analogously to \Cref{e-iso-more}, \Cref{ANB-c} extends as follows: for any pro-(finite \'{e}tale) affinoid perfectoid $\Delta'$-cover
\[
\tst \Spa(R_\infty'[\f{1}{p}], R_\infty') \ra \Spa(R[\f{1}{p}], R) \cong \fX_C^\ad \qq \text{that refines} \qq \fX_{C,\, \infty}^\ad \ra \fX_C^\ad
\]
subject to the same conditions as in \Cref{e-iso-more}, the edge map $e'$ defined analogously to \eqref{another-e} induces the isomorphism
\[
L\eta_{(\mu)}(e') \colon  L\eta_{(\mu)}(R\Gamma_\cont(\Delta', \bA_\Inf(R'_\infty))) \isomto L\eta_{(\mu)}(R\Gamma_\proet(\fX_C^\ad, \bA_{\Inf,\, \fX^\ad_C})).
\]
Indeed, as in \Cref{e-iso-more}, by the almost purity theorem and the octahedral axiom, $[\fm^\flat]A_\Inf$ kills the cohomology modules of the cone of the map $e_0\colon R\Gamma_\cont(\Delta, \bA_\Inf(R_\infty)) \ra R\Gamma_\cont(\Delta', \bA_\Inf(R'_\infty))$ and, by \cite{BS15}*{3.4.4 and 3.4.14}, these modules are derived $p$-adically complete; thus, by \Cref{pump-zero}, even $W(\fm^\flat)$ kills them, to the effect that \Cref{Bha-lemma} applies to the map $e_0$ and proves the claim.
\erem

As a final goal of \S\ref{local-analysis}, we wish to show in \Cref{more-ANB} that even the maps $L\eta_{(\mu)}(e \wh{\tensor}^\bL_{A_\Inf} A_\cris^{(m)})$ are isomorphisms for $A_\Inf$-algebras $A_\cris^{(m)}$ reviewed in \S\ref{Acris-m} below. 
This extension of \Cref{ANB-c} will be important for relating $A\Omega_\fX$ to logarithmic crystalline cohomology in \S\ref{section-Acris}. Our analysis of $L\eta_{(\mu)}(e \wh{\tensor}^\bL_{A_\Inf} A_\cris^{(m)})$ will use the following further consequences of the proof of \Cref{coho-Ainf}.

\bpp[The decomposition of $N_\infty$] \lab{N-dec-sec}
For $m \ge 0$, let $N_m^\square$ be the $(p, \mu)$-adically completed direct sum of those summands $A_\Inf \cdot X_0^{a_0}\cdots X_d^{a_d}$ that contribute to $\bA_\Inf(R_\infty^\square)$ in \S\ref{Ainf} for which $m$ is the smallest nonnegative integer such that $p^m a_j \in \bZ$ for all $j$ (equivalently, in the notation of \eqref{bj-def}, such that $p^mb_j \in \bZ$ for all $j$). For varying $m > 0$, the $A(R^\square)$-modules $N_m^\square$ and the $A(R)$-modules $N_m \ce N_m^\square \wh{\tensor}_{A(R^\square)} A(R)$ comprise the $(p, \mu)$-adically completed direct sum decompositions
\be \lab{N-dec}
\tst N_\infty^\square \cong \wh{\bigoplus}_{m > 0}\, N_m^\square \qq \text{and} \qq N_\infty \cong \wh{\bigoplus}_{m > 0}\, N_m.
\ee
For a fixed $i$, \Cref{Koszul} and \eqref{Smu-concl}--\eqref{S-bc} imply that
\be \lab{Nm-box-dec}
\tst H^i_\cont(\Delta, N_m^\square/(\mu, p^n)) \simeq \bigoplus_{I'} A_\Inf/(\varphi^{-m}(\mu), p^n) \q \text{for some set $I'$ and every $n > 0$.}
\ee
\epp

\bcor \lab{coho-Ainf-cor}
For all $i$ and $n, m \ge 0$,
\[
H^i_\cont(\Delta, N_m/(\mu, p^n)) \q \text{is killed by $\varphi^{-m}(\mu)$} \q \text{and} \q \text{is a flat $A_\Inf/(\varphi^{-m}(\mu), p^n)$-module.}
\]
\ecor

\bpf
If $R = R^\square$, then \eqref{Nm-box-dec} suffices. In addition, by Lazard's theorem, $A(R)/(\mu, p^n)$ is a filtered direct limit of finite free $A(R^\square)/(\mu, p^n)$-modules. Thus, the general case of the claim follows by using \eqref{Hi-N-bc} and its analogue for $N_0$ and $N_0^\square$. 
\epf

We wish to supplement \Cref{coho-Ainf} with \Cref{Ninfty-coho} that analyzes the cohomology of $N_\infty$ without reducing modulo $\mu$. Its proof will use the following base change result for $L\eta$.

\blem[\cite{Bha16}*{5.16}] \lab{Bha-BC}
For a ring $A$, elements $f, g \in A$ with $g$ a nonzerodivisor, and a $K \in D(A)$, if the modules $H^i(K \tensor^\bL_A A/f)$ have no nonzero $g$-torsion, then the natural map
\[
L\eta_{(f)}(K) \tensor^\bL_A A/g  \ra L\eta_{(\ov{f})} (K \tensor^\bL_A A/g), \qq \text{where $\ov{f}$ denotes the image of $f$ in $A/g$,}
\]
is an isomorphism. \QED
\elem

\bprop \lab{Ninfty-coho}
The element $\mu$ kills every $H^i_\cont(\Delta, N_\infty)$.
\eprop

\bpf
Let $\delta_1, \dotsc, \delta_d$ be the free generators of $\Delta$ fixed in \S\ref{Rinfty}. By \Cref{Koszul}, we need to prove that 
\be \lab{Leta-N-zero}
L\eta_{(\mu)}(K_{N_\infty}(\delta_1 -1, \dotsc, \delta_d -1)) \cong 0.
\ee
The key point, with which we start, is to prove the vanishing \eqref{Leta-N-zero} modulo $\varphi(\xi)$. The isomorphism 
\[
K_{N_\infty}(\delta_1 -1, \dotsc, \delta_d -1) \tensor^\bL_{A_\Inf} A_\Inf/\mu \cong K_{N_\infty/\mu}(\delta_1 -1, \dotsc, \delta_d -1),
\]
\Cref{Koszul}, and \Cref{coho-Ainf} show that the cohomology of $K_{N_\infty}(\delta_1 -1, \dotsc, \delta_d -1) \tensor^\bL_{A_\Inf} A_\Inf/\mu$ is $p$-torsion free. Therefore, \Cref{Bha-BC} supplies the identification
\be \lab{LKN-base-change}
L\eta_{(\mu)}(K_{N_\infty}(\delta_1 -1, \dotsc, \delta_d -1)) \tensor^{\bL}_{A_\Inf} A_\Inf/\varphi(\xi) \cong L\eta_{(\zeta_p - 1)}(K_{N_\infty/\varphi(\xi)}(\delta_1 -1, \dotsc, \delta_d -1)).
\ee
The inverse Frobenius $\varphi\i$ maps $N_\infty^\square$ isomorphically onto a direct summand of $N_\infty^\square$, so it maps $N_\infty$ isomorphically onto a direct summand of $N_\infty$. Thus, $\varphi\i$ maps $N_\infty/\varphi(\xi)$ isomorphically onto a direct summand of $N_\Infty/\xi \cong M_\infty$ (see \eqref{dec-dec}). In particular, by \Cref{Koszul} and \Cref{Rinfty-clean}, $\zeta_p - 1$ kills the cohomology of $K_{N_\infty/\varphi(\xi)}(\delta_1 -1, \dotsc, \delta_d -1)$, so both sides of \eqref{LKN-base-change} are acyclic.

Since $N_\infty$ is $(p, \mu)$-adically complete, it is also $\varphi(\xi)$-adically complete (see \cite{SP}*{\href{http://stacks.math.columbia.edu/tag/090T}{090T}}). Thus, $K_{N_\infty}(\delta_1 -1, \dotsc, \delta_d -1)$ is derived $\varphi(\xi)$-adically complete, and \cite{BMS16}*{6.19} implies the same for $L\eta_{(\mu)}(K_{N_\infty}(\delta_1 -1, \dotsc, \delta_d -1))$. The established acyclicity of the left side of \eqref{LKN-base-change} therefore implies the desired vanishing \eqref{Leta-N-zero}.
\epf

\bpp[The $A_\Inf$-algebras $A_\cris^{(m)}$] \lab{Acris-m}
For $m \in \bZ_{\ge 1}$, we let $A_\cris^{(m)}$ be the $p$-adic completion of the $A_\Inf$-subalgebra $A_\cris^{0,\, (m)}$ of $A_\Inf[\f{1}{p}]$ generated by the elements $\f{\xi^s}{s!}$ with $s \le m$.
In particular, $A_\cris^{(m)} \cong A_\Inf$ for $m < p$. In contrast, if $m \ge p$, then, since $\f{\mu^{p}}{p!} \in A_\cris^{(m)}$, the $p$-adic and $(p, \mu)$-adic topologies of $A_\cris^{(m)}$ agree. By its definition, $A_\cris^{(m)}$ is $p$-torsion free; in fact, although we will not use this, \Cref{Acrism-inj} below implies that $A_\cris^{(m)}$ is even a domain. The map $\theta$ of \eqref{theta-intro} extends to~$A_\cris^{(m)}$:
\be \lab{theta-Acrism}
\theta \colon A_\cris^{(m)} \surjects \cO_C.
\ee

Due to the ``finite type nature'' of the $A_\Inf$-algebra $A_\cris^{(m)}$, more precisely, due to \cite{BMS16}*{12.8 (ii)}, the systems of ideals
\[
 ( p^n A_\cris^{(m)} )_{n \ge 1} \q \text{and} \q (\{x \in A_\cris^{(m)}\  |\  \mu x \in  p^n A_\cris^{(m)} \})_{n \ge 1} \qq \text{of $A_\cris^{(m)}$ are intertwined.}
\]
Equivalently, 
\be \lab{red-vanishes}
 \text{for every} \q n \ge 1, \q \text{the map} \q (A_\cris^{(m)}/p^{n'})[\mu] \ra A_\cris^{(m)}/p^n \q \text{vanishes for large} \q n' > n.
\ee
Therefore, by taking the inverse limit over $n$ of the sequences 
\be \lab{mu-SES-2}
0 \ra (A_\cris^{(m)}/p^n)[\mu] \ra A_\cris^{(m)}/p^n \xra{\mu} A_\cris^{(m)}/p^n \ra A_\cris^{(m)}/(\mu, p^n) \ra 0,
\ee
we conclude that
\be \lab{cris-m-mu}
A_\cris^{(m)} \q \text{is $\mu$-torsion free} \qq \text{and} \qq A_\cris^{(m)}/\mu \q \text{is $p$-adically complete.}
\ee

The Frobenius automorphism of $A_\Inf$ preserves the subring $A_\cris^{0,\,(m)} \subset A_\Inf[\f{1}{p}]$: indeed, for $m \ge p$, since $\xi = \sum_{i = 0}^{p - 1} [\eps^{i/p}]$ and $\xi^p \in pA_\cris^{0,\, (m)}$, we have $\varphi(\xi) = \sum_{i = 0}^{p - 1} [\eps^i]$ and $\varphi(\xi) \in pA_\cris^{0,\, (m)}$. Thus, the Frobenius induces a ring endomorphism
\[
\varphi \colon A_\cris^{(m)} \ra A_\cris^{(m)},
\]
which, via the map $\theta$, intertwines the absolute Frobenius of $\cO_C/p$ (compare with \eqref{theta-intro}). 
\epp

\bpp[The $A(R)$-algebras $A_\cris^{(m)}(R)$] \lab{AcrisR-m}
The ``relative version'' of $A_\cris^{(m)}$ (resp.,~a ``highly ramified cover'' of this relative version) is the $A(R)$-algebra (resp.,~$\bA_\Inf(R_\infty)$-algebra)
\[
A_{\cris}^{(m)}(R) \ce A(R)\wh{\otimes}_{A_\Inf} A_{\cris}^{(m)} \qq \text{(resp.,} \qq \bA_{\cris}^{(m)}(R_\infty) \ce \bA_{\Inf}(R_\infty)\wh{\otimes}_{A_\Inf}A_{\cris}^{(m)}),
\]
where the completion is $(p, \mu)$-adic (equivalently, $p$-adic if $m \ge p$). In the case $m < p$, due to \Cref{A-Apr} and \S\ref{Acris-m}, we have $A_\cris^{(m)}(R) \cong A(R)$ and $\bA_\cris^{(m)}(R_\infty) \cong \bA_\Inf(R_\infty)$. 

Due to the decomposition \eqref{dec-dec}, the subring $A_{\cris}^{(m)}(R) \subset \bA_{\cris}^{(m)}(R_\infty)$ is an $A_{\cris}^{(m)}(R)$-module direct summand. Explicitly, the decomposition of $\bA_\Inf(R_\infty^\square)$ described in \S\ref{Ainf} gives the decomposition
\be \lab{Acrism-dec}
\tst \bA_\cris^{(m)}(R_\infty^\square) \cong \wh{\bigoplus}_{\substack{ (a_0, \ldots, a_d) \in (\bZ[\f{1}{p}]_{\ge 0})^{\oplus (r + 1)} \oplus (\bZ[\f{1}{p}])^{\oplus(d - r)}, \\ \text{$a_{j} = 0$ for some $0 \le j \le r$}}} \ \ A_\cris^{(m)} \cdot X_0^{a_0}\cdots X_d^{a_d},
\ee
where the completion is $(p, \mu)$-adic (equivalently, $p$-adic if $m \ge p$), and $\bA_{\cris}^{(m)}(R_\infty)$ is $(p, \mu)$-adically formally \'{e}tale over $\bA_{\cris}^{(m)}(R_\infty^\square)$ (see \S\ref{Ainf}). In particular, \eqref{red-vanishes} holds with $A_\cris^{(m)}$ replaced by $\bA_{\cris}^{(m)}(R^\square_\infty)$, and hence also by $\bA_{\cris}^{(m)}(R_\infty)$. Consequently, the generalization of \eqref{cris-m-mu} holds, too:
\be \lab{Cris-m-mu}
\bA_\cris^{(m)}(R_\infty) \q \text{is $\mu$-torsion free} \qq \text{and} \qq \bA_\cris^{(m)}(R_\infty)/\mu \q \text{is $p$-adically complete.}
\ee
In addition, by \eqref{Acrism-dec} and the formal \'{e}taleness, each $\bA_\cris^{(m)}(R_\infty)$ is $p$-torsion free.
By \S\ref{Ainf} and \S\ref{Acris-m}, the rings $A_{\cris}^{(m)}(R)$ and $\bA_{\cris}^{(m)}(R_\infty)$ come equipped with compatible $A_\cris^{(m)}$-semilinear Frobenius endomorphisms that are compatible as $m$ varies. 

The group $\Delta$ acts continuously, Frobenius-equivariantly, and $A_\cris^{(m)}$-linearly on $A_{\cris}^{(m)}(R)$ and $\bA_{\cris}^{(m)}(R_\infty)$. For each $\delta \in \Delta$, the $A_\Inf$-module endomorphism $\f{\delta - 1}{\mu}$ of $A(R)$ induces an $A_\cris^{(m)}$-module endomorphism $\f{\delta - 1}{\mu}$ of $A_\cris^{(m)}(R)$ that satisfies $\delta = 1 + \mu \cdot \f{\delta - 1}{\mu}$. In particular, $\Delta$ acts trivially on $A_\cris^{(m)}(R)/\mu$. 
\epp

\bpp[The $A_\cris^{(m)}$-base change of the edge map] \lab{Acrism-edge}
Since $A_\cris^{(m)} \cong A_\Inf$ for $m < p$, for the sake of analyzing the map $e \wh{\tensor}^\bL_{A_\Inf} A_\cris^{(m)}$, where $e$ is as in \eqref{another-e}, we suppose that $m \ge p$. Then, for each $n > 0$, we have $A_\cris^{(m)}/p^n \cong A_\cris^{(m)}/(p^n, \mu^{n'})$ for every large enough $n' > 0$ (see \S\ref{Acris-m}). Consequently, since $(p^n, \mu^{n'})$ is an $\bA_\Inf(R_\infty)$-regular sequence with $\bA_\Inf(R_\infty)/(p^n, \mu^{n'})$ flat over $A_\Inf/(p^n, \mu^{n'})$ (see \Cref{A-Apr}), the projection formula \cite{SP}*{\href{http://stacks.math.columbia.edu/tag/0944}{0944}} and \Cref{Koszul} imply that 
\[
\tst R\Gamma_\cont(\Delta, \bA_\Inf(R_\infty)) \wh{\tensor}^\bL_{A_\Inf} A_\cris^{(m)} 
\cong R\Gamma_\cont(\Delta, \bA_\cris^{(m)}(R_\infty)).
\]
Consequently, the edge map $e$ defined in \eqref{another-e} gives rise to the map
\be \lab{e-Acrism}
e \wh{\tensor}^\bL_{A_\Inf} A_\cris^{(m)} \colon R\Gamma_\cont(\Delta, \bA_\cris^{(m)}(R_\infty)) \ra R\Gamma_\proet(\fX^\ad_{C}, \bA_{\Inf}) \wh{\tensor}^\bL_{A_\Inf} A_\cris^{(m)}.
\ee
Since $[\fm^\flat]$ kills each $H^i(\Cone(e))$ (see \S\ref{e-reedged}) and $[\fm^\flat]^2 = [\fm^\flat]$, by using a free $A_\Inf$-module resolution of $A_\cris^{(m)}/p^n$ and the definition \cite{SP}*{\href{http://stacks.math.columbia.edu/tag/064M}{064M}}, we see that $[\fm^\flat]$ also kills each $H^i(\Cone(e) \tensor^\bL_{A_\Inf} A_\cris^{(m)}/p^n)$. Consequently, by \cite{SP}*{\href{http://stacks.math.columbia.edu/tag/0D6K}{0D6K}}, the ideal $[\fm^\flat]A_\Inf$ kills each $H^i(\Cone(e) \wh{\tensor}^\bL_{A_\Inf} A_\cris^{(m)})$, to the effect that, by \Cref{pump-zero} (and \cite{BS15}*{3.4.4 and 3.4.14}), so does $W(\fm^\flat)$. 
In conclusion,
\be \lab{Wflat-kills}
W(\fm^\flat) \q \text{kills the cohomology modules of} \q \Cone(e \wh{\tensor}^\bL_{A_\Inf} A_\cris^{(m)}) \simeq \Cone(e) \wh{\tensor}^\bL_{A_\Inf} A_\cris^{(m)}.
\ee
\epp

By applying \Cref{Bha-lemma}, we will show in \Cref{more-ANB} that $L\eta_{(\mu)}(e \wh{\tensor}^\bL_{A_\Inf} A_\cris^{(m)})$ is an isomorphism. Thus, we need to know that the $A_\Inf$-modules $H^i_\cont(\Delta, \bA_{\cris}^{(m)}(R_\infty)/\mu)$ have no nonzero $W(\fm^\flat)$-torsion (compare with \Cref{coho-Ainf} for $\bA_\Inf(R_\infty)/\mu$). The following result is a step in that direction:

\bprop \lab{AcrisR-clean}
Each $\bA_\cris^{(m)}(R_\infty)/(\mu, p^n)$ and also $\bA_\cris^{(m)}(R_\infty)/\mu$ have no nonzero $W(\fm^\flat)$-torsion.
\eprop

\bpf
By the $p$-adic completeness of $\bA_\cris^{(m)}(R_\infty)/\mu$ (see \eqref{Cris-m-mu}), we may focus on $\bA_\cris^{(m)}(R_\infty)/(\mu, p^n)$. The argument for the latter is similar to that of \cite{BMS16}*{12.8~(iii)} and uses approximation by Noetherian rings. Namely, by the $(p, \varphi\i(\mu))$-adic completeness of $A_\Inf$, the assignment 
\be \lab{noeth-map}
T \mapsto [\eps]^{1/p} - 1 \qq \text{defines a $\bZ_p$-algebra morphism} \qq \bZ_p\llb T\rrb \ra A_\Inf.
\ee
By \cite{BMS16}*{4.31}, this makes $A_\Inf$ a faithfully flat $\bZ_p\llb T\rrb$-algebra. Thus, letting $M$ be the mod $((T + 1)^p - 1, p^n)$ reduction of the $\bZ_p\llb T\rrb$-subalgebra of $\bZ_p\llb T\rrb[\f{1}{p}]$ generated by the $\f{1}{s!}(\sum_{i = 0}^{p - 1} (T + 1)^i)^s$ with $s \le m$, we have the identification
\[
\bA_\cris^{(m)}(R_\infty)/(\mu, p^n) \cong M \tensor_{\bZ_p\llb T \rrb/((T + 1)^p - 1,\, p^n)} \bA_\Inf(R_\infty)/(\mu, p^n).
\]
The $(\bZ_p\llb T \rrb/((T + 1)^p - 1, p^n))$-flatness of $\bA_\Inf(R_\infty)/(\mu, p^n)$ ensures that the $\varphi\i(\mu)$-torsion submodule of $\bA_\cris^{(m)}(R_\infty)/(\mu, p^n)$ is the base change of the $T$-torsion submodule $M[T] \subset M$. Consequently, since $\varphi\i(\mu) \in W(\fm^\flat)$, the consideration of the $p$-adic filtration of $M[T]$ reduces us to proving that 
\[
\bF_p \tensor_{\bZ_p\llb T \rrb/((T + 1)^p - 1,\, p^n)} \bA_\Inf(R_\infty)/(\mu, p^n) \cong R_\infty^\flat/\varphi\i(\mu) \qq \text{has no nonzero $\fm^\flat$-torsion,}
\]
which follows from \Cref{Rflat-clean}.
\epf

To relate $H^i_\cont(\Delta, \bA_{\cris}^{(m)}(R_\infty)/\mu)$ to $H^i_\cont(\Delta, \bA_{\Inf}(R_\infty)/\mu)$ in \Cref{Acrism-coho}, we will use the following general result about exactness properties of $p$-adically completed tensor products. For concreteness, we state it for $A_\Inf$ and its algebra $A_\cris^{(m)}$, but the proof is not specific to these choices.

\blem \lab{Tor-completed}
For a fixed $m \ge p$, consider the following condition on an $A_\Inf$-module $L$\ucolon
\be \tag{$\star$} \lab{ML-cond}
\text{for all $j > 0$,} \q \{\Tor_j^{A_\Inf} (L, A_\cris^{(m)}/p^n)\}_{n > 0} \q  \text{is Mittag--Leffler with vanishing eventual images,}
\ee
which means that for every $j, n$, the map $\Tor_j^{A_\Inf} (L, A_\cris^{(m)}/p^{n'}) \ra \Tor_j^{A_\Inf} (L, A_\cris^{(m)}/p^n)$ vanishes for some $n' > n$. For a bounded complex  
\[
M^\bullet = \ldots \ra M^i \xra{d^i} M^{i + 1} \ra \ldots
\]
of $A_\Inf$-modules, if each $M^i$ and each $H^i(M^\bullet)$ satisfy \eqref{ML-cond}, then, for every $i$,  we have
\be \lab{TC-eq}
\tst H^i(M^\bullet\wh{\otimes}_{A_\Inf}A_\cris^{(m)}) \cong \varprojlim_n (H^i(M^\bullet \tensor_{A_\Inf} A_\cris^{(m)}/p^n)) \cong  H^i(M^\bullet)\wh{\otimes}_{A_\Inf}A_\cris^{(m)}.
\ee
\elem

\bpf
For an inverse system $\{ 0 \ra I'_n \ra I_n \ra I''_n \ra 0\}_{n > 0}$ of short exact sequences of abelian groups, $\{ I_n\}_{n > 0}$ is Mittag--Leffler with vanishing eventual images if and only if so are both $\{ I'_n \}_{n > 0}$ and $\{ I''_n \}_{n > 0}$. Therefore, the short exact sequences 
\be \lab{TC-eq-2}
 0 \ra \Ker(d^i) \ra M^i \ra \im(d^i) \ra 0 \qq \text{and} \qq 0 \ra \im(d^{i - 1}) \ra \Ker(d^i) \ra H^i(M^\bullet) \ra 0
\ee
imply, by descending induction on $i$, that each $\Ker(d^i)$ and each $\im(d^i)$ satisfy \eqref{ML-cond}. Consequently, these sequences stay short exact after applying $- \wh{\tensor}_{A_\Inf} A_\cris^{(m)}$, to the effect that the flanking terms of \eqref{TC-eq} get identified. By construction, this identification is compatible with the canonical maps to $\varprojlim_n \p{H^i(M^\bullet \tensor_{A_\Inf} A_\cris^{(m)}/p^n)}$, so it remains to establish the second identification in \eqref{TC-eq}. 

By \cite{SP}*{\href{http://stacks.math.columbia.edu/tag/0662}{0662} and \href{http://stacks.math.columbia.edu/tag/0130}{0130}},  
the spectral sequences associated to a double complex give the following spectral sequences that converge to $H^{i + j}(M^\bullet \tensor^\bL_{A_\Inf} A_\cris^{(m)}/p^n)$:
\[
\leftexp{(n)}E_2^{ij} = H^i(H^j(M^\bullet) \tensor^\bL_{A_\Inf} A_\cris^{(m)}/p^n) \qq \text{and} \qq \leftexp{(n) \prime}E_1^{ij} = H^j(M^i \tensor^\bL_{A_\Inf} A_\cris^{(m)}/p^n),
\]
where the differential on the $\leftexp{(n) \prime}E_1$-page is $H^j(d^i \tensor^\bL_{A_\Inf} A_\cris^{(m)}/p^n)$. As $n$ varies, both families of spectral sequences form  inverse systems. Moreover, by assumption, the systems $\{\leftexp{(n)}E_2^{ij} \}_{n > 0}$ with $i \neq 0$ and $\{\leftexp{(n) \prime}E_1^{ij}\}_{n > 0}$ with $j \neq 0$ are Mittag--Leffler with vanishing eventual images. This persists to the subsequent pages: namely, by the first sentence of the proof, to the systems $\{\leftexp{(n)}E_s^{ij} \}_{n > 0}$ with $i \neq 0$ and $\{\leftexp{(n) \prime}E_s^{ij}\}_{n > 0}$ with $j \neq 0$ for any $s \le \infty$. Consequently, the edge maps 
\[
H^{i}(M^\bullet) \tensor A_\cris^{(m)}/p^n \ra H^{i}(M^\bullet \tensor^\bL A_\cris^{(m)}/p^n) \q \text{and} \q H^{i}(M^\bullet \tensor^\bL A_\cris^{(m)}/p^n) \ra H^{i}(M^\bullet \tensor A_\cris^{(m)}/p^n)
\]
become isomorphisms after applying the functor $\varprojlim_n$. It remains to note that then so does their composition, which is the canonical map $H^{i}(M^\bullet) \tensor_{A_\Inf} A_\cris^{(m)}/p^n \ra H^{i}(M^\bullet \tensor_{A_\Inf} A_\cris^{(m)}/p^n)$.
\epf

To make \Cref{Tor-completed} practical to use, we now establish its condition \eqref{ML-cond} in several key cases.

\blem \lab{ML-crit}
For a fixed $m \ge p$, the condition \eqref{ML-cond} holds in any of the following cases\ucolon
\benumr
\item \lab{MLC-ii}
for any $n, n' > 0$, the sequence $(p^n, \mu^{n'})$ is $L$-regular and $L/(p^n, \mu^{n'})$ is $A_\Inf/(p^n, \mu^{n'})$-flat\uscolon

\item \lab{MLC-iii}
the module $L$ has no nonzero $p$-torsion and each $L/p^n$ is a filtered direct limit of direct sums of $A_\Inf$-modules of the form $A_\Inf/(\varphi^{-s}(\mu), p^n)$ for variable $s \ge 0$.

\eenum
Thus,  \eqref{ML-cond} holds for $\bA_\Inf(R_\Infty)$ and $\bA_\Inf(R_\Infty)/\mu$, and for each $H^i_\cont(\Delta, N_\infty)$ and $H^i_\cont(\Delta, \bA_\Inf(R_\Infty)/\mu)$.
\elem

\bpf

If \ref{MLC-ii} holds, then, by the regular sequence aspect, $L \tensor^\bL_{A_\Inf} A_\Inf/(p^n, \mu^{n'}) \cong L/(p^n, \mu^{n'})$, so, by the flatness aspect, $L \tensor^\bL_{A_\Inf} A_\cris^{(m)}/p^n$ is concentrated in degree $0$. Thus, in the case \ref{MLC-ii}, the inverse systems in \eqref{ML-cond} vanish termwise.

If \ref{MLC-iii} holds, then each $L \tensor^\bL_{A_\Inf} A_\Inf/p^n$ is concentrated in degree $0$, so
\be \lab{EZ-temp-system}
\{ \Tor_j^{A_\Inf}(L, A_\cris^{(m)}/p^{n}) \}_{n > 0} \cong  \{ \Tor_j^{A_\Inf/p^{n}}(L/p^n, A_\cris^{(m)}/p^{n}) \}_{n > 0}
\ee
for every $j \ge 0$. In addition, since $\varphi^{-s}(\mu) \mid \mu$ for $s \ge 0$ and each $A_\Inf/p^n$ has no nonzero $\mu$-torsion, the assumption on $L/p^n$ in \ref{MLC-iii} ensures that the right side of \eqref{EZ-temp-system} vanishes termwise for $j > 1$. In contrast, for $j = 1$ and every $n > 0$, there is an $n' > n$ such that the transition map between positions $n'$ and $n$ in the right side system of \eqref{EZ-temp-system} vanishes: this follows from the identification
\[
\Tor_1^{A_\Inf/p^{n'}}(A_\Inf/(\varphi^{-s}(\mu), p^{n'}), A_\cris^{(m)}/p^{n'}) \cong (A_\cris^{(m)}/p^{n'})[\varphi^{-s}(\mu)]
\]
and \eqref{red-vanishes}. Consequently, \ref{MLC-iii} implies \eqref{ML-cond}, as claimed.

By \Cref{A-Apr}, \ref{MLC-ii} holds for $\bA_\Inf(R_\Infty)$ and then, by Lazard's theorem, \ref{MLC-iii} holds for $\bA_\Inf(R_\Infty)/\mu$. Likewise, \Cref{coho-Ainf}, \Cref{coho-Ainf-cor}, and Lazard's theorem imply that \ref{MLC-iii} holds for each $H^i_\cont(\Delta, \bA_\Inf(R_\Infty)/\mu)$. By \Cref{Koszul}, $H^i_\cont(\Delta, N_\infty)$ vanishes for large $i$ and, by \Cref{Ninfty-coho}, we have the short exact sequences
\[
 0 \ra H^i_\cont(\Delta, N_\infty)\ra H^i_\cont(\Delta, N_\infty/\mu)\ra H^{i+1}_\cont(\Delta, N_\infty)\ra 0.
\]
Therefore, due to the first sentence of the proof of \Cref{Tor-completed}, descending induction on $i$ shows  that \eqref{ML-cond} for $H^i_\cont(\Delta, \bA_\Inf(R_\Infty)/\mu)$ implies  \eqref{ML-cond} for $H^i_\cont(\Delta, N_\infty)$. 
\epf

Thanks to \Cref{ML-crit}, we may draw the following concrete consequences from \Cref{Tor-completed}.

\bprop \lab{N-cris-coho}
For every $m \ge p$ and $i \in \bZ$, we have the identifications
\be \lab{NCC-eq}
\tst H^i_\cont(\Delta, N_\infty \wh{\tensor}_{A_\Inf} A_\cris^{(m)}) \cong \varprojlim_n (H^i_\cont(\Delta, N_\infty \tensor_{A_\Inf} A_\cris^{(m)}/p^n)) \cong  H^i_\cont(\Delta, N_\infty) \wh{\tensor}_{A_\Inf} A_\cris^{(m)}.
\ee
In particular, $\mu$ kills every $H^i_\cont(\Delta, N_\infty \wh{\tensor}_{A_\Inf} A_\cris^{(m)})$.
\eprop

\bpf
By \Cref{Koszul}, the Koszul complex $M^\bullet$ of $N_\infty$ with respect to $\delta_1, \dotsc, \delta_d$ satisfies
\[
H^i(M^\bullet) \cong H^i_\cont(\Delta, N_\infty)  \qq \text{and} \qq H^i(M^\bullet \wh{\tensor} A_\cris^{(m)}) \cong  H^i_\cont(\Delta, N_\infty \wh{\tensor}_{A_\Inf} A_\cris^{(m)}),
\]
as well as $H^i(M^\bullet \tensor A_\cris^{(m)}/p^n) \cong H^i_\cont(\Delta, N_\infty \tensor A_\cris^{(m)}/p^n)$ for every $n > 0$. Moreover, by \Cref{ML-crit}, each $M^i$ and each $H^i(M^\bullet)$ satisfy \eqref{ML-cond}. Thus, \eqref{NCC-eq} is a special case of \eqref{TC-eq}. By \Cref{Ninfty-coho}, $\mu$ kills every $H^i_\cont(\Delta, N_\infty)$, so, by \eqref{NCC-eq}, it also kills every $H^i_\cont(\Delta, N_\infty \wh{\tensor}_{A_\Inf} A_\cris^{(m)})$.
\epf

\bprop \lab{Acrism-coho}
For every $m \ge p$ and $i \in \bZ$, we have the identifications
\[
\tst H^i_\cont(\Delta, \bA_\cris^{(m)}(R_\infty)/\mu) \cong \varprojlim_n (H^i_\cont(\Delta, \bA_\cris^{(m)}(R_\infty)/(\mu, p^n))) \cong H^i_\cont(\Delta, \bA_\Inf(R_\infty)/\mu) \wh{\tensor}_{A_\Inf} A_\cris^{(m)}.
\]
Moreover, the $A_\Inf$-module $H^i_\cont(\Delta, \bA_\cris^{(m)}(R_\infty)/\mu)$ has no nonzero $W(\fm^\flat)$-torsion.
\eprop

\bpf
Similarly to the proof of \Cref{N-cris-coho}, \Cref{Tor-completed} applies to the Koszul complex of $\bA_\Inf(R_\infty)/\mu$ and, due to \eqref{Cris-m-mu}, gives the identifications. Thus, it suffices to show that each 
\[
H^i_\cont(\Delta, \bA_\Inf(R_\infty)/\mu) \tensor_{A_\Inf} A_\cris^{(m)}/p^n \overset{\eqref{ANB-red}}{\cong} H^i_\cont(\Delta, \bA_\Inf(R_\infty)/(\mu, p^n)) \tensor_{A_\Inf/p^n} A_\cris^{(m)}/p^n
\]
has no nonzero $W(\fm^\flat)$-torsion. Since $\Delta$ acts trivially on $A(R)/(\mu, p^n)$, \Cref{Koszul} and \Cref{AcrisR-clean} imply that each $H^i_\cont(\Delta, A(R)/(\mu, p^n)) \tensor_{A_\Inf/p^n} A_\cris^{(m)}/p^n$ has no nonzero $W(\fm^\flat)$-torsion. Consequently, due to the decomposition \eqref{N-dec}, it suffices to show that for $j > 0$, the module
\[
H^i_\cont(\Delta, N_j/(\mu, p^n)) \tensor_{A_\Inf/p^n} A_\cris^{(m)}/p^n \overset{\ref{coho-Ainf-cor}}{\cong} H^i_\cont(\Delta, N_j/(\mu, p^n)) \tensor_{A_\Inf/(\varphi^{-j}(\mu),\, p^n)} A_\cris^{(m)}/(\varphi^{-j}(\mu), p^n)
\]
 has no nonzero $W(\fm^\flat)$-torsion. For this, similarly to the proof of \Cref{AcrisR-clean}, we will approximate by Noetherian rings. More precisely, similarly to \eqref{noeth-map}, the assignment 
\[
T \mapsto [\eps]^{1/p^j} - 1 \qq \text{defines a $\bZ_p$-algebra morphism} \qq \bZ_p\llb T \rrb \ra A_\Inf,
\]
for which $A_\Inf$ is $\bZ\llb T\rrb$-flat. In terms of this morphism, the $A_\Inf$-algebra $A_\cris^{(m)}/(\varphi^{-j}(\mu), p^n)$ is the $A_\Inf/(\varphi^{-j}(\mu), p^n)$-base change of the mod $(T, p^n)$ reduction $M$ of the $\bZ_p\llb T \rrb$-subalgebra of $\bZ_p\llb T \rrb[\f{1}{p}]$ generated by the elements $\f{1}{s!}(\sum_{i = 0}^{p - 1} (T + 1)^{p^{j - 1} \cdot i})^s$ with $s \le m$. Consequently, we need to show that
\[
H^i_\cont(\Delta, N_j/(\mu, p^n)) \tensor_{\bZ_p\llb T\rrb/(T,\, p^n)} M
\]
has no nonzero $W(\fm^\flat)$-torsion. By \Cref{coho-Ainf-cor}, the module $H^i_\cont(\Delta, N_j/(\mu, p^n))$ is $\bZ_p\llb T\rrb/(T,\, p^n)$-flat. Thus, by $p$-adically filtering $M$, we reduce to showing that $H^i_\cont(\Delta, N_j/(\mu, p^n))/p$ has no nonzero $W(\fm^\flat)$-torsion. This, in turn, follows from \Cref{coho-Ainf} and \Cref{Rflat-clean}.
\epf

With \Cref{Acrism-coho} in hand, we are ready for the promised claim about $L\eta_{(\mu)}(e \wh{\tensor}^\bL_{A_\Inf} A_\cris^{(m)})$:

\bthm \lab{more-ANB}
For each $m \ge p$, the map $e \wh{\tensor}^\bL_{A_\Inf} A_\cris^{(m)}$ from \eqref{e-Acrism} induces the isomorphism
\[
\tst L\eta_{(\mu)}(e \wh{\tensor}^\bL_{A_\Inf} A_\cris^{(m)}) \colon L\eta_{(\mu)}(R\Gamma_\cont(\Delta, \bA_{\cris}^{(m)}(R_\infty))) \isomto L\eta_{(\mu)}(R\Gamma_\proet(\fX_C^\ad, \bA_{\Inf}) \wh{\tensor}^\bL_{A_\Inf} A_\cris^{(m)}).
\]
\ethm

\bpf
By \eqref{Wflat-kills}, the ideal $W(\fm^\flat) \subset A_\Inf$ kills the cohomology of $\Cone(e \wh{\tensor}^\bL_{A_\Inf} A_\cris^{(m)})$. By \Cref{Acrism-coho} (and the projection formula \cite{SP}*{\href{http://stacks.math.columbia.edu/tag/0944}{0944}} with \eqref{Cris-m-mu}), the cohomology modules of
\[
R\Gamma_\cont(\Delta, \bA_{\cris}^{(m)}(R_\infty)) \tensor^\bL_{A_\Inf} A_\Inf/\mu
\]
 have no nonzero $W(\fm^\flat)$-torsion. Thus, \Cref{Bha-lemma} applies and gives the desired conclusion.
\epf

\brem \lab{more-ANB-more}
Analogously to \Cref{ANB-c-more}, we may extend \Cref{more-ANB} to any affinoid perfectoid $\Delta'$-cover that refines $\fX_{C,\, \infty}^\ad \ra \fX_C^\ad$ and is subject to the same conditions as in \Cref{e-iso-more}: more precisely, with the notation used there, we have
\[
\tst L\eta_{(\mu)}(e' \wh{\tensor}^\bL_{A_\Inf} A_\cris^{(m)}) \colon L\eta_{(\mu)}(R\Gamma_\cont(\Delta', \bA_{\cris}^{(m)}(R'_\infty))) \isomto L\eta_{(\mu)}(R\Gamma_\proet(\fX_C^\ad, \bA_{\Inf}) \wh{\tensor}^\bL_{A_\Inf} A_\cris^{(m)}),
\]
where $\bA_{\cris}^{(m)}(R'_\infty) \ce \bA_{\Inf}(R'_\infty) \wh{\tensor}_{A_\Inf} A_{\cris}^{(m)}$. Indeed, as there (see also \S\ref{Acrism-edge}), the ideal $W(\fm^\flat)$ kills the cohomology of the cone of the map $R\Gamma_\cont(\Delta, \bA_{\cris}^{(m)}(R_\infty)) \ra R\Gamma_\cont(\Delta', \bA_{\cris}^{(m)}(R'_\infty))$, so \Cref{Bha-lemma} applies to this map and gives the claim.
\erem

}

%% file: AOmega.tex

\section{The de Rham specialization of $A\Omega_\fX$} \lab{AOmegaX}

\ready{
With the local analysis of \S\ref{local-analysis} at our disposal, we turn to relating $A\Omega_\fX$ to the logarithmic de Rham complex of $\fX$ in \Cref{dR-spec}. The key steps for this are the identification and the analysis of the Hodge--Tate specialization of $A\Omega_\fX$ in \Cref{HT-iso,Hi-OmegaX}. These steps were also used in the smooth case in \cite{BMS16}*{\S8 and \S9} but, due to the difficulties mentioned in the beginning of \S\ref{local-analysis}, we carry them out differently. Namely, we rely on the analysis of group cohomology presented in \S\ref{local-analysis} and, in the identification step, we use \Cref{Bha-BC} (which comes from \cite{Bha16}). Nevertheless, similarly to \cite{BMS16}*{\S9.2}, we will take advantage of the following formalism of presheaves.

\bpp[The presheaf version $A\Omega_\fX^\psh$] \lab{psh-ver}
In addition to the \'{e}tale site $\fX_\et$, we consider the site $\fX_\et^\psh$ whose objects are those connected 
affine opens of $\fX_\et$ that have an \'{e}tale coordinate map \eqref{sst-coord} and coverings are the isomorphisms. Thus, the topology of $\fX_\et^\psh$ is the coarsest possible and any presheaf is already a sheaf. Since the objects of $\fX_\et^\psh$ form a basis of $\fX_\et$, there is a morphism of~topoi
\[
(\phi\i, \phi_*) \colon \fX_\et \ra \fX_\et^\psh
\]
for which $\phi_*$ is given by restricting sheaves on $\fX_\et$ to $\fX_\et^\psh$ and $\phi\i$ is given by sheafifying. In particular, since any sheaf is the sheafification of its associated presheaf, $\phi\i \circ \phi_* \cong \id$. We let 
\[
\nu^\psh \ce \phi \circ \nu \colon (\fX^\ad_C)_\proet \ra \fX^\psh_\et
\]
be the indicated composition of morphisms of topoi (with $\nu$ defined in \eqref{nu-def}) and set
\be\lab{nu-psh-intro}
 A\Omega_\fX^\psh \ce L\eta_{(\mu)}(R \nu^\psh_*(\bA_{\Inf,\, \fX^\ad_C})) \in D^{\ge 0}(\fX_\et^\psh, A_\Inf).
\ee
Since $L\eta$ commutes with pullback under flat morphisms of ringed topoi (see \cite{BMS16}*{6.14}),
\be \lab{phi-AO-pullb}
\phi\i(A \Omega_\fX^\psh) \cong A\Omega_{\fX}.
\ee
Moreover, $A \Omega_\fX^\psh$ may be described explicitly: for every object $\fU$ of $\fX_\et^\psh$, we have
\be \lab{psh-ver-sections}
R \Gamma(\fU, A\Omega_\fX^\psh) \cong L\eta_{(\mu)}(R\Gamma((\fU^\ad_C)_\proet, \bA_{\Inf,\, \fU^\ad_C})).
\ee
In particular, since, by \cite{BMS16}*{6.19}, the functor $L\eta$ preserves derived completeness when used in the context of a \emph{replete} topos (such as that of sets), we see from \eqref{psh-ver-sections} that $A\Omega_\fX^\psh$ is derived $\xi$-adically (and also $\varphi(\xi)$-adically) complete (compare with \Cref{der-comp} below).
\epp

Armed with the formalism of \S\ref{psh-ver}, we now identify the Hodge--Tate specialization of $A\Omega_\fX$.

\bthm \lab{HT-iso}
We have the identification
\be \lab{HT-iso-eq}
\tst A\Omega_\fX \tensor^\bL_{A_\Inf,\, \theta \circ \varphi\i} \cO_C \isomto L\eta_{(\zeta_p - 1)}(R\nu_*(\wh{\cO}_{\fX_C^\ad}^+)),
\ee
where in the target $L\eta$ is with respect to the ideal sheaf $(\zeta_p - 1)\cO_{\fX,\, \et} \subset \cO_{\fX,\, \et}$. If the coordinate morphisms \eqref{sst-coord} exist Zariski locally on $\fX$, then \eqref{HT-iso-eq} also holds for $A\Omega_{\fX_\Zar}$ \up{defined in \eqref{AO-Zar-def}}.

\ethm

\bpf
The kernel of $\theta_{\fX_C^\ad} \circ \varphi\i \colon \bA_{\Inf,\, \fX_C^\ad} \surjects \wh{\cO}_{\fX_C^\ad}^+$ is generated by the nonzero divisor $\varphi(\xi)$ (see \S\ref{AOX-def}), so the projection formula \cite{SP}*{\href{http://stacks.math.columbia.edu/tag/0944}{0944}} provides the identification
\[
R\nu_* (\bA_{\Inf,\, \fX_C^\ad})  \tensor^{\bL}_{A_\Inf,\, \theta \circ \varphi\i} \cO_C \cong R\nu_* (\wh{\cO}_{\fX_C^\ad}^+).
\]
Since $(\theta \circ \varphi\i)(\mu) = \zeta_p - 1$, this induces the map \eqref{HT-iso-eq} and, likewise, also  its presheaf version
\be \lab{HT-iso-psh}
\tst A\Omega_\fX^\psh \tensor^\bL_{A_\Inf,\, \theta \circ \varphi\i} \cO_C \ra L\eta_{(\zeta_p - 1)}(R\phi_*(R\nu_*(\wh{\cO}_{\fX_C^\ad}^+))).
\ee
Due to \eqref{phi-AO-pullb}, $\phi\i$ brings \eqref{HT-iso-psh} to \eqref{HT-iso-eq}, so we seek to show that \eqref{HT-iso-psh} is an isomorphism. 

For every object $\fU \cong \Spf(R)$ of $\fX_\et^\psh$ equipped with an \'{e}tale morphism as in \eqref{sst-coord}, the discussion and the notation of \S\ref{local-analysis} apply. In particular, \Cref{coho-Ainf} and \eqref{proj-proj} ensure that the cohomology of $R\Gamma_{\cont}(\Delta, \bA_{\Inf}(R_\Infty)) \tensor^{\bL}_{A_\Inf} A_\Inf/\mu$ is $p$-torsion free. Thus, since $\varphi(\xi) \equiv p \bmod (\mu)$ (see \S\ref{Ainf-not}), \Cref{Bha-BC}  implies that
\[
 L\eta_{(\mu)}(R\Gamma_{\cont}(\Delta, \bA_{\Inf}(R_\Infty))) \tensor^\bL_{A_\Inf,\, \theta \circ \varphi\i} \cO_C \isomto  L\eta_{(\zeta_p - 1)}(R\Gamma_{\cont}(\Delta, R_\Infty)). 
\]
Since the edge maps \eqref{e-ebox} and \eqref{another-e} are compatible, \Cref{e-iso,ANB-c} then imply that
\[
 L\eta_{(\mu)}(R \Gamma((\fU^\ad_C)_\proet, \bA_{\Inf})) \tensor^\bL_{A_\Inf,\, \theta \circ \varphi\i} \cO_C \isomto L\eta_{(\zeta_p - 1)}(R\Gamma((\fU^\ad_C)_\proet, \wh{\cO}^+)).
\]
Consequently, \eqref{HT-iso-psh} is an isomorphism on every $\fU$, as desired.
\epf

\bpp[The object $\wt{\Omega}_\fX$] \lab{OmegaX-def}
To proceed further, we need to analyze the right side of \eqref{HT-iso-eq}, namely,
\be \lab{wt-Omega-def}
\wt{\Omega}_{\fX} \ce L\eta_{(\zeta_p - 1)}(R\nu_*(\wh{\cO}_{\fX_C^\ad}^+)) \in D^{\ge 0}(\cO_{\fX,\, \et}),
\ee
where, as in \Cref{HT-iso}, the  functor $L\eta$ is formed with respect to the ideal sheaf $(\zeta_p - 1)\cO_{\fX,\, \et}$.
\epp

\bprop \lab{Hi-loc-free}
For $i \ge 0$, the $\cO_{\fX,\, \et}$-module $H^i(\wt{\Omega}_\fX)$ is locally free of rank ${ \dim_x(\fX_k) \choose i }$ at a variable closed point $x$ of $\fX_k$ \up{in particular, each $H^i(\wt{\Omega}_\fX)/p^n$ is a quasi-coherent $\cO_{\fX,\, \et}/p^n$-module}. Moreover,
\be \lab{H0-triv-eq}
\nu^\sharp\colon \cO_{\fX,\, \et} \isomto \nu_* (\wh{\cO}_{\fX_C^\ad}^+), \qq\text{so that} \qq 
H^0(\wt{\Omega}_\fX) \cong \cO_{\fX,\, \et}.
\ee
\eprop

\bpf
The claims are \'{e}tale local (see \cite{SP}*{\href{http://stacks.math.columbia.edu/tag/058S}{058S}}), so we assume  that $\fX = \Spf( R)$, that $\fX$ is connected, and that there is an \'{e}tale $\Spf(\cO_C)$-morphism as in \eqref{sst-coord}:
\be \lab{et-coord-2}
\fX = \Spf( R) \ra \Spf( R^{\square}) \equalscolon \fX^\square \q \text{with} \q R^\square \ce \cO_C\{t_0, \dotsc, t_r, t_{r + 1}^{\pm 1}, \dotsc,  t_d^{\pm 1}\}/(t_0\cdots t_r - p^q),
\ee
so that the discussion and the notation of \S\ref{local-analysis} apply. In particular, since $R$ is $R^\square$-flat (see \S\ref{fir-setup}) and $\Delta$ acts trivially on $R^\square$ and $R$, \Cref{Koszul} and \Cref{Rinfty-clean} imply that
\be \lab{HLF-eq}
\tst R^{\oplus {d \choose i}} \cong H^i_\cont(\Delta, R^\square) \tensor_{R^\square} R \cong \f{H^i_\cont(\Delta,\, R_\infty^\square)}{H^i_\cont(\Delta,\, R_\infty^\square)[\zeta_p - 1]} \tensor_{R^\square} R \isomto \f{H^i_\cont(\Delta,\, R_\infty)}{H^i_\cont(\Delta,\, R_\infty)[\zeta_p - 1]}.
\ee
Thus,  since the edge maps $e$ of \eqref{e-ebox} are compatible for $R$ and $R^\square$, \Cref{e-iso} shows that
\be \lab{prev-iso}
\tst \q \f{H^i((\fX^\square)_C^\ad,\, \wh{\cO}^+)}{H^i((\fX^\square)_C^\ad,\, \wh{\cO}^+)[\zeta_p - 1]} \tensor_{R^\square} R \isomto \f{H^i(\fX_C^\ad,\, \wh{\cO}^+)}{H^i(\fX_C^\ad,\, \wh{\cO}^+)[\zeta_p - 1]}
\ee
is an isomorphism of free $R$-modules of rank ${d \choose i}$. Consequently, 
\be \lab{fol-iso}
\tst \f{H^i((\fX^\square)_C^\ad,\, \wh{\cO}^+)}{H^i((\fX^\square)_C^\ad,\, \wh{\cO}^+)[\zeta_p - 1]} \tensor_{R^\square} \cO_{\Spf (R),\, \et} \isomto \f{R^i \nu_*(\wh{\cO}^+)}{\p{R^i \nu_*(\wh{\cO}^+)}[\zeta_p - 1]} \cong H^i(\wt{\Omega}_\fX),
\ee
to the effect that $H^i(\wt{\Omega}_\fX)$ is free of rank ${d \choose i}$, as desired. For \eqref{H0-triv-eq}, by \S\ref{alm-pur}, we need to show that $R \isomto (R_\infty)^\Delta$. This map is an inclusion of a direct summand whose complementary summand $M_\infty^\Delta$ is both $p$-torsion free and, by \Cref{Rinfty-clean}, killed by $\zeta_p - 1$, so the claim follows.
\epf

\brem \lab{desheafify}
The proof of \Cref{Hi-loc-free}, specifically, \eqref{prev-iso} and \eqref{fol-iso}, shows that if $\fX$ is affine, connected, and admits a coordinate map as in \eqref{sst-coord}, then the presheaf assigning $\f{H^i(\fX'^\ad_C,\, \wh{\cO}^+)}{H^i(\fX'^\ad_C,\, \wh{\cO}^+)[\zeta_p - 1]}$ to a variable $\fX$-\'{e}tale affine $\fX'$ is already a sheaf. In particular, if the coordinate maps \eqref{sst-coord} exist Zariski locally on $\fX$ (for instance, if $\fX$ is $\cO_C$-smooth or arises as in \eqref{fX-sst} from a strictly semistable $\cX$), then the sheaves $H^i(\wt{\Omega}_{\fX})$ may be computed using the Zariski topology: more precisely, then the object $\wt{\Omega}_{\fX_\Zar}$ defined by the formula \eqref{wt-Omega-def} using the Zariski topology of $\fX$ satisfies
\be \lab{Omega-Zar-et}
H^i(\wt{\Omega}_{\fX_\Zar}) \isomto (H^i(\wt{\Omega}_{\fX}))|_{\fX_\Zar} \qq \text{for every $i$.}
\ee
\erem

\bcor \lab{der-comp}
The object $A\Omega_\fX$ is derived $\xi$-adically complete and 
\be \lab{der-comp-eq}
A\Omega_{\fX}^\psh \isomto R\phi_* (A\Omega_\fX) \overset{\eqref{phi-AO-pullb}}{\cong} R\phi_*(\phi\i(A\Omega_\fX^\psh)).
\ee
\ecor

\bpf 
For the derived $\xi$-adic completeness, since $\phi\i \circ R\phi_* \cong \id$, it suffices to show that the map
\[
\tst A\Omega_\fX \ra R\lim_n (A\Omega_\fX \tensor^\bL_{A_\Inf} A_\Inf/\xi^n)
\]
becomes an isomorphism after applying $R\phi_*$. Thus, since $A\Omega_{\fX}^\psh$ is derived $\xi$-adically complete (see \S\ref{psh-ver}), it suffices to establish the adjunction isomorphism \eqref{der-comp-eq}.  For this, by the definition of $\fX_\et^\psh$ given in \S\ref{psh-ver}, we may assume that $\fX$ is affine, connected, and admits an \'{e}tale morphism \eqref{sst-coord}. In addition, since $A\Omega_\fX^\psh$ is derived $\varphi(\xi)$-adically complete, the $\fX_\et^\psh$-analogue of \cite{BMS16}*{9.15} reduces us to proving that
\[
A\Omega_{\fX}^\psh \tensor^{\bL}_{A_\Inf} A_\Inf/(\varphi(\xi)^n) \isomto  R\phi_*(\phi\i(A\Omega_{\fX}^\psh \tensor^{\bL}_{A_\Inf} A_\Inf/(\varphi(\xi)^n))).
\]
By the five lemma, we may assume that $n = 1$ and, by the proof of \Cref{HT-iso}, 
\[
A\Omega_{\fX}^\psh \tensor^{\bL}_{A_\Inf} A_\Inf/(\varphi(\xi)) \cong L\eta_{(\zeta_p - 1)}(R\phi_*(R\nu_*(\wh{\cO}_{\fX_C^\ad}^+))) \equalscolon \Omega_\fX^\psh.
\]
It remains to recall from \Cref{desheafify} that the cohomology presheaves of $\Omega_\fX^\psh$ are in fact sheaves.
\epf

Our next task is to identify the vector bundles $H^i(\wt{\Omega}_\fX)$ with the twists of the bundles given by logarithmic differentials (see \Cref{Hi-OmegaX}). For this, in \Cref{cup-iso}, we first express $H^i(\wt{\Omega}_\fX)$ as $\bigwedge^i H^1(\wt{\Omega}_\fX)$, and then, in \eqref{omg-comp-map}, construct a map that relates $H^1(\wt{\Omega}_\fX)$ to K\"{a}hler differentials.

\bpp[The cup product maps]
By the same arguments as in \cite{SP}*{\href{http://stacks.math.columbia.edu/tag/068G}{068G}}, there are product maps
\[
R^j\nu_*(\wh{\cO}^+) \tensor_{\cO_{\fX,\, \et}} R^{j'}\nu_*(\wh{\cO}^+) \xra{- \cup -} H^{j + j'}(R\nu_*(\wh{\cO}^+) \tensor^{\bL}_{\cO_{\fX,\, \et}} R\nu_*(\wh{\cO}^+))
\]
that satisfy $x \cup y = (-1)^{jj'} y\cup x$ (see \cite{SP}*{\href{http://stacks.math.columbia.edu/tag/0BYI}{0BYI}}). 
By \cite{SP}*{\href{http://stacks.math.columbia.edu/tag/0B6C}{0B6C}}, 
there is a cup product map
\[
R\nu_*(\wh{\cO}^+) \tensor^{\bL}_{\cO_{\fX,\, \et}} R\nu_*(\wh{\cO}^+) \ra R\nu_*(\wh{\cO}^+).
\]
These maps  combine to give the ``cup product map'' (where the tensor product is over $\cO_{\fX,\, \et}$)
\be  \lab{cup-map}
\tst \bigotimes_{s = 1}^i R^1 \nu_*(\wh{\cO}^+) \ra R^i \nu_*(\wh{\cO}^+) \q\ \ \text{for each} \q i >0.
\ee
\epp

\bprop \lab{cup-iso}
For each $i > 0$, the map \eqref{cup-map} induces the isomorphism
\be \lab{Hi-cup-prod}
\tst \bigwedge^i\p{ \f{R^1 \nu_*(\wh{\cO}^+)}{R^1 \nu_*(\wh{\cO}^+)[\zeta_p - 1]} } \cong \bigwedge^i H^1(\wt{\Omega}_\fX) \isomto H^i(\wt{\Omega}_\fX) \cong \f{R^i \nu_*(\wh{\cO}^+)}{R^i \nu_*(\wh{\cO}^+)[\zeta_p - 1]}.
\ee
\eprop

\bpf
By \Cref{Hi-loc-free}, each $H^i(\wt{\Omega}_\fX)$ has no nontrivial $2$-torsion, so the antisymmetry of the map \eqref{cup-map} in each pair of variables indeed induces the $\cO_{\fX,\, \et}$-module map \eqref{Hi-cup-prod}. For the isomorphism claim, we may work \'{e}tale locally, so we put ourselves in the situation \eqref{et-coord-2}. The edge maps
\[
e\colon H^i_\cont(\Delta, R_\infty) \ra H^i(\fX_C^\ad, \wh{\cO}^+)
\]
of \eqref{e-ebox} are compatible with cup products: to check this, one identifies $H^i(\fX_C^\ad, \wh{\cO}^+)$ with the direct limit of the $i\th$ \v{C}ech cohomology groups of $\wh{\cO}^+$ with respect to a variable pro\'{e}tale hypercovering of $\fX_C^\ad$ (see \cite{SP}*{\href{http://stacks.math.columbia.edu/tag/01H0}{01H0}}) and uses the hypercovering construction of the cup product (see \cite{SP}*{\href{http://stacks.math.columbia.edu/tag/01FP}{01FP}}). Due to \Cref{e-iso} and \eqref{HLF-eq}, it then remains to argue that via the cup product the identification
\[
\tst H^1_\cont(\Delta, R)  \overset{\text{\ref{Koszul}}}{\cong}  R^d \qq \text{induces} \qq  H^i_\cont(\Delta, R)  \overset{\text{\ref{Koszul}}}{\cong} \bigwedge^i (R^d),\qq 
\]
which follows from \cite{BMS16}*{7.3 and 7.5}.
\epf

To relate $H^1(\wt{\Omega}_\fX)$ to K\"{a}hler differentials, we now review the needed material on cotangent complexes.

\bpp[The completed cotangent complex $\wh{\bL}_{\wh{\cO}^+/\bZ_p}$] \lab{p-cot-cx}
Affinoid perfectoids form a basis of $(\fX^\ad_C)_\proet$ (see \cite{Sch13}*{4.7}). Therefore, \cite{BMS16}*{3.14} ensures that for the sheaf of rings $\wh{\cO}_{\fX^\ad_C}^+$, the cotangent complex $\bL_{\wh{\cO}^+/\cO_C}$, whose terms are $\wh{\cO}_{\fX^\ad_C}^+$-flat and which gives an object of $D^{\le 0}(\wh{\cO}_{\fX^\ad_C}^+)$, satisfies 
\[
\bL_{\wh{\cO}^+/\cO_C} \tensor^\bL_{\bZ} \bZ/p\bZ \cong 0, \qq \text{and hence also} \qq \wh{\bL}_{\wh{\cO}^+/\cO_C} \cong 0.
\]
Consequently, the derived $p$-adic completion turns the canonical morphism
\[
\bL_{\cO_C/\bZ_p} \tensor_{\cO_C} \wh{\cO}^+_{\fX^\ad_C} \ra \bL_{\wh{\cO}^+/\bZ_p} \qq \text{into an isomorphism} \qq (\bL_{\cO_C/\bZ_p} \tensor_{\cO_C} \wh{\cO}^+_{\fX^\ad_C})\wh{\ \ } \isomto \wh{\bL}_{\wh{\cO}^+/\bZ_p}
\]
in the derived category. By \cite{GR03}*{6.5.12 (ii)}, the complex $\bL_{\cO_C/\bZ_p}$ is quasi-isomorphic to $\Omega^1_{\cO_C/\bZ_p}$ placed in degree $0$. The $p$-divisibility of $\Omega^1_{\cO_C/\bZ_p}$ then ensures that for every $n > 0$ we have
\[
\bL_{\cO_C/\bZ_p} \tensor_{\cO_C}^\bL (\wh{\cO}^+/p^n\wh{\cO}^+) \cong (\Omega^1_{\cO_C/\bZ_p}[p^n] \tensor_{\cO_C} \wh{\cO}^+)[1] \overset{\text{\cite{Sch13}*{4.2 (iii)}}}{\cong} (\Omega^1_{\cO_C/\bZ_p}[p^n] \tensor_{\cO_C} (\cO^+/p^n\cO^+))[1],
\]
where $\cO^+$ abbreviates the integral structure sheaf $\cO^+_{\fX^\ad_C}$. Moreover, by \cite{Fon82}*{Thm.~1$\pr$~(ii)},\footnote{For passage from $\Omega^1_{\ov{\bZ}_p/\bZ_p}$ of \emph{loc.~cit.}~to $\Omega^1_{\cO_C/\bZ_p}$, one may use \cite{GR03}*{6.5.20 (i)} to conclude that $\Omega^1_{\cO_C/\ov{\bZ}_p}[p] = 0$.} 
\[
\tst \cO_C\{1\} \ce \varprojlim_{n,\, y \mapsto py} \p{\Omega^1_{\cO_C/\bZ_p}[p^n]} \qq \text{is a free $\cO_C$-module of rank $1$.}
\]
In conclusion, letting $\{ 1\}$ abbreviate the $\cO_C$-tensor product with $\cO_C\{1\}$, we obtain an isomorphism
\be \lab{p-cot-cx-eq}
(\bL_{\cO_C/\bZ_p} \tensor_{\cO_C} \wh{\cO}^+_{\fX^\ad_C})\wh{\ \ } \cong (\wh{\cO}^+_{\fX^\ad_C}\{1\})[1], \ \ \text{and hence also} \ \ \wh{\bL}_{\wh{\cO}^+/\bZ_p} \cong (\wh{\cO}^+_{\fX^\ad_C}\{1\})[1], \ \, \text{in} \ D(\wh{\cO}^+_{\fX^\ad_C}).
\ee
\epp

\bpp[The relation between $\wt{\Omega}_\fX$ and K\"{a}hler differentials] \lab{Omega1-map}
The functoriality of the cotangent complex supplies the pullback morphism 
\be \lab{cot-pullb}
\wh{\bL}_{\cO_{\fX,\, \et}/\bZ_p} \ra R\nu_* (\wh{\bL}_{\wh{\cO}^+/\bZ_p})  \overset{\eqref{p-cot-cx-eq}}{\cong} (R\nu_*(\wh{\cO}^+_{\fX^\ad_C}\{1\}))[1].
\ee
To explicate its source, we note that, as in \S\ref{p-cot-cx}, the explicit description of $\bL_{\cO_C/\bZ_p}$ gives
\[
(\bL_{\cO_C/\bZ_p} \tensor_{\cO_C} \cO_{\fX,\, \et})\wh{\ \ } \cong (\cO_{\fX,\, \et} \{1\})[1], \q \text{so} \q H^0(\wh{\bL}_{\cO_{\fX,\, \et}/\bZ_p}) \cong H^0(\wh{\bL}_{\cO_{\fX,\, \et}/\cO_C}).
\]
Moreover, the short exact sequence \cite{SP}*{\href{http://stacks.math.columbia.edu/tag/0D6K}{0D6K}} leads to the identification $H^0(\wh{\bL}_{\cO_{\fX,\, \et}/\cO_C}) \cong \Omega_{\fX/\cO_C}^1$ (the $R^1\lim$ term vanishes due to the description \cite{Ill71}*{III.3.2.7}: each $\fX_{\cO_C/p^n}$ is a local complete intersection over $\cO_C/p^n$ and, as may be seen using \eqref{sst-coord}, no nonzero local section of a vector bundle on $\fX_{\cO_C/p^n}$ vanishes on $\fX^\sm_{\cO_C/p^n}$). By \cite{Ill71}*{III.3.1.2}, over $\fX^\sm$, this identification gives a quasi-isomorphism~between
 \[
 \wh{\bL}_{\cO_{\fX^\sm,\, \et}/\cO_C}  \qq \text{and} \qq \Omega^1_{\fX^\sm/\cO_C} \text{ placed in degree $0$.}
 \]
Consequently, by applying $H^0(-)$ to the map \eqref{cot-pullb} and twisting by $\cO_C\{-1\}$ we obtain the first map in the following composition of $\cO_{\fX,\, \et}$-module morphisms:
\be \lab{omg-comp-map}
\tst \Omega^1_{\fX/\cO_C}\{-1\} \ra R^1\nu_*(\wh{\cO}^+_{\fX^\ad_C}) \surjects \f{R^1\nu_*(\wh{\cO}^+)}{(R^1\nu_*(\wh{\cO}^+))[\zeta_p - 1]} \cong H^1(\wt{\Omega}_\fX).
\ee
By \cite{BMS16}*{8.15 and its proof}, the restriction of this composition to $\fX^\sm$ is an isomorphism onto $((\zeta_p - 1) \cdot H^1(\wt{\Omega}_\fX))\vert_{\fX^\sm}$. Moreover, by \Cref{Hi-loc-free}, the $\cO_{\fX,\, \et}$-module $H^1(\wt{\Omega}_\fX)$ is a vector bundle, so it has no nonzero $(\zeta_p - 1)$-torsion and $(H^1(\wt{\Omega}_\fX))/(\zeta_p - 1)$ has no nonzero local sections that vanish on $\fX_{\cO_C/(\zeta_p - 1)}^\sm$. In conclusion, we may divide the composition \eqref{omg-comp-map} by $\zeta_p - 1$ to obtain a map 
\be \lab{over-sm-id}
\Omega^1_{\fX/\cO_C}\{-1\} \ra H^1(\wt{\Omega}_\fX)  \qq \text{that is an isomorphism over $\fX^\sm$.}
\ee
\epp

\bthm \lab{Hi-OmegaX}
The restriction of the map \eqref{over-sm-id} to $\fX^\sm$ extends uniquely to an $\cO_{\fX,\, \et}$-isomorphism
\be \lab{unique-extn}
\Omega^1_{\fX/\cO_C,\, \log} \{ -1\} \cong H^1(\wt{\Omega}_\fX),
\ee
which, by {\upshape\eqref{H0-triv-eq}} and Proposition {\upshape\ref{cup-iso}}, induces an $\cO_{\fX,\, \et}$-module identification
\be \lab{Hi-OmegaX-eq}
\Omega^i_{\fX/\cO_C,\, \log} \{ -i\} \cong H^i(\wt{\Omega}_\fX) \qq \text{for every $i \ge 0$.}
\ee
\ethm

The proof of \Cref{Hi-OmegaX} will use the formal GAGA and Grothendieck existence theorems. The Noetherian cases of these theorems proved in \cite{EGAIII1}*{\S5} have been extended to suitable non-Noetherian settings by K.~Fujiwara and F.~Kato (with important inputs due to O.~Gabber). The following theorem summarizes the relevant to our aims special case of this extension.

\bthm[Fujiwara--Kato] \lab{FK-input}
For a valuation ring $V$ of height $1$, a nonzero nonunit $a \in V$ such that $V$ is $a$-adically complete, and a proper, finitely presented $V$-scheme $Y$, 
the functor
\be \lab{FK-functor}
\cF \mapsto (\cF/a^n\cF)_{n > 0}
\ee
is an equivalence from the category of finitely presented $\cO_Y$-modules $\cF$ to that of sequences $(\cF_n)_{n > 0}$ of finitely presented $\cO_{Y_{V/a^n}}$-modules $\cF_n$ equipped with isomorphisms $\cF_{n + 1}|_{Y_{V/a^n}} \simeq \cF_n$.
\ethm

\bpf
The claim is a special case of \cite{FK14}*{I.10.1.2}. In order to explain why \emph{loc.~cit.}~implies our assertion, we first reinterpret our source and target categories.

By a result of Gabber \cite{FK14}*{0.9.2.7}, the ring $V$ is ``$a$-adically topologically universally adhesive,'' so, by \cite{FK14}*{0.8.5.25~(2)}, it is also ``topologically universally coherent with respect to $(a)$.'' In particular, by \cite{FK14}*{0.8.5.24}, every finitely presented $V$-algebra is a coherent ring, and hence, by \cite{FK14}*{0.5.1.2}, the $\cO_Y$-module $\cO_Y$ is coherent (in the sense of \cite{FK14}*{0.4.1.4 (2)} or \cite{EGAI}*{0.5.3.1}). 
In particular, by \cite{FK14}*{0.4.1.8}, an $\cO_Y$-module $\cF$ is finitely presented if and only if $\cF$ is coherent, and likewise for $\cO_{Y_{V/a^n}}$-modules for $n > 0$.

By \cite{FK14}*{0.8.4.2 and 0.8.5.19 (3)}, the formal $a$-adic completion $\wh{Y}$ of $Y$ is covered by open affines whose coordinate rings are ``topologically universally adhesive'' and hence, by \cite{FK14}*{0.8.5.18}, also ``topologically universally Noetherian outside $(a)$.'' In particular, by \cite{FK14}*{I.2.1.1~(1) and I.2.1.7}, the topological ring $V$ is ``topologically universally rigid-Noetherian'' and the formal scheme $\wh{Y}$ is ``universally rigid-Noetherian.'' In addition, by \cite{FK14}*{0.8.4.5}, the formal scheme $\wh{Y}$ is locally of finite presentation over $\Spf(V)$. Thus, \cite{FK14}*{I.7.2.2} applied with $A = V$ and \cite{FK14}*{I.7.2.1} imply that $\wh{Y}$ is ``universally cohesive.'' 
Then, by \cite{FK14}*{I.7.2.4 and I.3.4.1}, the functor $(\cF_n) \mapsto \varprojlim \cF_n$ is an equivalence from the target category of \eqref{FK-functor} to the category of coherent $\cO_{\wh{Y}}$-modules. 

In conclusion, our claim is that the quasi-coherent pullback $i^*$ along the morphism $i\colon \wh{Y} \ra Y$ of locally ringed spaces induces an equivalence between the category of coherent $\cO_Y$-modules and that of coherent $\cO_{\wh{Y}}$-modules. This is a special case of \cite{FK14}*{I.10.1.2} (see also \cite{FK14}*{I.\S9.1}).
\epf

\brems
\remi \lab{loc-free-lift}
In \Cref{FK-input}, if each $\cF_n$ is locally free, then the $\cO_Y$-module $\cF$ that algebraizes the sequence $(\cF_n)_{n > 0}$ is also locally free. Indeed, it is enough to argue that the stalks of $\cF$ at the points of $Y_{V/a}$ are flat, so, since $i$ is flat by \cite{FK14}*{I.1.4.7~(2), 0.8.5.8 (2), 0.8.5.17}, it suffices to note that the $\cO_{\wh{Y}}$-module $i^*\cF \cong \varprojlim \cF_n$ is locally free because the  Nakayama lemma ensures that $\cF_{n + 1}$ is locally trivialized by any lifts of local sections that trivialize $\cF_n$.

\remi \lab{i-flat}
Remark \ref{loc-free-lift} and the proof of \Cref{FK-input} also show that $i$ is flat and that the functor $(\cF_n) \mapsto \varprojlim \cF_n$ is an equivalence to the category of finitely presented $\cO_{\wh{Y}}$-modules.
\erems

\bpp[Proof of Theorem {\upshape\ref{Hi-OmegaX}}]
As we observed in \S\ref{Omega1-map}, no nonzero local section of a vector bundle on $\fX$ vanishes on $\fX^\sm$. Thus, the desired isomorphism \eqref{unique-extn} is unique if it exists. Consequently, we may assume that $\fX = \Spf\p{\cO_C\{t_0, \dotsc, t_r, t_{r + 1}^{\pm 1}, \dotsc,  t_d^{\pm 1}\}/(t_0\cdots t_r - p^q) }$ with $r$, $d$, and $q$ as in \eqref{sst-coord}. In this case, $\fX$ is an open subscheme of the formal $p$-adic completion of some \emph{proper}, flat $\ov{W(k)}$-scheme $\cX$ that Zariski locally has \'{e}tale ``coordinate morphisms'' as in \eqref{sst-approx} with $\cO$ there replaced by $\ov{W(k)}$. 
Thus, finally, we may drop the previous assumptions and assume instead that $\fX = \wh{\cX}$ with $\cX$ as above. We equip $\cX$ with the log structure $\cO_{\cX} \cap (\cO_{\cX}[\f{1}{p}])^\times$, so that $\cX$ is log smooth over $\ov{W(k)}$ (see \S\ref{log-str-def}, especially, \Cref{log-claim}) and the map $\fX \ra \cX$ of log ringed \'{e}tale sites is strict (see \Cref{log-comp}). By \Cref{FK-input}, the map \eqref{over-sm-id} algebraizes to an $\cO_{\cX}$-module map 
\[
f\colon \Omega^1_{\cX/\cO_C}\{ - 1\} \ra \cH.
\]
By \Cref{Hi-loc-free} and Remark \ref{loc-free-lift}, the $\cO_\cX$-module $\cH$ is locally free. By \eqref{over-sm-id} and the Nakayama lemma, $f$ is surjective at every point of $\cX_k^\sm$.
\epp

\bcl \lab{adic-GG}
There is an isomorphism $\cH_C \simeq \Omega^1_{\cX_C/C}$.
\ecl

\bpf
By the adic GAGA (see \cite{Sch13}*{9.1~(i)}), it suffices to find an analogous isomorphism after pullback to $(\cX_C)^\ad \cong \fX_C^\ad$. On the one hand, such a pullback of $\cH_C$ is isomorphic to $(R^1\nu_* (\wh{\cO}_{\fX_C^\ad}^+))[\f{1}{p}]$. On the other, \cite{Sch13b}*{3.23--3.24 and their proofs} supply an isomorphism between $(R^1\nu_*( \wh{\cO}_{\fX_C^\ad}^+))[\f{1}{p}]$ and the pullback of $\Omega^1_{\cX_C/C}$ to $(\cX_C)^\ad$.
\epf

\Cref{adic-GG} ensures that $f_C$ is a generically surjective morphism between isomorphic vector bundles on $\cX_C$. Since $\cX_C$ is proper and smooth, every global section of the structure sheaf of each connected component of $\cX_C$ is constant, so $\det(f_C)$ is an isomorphism, and hence $f_C$ is also an isomorphism. In conclusion, $f|_{\cX^\sm}$ is a surjection between vector bundles of the same rank, so 
\be \lab{f-iso}
f|_{\cX^\sm}\colon \Omega^1_{\cX^\sm/\cO_C}\{-1\} \isomto \cH|_{\cX^\sm}.
\ee
Since $\cX \setminus \cX^\sm$ is of codimension $\ge 2$ in $\cX$, limit arguments and \cite{EGAIV2}*{5.10.5} 
ensure that $\cH$ is the unique vector bundle extension of $\cH|_{\cX^\sm}$ to $\cX$. The isomorphism \eqref{f-iso} then leads to an isomorphism $\Omega^1_{\cX/\cO_C,\, \log}\{ -1 \} \simeq \cH$ whose formal $p$-adic completion gives the desired \eqref{unique-extn}.
\QED

\brem \lab{Hi-id-Zar-et}
If the coordinate morphisms \eqref{sst-coord} exist Zariski locally on $\fX$, then, by \eqref{Omega-Zar-et}, the identifications of \Cref{Hi-OmegaX} hold already for the Zariski topology; more precisely, then
\[
H^i(\wt{\Omega}_{\fX_\Zar}) \cong \Omega^i_{\fX/\cO_C,\, \log} \{ -i\} \qq \text{as $\cO_{\fX_\Zar}$-modules for every $i \ge 0$.}
\]
\erem

We are ready to relate the de Rham specialization of $A\Omega_\fX$ to differential forms by combining the results above with the  argument from the proof of \cite{BMS16}*{14.1}.

\bthm \lab{dR-spec}
There is an identification
\be \lab{dR-spec-id}
A\Omega_\fX \tensor^\bL_{A_\Inf,\, \theta} \cO_C \cong \Omega^\bullet_{\fX/\cO_C,\, \log}.
\ee
If the coordinate morphisms \eqref{sst-coord} exist Zariski locally on $\fX$, then \eqref{dR-spec-id} also holds for $A\Omega_{\fX_\Zar}$.
\ethm

\bpf
Since $\varphi(\mu) = \varphi(\xi) \mu$ (see \S\ref{Ainf-not}), \cite{BMS16}*{6.11} gives the second identification in
\[
A\Omega_\fX \tensor^\bL_{A_\Inf,\, \theta} \cO_C \cong A\Omega_\fX \tensor^\bL_{A_\Inf,\, \varphi} A_\Inf  \tensor^\bL_{A_\Inf,\, \theta \circ \varphi\i} \cO_C \cong (L\eta_{(\varphi(\xi))}(A \Omega_\fX)) \tensor^\bL_{A_\Inf,\, \theta \circ \varphi\i} \cO_C.
\]
By \cite{BMS16}*{6.12}, since $A_\Inf/(\varphi(\xi)) \cong \cO_C$ via $\theta \circ \varphi\i$, the object $(L\eta_{(\varphi(\xi))}(A \Omega_\fX)) \tensor^\bL_{A_\Inf,\, \theta \circ \varphi\i} \cO_C$ is identified with the complex whose $i$-th degree term is 
\[
\tst H^i(A\Omega_\fX \tensor^\bL_{A_\Inf,\, \theta \circ \varphi\i} \cO_C) \tensor_{\cO_C} \p{\f{\Ker (\theta \circ \varphi\i)}{(\Ker (\theta \circ \varphi\i))^2}}^{\tensor i} \overset{\text{\eqref{HT-iso-eq}}}\cong H^i(\wt{\Omega}_\fX) \tensor_{\cO_C} \p{\f{\Ker(\theta \circ \varphi\i)}{(\Ker(\theta \circ \varphi\i))^2}}^{\tensor i} 
\]
and the differentials are given by Bockstein homomorphisms. 

Since $\cO_C^\flat$ is perfect, $\wh{\bL}_{A_\Inf/\bZ_p} \cong 0$. Moreover, \eqref{p-cot-cx-eq} applied with $\fX = \Spf(\cO_C)$ implies that $\wh{\bL}_{\cO_C/\bZ_p} \cong (\cO_C\{1\})[1]$. Thus, $\wh{\bL}_{\cO_C/A_\Inf} \cong (\cO_C\{1\})[1]$, where $\cO_C$ is an $A_\Inf$-algebra via $\theta \circ \varphi\i$. In particular, due to \cite{Ill71}*{III.3.2.4 (iii)}, we have $\f{\Ker(\theta \circ \varphi\i)}{(\Ker(\theta \circ \varphi\i))^2} \cong \cO_C\{1\}$. 

In conclusion, by \eqref{Hi-OmegaX-eq} and the preceding discussion, $A\Omega_\fX \tensor^\bL_{A_\Inf,\, \theta} \cO_C$ is identified with the complex whose $i$-th degree term is $\Omega_{\fX/\cO_C,\, \log}^i$ and the differentials are certain Bockstein homomorphisms. Each $\Omega_{\fX/\cO_C,\, \log}^i$ is a vector bundle, so the agreement of the Bockstein differentials with those of $\Omega^\bullet_{\fX/\cO_C,\, \log}$ may be checked over $\fX^\sm$ (compare with the argument for \eqref{over-sm-id}),  where it follows from \cite{BMS16}*{14.1~(ii)} (or \cite{Bha16}*{proof of Prop.~7.9}).

Due to \Cref{Hi-id-Zar-et}, the proof for $A\Omega_{\fX_\Zar}$ is the same.
\epf

\bcor \lab{G-dR-spec}
The de Rham specialization of $R\Gamma(\fX_\et, A\Omega_\fX)$ may be identified as follows\ucolon
\be \lab{G-dR-spec-eq}
R\Gamma(\fX_\et, A\Omega_\fX) \tensor^\bL_{A_\Inf,\, \theta} \cO_C \cong  R\Gamma_{\log\dR}(\fX/\cO_C). 
\ee
\ecor

\bpf
The claim follows from \Cref{dR-spec} and the projection formula \cite{SP}*{\href{http://stacks.math.columbia.edu/tag/0944}{0944}}.
\epf

\brem \lab{comp-sch}
In the case when $\fX \cong \wh{\cX}$ for a proper, flat  $\ov{W(k)}$-scheme $\cX$ that \'{e}tale locally has \'{e}tale coordinate morphisms \eqref{sst-approx} with $\cO$ there replaced by $\ov{W(k)}$, we have the further identification
\[
R\Gamma(\cX_\et, \Omega^\bullet_{\cX/\ov{W(k)},\, \log}) \tensor^\bL_{\ov{W(k)}} \cO_C \isomto R\Gamma(\fX_\et, \Omega^\bullet_{\fX/\cO_C,\, \log}) = R\Gamma_{\log\dR}(\fX/\cO_C),
\]
where $\cX$ is endowed with the log structure $\cO_{\cX,\,\et} \cap (\cO_{\cX,\,\et}[\f{1}{p}])^\times$ (whose pullback to $\fX$ is the log structure $\cO_{\fX,\,\et} \cap (\cO_{\fX,\,\et}[\f{1}{p}])^\times$ of $\fX$, see \Cref{log-comp}) and $\ov{W(k)}$ is endowed with the log structure associated to $\ov{W(k)} \setminus \{ 0\} \hra \ov{W(k)}$. Indeed, the pullback map between the $E_1$-spectral sequences
\[
\ba
H^j(\cX_{\cO_C}, \Omega^i_{\cX_{\cO_C}/\cO_C,\, \log}) &\Rightarrow H^{i + j}(R\Gamma(\cX_\et, \Omega^\bullet_{\cX_{\cO_C}/\cO_C,\, \log})), \\
 H^j(\fX, \Omega^i_{\fX/\cO_C,\, \log}) &\Rightarrow H^{i + j}(R\Gamma_{\log\dR}(\fX/\cO_C))
 \ea
\]
is an isomorphism because, by the Grothendieck finiteness and comparison theorems \cite{EGAIII1}*{3.2.1 and 4.1.7} (combined with limit arguments, which use \Cref{log-claim} and the fact that $\cX$ is necessarily finitely presented, see \cite{SP}*{\href{http://stacks.math.columbia.edu/tag/053E}{053E}}; alternatively, directly by \cite{FK14}*{I.9.2.1}), 
\[
H^j(\cX_{\cO_C}, \Omega^i_{\cX_{\cO_C}/\cO_C,\, \log}) \isomto  H^j(\fX, \Omega^i_{\fX/\cO_C,\, \log}) \qq \text{for all $i, j$.}
\]
\erem

\bcor \lab{perfect}
If $\fX$ is proper over $\cO_C$, then $R\Gamma(\fX_\et, A\Omega_\fX)$ is a perfect object of $D^{\ge 0}(A_\Inf)${\upshape;} in other words, then $R\Gamma(\fX_\et, A\Omega_\fX)$ is quasi-isomorphic to a bounded complex of finite free $A_\Inf$-modules.
\ecor

\bpf
By the Grothendieck finiteness theorem \cite{Ull95}*{5.3} 
and the spectral sequence as in \Cref{comp-sch}, the $\cO_C$-modules $H^j(R\Gamma_{\log\dR}(\fX/\cO_C))$ are finitely presented, and hence also perfect (see \cite{SP}*{\href{http://stacks.math.columbia.edu/tag/0ASP}{0ASP}}). Thus, by \Cref{G-dR-spec} and \cite{SP}*{\href{http://stacks.math.columbia.edu/tag/066U}{066U}}, the object $R\Gamma(\fX_\et, A\Omega_\fX) \tensor_{A_\Inf}^\bL A_\Inf/(\xi)$ of $D^{\ge 0}(\cO_C)$ is perfect. Moreover, by \Cref{der-comp}, the object $R\Gamma(\fX_\et, A\Omega_\fX)$ is derived $\xi$-adically complete. Therefore, by \cite{SP}*{\href{http://stacks.math.columbia.edu/tag/09AW}{09AW}}, it is perfect as well, as desired.
\epf

We close the section by comparing $R\Gamma(\fX_\et, A\Omega_\fX)$ to its analogue defined using the Zariski topology.

\bcor \lab{Zar-et-comp}
If the coordinate morphisms \eqref{sst-coord} exist Zariski locally on $\fX$, then $R\Gamma(\fX_\et, A\Omega_\fX)$ may be computed using the Zariski topology of $\fX${\upshape;} more precisely, then 
\be \lab{ZEC-map}
R\Gamma(\fX_\Zar, A\Omega_{\fX_\Zar}) \isomto R\Gamma(\fX_\et, A\Omega_\fX).
\ee
\ecor

\bpf
By \Cref{dR-spec} and \Cref{G-dR-spec}, the reduction of \eqref{ZEC-map} modulo $\xi$ is identified with
\[
R\Gamma(\fX_\Zar, \Omega_{\fX/\cO_C,\, \log}^{\bullet}) \isomto R\Gamma(\fX_\et, \Omega_{\fX/\cO_C,\, \log}^{\bullet}),
\]
and hence is an isomorphism as indicated. 
Thus, since, by \Cref{der-comp} (and its Zariski analogue), $R\Gamma(\fX_\Zar, A\Omega_{\fX_\Zar})$ and $R\Gamma(\fX_\et, A\Omega_\fX)$ are derived $\xi$-adic complete, \eqref{ZEC-map} is an isomorphism.
\epf

\beg \lab{Zar-et-comp-eg}
By \S\ref{fir-setup}, \Cref{Zar-et-comp} applies to any $\cO_C$-smooth $\fX$ and, more generally, to any $\fX$ that Zariski locally arises from a strictly semistable scheme defined over a discrete valuation ring.
\eeg
}

%% file: Acris.tex

\section{The absolute crystalline comparison isomorphism \nopunct} \lab{section-Acris}

\ready{
In \Cref{dR-spec}, we identified the $\cO_C$-base change along $\theta$ of the object $A\Omega_\fX$ with $\Omega_{\fX/\cO_C,\,\log}^\bullet$. The goal of this section is to similarly identify the $A_\cris$-base change of $A\Omega_\fX$ with an object that computes the logarithmic crystalline (that is, Hyodo--Kato) cohomology of $\fX_{\cO_C/p}$ over $A_\cris$ (see \Cref{abs-crys-comp}). This is more general because, on the one hand, $\theta$ factors through the map $A_\Inf \ra A_\cris$, while, on the other, $\Omega_{\fX/\cO_C,\,\log}^\bullet$ computes the log crystalline cohomology of $\fX_{\cO_C/p}$ over $\cO_C$. In fact, even the map $A_\Inf \surjects A_\Inf/\mu$ factors through $A_\Inf \ra A_\cris$, 
so the identification of the $A_\cris$-base change of $A\Omega_\fX$ will capture the entire $\mu = 0$ locus of $A_\Inf$ (in contrast, the comparison with the $p$-adic \'{e}tale cohomology captured the $\mu \neq 0$ locus, see \Cref{RG-et-id}).

In comparison to the case when $\fX$ is smooth treated in \cite{BMS16}*{\S12}, controlling the interaction of the functor $L\eta_{(\mu)}$ with the relevant base changes seems more subtle. To overcome this, we resort to the analysis of continuous group cohomology carried out in \S\ref{local-analysis}. Another major complication is the presence of log structures. Specifically, not knowing the existence of logarithmic divided power envelopes of certain nonexact logarithmic closed immersions in mixed characteristic, we are forced to devise slightly indirect arguments when analyzing the relevant divided power envelopes. For this, we rely on the results and arguments from \cite{Kat89} and \cite{Bei13a};\footnote{\revise{We are citing the post-publication arXiv version of the article, which slightly differs from the published version.}} the latter reference is especially useful for us because some log structures that we use are not coherent (only quasi-coherent).

\bpp[The ring $A_\cris$] \lab{Acris-ring}
Using the generator $\xi$ of the kernel of $\theta \colon A_\Inf \surjects \cO_C$, we let $A_\cris^0$ be the $A_\Inf$-subalgebra of $A_\Inf[\f{1}{p}]$ generated by the divided powers $\f{\xi^n}{n!}$ for $n \ge 1$. The induced map $\theta \colon A_\cris^0 \surjects \cO_C$ identifies $A_\cris^0$ with the divided power envelope of $\theta\colon A_\Inf \surjects \cO_C/p$ over $(\bZ_p, p\bZ_p)$ equipped with the unique divided powers on $p\bZ_p$, see \cite{Tsu99}*{A2.8}. Since $\theta(\mu) = 0$, we have $\mu^p \in pA_\cris^0$, so the $p$-adic topology of $A_\cris^0$ agrees with the $(p, \mu)$-adic topology. We set
\[
A_\cris \ce (A_\cris^0)\wh{\ \,}, \qq \text{where the completion is $p$-adic (equivalently, $(p, \mu)$-adic).}
\]
The induced map $\theta \colon A_\cris \surjects \cO_C$ identifies $A_\cris$ with the initial $p$-adically complete divided power thickening of $\cO_C$ over $\bZ_p$ (see \cite{Tsu99}*{A1.3 and A1.5}). By \Cref{Acrism-inj} below (or by \cite{Tsu99}*{A2.13} and \cite{Bri06}*{2.33}), the map $A_\cris^0 \ra A_\cris$ is an injection into an integral domain.

Analogously to \S\ref{Acris-m}, the ring $A_\cris$ comes equipped with the Frobenius endomorphism $\varphi$ that intertwines the absolute Frobenius endomorphism of $\cO_C/p$ via the map $\theta$. The identification
\be \lab{Acris-m-eq}
\tst A_\cris \cong (\varinjlim_m A_\cris^{(m)})\wh{\ \ }, \qq \text{which results from the evident} \qq A_\cris^0 \cong \varinjlim_m A_\cris^{0,\,(m)},
\ee
is Frobenius equivariant and compatible with the maps $\theta$.
\epp

\bpp[The log structure on $A_\cris$] \lab{Acris-log}
For each $n  > 0$, the ring $A_\cris/p^n$ is a divided power thickening of $\cO_C/p$ over $\bZ/p^n$. Therefore, by \cite{Bei13a}*{\S1.17, Lemma}, every quasi-coherent, integral log structure $\cN$ on $\cO_C/p$ for which $\cN/(\cO_C/p)^\times$ is uniquely $p$-divisible lifts uniquely to a quasi-coherent, integral log structure on $A_\cris/p^n$. Thus, letting $\cN$ be the default log structure \S\ref{log-str-def}~\ref{OC-log-str} on $\cO_C/p$, for which $\cN/(\cO_C/p)^\times \cong \bQ_{\ge 0}$, we obtain compatible, quasi-coherent, integral log structures on the rings $A_\cris/p^n$, to the effect that each $A_\cris/p^n$ becomes a log PD thickening of $\cO_C/p$. Explicitly, these log structures are the pullbacks of the log structure on $A_\cris$ associated to the prelog structure
\be \lab{AL-1}
\cO_C^\flat \setminus \{ 0 \} \ra A_\cris, \qq x \mapsto [x].
\ee
In what follows, we always equip 
\begin{itemize}
\item
each $A_\cris/p^n$, as well as $A_\cris$, with the log structure described above;

\item
each $\bZ/p^n\bZ$ with the standard divided powers on $p\bZ/p^n\bZ$ and the trivial log structure.
\end{itemize}

For every divided power thickening $\wt{Z}$ over $\bZ/p^n\bZ$ of an $\cO_C/p$-scheme $Z$, the map $Z \xra{z} \Spec(\cO_C/p)$ extends uniquely to a PD map $\wt{Z} \xra{\wt{z}} \Spec(A_\cris/p^n)$ (see the proof of \cite{Tsu99}*{A1.5}). If, in addition, $\wt{Z}$ is equipped with a quasi-coherent, integral log structure for which $z$ is enhanced to a map $z^\sharp$ of log schemes, then, by \cite{Bei13a}*{\S1.17, Exercise}, the map $z^\sharp$ extends uniquely to a PD map $\wt{z}^\sharp \colon \wt{Z} \ra \Spec(A_\cris/p^n)$ of log schemes.
\epp

\bpp [The absolute crystalline cohomology of $\fX_{\cO_C/p}$] \lab{abs-crys-log}
We let
\[
(\fX_{\cO_C/p}/\bZ_p )_{\log\cris}
\]
 be the log crystalline site of $\fX_{\cO_C/p}$ over $\bZ_p$ defined as in \cite{Bei13a}*{\S1.12}: the objects are the  \'{e}tale $\fX_{\cO_C/p}$-schemes $Z$ equipped with a divided power thickening $\wt{Z}$ over some $\bZ/p^n\bZ$ such that $\wt{Z}$ is, in turn, equipped with a quasi-coherent, integral log structure whose pullback to $Z$ is identified with the pullback of the log structure of $\fX_{\cO_C/p}$ (which is defined in \S\ref{log-str-def}~\ref{X-log-str}); the coverings are the jointly surjective \'{e}tale log PD morphisms. The universal property of $A_\cris$ reviewed in the last paragraph of \S\ref{Acris-log} gives the following identification of sites:
\[
(\fX_{\cO_C/p}/\bZ_p )_{\log\cris} \cong (\fX_{\cO_C/p}/A_\cris)_{\log\cris}, 
\]
where $(\fX_{\cO_C/p}/A_\cris)_{\log\cris}$ is the log crystalline site of $\fX_{\cO_C/p}$ over $A_\cris$ defined analogously to the site $(\fX_{\cO_C/p}/\bZ_p )_{\log\cris}$ reviewed above (simply replace $\bZ/p^n\bZ$ by $A_\cris/p^n$).  The \emph{absolute logarithmic crystalline cohomology} of $\fX_{\cO_C/p}$ is the cohomology of the structure sheaf:
\[
R\Gamma_{\log\cris} (\fX_{\cO_C/p}/A_\cris)\ce R\Gamma ((\fX_{\cO_C/p}/ A_\cris)_{\log\cris},\, \cO_{\fX_{\cO_C/p}/A_{\cris}}).
\]
We consider the morphism of topoi
\[
u \colon (\fX_{\cO_C/p}/ A_{\cris})_{\log\cris} \ra (\fX_{\cO_C/p})_{\et} \cong \fX_{\et}
\]
that ``forgets the thickenings $\wt{Z}$'' (see \cite{Bei13a}*{\S1.5}), and we use it to obtain the identification
\[
R\Gamma_{\log\cris} (\fX_{\cO_C/p}/A_\cris) \cong R\Gamma (\fX_{\et}, Ru_{*} (\cO_{\fX_{\cO_C/p}/A_\cris})).
\]
By the functoriality discussed in \cite{Bei13a}*{\S1.5, Corollary}, the absolute Frobenius of $\fX_{\cO_C/p}$ (which is the multiplication by $p$ on log structures) and the Frobenius of $A_\cris$ induce the $A_\cris$-semilinear Frobenius endomorphisms of $Ru_{*} (\cO_{\fX_{\cO_C/p}/A_\cris})$ and $R\Gamma_{\log\cris} (\fX_{\cO_C/p}/A_\cris)$.
\epp

The main goal of this section is the following identification of the $A_\cris$-base change of $A\Omega_\fX$.

\bthm \label{abs-crys-comp}
There is a Frobenius-equivariant identification
\be \lab{abs-crys-comp-eq}
A\Omega_{\fX}\wh{\otimes}^{\bL}_{A_\Inf}A_{\cris} \cong Ru_{*} (\cO_{\fX_{\cO_C/p}/A_{\cris}}), 
\ee
where the Frobenii result from those discussed in {\upshape\S\ref{Ainf-not}}, {\upshape\S\ref{AOX-def}}, {\upshape\S\ref{Acris-ring}}, and {\upshape\S\ref{abs-crys-log}} and, consistently with the notation \eqref{der-comp-def}, we have $A\Omega_{\fX}\wh{\otimes}^{\bL}_{A_\Inf}A_{\cris} = R\lim_n (A\Omega_{\fX} \otimes^{\bL}_{A_\Inf}A_{\cris}/p^n)$.
\ethm

We will first prove a version of \Cref{abs-crys-comp} in the presence of fixed semistable coordinates. We will then complete the proof by using ``all possible coordinates''  to globalize the argument. This overall strategy is similar to the one used in \cite{BMS16}*{\S12} in the smooth case.

\bpp[The local setup] \lab{local-setup-sec}
For the local argument, we assume until \S\ref{all-possible-coordinates} that $\fX = \Spf(R)$, that $\fX$ is connected, and that for some $0 \le r \le d$ and $q \in \bQ_{> 0}$ there is an \'{e}tale $\cO_C$-morphism
\be \lab{local-setup}
\fX = \Spf(R) \ra \Spf(R^\square) \qq \text{with} \qq R^\square = \cO_C\{t_0, \dotsc, t_r, t_{r + 1}^{\pm 1}, \dotsc,  t_d^{\pm 1}\}/(t_0\cdots t_r - p^q).
\ee
We use the rings $R_\infty^\square$ and $R_\infty$, the group $\Delta$, and its generators $\delta_i$ introduced in \S\ref{Rinfty}, the rings $\bA_\Inf(R_\infty^\square)$, $\bA_\Inf(R_\infty)$, $A(R^\square)$, and $A(R)$ and the modules $N_\infty^\square$ and $N_\infty$ introduced in \S\ref{Ainf}, the rings $\bA_\cris^{(m)}(R_\infty)$ and $A_\cris^{(m)}(R)$ introduced in \S\ref{AcrisR-m}, and the object $A\Omega_\fX^\psh$ introduced in \S\ref{psh-ver}.
\epp

Roughly speaking, with the coordinates above, we will access the right side of \eqref{abs-crys-comp-eq} through the logarithmic de Rham complex of an explicit log smooth lift $\Spf(A_\cris(R))$ over $\Spf(A_\cris)$ of $\fX_{\cO_C/p}$ over $\Spec(\cO_C/p)$ (see \Cref{RHS-rewrite}). This complex may be made explicit by expressing its differentials in terms of the $\Delta$-action on $A_\cris(R)$ (see \Cref{delta-exp}). In contrast, results from \S\ref{local-analysis}, namely, \Cref{ANB-c} and \eqref{Leta-N-zero}, make the left side of \eqref{abs-crys-comp-eq} explicit. Once both sides are explicit, one identifies them and establishes (the presheaf version of) the local case of \Cref{abs-crys-comp}.

However, this relatively short local proof, whose detailed version in the good reduction case is given in \cite{BMS16}*{12.5}, is ill-suited for globalizing. This is so because it appears difficult to extend the implicit exchange of the order of the functors $L\eta_{(\mu)}$ and $- \wh{\tensor}^\bL_{A_\Inf} A_\cris$ in this argument to general perfectoid covers that appear in the ``all possible coordinates'' technique. For instance, one may attempt to use the almost purity theorem and \Cref{Bha-lemma} to reduce such commutativity to the ``base case'' of $R_\infty$, but this requires understanding the $W(\fm^\flat)$-torsion in the groups $H^i_\cont(\Delta, (\bA_\Inf(R_\infty) \wh{\tensor}_{A_\Inf} A_\cris)/\mu)$ that seem difficult to access due to pathologies of the ring $A_\cris/\mu$.

Similarly to \cite{BMS16}*{\S12.2}, to overcome this difficulty we will use the rings $A_\cris^{(m)}$ reviewed in \S\ref{Acris-m} that retain better finite type properties over $A_\Inf$ than $A_\cris$. In particular, we use the work of \S\ref{local-analysis} to commute the functors $L\eta_{(\mu)}$ and $- \wh{\tensor}^\bL_{A_\Inf} A_\cris^{(m)}$ in the following proposition:

\bprop \label{L-cris-commute}
In the local setting of \uS\uref{local-setup-sec}, for every $m \ge p$, we have
\be \lab{LCC-eq}
L\eta_{(\mu)} (R\Gamma_\proet(\fX^\ad_C, \bA_{\Inf})) \wh{\tensor}^\bL_{A_\Inf} A_\cris^{(m)} \isomto L\eta_{(\mu)}(R\Gamma_\proet(\fX^\ad_C, \bA_{\Inf}) \wh{\tensor}^\bL_{A_\Inf} A_\cris^{(m)}).
\ee
\eprop

\bpf  
The map \eqref{LCC-eq} exists because its target is derived $p$-adically complete (see \cite{BMS16}*{6.19}). Moreover, by \Cref{ANB-c,more-ANB}, it suffices to prove that
\[
L\eta_{(\mu)} (R\Gamma_\cont(\Delta, \bA_{\Inf}(R_\infty))) \wh{\tensor}^\bL_{A_\Inf} A_\cris^{(m)} \isomto L\eta_{(\mu)}(R\Gamma_\cont(\Delta, \bA_\cris^{(m)}(R_\infty))).
\]
By \Cref{Ninfty-coho,N-cris-coho}, the ``nonintegral'' part $N_\infty$ does not contribute, so it suffices to show:
\be \lab{LCC-rem-3}
L\eta_{(\mu)} (R\Gamma_\cont(\Delta, A(R))) \wh{\tensor}^\bL_{A_\Inf} A_\cris^{(m)} \isomto L\eta_{(\mu)}(R\Gamma_\cont(\Delta, A_\cris^{(m)}(R))).
\ee
In turn, \eqref{LCC-rem-3} follows from the triviality of the $\Delta$-action on $A(R)/\mu$ and $ A_\cris^{(m)}(R)/\mu$ (see \S\ref{Ainf} and \S\ref{AcrisR-m}): namely, due to \Cref{Koszul} and this triviality, the left (resp.,~right) side of \eqref{LCC-rem-3} becomes
\[
\tst K_{A(R)}(\frac{\delta_1 -1}{\mu}, \dots, \frac{\delta_d -1}{\mu}) \wh{\tensor}_{A_\Inf} A_\cris^{(m)} \qq \text{(resp.,} \qq K_{A_\cris^{(m)}(R)}(\frac{\delta_1 -1}{\mu}, \dots, \frac{\delta_d -1}{\mu})),
\]
where the completed tensor product is nonderived (that is, termwise) because each $p^n$, $\mu^{n'}$ is an $A(R)$-regular sequence with $A(R)/(p^n, \mu^{n'})$ flat over $A_\Inf/(p^n, \mu^{n'})$ (see \Cref{A-Apr}); the two Koszul complexes may then be identified termwise (see \S\ref{AcrisR-m}).
\epf

Continuing to work in the local setting, we now express the (presheaf version of the) left side of \eqref{abs-crys-comp-eq} in the form that will be convenient for the ``all possible coordinates'' technique.

\bcor \label{LHS-rewrite}
In the local setting of \uS\uref{local-setup-sec}, there is a Frobenius-equivariant identification
\be \lab{LHS-rewrite-eq}
\tst R\Gamma(\fX_\et^\psh, A\Omega_\fX^\psh) \wh{\tensor}^\bL_{A_\Inf} A_\cris \cong
\p{\varinjlim_m \p{ \eta_{(\mu)} \p{K_{\bA_{\cris}^{(m)}(R_\infty)}(\delta_1 -1, \dotsc, \delta_d -1)}}}^{\wh{\ \ \   }}
\ee
\up{see \eqref{nu-psh-intro} for $A\Omega_\fX^\psh$} where, on the right side, the direct limit and the $p$-adic completion are~termwise.
\ecor

\bpf 
The $\Delta$-equivariant Frobenii of the rings $\bA_\cris^{(m)}(R_\infty)$ are compatible as $m$ varies (see \S\ref{AcrisR-m}), so, due to the divisibility $\mu \mid \varphi(\mu)$, they induce the Frobenius on the right side of \eqref{LHS-rewrite-eq}. \Cref{L-cris-commute} and \Cref{more-ANB} give the Frobenius-equivariant identification
\[
\tst R\Gamma(\fX_\et^\psh, A\Omega_\fX^\psh) \wh{\tensor}^\bL_{A_\Inf} A_\cris^{(m)} \cong
 \eta_{(\mu)} \p{K_{\bA_{\cris}^{(m)}(R_\infty)}(\delta_1 -1, \dotsc, \delta_d -1)},
\]
so it remains to pass to the direct limit and to form the $p$-adic completion.
\epf

We turn to the right side of \eqref{abs-crys-comp-eq} and begin by constructing a log smooth lift $A_\cris(R)$ of $R/p$.

\bpp[The ring $A_\cris(R)$] \lab{Acris-R}
\ready{The ``relative version'' of $A_\cris$ (resp.,~a ``highly ramified cover'' of this relative version) is the $A(R)$-algebra (resp.,~$\bA_\Inf(R_\infty)$-algebra)
\[
A_{\cris}(R) \ce A(R)\wh{\otimes}_{A_\Inf} A_{\cris} \qq \text{(resp.,} \qq \bA_{\cris}(R_\infty) \ce \bA_\Inf(R_\infty)\wh{\otimes}_{A_\Inf} A_{\cris}),
\]
where the completion is $p$-adic (equivalently, $(p, \mu)$-adic, see \S\ref{Acris-ring}). Due to the decomposition \eqref{dec-dec}, the subring $A_\cris(R)\subset \bA_{\cris}(R_\infty)$ is an $A_\cris(R)$-module direct summand. The maps $\theta$ from \S\ref{Ainf} and \S\ref{Acris-ring} induce compatible surjections 
\[
\theta\colon A_\cris(R) \surjects R \qq \text{and} \qq \theta\colon \bA_{\cris}(R_\infty) \surjects R_\infty.
\]
We let $\bA_\cris^0(R_\infty)$ be the $\bA_\Inf(R_\infty)$-subalgebra of $\bA_\Inf(R_\infty)[\f{1}{p}]$ generated by the elements $\f{\xi^n}{n!}$ for $n \ge 1$. By \cite{Tsu99}*{proof of A2.8}, 
letting $\bA_\Inf(R_\infty)[\f{T^n}{n!}]_{n \ge 1}$ denote the divided power polynomial algebra over $\bA_\Inf(R_\infty)$ in one variable, we have
\[
\tst \bA_\cris^0(R_\infty) \cong (\bA_\Inf(R_\infty)[\f{T^n}{n!}]_{n \ge 1})/(T - \xi), \qq \text{so also} \qq \bA_\cris^0(R_\infty) \cong \bA_\Inf(R_\infty) \tensor_{A_\Inf} A_\cris^0.
\]
Consequently, since $\xi$ generates $\Ker(\theta) \subset \bA_\Inf(R_\infty)$, the ring  $\bA_\cris^0(R_\infty)$ is identified with the divided power envelope of $(\bA_\Inf(R_\infty), \Ker(\theta) + p\bA_\Inf(R_\infty))$ over $(\bZ_p, p\bZ_p)$. By the previous display, 
\[
\bA_{\cris}(R_\infty) \cong (\bA_\cris^0(R_\infty))\wh{\ \ }.
\]
By \S\ref{Ainf}, the ring $\bA_{\cris}(R_\infty)$ (resp.,~$A_\cris(R)$) is $p$-adically formally \'{e}tale as an $\bA_{\cris}(R_\infty^\square)$-algebra (resp.,~$A_\cris(R^\square)$-algebra) and $p$-adically formally flat as an $A_\cris$-algebra. In particular, $\bA_{\cris}(R_\infty)$ inherits $p$-torsion freeness from $A_\cris$. Moreover, even though we will not use this, $\bA_{\cris}(R_\infty)$ is also $\mu$-torsion free, as follows from \Cref{Acrism-inj} below (contrast this with an argument for \eqref{Cris-m-mu}).

The rings $A_{\cris}(R)$ and $\bA_{\cris}(R_\infty)$ come equipped with $A_\cris$-semilinear Frobenius endomorphisms that are compatible with their counterparts for $A_{\cris}^{(m)}(R)$ and $\bA_{\cris}^{(m)}(R_\infty)$ discussed in \S\ref{AcrisR-m}. The profinite group $\Delta$ acts continuously, Frobenius-equivariantly, and $A_\cris$-linearly on $A_\cris(R)$ and $\bA_\cris(R_\infty)$. As in \S\ref{AcrisR-m}, the induced $\Delta$-action on $A_\cris(R)/\mu$ is trivial.}
\epp

We will endow $A_\cris(R)$ with a log structure, which will in fact come from $A(R)$.

\bpp[The log structure on $A(R)$] \lab{log-A(R)}
\ready{Provisionally, we consider the (fine) log structures on $A_\Inf$  and $A(R)$ associated to the prelog structures
\[
\bN_{\ge 0} \xra{a\, \mapsto\, [(p^{1/p^\infty})^q]^a} A_\Inf \qq \text{and} \qq \bN_{\ge 0}^{r + 1} \xra{(a_i)\, \mapsto\, \prod X_i^{a_i}} A(R).
\]
Then, under the diagonal map $\bN_{\ge 0} \ra \bN_{\ge 0}^{r + 1}$, the ring $A(R)$ is a $(p, \mu)$-adically formally log smooth $A_\Inf$-algebra (see \eqref{AR-box} and \cite{Kat89}*{3.5--3.6}).  
To eliminate the dependence on $q$, we always, unless noted otherwise, equip $A_\Inf$ with the log structure associated to the prelog structure
\be \lab{Ainf-log-str}
\cO_C^\flat \setminus \{ 0 \} \ra A_\Inf, \qq x \mapsto [x].
\ee
Likewise, we always equip $A(R)$ with the log structure that is the base change of the fine log structure on $A(R)$ described above along the ``change of log structure'' self-map of $A_\Inf$ determined by $\bN_{\ge 0} \xra{a\, \mapsto\, ((p^{1/p^\infty})^q)^a} \cO_C^\flat \setminus \{ 0 \}$. Explicitly, this log structure is associated to the prelog structure
\be \lab{log-str-A(R)}
\tst  \bN_{\ge 0}^{r + 1} \bigsqcup_{\bN_{\ge 0}} (\cO_C^\flat \setminus \{ 0\}) \ra A(R)
\ee
that embeds $\bN_{\ge 0}$ diagonally into $\bN_{\ge 0}^{r + 1}$, sends an $a \in \bN_{\ge 0}$ to $((p^{1/p^\infty})^q)^a$, and sends the $i\th$ standard basis vector of $\bN_{\ge 0}^{r + 1}$ (resp.,~an $x \in \cO_C^\flat \setminus \{ 0 \}$) to $X_i$ (resp.,~to $[x]$).

These latter ``default'' log structures on $A_\Inf$ and $A(R)$ are quasi-coherent and integral 
and, by base change, with them $A(R)$ is $(p, \mu)$-adically formally log smooth over $A_\Inf$. In fact, via the map $\theta$, the ring $A(R)$ over $A_\Inf$ becomes a $(p, \mu)$-adically formally log smooth thickening of $R/p$ over $\cO_C/p$ (where $R/p$ is endowed with the log structure discussed in \S\ref{log-str-def}). 

The Frobenii of $A_\Inf$ and $A(R)$ extend to the log structures by letting them act as multiplication by $p$ on $\bN_{\ge 0}^{r + 1}$ and $\bN_{\ge 0}$ and as the $p$-th power map on $\cO_C^\flat \setminus \{ 0\}$. Consequently, the Frobenius of the log $A_\Inf$-algebra $A(R)$ lifts the absolute Frobenius of the log $\cO_C/p$-algebra $R/p$.
}

The Frobenius-equivariant $\Delta$-action on the $A_\Inf$-algebra $A(R)$ (see \S\ref{Ainf}) extends to a Frobenius-equivariant $\Delta$-action on the log $A_\Inf$-scheme $\Spec(A(R))$: indeed, a $\delta \in \Delta$ sends each $X_i$ with $0 \le i \le r$ to $u_{\delta,\, i} \cdot X_i$ for some Teichm\"{u}ller unit $u_{\delta,\, i} \in A(R)^\times$ (see \S\ref{Ainf}) and the prelog structures
\[
\bN_{\ge 0}^{r + 1} \xra{(a_i)\, \mapsto\, \prod X_i^{a_i}} A(R) \qq \text{and} \qq \bN_{\ge 0}^{r + 1} \xra{(a_i)\, \mapsto\, \prod (u_{\delta,\, i} \cdot X_i)^{a_i}} A(R)
\]
determine the same log structure on $\Spec(A(R))$, namely, the one associated to the prelog structure
\[
\bZ^{r + 1} \times \bN_{\ge 0}^{r + 1} \xra{((z_i),\, (a_i))\, \mapsto\, \prod u_{\delta,\, i}^{z_i} \cdot \prod X_i^{a_i}} A(R) .
\]
\epp

\bpp[The logarithmic de Rham complex] \lab{log-dR-cx}
With a slight abuse of notation, we let
\[
\tst \Omega^{\bullet}_{A(R)/A_\Inf,\, \log}
\]
be the (global section complex of the) logarithmic de Rham complex of $\Spf(A(R))$ over $\Spf( A_\Inf)$. More precisely, $\Omega^{\bullet}_{A(R)/A_\Inf,\, \log}$ is the (termwise) inverse limit over $n, n' > 0$ of the logarithmic de Rham complexes of $A(R)/(p^n, \mu^{n'})$ over $A_\Inf/(p^n, \mu^{n'})$ (described, for instance, in \cite{Ogu18}*{V.2.1.1}). Due to the formal log smoothness of $A(R)$ over $A_\Inf$, each $\Omega^i_{A(R)/A_\Inf,\, \log}$ is a free $A(R)$-module: indeed, the logarithmic differentials 
\[
d\log(X_1), \dotsc, d\log(X_d)
\]
form an $A(R)$-basis of $\Omega^1_{A(R)/A_\Inf,\, \log}$. We let 
\be \lab{derivations}
\tst \f{\partial}{\partial\log(X_i)} \colon A(R) \ra A(R) \qq \text{for} \qq  i = 1, \ldots, d
\ee
denote the dual basis of log $A_\Inf$-derivations (we do not notationally explicate the accompanying homomorphisms from the log structure to $A(R)$). These satisfy the following explicit formulas that are derived using the relation $d\log(X_0) + \cdots + d\log(X_r) = 0$:
\be \lab{log-expl}
\tst \f{\partial}{\partial\log(X_i)} (X_j) = \begin{cases}  0, \q\ \ \text{if $0 < j \neq i$,} \\ X_i, \q \text{if $j = i$,}\end{cases} \q \text{and} \qq   \f{\partial}{\partial\log(X_i)} (X_0) = \begin{cases} -X_0, \q \text{if $0 < i \le r$,} \\ 0, \q \ \ \, \ \ \text{if $r < i$.} \end{cases} 
\ee
The $\f{\partial}{\partial\log(X_i)}$ also define an isomorphism $\Omega^1_{A(R)/A_\Inf,\, \log} \cong A(R)^{\oplus d}$, which extends to an isomorphism
\be \lab{log-Koszul}
\tst \Omega^{\bullet}_{A(R)/A_\Inf,\, \log} \cong K_{A(R)}\p{\frac{\partial}{\partial\log(X_1)}, \dots, \frac{\partial}{\partial \log(X_d)}}
\ee
that may be considered canonical because its construction only uses data determined by the local coordinate map \eqref{local-setup}.  The Frobenius of the log $A_\Inf$-algebra $A(R)$ multiplies each $d\log(X_i)$ by $p$, so its effect on the right side of \eqref{log-Koszul} is given in each degree $j$ by $p^j$ times the Frobenius of $A(R)$.
\epp

\bpp[The log structure on $A_\cris(R)$] \lab{log-Acris}
\ready{
Unless specified otherwise, we equip the $A(R)$-algebras $A_\cris(R)$ and 
$A_\cris(R)/p^n$ 
for $n > 0$ with the pullback of the log structure on $A(R)$ determined by \eqref{log-str-A(R)}. Thus, since the log structures \eqref{Ainf-log-str} on $A_\Inf$ and \eqref{AL-1} on $A_\cris$ agree, $A_\cris(R)$ is $p$-adically formally log smooth over $A_\cris$. 
Letting the completion be $p$-adic, we set 
\[
\tst \Omega^{\bullet}_{A_\cris(R)/A_\cris,\, \log} \ce \Omega^{\bullet}_{A(R)/A_\Inf,\, \log} \wh{\tensor}_{A_\Inf} A_\cris, 
\]
which is the (global sections of the) logarithmic de Rham complex of $\Spf(A_\cris(R))$ 
over $\Spf(A_\cris)$. 


We use the $p$-adic completeness of $A_\cris(R)$ and its $p$-adic formal flatness over $A_\cris$ (see \S\ref{Acris-R}) to extend the divided power structure of $A_\cris$ to $A_\cris(R)$ (see \cite{SP}*{\href{http://stacks.math.columbia.edu/tag/07H1}{07H1}}). In effect, $A_\cris(R)$ over $A_\cris$ becomes a $p$-adically formally log smooth log PD thickening of $R/p$ over $\cO_C/p$ (compare with~\S\ref{log-A(R)}). 
}
\epp

Through results of \cite{Bei13a}, the following lemma will be key for relating the right side of \eqref{abs-crys-comp-eq} to the logarithmic de Rham cohomology of $\Spf(A_\cris(R))$ over $\Spf(A_\cris)$ in \Cref{RHS-rewrite}.

\blem \lab{PD-smooth}
For each $n \ge 1$, the log smooth log PD thickening $A_\cris (R)/p^n$ over $A_\cris /p^n$ of $R/p$ over $\cO_C/p$ is PD smooth in the sense of \cite{Bei13a}*{\S1.4} \up{see the proof for the definition}.
\elem

\bpf
The PD smoothness is the claim that for every log PD thickening $U \hra \wt{U}$ over the log PD scheme $A_\cris/p^n$ such that $U$ is affine and the log structure of $\wt{U}$ (and hence also of $U$) is integral and quasi-coherent, the indicated diagonal  log PD morphism exists in every commutative diagram
\[
\xymatrix{
U \ar[r] \ar@{^(->}[d] & \Spec(R/p) \ar@{^(->}[r] & \Spec(A_\cris (R)/p^n) \ar[d]^-{\text{log PD}} \\
\wt{U} \ar[rr]_-{\text{log PD}} \ar@{-->}[rru]^-{\text{log PD}} && \Spec(A_\cris/p^n)
}
\]
of log schemes and log (or log PD where indicated) scheme  morphisms over $A_\cris/p^n$  (see \emph{loc.~cit.}).

This sought property of $A_\cris(R)/p^n$ is invariant under base change that changes the log structure on $A_\cris/p^n$, so we may assume that $A_\cris/p^n$ and $A_\cris(R)/p^n$ are instead equipped with the pullbacks of the ``provisional'' fine log structures defined in \S\ref{log-A(R)}. Moreover, since the PD structure of $A_\cris(R)/p^n$ is extended from $A_\cris/p^n$, the log PD thickening $\Spec(R/p) \hra \Spec(A_\cris(R)/p^n)$ over $A_\cris/p^n$ is its own log PD-envelope over $A_\cris/p^n$ (in the sense of \cite{Bei13a}*{\S1.3}). Thus, the log smoothness of $A_\cris(R)/p^n$ over $A_\cris/p^n$ and \cite{Bei13a}*{\S1.4, Remarks (ii)} give the claimed PD smoothness.
\epf


\bprop \label{RHS-rewrite}
In the local setting of \uS\uref{local-setup-sec}, letting $\frac{\partial}{\partial \log(X_i)}\colon A_\cris^{(m)}(R) \ra A_\cris^{(m)}(R)$ denote the $A_\cris^{(m)}$-derivations induced from \eqref{derivations} by base change, we have Frobenius-equivariant identifications
\[
\tst R\Gamma_{\log\cris}(\cO_{\fX_{\cO_C/p}/A_{\cris}}) \cong \Omega^{\bullet}_{A_{\cris}(R)/A_\cris,\, \log} \overset{\eqref{log-Koszul}}{\cong} \p{\varinjlim_{m \ge p} \p{K_{A_\cris^{(m)}(R)}\p{\frac{\partial}{\partial\log(X_1)}, \dots, \frac{\partial}{\partial \log(X_d)}}}}\wh{\ \ }
\]
\up{the Frobenius action on the last term is analogous to the one described after \eqref{log-Koszul}}.
\eprop

\bpf
By \Cref{PD-smooth}, each $A_\cris(R)/p^n$ over $A_\cris/p^n$ is a PD smooth thickening of $R/p$ over $\cO_C/p$, so \cite{Bei13a}*{(1.8.1)} gives the Frobenius-equivariant identification\footnote{\ready{\emph{Loc.~cit.}~uses the logarithmic PD de Rham complex, that is, the quotient of $\Omega^{\bullet}_{\Spf(A_{\cris}(R))/\Spf(A_\cris),\, \log}$ by the PD relations $d(u^{[m]}) = u^{[m - 1]}du$, see \cite{Bei13a}*{\S1.7}. In our situation, there is no difference: since the PD structure of $A_\cris(R)/p^n$ is extended from the base $A_\cris/p^n$, the PD relations hold already in $\Omega^{\bullet}_{\Spf(A_{\cris}(R))/\Spf(A_\cris),\, \log}$. }}
\[
R\Gamma_{\log\cris}(\cO_{\fX_{\cO_C/p}/A_{\cris}}) \cong R\Gamma_\et(\Spf(A_\cris(R)), \Omega^\bullet_{\Spf(A_\cris(R))/\Spf(A_\cris),\, \log}).
\]
Since the sheaves $\Omega^i_{\Spf(A_\cris(R))/\Spf(A_\cris),\, \log}$ are locally free and, in particular, quasi-coherent, they are acyclic for $\Gamma_\et(\Spf(A_\cris(R)), -)$ (see \cite{FK14}*{I.1.1.23~(2)}
), so we have
\[
R\Gamma_\et(\Spf(A_\cris(R)), \Omega^\bullet_{\Spf(A_\cris(R))/\Spf(A_\cris),\, \log}) \cong \Gamma_\et(\Spf(A_\cris(R)), \Omega^\bullet_{\Spf(A_\cris(R))/\Spf(A_\cris),\, \log}).
\]
It remains to observe that the latter complex is identified with $\Omega^{\bullet}_{A_{\cris}(R)/A_\cris,\, \log}$.
\epf

Having expressed the presheaf versions of both sides of \eqref{abs-crys-comp-eq} in the desired forms in \Cref{LHS-rewrite} and \Cref{RHS-rewrite}, we now seek to exhibit an isomorphism between them in \Cref{local-isom}. 

\bpp[The element $\log(\sqp{\eps})$] \lab{log-eps}
Fix an $m \ge p^2$. By the proof of \cite{BMS16}*{12.2},\footnote{\ready{The argument is as follows. Since $p, \mu\xi$ is an $A_\Inf$-regular sequence, $\mu^p - \mu\xi^p \in p\mu\xi A_\Inf$, so $\f{\mu^{p - 1}}{p} = \f{\xi^p}{p} + \xi a$ with $a \in A_\Inf$. Thus, since $\f{(p^2)!}{p^p} \in p\bZ$, we have $\p{\f{\mu^{p - 1}}{p}}^p \in pA_\cris^{(m)}$, so $\f{\mu^{p - 1}}{p}$ is topologically nilpotent in $A_\cris^{(m)}$. In effect, since $\f{1}{(n + 1)!} p^{\lfloor \f{n}{p - 1} \rfloor} \in \bZ_p$, the elements $\f{\mu^n}{(n + 1)!}$ tend to $0$ in the $p$-adic topology of $A_\cris^{(m)}$ and are topologically nilpotent.}}
each $\f{\mu^n}{(n + 1)!} \in A_\cris^{(m)}[\f{1}{p}]$ with $n \ge 1$ lies in $A_\cris^{(m)}$, is $p$-adically topologically nilpotent there, and $p$-adically tends to $0$ as $n \ra \infty$.  
Consequently, recalling that $\mu = [\eps] - 1$, we may define
\[
\tst \log([\eps]) \ce \mu - \f{\mu^2}{2} + \f{\mu^3}{3} - \ldots  \qq \text{in} \qq A_\cris^{(m)},
\]
so that the Frobenius maps $\log(\sqp{\eps})$ to $p \cdot \log([\eps])$. By \emph{loc.~cit.},\footnote{\ready{The argument is as follows. By the previous footnote, $\sum_{n \ge p} \f{(-1)^n\mu^n}{n + 1}$ lies in $pA_\cris^{(m)}$. Thus, since each $\f{\mu^n}{n + 1}$ with $0 < n < p$ is topologically nilpotent in $A_\cris^{(m)}$, so is $\sum_{n \ge 1}  \f{(-1)^n\mu^n}{n + 1}$. In conclusion, $\f{\log([\eps])}{\mu}$ is a unit in $A_\cris^{(m)}$.}}~the elements $\log([\eps])$ and $\mu$ are unit multiples of each other in $A_\cris^{(m)}$, so $\f{(\log([\eps]))^n}{\mu \cdot n!}$ lies in $A_\cris^{(m)}$, is topologically nilpotent if $n > 1$, and $p$-adically tends to $0$ in $A_\cris^{(m)}$ as $n \ra \infty$. 
\epp

The following lemma describes the $\Delta$-action on $A_\cris^{(m)}(R)$ in terms of the derivations $\f{\partial}{\partial\log(X_i)}$ induced on $A_\cris^{(m)}(R)$ by base change from the derivations \eqref{derivations}.

\blem \lab{delta-exp}
For $m \ge p^2$, a $\delta_i \in \Delta$ with $i = 1, \dotsc, d$ \up{see {\upshape\S\ref{Rinfty}}} acts on $A_\cris^{(m)}(R)$ as the series
\be \lab{exp-def}
\tst \exp(\log([\eps]) \cdot \f{\del}{\del \log(X_i)}) \ce \sum_{n \ge 0} \f{(\log([\eps]))^n}{n!} (\f{\del}{\del \log(X_i)})^n.
\ee
In particular, for $m$ and $i$ as above, we have the following description of the ``$q$-derivative'' $\f{\delta_i - 1}{\mu}${\upshape:}
\be \lab{del-id}
\tst \f{\delta_i - 1}{\mu} = \frac{\partial}{\partial\log (X_i)} \cdot \p{\sum_{n\geq 1}\frac{(\log([\eps]))^n}{\mu \cdot n!}(\frac{\partial}{\partial\log (X_i)})^{n - 1}} \qq \text{as maps} \qq A_\cris^{(m)}(R) \ra A_\cris^{(m)}(R),
\ee
where the parenthetical factor defines an $A_\cris^{(m)}$-linear additive automorphism of $A_\cris^{(m)}(R)$.
\elem

\bpf
The argument is similar to that of \cite{BMS16}*{12.4}. Firstly, $\f{(\log([\eps]))^n}{n!}$ tends to $0$ in the $p$-adic topology of $A_\cris^{(m)}$ (see \S\ref{log-eps}), so the series \eqref{exp-def} does define an $A_\cris^{(m)}$-linear additive endomorphism of $A_\cris^{(m)}(R)$. This endomorphism is also multiplicative because, by the Leibniz rule, 
\[
\tst \f{(\log([\eps]))^n}{n!}(\f{\del}{\del \log(X_i)})^n(ab) = \sum_{j = 0}^n\p{ \f{(\log([\eps]))^j}{j!} (\f{\del}{\del \log(X_i)})^j(a) \cdot \f{(\log([\eps]))^{n - j}}{(n - j)!} (\f{\del}{\del \log(X_i)})^{n - j}(b)}.
\]
Therefore, in the case $R = R^\square$, the desired equality
\be \lab{delta-exp-eq}
\tst \delta_i = \exp(\log([\eps]) \cdot \f{\del}{\del \log(X_i)}) \qq \text{of endomorphisms} \qq A_\cris^{(m)}(R^\square) \ra A_\cris^{(m)}(R^\square)
\ee
 follows by noting that both of its sides agree on every $X_j$: indeed, due to the formulas \eqref{log-expl}, they send $X_i$ to $[\eps]X_i$, fix each $X_j$ with $0 < j \neq i$, and send $X_0$ to $[\eps\i]X_0$ if $i \le r$ and to $X_0$ if $r < i$. 

In the general case, since $\mu$, and hence also $\xi$, divides each $\f{(\log([\eps]))^n}{n!}$ with $n \ge 1$ (see \S\ref{log-eps}), both sides of the equality \eqref{delta-exp-eq} induce the trivial action modulo $(p, \xi)$ (see \S\ref{AcrisR-m}). Therefore, due to the formal \'{e}taleness of $A_\cris^{(m)}(R)$ over $A_\cris^{(m)}(R^\square)$ and the settled $R = R^\square$ case, the sides agree. 

Since $A_\cris^{(m)}(R)$ is $\mu$-torsion free (see \eqref{Cris-m-mu}) and $\mu \mid \frac{(\log([\eps]))^n}{n!}$ in $A_\cris^{(m)}$, the equality \eqref{del-id} follows from \eqref{delta-exp-eq}. Since $\frac{(\log([\eps]))^n}{\mu \cdot n!}$ is a unit for $n = 1$, is topologically nilpotent for $n > 1$ (see \S\ref{log-eps}), and $p$-adically tends to $0$ as $n \ra \infty$, the parenthetical factor of \eqref{del-id} is indeed an automorphism.
\epf

We are ready to settle the (presheaf version of the) local case of \Cref{abs-crys-comp}.

\bprop \label{local-isom}
In the local setting of \uS\uref{local-setup-sec}, for $m \ge p^2$ and $i = 1, \ldots, d$, the morphism
\be \lab{local-isom-1}
\tst
 \p{A_\cris^{(m)}(R) \xra{\frac{\partial}{\partial\log(X_i)}} A_\cris^{(m)}(R)} \xra{\p{\id,\,\, \sum_{n\geq 1}\frac{(\log([\eps]))^n}{ n!}(\frac{\partial}{\partial\log (X_i)})^{n - 1}}}  \p{A_\cris^{(m)}(R) \xra{\delta_i - 1} A_\cris^{(m)}(R)}
\ee
of complexes in degrees $0$ and $1$ is Frobenius equivariant, granted that the usual Frobenius action on the copy of $A_\cris^{(m)}(R)$ in degree $1$ of the source is multiplied by $p$ \up{compare with the description after \eqref{log-Koszul}}. For $m \ge p^2$, these morphisms induce a Frobenius-equivariant quasi-isomorphism 
\be \lab{local-isom-3}
\tst K_{A_\cris^{(m)}(R)}\p{\frac{\partial}{\partial\log(X_1)}, \dots, \frac{\partial}{\partial \log(X_d)}} \isomto \eta_{(\mu)}\p{K_{\bA_\cris^{(m)}(R_\infty)}\p{\delta_1 - 1, \dots, \delta_d - 1}},
\ee
which, as $m$ varies, induces the Frobenius-equivariant identification \up{a local version of~\eqref{abs-crys-comp-eq}}{\upshape:}
\be \lab{local-isom-2}
 R\Gamma_{\log\cris}(\cO_{\fX_{\cO_C/p}/A_{\cris}}) \cong R\Gamma(\fX_\et^\psh, A\Omega_\fX^\psh) \wh{\tensor}^\bL_{A_\Inf} A_\cris.
\ee
\eprop

\bpf
The Frobenius-equivariance of \eqref{local-isom-1} follows from the equations  
\[
\tst \frac{\partial}{\partial\log (X_i)} \circ \varphi = p \cdot \p{\varphi \circ \frac{\partial}{\partial\log (X_i)}} \qq \text{and} \qq \varphi(\log([\eps])) = p \cdot \log([\eps])
\]
(see \S\ref{log-dR-cx} and \S\ref{log-eps}). Since $\Delta$ acts trivially on $A_\cris^{(m)}(R)/\mu$ (see \S\ref{AcrisR-m}), the subcomplex 
\[
\eta_{(\mu)}\p{K_{A_\cris^{(m)}(R)}\p{\delta_1 - 1, \dots, \delta_d - 1}} \subset K_{A_\cris^{(m)}(R)}\p{\delta_1 - 1, \dots, \delta_d - 1}
\]
is obtained by letting its $j$-th term be the submodule of the $j$-th term of $K_{A_\cris^{(m)}(R)}\p{\delta_1 - 1, \dots, \delta_d - 1}$ comprised of the $\mu^j$-multiples (see \eqref{Koszul-def} and \eqref{eta-def});  since $\mu \mid \varphi(\mu)$, this subcomplex is Frobenius-stable. Thus, \Cref{delta-exp} implies that the morphisms \eqref{local-isom-1} induce an isomorphism
\be \lab{local-isom-4}
\tst K_{A_\cris^{(m)}(R)}\p{\frac{\partial}{\partial\log(X_1)}, \dots, \frac{\partial}{\partial \log(X_d)}} \isomto \eta_{(\mu)}\p{K_{A_\cris^{(m)}(R)}\p{\delta_1 - 1, \dots, \delta_d - 1}}.
\ee
\Cref{N-cris-coho} (with \Cref{Koszul}) implies that the natural inclusion of the target of \eqref{local-isom-4} into the target of \eqref{local-isom-3} is a quasi-isomorphism, and \eqref{local-isom-3} follows. The maps \eqref{local-isom-3} are compatible as $m$ varies, so, by passing to their limit over $m$, forming the termwise $p$-adic completions, and applying \Cref{LHS-rewrite} and \Cref{RHS-rewrite}, we obtain the desired identification \eqref{local-isom-2}.
\epf

\Cref{local-isom} concludes the ``single coordinate patch'' part of the proof of \Cref{abs-crys-comp}, so we turn to the ``all possible coordinates'' technique that will globalize the argument. For this, the key steps are, for a small enough affine $\fX$, to build in \S\ref{functorial-AO} a \emph{functorial} in $\fX$ explicit complex that computes the presheaf version of the left side of \eqref{abs-crys-comp-eq}, to then build in \S\ref{functorial-cris} such a complex for the right side of \eqref{abs-crys-comp-eq}, and, finally, to build in \S\ref{comp-map} and \Cref{alphaR-qiso} a natural quasi-isomorphism between these complexes. Each of these steps will use our work in the setting of \S\ref{local-setup-sec} discussed so far.

\bpp [More general coordinates] \label{all-possible-coordinates}
\ready{
Continuing to work locally, we now assume until the final part of the proof of \Cref{abs-crys-comp} given in \S\ref{Acris-main-pf} that $\fX = \Spf(R)$ is affine and nonempty, that every two irreducible components of $\Spec(R \tensor_{\cO_C} k)$ meet (so that $\fX$ is connected), 
and that we~have
\begin{itemize}
\item
a finite set $\Sigma$ that indexes the coordinates of the formal $\cO_C$-torus
\[
R^\square_\Sigma \ce \cO_C\{ t_\sigma^{\pm 1} \, | \,  \sigma \in \Sigma \};
\]

\item
a nonempty finite set $\Lambda$ and, for each $\gL \in \Lambda$, an $\cO_C$-algebra 
\[
\qq R^{\square}_\lambda \ce \cO_C\{t_{\lambda,\, 0}, \dotsc, t_{\lambda,\, r_\lambda}, t_{\lambda,\, r_\lambda + 1}^{\pm 1}, \dotsc,  t_{\lambda,\, d}^{\pm 1}\}/(t_{\lambda,\, 0}\cdots t_{\lambda,\,r_\lambda} - p^{q_\lambda}) \q \text{with} \q q_\gL \in \bQ_{> 0};
\]

\item
a closed immersion 
\be \lab{more-coord}
\tst \fX = \Spf(R) \ra \Spf( R^\square_\Sigma) \times \prod_{\lambda\in\Lambda} \Spf(R^{\square}_\lambda), 
\ee
where the products are formed over $\Spf( \cO_C)$, subject to the requirements that already 
\be \lab{Sigma-ci}
\fX = \Spf (R)\ra \Spf (R^\square_\Sigma) \qq \text{is a closed immersion}
\ee
and, for each $\gL \in \Lambda$, the induced map 
\be \lab{lambda-proj}
\qq \fX = \Spf (R)\ra \Spf( R^{\square}_\gL) \qq \text{is \'{e}tale.}
\ee
\end{itemize}
By \eqref{lambda-proj}, for each $\lambda \in \Lambda$, the irreducible components of $\Spec(R \tensor_{\cO_C} k)$ are \emph{a priori} identified with the connected components of $\bigsqcup_i \Spec((R \tensor_{\cO_C} k)/(t_{\lambda,\, i}))$. Thus, our assumption on $\Spec(R \tensor_{\cO_C} k)$ implies that each irreducible component of $\Spec(R \tensor_{\cO_C} k)$ is cut out by a unique $t_{\lambda,\, i}$ with $0 \le i \le r_\gL$.

By \S\ref{fir-setup}, if $R \tensor_{\cO_C} k$ is not $k$-smooth, then $R$ determines $q_\lambda$, which therefore does not depend on $\lambda$. On the other hand, if $R \tensor_{\cO_C} k$ is $k$-smooth, then $q_\lambda$ may depend on $\lambda$. This, together with the possibility that $r_\lambda > 0$, complicates matters in the ``simpler'' smooth case but is crucial to allow in order for the eventual ``all possible coordinates'' constructions to be functorial in $R$.

For any $\fX$, the data above exist on a basis for $\fX_\et$: indeed, a coordinate map \eqref{lambda-proj} exists \'{e}tale locally on $\fX$ (see \S\ref{fir-setup}), and then $R$ is the $p$-adic completion of a finite type $\cO_C$-algebra, so the Zariski topology of $\Spf(R)$ has a basis whose elements embed into some (variable) $\wh{\bG}_m^\Sigma$.}

Each \eqref{lambda-proj} is an instance of the local setting of \uS\uref{local-setup-sec}, so the discussion between \S\ref{local-setup-sec} and the present section applies to it. Another instance is the identity map $\Spf( R_\Sigma^\square) \xra{=} \Spf (R_\Sigma^\square)$ (with $r = 0$ and $d = \# \Sigma$), so the indicated discussion also applies to the ring $R_\Sigma^\square$ in place of $R^\square$.
\epp

Our first aim in this setup is to reexpress the (presheaf version of the) left side of \eqref{abs-crys-comp-eq} in \S\ref{functorial-AO}.

\bpp[The perfectoid $R_{\Sigma,\, \Lambda,\, \infty}$] \lab{more-cover}
\ready{For each $\lambda \in \Lambda$, we set
\[
\tst \Delta_\lambda \ce \left\{ (\eps_0, \dotsc, \eps_d) \in \p{\varprojlim_{m \ge 0} \p{\mu_{p^m}(\cO_C)}}^{\oplus (d + 1)}\, \Big\vert\, \eps_0\cdots \eps_{r_\lambda} = 1 \right\} \simeq \bZ_p^{\oplus d}
\]
and let 
\[
\tst \Spa(R_{\lambda,\, \infty}[\f{1}{p}], R_{\lambda,\, \infty}) \ra \Spa(R[\f{1}{p}], R) \q \text{and} \q \Spa(R^\square_{\lambda,\, \infty}[\f{1}{p}], R^\square_{\lambda,\, \infty}) \ra \Spa(R^\square_{\lambda}[\f{1}{p}], R^\square_{\lambda})
\]
be the affinoid perfectoid pro-(finite \'{e}tale) $\Delta_\lambda$-covers defined as in \S\ref{Rinfty} using the coordinate map $\Spf (R) \xra{\eqref{lambda-proj}} \Spf( R^{\square}_\gL)$. Similarly, we set
\[
\tst \Delta_\Sigma \ce \p{\varprojlim_{m \ge 0} (\mu_{p^m}(\cO_C))}^{\Sigma} \simeq \bZ_p^{\Sigma}
\]
and let 
\[
\tst \Spa(R^\square_{\Sigma,\, \infty}[\f{1}{p}], R^\square_{\Sigma,\, \infty}) \ra \Spa(R^\square_{\Sigma}[\f{1}{p}], R^\square_{\Sigma})
\]
be the affinoid perfectoid pro-(finite \'{e}tale) $\Delta_\Sigma$-cover defined as in \S\ref{Rinfty} using the coordinate map $\Spf (R_\Sigma^\square) \xra{=} \Spf (R_\Sigma^\square)$, so that, explicitly,
\[
\tst R^\square_{\Sigma,\, \infty} \cong \p{\varinjlim_{m \ge 0} (\cO_C\{ t_\sigma^{\pm 1/p^m} \, | \,  \sigma \in \Sigma \})}\widehat{\ \ }.
\]
} By forming products over $\Spa(\cO_C[\f{1}{p}], \cO_C)$ and setting
\[
\tst \Delta_{\Sigma,\, \Lambda} \ce \Delta_\Sigma \times \prod_{\gL \in \Lambda} \Delta_\gL,
\]
we obtain the affinoid perfectoid pro-(finite \'{e}tale) $\Delta_{\Sigma,\, \Lambda}$-cover
\[
\tst \Spa(R^\square_{\Sigma,\, \infty}[\f{1}{p}], R^\square_{\Sigma,\, \infty}) \times \prod_{\gL \in \Lambda} \Spa(R^\square_{\lambda,\, \infty}[\f{1}{p}], R^\square_{\lambda,\, \infty})  \ra \Spa(R^\square_{\Sigma}[\f{1}{p}], R^\square_{\Sigma}) \times \prod_{\gL \in \Lambda} \Spa(R^\square_{\lambda}[\f{1}{p}], R^\square_{\lambda}),
\]
which we abbreviate as
\[
\tst \Spa(R^\square_{\Sigma,\, \Lambda,\, \infty}[\f{1}{p}], R^\square_{\Sigma,\, \Lambda,\, \infty}) \ra \Spa(R^\square_{\Sigma,\, \Lambda}[\f{1}{p}], R^\square_{\Sigma,\, \Lambda}).
\]
Its base change along the generic fiber of \eqref{more-coord} is the pro-(finite \'{e}tale)~$\Delta_{\Sigma,\, \Lambda}$-cover
\be \lab{big-cover-def}
\tst \Spa(R_{\Sigma,\, \Lambda,\, \infty}[\f{1}{p}], R_{\Sigma,\, \Lambda,\, \infty}) \ra \Spa(R[\f{1}{p}], R),
\ee
which contains each $\Spa(R_{\lambda,\, \infty}[\f{1}{p}], R_{\lambda,\, \infty}) \ra \Spa(R[\f{1}{p}], R)$ as a subcover. Thus, by the almost purity theorem \cite{Sch12}*{7.9~(iii)}, the $\cO_C$-algebra $R_{\Sigma,\, \Lambda,\, \infty}$ defined by \eqref{big-cover-def} is perfectoid (the notions of `perfectoid' used here and in \cite{Sch12} agree by \cite{BMS16}*{3.20}). 

The topological generators for $\Delta_\Sigma$ and $\Delta_\lambda$ fixed in \S\ref{Rinfty} are
\[
\delta_\sigma \ce (1, \ldots, 1, \eps, 1, \ldots, 1) \q \text{for $\sigma \in \Sigma$,} \qq \text{where the $\sigma$-th entry is nonidentity,}
\]
and
\[
\ba
\q \delta_{\lambda,\, i} &\ce (\eps\i, 1, \dotsc, 1, \eps, 1, \dotsc, 1) \q \text{for $i = 1, \ldots, r_\gL$,} \q \text{where the $0$-th and $i$-th entries are nonidentity;} \\
\q \delta_{\lambda,\, i} &\ce (1, \dotsc, 1, \eps, 1, \dotsc, 1) \q \text{for $i = r_\gL + 1, \ldots, d$,} \q\, \text{where the $i$-th entry is nonidentity.} 
\ea
\]
Jointly, the $\delta_\sigma$'s and the $\delta_{\lambda,\, i}$'s topologically freely generate $\Delta_{\Sigma,\, \Lambda}$. 
\epp

\bpp[The rings $\bA_{\Inf}(R_{\Sigma,\, \Lambda,\, \infty})$ and $\bA_{\cris}^{(m)}(R_{\Sigma,\, \Lambda,\, \infty})$] \lab{more-Ainf}
\ready{Similarly to \S\ref{Ainf}, we set
\[
\bA_{\Inf}(R_{\Sigma,\, \Lambda,\, \infty}) \ce W(R_{\Sigma,\, \Lambda,\, \infty}^\flat).
\]
By \Cref{A-Apr}, for $n, n' > 0$, the sequence $(p^n, \mu^{n'})$ is $\bA_{\Inf}(R_{\Sigma,\, \Lambda,\, \infty})$-regular, $\bA_{\Inf}(R_{\Sigma,\, \Lambda,\, \infty})/(p^n, \mu^{n'})$ is $A_\Inf/(p^n, \mu^{n'})$-flat, and $\bA_{\Inf}(R_{\Sigma,\, \Lambda,\, \infty})/\mu$ is $p$-adically complete. As in \S\ref{Ainf}, we have the surjection
\be \lab{more-theta}
\theta\colon \bA_{\Inf}(R_{\Sigma,\, \Lambda,\, \infty}) \surjects R_{\Sigma,\, \Lambda,\, \infty}
\ee
that intertwines the Witt vector Frobenius of $\bA_{\Inf}(R_{\Sigma,\, \Lambda,\, \infty})$ with the absolute Frobenius of $R_{\Sigma,\, \Lambda,\, \infty}/p$ and whose kernel is generated by the regular element $\xi$.
To fix further notation, we let 
\be \lab{A-box-def} \ba
A(R^\square_\Sigma) &\cong A_\Inf\{X_{\sigma}^{\pm 1}\, \vert\, \sigma \in \Sigma\}, \\ A(R^\square_\lambda) &\cong A_\Inf\{X_{\lambda,\, 0}, \dotsc, X_{\lambda,\, r_\lambda}, X_{\lambda,\, r_\lambda + 1}^{\pm 1}, \dotsc,  X_{\lambda,\, d}^{\pm 1}\}/(X_{\lambda,\, 0}\cdots X_{\lambda,\, r_\lambda} - [(p^{1/p^\infty})^{q_\lambda}])
\ea
\ee
be the isomorphisms \eqref{AR-box} for $R_\Sigma^\square$ and $R_\lambda^\square$.
Similarly to \S\ref{AcrisR-m}, for an $m\in \bZ_{ \ge 1}$, we~set
\be \lab{more-Acrism}
\bA_\cris^{(m)}(R_{\Sigma,\, \Lambda,\, \infty}) \ce \bA_{\Inf}(R_{\Sigma,\, \Lambda,\, \infty}) \wh{\tensor}_{A_\Inf} A_\cris^{(m)},
\ee
where the completion is $(p, \mu)$-adic (equivalently, $p$-adic if $m \ge p$). Since $\bA_{\Inf}(R_{\Sigma,\, \Lambda,\, \infty})$ is $(p, \mu)$-adically formally flat over $A_\Inf$, the ring $\bA_\cris^{(m)}(R_{\Sigma,\, \Lambda,\, \infty})$ inherits $p$-torsion freeness from $A_\cris^{(m)}$. By also using the short exact sequences \eqref{mu-SES-2} and the vanishing \eqref{red-vanishes}, we see that 
\[
\bA_\cris^{(m)}(R_{\Sigma,\, \Lambda,\, \infty}) \q  \text{is $\mu$-torsion free} \qq \text{and}\qq  \bA_\cris^{(m)}(R_{\Sigma,\, \Lambda,\, \infty})/\mu \q \text{is $p$-adically complete.}
\]
As in \S\ref{AcrisR-m}, the rings $\bA_\cris^{(m)}(R_{\Sigma,\, \Lambda,\, \infty})$ come equipped with $A_\cris^{(m)}$-semilinear Frobenius endomorphisms that are compatible as $m$ varies. The maps \eqref{theta-Acrism} and \eqref{more-theta} give rise to the surjection
\be \lab{AAcrism-theta} 
\theta\colon \bA_\cris^{(m)}(R_{\Sigma,\, \Lambda,\, \infty}) \surjects R_{\Sigma,\, \Lambda,\, \infty}.
\ee
The actions of the profinite group $\Delta_{\Sigma,\, \Lambda}$ on $\bA_{\Inf}(R_{\Sigma,\, \Lambda,\, \infty})$ and $\bA_\cris^{(m)}(R_{\Sigma,\, \Lambda,\, \infty})$ are compatible, continuous, and Frobenius-equivariant. }
\epp

The following consequence of \Cref{more-ANB-more} will help us build  a desired functorial complex in \S\ref{functorial-AO}.

\bprop \label{crys-limit}
In the local setting of {\upshape\S\ref{all-possible-coordinates}}, for every $m \ge p$, the analogue for $R_{\Sigma,\,\Lambda,\,\infty}$ of the edge map \eqref{another-e} induces the Frobenius-equivariant identification
\be \lab{CL-eq}
\eta_{(\mu)} \p{ K_{\bA_{\cris}^{(m)}(R_{\Sigma,\, \Lambda,\, \infty})}((\delta_\sigma - 1)_{\sigma\in\Sigma}, (\delta_{\lambda,\, i}-1)_{\lambda\in\Lambda,\, 1 \le i \le  d})} \isomto R\Gamma(\fX_\et^\psh, A\Omega_\fX^\psh) \wh{\tensor}^\bL_{A_\Inf} A_\cris^{(m)}.
\ee
In particular, we have the following Frobenius-equivariant identification in the derived category{\upshape:}
\[
\tst \p{\varinjlim_{m} \p{\eta_{(\mu)} \p{ K_{\bA_{\cris}^{(m)}(R_{\Sigma,\, \Lambda,\, \infty})}((\delta_\sigma - 1)_{\sigma\in\Sigma}, (\delta_{\lambda,\, i}-1)_{\lambda\in\Lambda,\, 1 \le i \le  d})}}}\wh{\ } \isomto R\Gamma(\fX_\et^\psh, A\Omega_\fX^\psh) \wh{\tensor}^\bL_{A_\Inf} A_\cris,
\]
where the direct limit and the $p$-adic completion of the complexes in the source are formed termwise.
\eprop

\bpf
\Cref{L-cris-commute} gives the Frobenius-equivariant identification
\[
R\Gamma(\fX_\et^\psh, A\Omega_\fX^\psh) \wh{\tensor}^\bL_{A_\Inf} A_\cris^{(m)} \cong L\eta_{(\mu)}(R\Gamma_\proet(\fX^\ad_C, \bA_{\Inf}) \wh{\tensor}^\bL_{A_\Inf} A_\cris^{(m)}).
\]
Therefore, since the pro-(finite \'{e}tale) affinoid perfectoid $\Delta_{\Sigma,\, \Lambda}$-cover $\Spa(R_{\Sigma,\, \Lambda,\, \infty}[\f{1}{p}], R_{\Sigma,\, \Lambda,\, \infty})$ of $\Spa(R[\f{1}{p}], R)$ 
contains $\Spa(R_{\lambda,\, \infty}[\f{1}{p}], R_{\lambda,\, \infty})$ as a subcover (see \S\ref{more-cover}), \Cref{more-ANB-more} applies and (with \Cref{Koszul}) gives \eqref{CL-eq}. The remaining assertion follows: each $\bA_{\cris}^{(m)}(R_{\Sigma,\, \Lambda,\, \infty})$ is $p$-torsion free, so the termwise $p$-adic completion of the source there agrees with the derived $p$-adic completion.
\epf

\bpp[A functorial complex that computes $R\Gamma(\fX_\et^\psh, A\Omega_\fX^\psh)  \wh{\otimes}^{\bL}_{A_{\Inf}}A_{\cris}$] \lab{functorial-AO}
For a fixed $R$, we form the filtered direct limit over the closed immersions \eqref{more-coord} for varying $\Sigma$ and $\Lambda$ to build the complex 
\be \lab{FAO-colim}
\tst \varinjlim_{\Sigma,\, \Lambda}\p{ \p{\varinjlim_{m \ge p} \p{ \eta_{(\mu)} \p{ K_{\bA_{\cris}^{(m)}(R_{\Sigma,\, \Lambda,\, \infty})}((\delta_\sigma - 1)_{\sigma\in\Sigma}, (\delta_{\lambda,\, i}-1)_{\lambda\in\Lambda,\, 1 \le i \le  d})}}}\wh{\ }},
\ee
where the direct limits and the $p$-adic completion are termwise. By its construction, this complex comes equipped with an $A_\cris$-semilinear Frobenius endomorphism. The isomorphisms of \Cref{crys-limit} are compatible with enlarging $\Sigma$ and $\Lambda$, so they show that in the derived category the complex \eqref{FAO-colim} is canonically and Frobenius-equivariantly identified with
\[
R\Gamma(\fX_\et^\psh, A\Omega_\fX^\psh) \wh{\tensor}^\bL_{A_\Inf} A_\cris.
\]
Moreover, if $R'$ is a $p$-adically formally \'{e}tale $R$-algebra equipped with data as in \S\ref{all-possible-coordinates} for some sets $\Sigma'$ and $\Lambda'$, then the term indexed by $\Sigma$, $\Lambda$ (and by the closed immersion \eqref{more-coord}) of the direct limit \eqref{FAO-colim} maps to the term indexed by $\Sigma \cup \Sigma'$, $\Lambda \cup \Lambda'$ (and by a closed immersion of $\Spf (R')$) of the analogous direct limit for $R'$, compatibly with the transition maps in \eqref{FAO-colim} and with the Frobenius. Thus, the complex \eqref{FAO-colim} equipped with its Frobenius is functorial in $R$, and so is its identification with $R\Gamma(\fX_\et^\psh, A\Omega_\fX^\psh) \wh{\tensor}^\bL_{A_\Inf} A_\cris$.
\epp

Our next aim is to similarly reexpress the (presheaf version of the) right side of \eqref{abs-crys-comp-eq} in \S\ref{functorial-cris}. 

\bpp [The completed log PD envelope $D_{\Sigma,\, \Lambda}$] \label{log-PD}
\ready{By \S\ref{log-A(R)}, the maps $\theta\colon A(R_\lambda^\square) \surjects R_\lambda^\square$ of \eqref{AR-surj} are compatible with log structures. Thus, they give rise to a Frobenius-equivariant closed immersion
\be \lab{log-ci}
\tst \Spec(R/p) \hra \Spf(A(R_\Sigma^\square)) \times \prod_{\lambda \in \Lambda } \Spf(A(R_\lambda^\square)) \equalscolon \Spf(A_{\Sigma,\, \Lambda}^\square)
\ee
of $(p, \mu)$-adic formal log schemes, where the products are over the $(p, \mu)$-adic formal log scheme $\Spf(A_\Inf)$. By \cite{Kat89}*{4.1, 4.4}, for $n, n' > 0$, the quasi-coherent log structure of $\Spec(A_{\Sigma,\, \Lambda}^\square/(p^n, \mu^{n'}))$ and the log scheme map $\Spec(A_{\Sigma,\, \Lambda}^\square/(p^n, \mu^{n'})) \ra \Spec(A_\Inf/(p^n, \mu^{n'}))$ are integral.  
 
 For $n, n' > 0$, by \cite{Bei13a}*{1.3, Theorem}, the $A_\Inf/(p^n, \mu^{n'})$-base change of the closed immersion \eqref{log-ci} has a log PD envelope 
 \[
 \Spec(D_{\Sigma,\,\Lambda,\, n,\, n'}) \qq \text{over} \qq (\bZ/p^n\bZ, p\bZ/p^n\bZ),
 \]
 which, in particular, is a nil thickening of $\Spec(R/p)$, so is affine as indicated (see \cite{SP}*{\href{http://stacks.math.columbia.edu/tag/01ZT}{01ZT}}). In fact, $D_{\Sigma,\,\Lambda,\, n,\, n'}$ is supplied already by \cite{Kat89}*{5.4} because the closed immersion \eqref{log-ci} is the base change of a similar closed immersion of \emph{fine} formal log schemes along a  ``change of log structure'' self-map of $A_\Inf$ (see \S\ref{log-A(R)}).\footnote{\ready{The two references characterize the log PD envelope differently, but they give the same $\Spec(D_{\Sigma,\, \Lambda,\, n,\, n'})$, in essence because the image of any monoid morphism $M \ra M'$ with $M$ finitely generated is finitely generated. }}  

If $n'$ is large enough relative to $n$, so that $\mu^{n'} \in p^nA_\cris$, then, by \S\S\ref{Acris-ring}--\ref{Acris-log}, $\Spec(A_\cris/p^n)$ is identified with the log PD envelope of the exact log closed immersion $\Spec(\cO_C/p) \hra \Spec(A_\Inf/(p^n, \mu^{n'}))$ over $(\bZ/p^n\bZ, p\bZ/p^n\bZ)$. Thus, for such $n$, $n'$, the envelope $\Spec(D_{\Sigma,\,\Lambda,\, n,\, n'})$ comes equipped with a canonical log PD morphism to $\Spec(A_\cris/p^n)$ that identifies it with the log PD envelope of 
\[
\Spec(R/p) \hra \Spec(A_{\Sigma,\, \Lambda}^\square \tensor_{A_\Inf} A_\cris/p^n) \qq \text{over} \qq \Spec(\cO_C/p) \hra \Spec(A_\cris/p^n).
\]
Thus, letting $D_{\Sigma,\,\Lambda,\, n}$ be this log PD envelope, that is, the common $D_{\Sigma,\,\Lambda,\, n,\, n'}$ for large $n'$, we have $D_{\Sigma,\,\Lambda,\, n}/p^{n - 1} \cong D_{\Sigma,\,\Lambda,\, n - 1}$ for $n > 1$ and obtain a $p$-adic formal log $\Spf(A_\cris)$-scheme $\Spf(D_{\Sigma,\, \Lambda})$ that fits into a factorization
\be \lab{Acris-fact}
\tst \Spec(R/p) \hra \Spf(D_{\Sigma,\, \Lambda}) \ra  \Spf(A_\cris(R_\Sigma^\square)) \times \prod_{\lambda \in \Lambda } \Spf(A_\cris(R_\lambda^\square)) \equalscolon \Spf(A_{\Sigma,\, \Lambda,\, \cris}^\square),
\ee
where the products are formed over the $p$-adic formal log scheme $\Spf(A_\cris)$ and we have 
\[
A_{\Sigma,\, \Lambda,\, \cris}^\square \cong A_{\Sigma,\, \Lambda}^\square \wh{\tensor}_{A_\Inf} A_\cris.
\] 
By functoriality, $\Spf(D_{\Sigma,\, \Lambda})$ comes equipped with an $A_\cris$-semilinear Frobenius. In addition, since, for each $n > 0$, the ideal defining the exact closed immersion $\Spec(R/p) \hra \Spec(R/p^n)$ inherits divided powers from $\bZ/p^n$, the universal property of $D_{\Sigma,\, \Lambda}$ supplies the factorization
\be \lab{R-D-map}
\Spec(R/p) \hra \Spf(R) \hra \Spf(D_{\Sigma,\, \Lambda}) \q \text{over} \q \Spec(\cO_C/p) \hra \Spf(\cO_C) \hra \Spf(A_\cris).
\ee
The profinite group $\Delta_{\Sigma,\, \Lambda}$ acts continuously and Frobenius-equivariantly on $A_{\Sigma,\, \Lambda}^\square$ over $A_\Inf$ (see \S\ref{Ainf}) and, due to the last paragraph of \S\ref{log-A(R)}, this action extends to a $\Delta_{\Sigma,\, \Lambda}$-action on the $(p, \mu)$-adic formal log scheme $\Spf(A_{\Sigma,\, \Lambda}^\square)$. Thus, since the closed immersion \eqref{log-ci} is $\Delta_{\Sigma,\, \Lambda}$-equivariant, $\Delta_{\Sigma,\, \Lambda}$  acts $A_\cris$-linearly and Frobenius-equivariantly on each $D_{\Sigma,\, \Lambda,\, n}$ and also on $D_{\Sigma,\, \Lambda}$. }
\epp

The main practical deficiency of $D_{\Sigma,\,\Lambda}$ is its inexplicit nature, for instance, we do not know whether  $D_{\Sigma,\,\Lambda}$ is $p$-torsion free. In contrast, its utility for us manifests itself through the following proposition.

\bprop \lab{D-compute}
In the local setting of {\upshape\S\ref{all-possible-coordinates}}, the complex \up{where the inverse limit is termwise}
\[
\tst \Omega^\bullet_{D_{\Sigma,\, \Lambda}/A_\cris,\, \log,\, \PD} \ce \varprojlim_{n > 0} \p{\Omega^{\bullet}_{(A_{\Sigma,\, \Lambda,\, \cris}^\square/p^n)/(A_\cris/p^n),\, \log} \tensor_{A_{\Sigma,\, \Lambda,\, \cris}^\square/p^n} D_{\Sigma,\, \Lambda,\, n}}
\]
is canonically and Frobenius-equivariantly identified in the derived category as follows{\upshape:}
\be \lab{D-compute-id}
 R\Gamma_{\log\cris}(\cO_{\fX_{\cO_C/p}/A_{\cris}}) \cong \Omega^\bullet_{D_{\Sigma,\, \Lambda}/A_\cris,\, \log,\, \PD}.
\ee
Under this identification, the map
\be \lab{cris-dR-map}
R\Gamma_{\log\cris}(\cO_{\fX_{\cO_C/p}/A_{\cris}}) \ra R\Gamma_{\log\dR}(\fX/\cO_C) \ \ \, \text{is} \ \ \, \Omega^\bullet_{D_{\Sigma,\, \Lambda}/A_\cris,\, \log,\, \PD} \xra{\eqref{R-D-map}} \Omega^\bullet_{\Spf(R)/\cO_C,\,\log}.
\ee
In particular, we have a Frobenius-equivariant identification
\be \lab{D-Koszul}
\tst R\Gamma_{\log\cris}(\cO_{\fX_{\cO_C/p}/A_{\cris}}) \cong K_{D_{\Sigma,\, \Lambda}}\p{\p{\frac{\partial}{\partial\log(X_\sigma)}}_{\sigma\in\Sigma}, \p{\frac{\partial}{\partial\log(X_{\lambda,\, i})}}_{\lambda\in\Lambda,\, 1 \le i \le  d}}
\ee
where the $\frac{\partial}{\partial\log(X_\sigma)}$ \up{resp.,~$\frac{\partial}{\partial\log(X_{\lambda,\, i})}$} are as in \eqref{derivations} with $R_\Sigma^\square$ \up{resp.,~$R_\lambda^\square$} in place of $R$ and the Frobenius acts in degree $j$ on the right side as $p^j$ times the Frobenius of $D_{\Sigma,\, \Lambda}$ \up{compare with {\upshape\S\ref{log-dR-cx}}}.
\eprop

\bpf
By \S\ref{log-Acris}, each $A_{\Sigma,\, \Lambda,\, \cris}^\square/p^n$ is a log smooth thickening of $R/p$ over $A_\cris/p^n$. Therefore, by \cite{Bei13a}*{1.4, Remarks (ii)} (and the second paragraph of \S\ref{log-PD}), the log PD thickening $D_{\Sigma,\, \Lambda,\, n}$ of $R/p$ is PD smooth over $A_\cris/p^n$ (see the proof of \Cref{PD-smooth}). Thus, as in the proof of \Cref{RHS-rewrite} above, \cite{Bei13a}*{(1.8.1)} ensures that the logarithmic PD de Rham complex $\Omega^\bullet_{D_{\Sigma,\, \Lambda,\, n}/(A_\cris/p^n),\, \log,\, \PD}$ Frobenius-equivariantly computes $R\Gamma_{\log\cris}(\cO_{\fX_{\cO_C/p}/(A_{\cris}/p^n)})$. By \cite{Bei13a}*{1.7,~Exercises, (i)}, 
\be \lab{PD-dif-id}
\Omega^\bullet_{D_{\Sigma,\, \Lambda,\, n}/(A_\cris/p^n),\, \log,\, \PD} \cong \Omega^{\bullet}_{(A_{\Sigma,\, \Lambda,\, \cris}^\square/p^n)/(A_\cris/p^n),\, \log} \tensor_{A_{\Sigma,\, \Lambda,\, \cris}^\square/p^n} D_{\Sigma,\, \Lambda,\, n},
\ee
so \eqref{D-compute-id} follows. Then, since each $R/p^n$ is a log smooth log PD thickening of $R/p$ over $\cO_C/p^n$, analogous reasoning applies to $R\Gamma_{\log\cris}(\cO_{\fX_{\cO_C/p}/\cO_C}) \overset{\text{\cite{Bei13a}*{(1.8.1)}}}{\cong} R\Gamma_{\log\dR}(\fX/\cO_C)$ (compare with the proof of \Cref{RHS-rewrite}) and gives \eqref{cris-dR-map}.

Finally, the identification \eqref{D-Koszul} results from \eqref{D-compute-id} and the Frobenius-equivariant identifications 
\[
\tst \Omega^{\bullet}_{(A_{\Sigma,\, \Lambda,\, \cris}^\square/p^n)/(A_\cris/p^n),\, \log} \cong K_{A_{\Sigma,\, \Lambda,\, \cris}^\square/p^n}\p{\p{\frac{\partial}{\partial\log(X_\sigma)}}_{\sigma\in\Sigma}, \p{\frac{\partial}{\partial\log(X_{\lambda,\, i})}}_{\lambda\in\Lambda,\, 1 \le i \le  d}}
\]
supplied by \eqref{log-Koszul}.
\epf

\brem \lab{cris-dR-bc-rem}
By \cite{Bei13a}*{(1.11.1)}, the first map in \eqref{cris-dR-map} induces the identification 
\be \lab{cris-dR-bc}
R\Gamma_{\log\cris}(\cO_{\fX_{\cO_C/p}/A_{\cris}}) \tensor^{\bL}_{A_\cris} \cO_C/p \cong R\Gamma_{\log\dR}(\fX/\cO_C) \tensor^\bL_{\cO_C} \cO_C/p
\ee
in the derived category, so the same holds for the second map:
\[
\Omega^\bullet_{D_{\Sigma,\, \Lambda}/A_\cris,\, \log,\, \PD} \tensor^{\bL}_{A_\cris} \cO_C/p \cong  \Omega^\bullet_{(R/p)/(\cO_C/p),\, \log}.
\]
\erem

\ready{To make the identification \eqref{D-Koszul} analogous to the identification in \Cref{crys-limit}, we will express $D_{\Sigma,\, \Lambda}$ as a completed direct limit of rings $D_{\Sigma,\, \Lambda}^{(m)}$ that ``are generated by divided powers of degree at most $m$,'' see \eqref{D-Dm}. For this, we will build on the ideas from the proof of \cite{Kat89}*{4.10~(1)} to identify  $D_{\Sigma,\, \Lambda}$ with the $p$-adic completion of the (non log) divided power envelope of an \emph{exact} closed immersion in \Cref{D-no-log}.\footnote{\lab{no-env-foot}\ready{The arguments below would become more direct if we could ``uncomplete'' $D_{\Sigma,\, \Lambda}$ by constructing the log PD envelope of the (possibly nonexact) log closed immersion $\Spec(R/p) \hra \Spec(A_{\Sigma,\, \Lambda}^\square)$. Neither \cite{Kat89}*{5.4} nor \cite{Bei13a}*{1.3, Theorem} gives this hypothetical envelope because $p$ is not nilpotent in $A_{\Sigma,\, \Lambda}^\square$.}} This will also make $D_{\Sigma,\, \Lambda}$ more explicit and easier to analyze.}

\bpp[A chart for $A_{\Sigma,\, \Lambda}^\square$] \lab{chart-target}
\ready{To express $D_{\Sigma,\, \Lambda}$ as the $p$-adic completion of a usual (non log) divided power envelope, in \S\S\ref{chart-target}--\ref{conv-chart} we build a chart for the (fine version) of the log closed~immersion 
\be \lab{ci-chart}
\Spec(R/p) \hra \Spec(A_{\Sigma,\, \Lambda}^\square).
\ee
For this, we fix the unique $q \in \bQ_{> 0}$ for which
\[
\tst \bZ \cdot q = \Sigma_{\lambda \in \Lambda}\,  \bZ \cdot q_\lambda  \q \text{inside} \q \bQ,
\]
so that $\f{q_\lambda}{q} \in \bZ_{> 0}$ for every $\lambda$ (and even $q_\lambda = q$ in the case when $R \tensor_{\cO_C} k$ is not $k$-smooth, see \S\ref{all-possible-coordinates}). 
We endow $\cO_C/p$ and $A_\Inf$ with the compatible via $\theta$ fine log structures determined by 
\[
\bN_{\ge 0} \ra \cO_C/p \q \text{with} \q 1 \mapsto p^q \qq \text{and} \qq \bN_{\ge 0} \ra A_\Inf \q \text{with} \q 1 \mapsto [(p^{1/p^\infty})^q].
\]
For each $\lambda \in \Lambda$, we consider the submonoid
\[
\tst Q_\lambda \subset \f{q}{q_\lambda} \prod_{0 \le i \le r_\lambda} \bN_{\ge 0} \qq\text{generated by} \qq \prod_{0 \le i \le r_\lambda} \bN_{\ge 0} \q \text{and the diagonal} \q (\f{q}{q_\lambda}, \dotsc, \f{q}{q_\lambda}),
\]
so that the chart 
\[
\tst Q_\lambda \ra A(R_\lambda^\square) \qq \text{given by} \qq \prod_{0 \le i \le r_\lambda} \bN_{\ge 0} \xra{(n_i)\, \mapsto\, \prod X_{\lambda,\, i}^{n_i}} A(R_\lambda^\square) \q \text{and} \q (\f{q}{q_\lambda}, \dotsc, \f{q}{q_\lambda}) \mapsto [(p^{1/p^\infty})^q]
\]
makes $\Spec(A(R_\lambda^\square))$ a fine log $\Spec(A_\Inf)$-scheme.} We let 
\[
\tst Q \ce \p{\prod_{\lambda \in \Lambda} Q_\lambda} / ((\f{q}{q_{\lambda_1}}, \dotsc, \f{q}{q_{\lambda_1}}) = (\f{q}{q_{\lambda_2}}, \dotsc, \f{q}{q_{\lambda_2}}) )_{\lambda_1 \neq \lambda_2}
\]
 be the quotient monoid obtained by identifying the diagonal elements $(\f{q}{q_\lambda}, \dotsc, \f{q}{q_\lambda})$, so that the map
\[
Q \ra A_{\Sigma,\, \Lambda}^\square \qq \text{that results from the charts} \qq Q_\lambda \ra A(R_\lambda^\square)
\]
is a chart for the target $\Spec(A_{\Sigma,\, \Lambda}^\square)$ of a fine version of the log closed immersion \eqref{ci-chart}. In terms of this chart, the Frobenius of $A_{\Sigma,\, \Lambda}^\square$ multiplies each element of $Q$ by $p$ (see \S\ref{log-A(R)}).
\epp

\bpp[A convenient chart in the smooth case] \lab{conv-chart-sm}
\ready{We consider the case when $R \tensor_{\cO_C} k$ is $k$-smooth, so that for each $\lambda \in \Lambda$ there is a unique $0 \le  i_{\lambda} \le r_\lambda$ with $t_{\lambda,\, i_\lambda} \not\in R^\times$, and build the monoids
\be \lab{P-sm-def}
\tst P_{\lambda_0} \ce \p{\bN_{\ge 0} \times \prod_{0 \le i \le r_{\lambda_0},\, i \neq i_{\lambda_0} } \bZ } \times \prod_{\lambda \neq \lambda_0} \p{ \p{ \prod_{0 \le i \le r_{\lambda}} \bZ }/\bZ } \qq \text{for} \qq \lambda_0 \in \Lambda,
\ee
where each $\bZ$ by which we quotient is embedded diagonally. For each $(\lambda, i)$ with $0 \le i \le r_\lambda$, 
\be \lab{vli-def}
t_{\lambda,\, i} = (p^{q})^{n_{\lambda,\, i}}\cdot  v_{\lambda,\, i} \q \text{in} \q R \qq \text{for unique}  \q n_{\lambda,\, i} \in \bZ_{\ge 0} \q \text{and} \q v_{\lambda,\, i} \in R^\times;
\ee
explicitly, $n_{\lambda,\, i_\lambda} = \f{q_\lambda}{q}$ and $n_{\lambda,\, i} = 0$ for $i \neq i_\lambda$. In particular, $\prod_{0 \le i \le r_\lambda} v_{\lambda,\, i} = 1$ for each $\lambda$.  The map 
\[
P_{\lambda_0} \ra R/p \qq \text{given by} \qq \bN_{\ge 0} \ni 1 \mapsto p^q, \q \bZ_{(\lambda, i)} \ni 1 \mapsto v_{\lambda, i},
\]
where the subscript $(\lambda, i)$ indicates the factor $\bZ$ of \eqref{P-sm-def} being considered, is a chart for the source $\Spec(R/p)$ of a fine version of the log closed immersion \eqref{ci-chart}. In terms of this chart, the Frobenius of $R/p$ multiplies each element of $P_{\lambda_0}$ by $p$.

Due to \eqref{vli-def}, knowing the indices $i_\lambda$, we may evidently express the image of every generator of $Q$ under $Q \ra A_{\Sigma,\, \Lambda}^\square \surjects R/p$ in terms of the images of elements of $P_{\lambda_0}$ without knowing the ``values'' of these images. Thus, the log closed immersion \eqref{ci-chart} has a natural Frobenius-equivariant chart
\[
\tst Q \ra P_{\lambda_0} =  \p{\bN_{\ge 0} \times \prod_{0 \le i \le r_{\lambda_0},\, i \neq i_{\lambda_0} } \bZ } \times \prod_{\lambda \neq \lambda_0} \p{ \p{ \prod_{0 \le i \le r_{\lambda}} \bZ }/\bZ}
\]
that, for instance, sends $1 \in (\bN_{\ge 0})_{(\lambda_0,\, i_{\lambda_0})}$ to the element $(\f{q_{\lambda_0}}{q}, -1, \dotsc, -1)$ of $\bN_{\ge 0} \times \prod_{0 \le i \le r_{\lambda_0},\, i \neq i_{\lambda_0} } \bZ $, each $(\f{q}{q_\lambda}, \dotsc, \f{q}{q_\lambda})$ to $1 \in \bN_{\ge 0}$, each $1 \in (\bN_{\ge 0})_{(\lambda, i)}$ with $i \neq i_\lambda$ to $1 \in \bZ_{(\lambda, i)}$, etc. 

More precisely, the resulting $A_{\Sigma,\, \Lambda}^\square$-algebra 
\[
A_{\Sigma,\, \Lambda}^\square \tensor_{\bZ[Q]} \bZ[P_{\lambda_0}]
\]
comes equipped with an $A_{\Sigma,\, \Lambda}^\square$-semilinear Frobenius and is initial among the $A_{\Sigma,\, \Lambda}^\square$-algebras $B$ equipped with a unit $V_{\lambda,\, i} \in B^\times$ for each $(\lambda, i)$ with $0 \le i \le r_{\lambda}$ subject to the relations
\be \lab{chart-rels}
\tst  X_{\lambda,\, i} = [((p^{1/p^\infty})^q)^{n_{\lambda,\, i}}] \cdot V_{\lambda,\, i}, \qq \tst \prod_{0 \le i \le r_{\lambda}} V_{\lambda,\, i} = 1.
\ee
In particular, 
\be \lab{R-nat-alg-1}
R \q \text{is naturally an} \q  \text{$(A_{\Sigma,\, \Lambda}^\square \tensor_{\bZ[Q]} \bZ[P_{\lambda_0}])$-algebra} \q \text{(with $V_{\lambda,\, i} = v_{\lambda,\, i}$).}
\ee
A fine version of the log closed immersion \eqref{ci-chart} factors Frobenius-equivariantly as follows:
\be \lab{ci-fact-0}
\xymatrix@C=18pt{
\tst \Spec(R/p) \ar@{^(->}[r]^-{j_{\lambda_0}} & \Spec\p{A_{\Sigma,\, \Lambda}^\square \tensor_{\bZ[Q]} \bZ[P_{\lambda_0}]} \ar[r]^-{q_{\lambda_0}} & \Spec(A_{\Sigma,\, \Lambda}^\square),
}
\ee
where $\Spec\p{A_{\Sigma,\, \Lambda}^\square \tensor_{\bZ[Q]} \bZ[P_{\lambda_0}]}$ is equipped with the log structure determined by $P_{\lambda_0}$. By construction, $j_{\lambda_0}$ is an exact closed immersion and, by \cite{Kat89}*{3.5}, the projection $q_{\lambda_0}$ is log \'{e}tale.}

The relations \eqref{chart-rels} do not depend on the choice of $\lambda_0$, so neither does the factorization \eqref{ci-fact-0}. More precisely, for another $\lambda_0' \in \Lambda$, we have the a natural isomorphism over $Q$ of charts for $R/p$:
\be \lab{P-P-iso}
\tst P_{\lambda_0} \isomto  P_{\lambda'_0},
\ee
which gives rise to the vertical Frobenius-equivariant isomorphism in the commutative diagram
\be \lab{indep-0} \ba
\xymatrix@R=-4pt{
& \Spec\p{A_{\Sigma,\, \Lambda}^\square \tensor_{\bZ[Q]} \bZ[P_{\lambda_0}]} \ar[rd]^-{q_{\lambda_0}} & \\
\Spec(R/p) \ar@{^(->}[ru]^-{j_{\lambda_0}} \ar@{^(->}[rd]_-{j_{\lambda'_0}} & & \Spec(A_{\Sigma,\, \Lambda}^\square). \\
& \Spec\p{A_{\Sigma,\, \Lambda}^\square \tensor_{\bZ[Q]} \bZ[P_{\lambda'_0}]} \ar[uu]^(.4)*[@]-{\sim} \ar[ru]_-{q_{\lambda'_0}} & 
}
\ea \ee
\epp

\bpp[A convenient chart in the nonsmooth case] \lab{conv-chart}
\ready{We now consider the case when $R \tensor_{\cO_C} k$ is not $k$-smooth, so that $q_\lambda = q$ and $Q_\lambda \cong \prod_{0 \le i \le r_\lambda} \bN_{\ge 0}$ for every $\lambda \in \Lambda$. Letting $\Delta_\lambda \subset Q_\lambda$ be the diagonal copy of $\bN_{\ge 0}$, we can then describe the chart $Q$ for a fine version of $\Spec(A_{\Sigma,\, \Lambda}^\square)$ as follows:
\[
\tst Q \cong \p{\prod_{\lambda\in \Lambda}\p{ \prod_{0 \le i \le r_\lambda} \bN_{\ge 0}}}/\p{\Delta_{\lambda_1} = \Delta_{\lambda_2}}_{\lambda_1 \neq \lambda_2}.
\]
By \S\ref{all-possible-coordinates}, each $t_{\lambda,\, i}\not \in R^\times$ cuts out a unique irreducible component $\ov{\{y_{\lambda,\, i}\}}$ of $\Spec(R \tensor_{\cO_C} k)$. Its generic point $y_{\lambda,\, i}$ determines the ideal $(t_{\lambda,\, i}) \subset R$: indeed, $(p^q) \subset (t_{\lambda,\, i})$ in $R$ and the ideal $(t_{\lambda,\, i})/(p^q) \subset R/(p^q)$ is the kernel of the localization map $R/(p^q) \ra (R/(p^q))_{y_{\lambda,\, i}}$, as may be seen over $R_\gL^\square$. Conversely, for each generic point $y$ of $\Spec(R \tensor_{\cO_C} k)$ and $\lambda \in \Lambda$, a unique $t_{\lambda,\, i_{\lambda}(y)}$ with $0 \le i_{\lambda}(y) \le r_\lambda$ cuts out $\ov{\{y\}}$ (see \S\ref{all-possible-coordinates}). Consequently, for each $y$ and $\lambda, \lambda_0 \in \Lambda$,
\be \lab{a-def}
t_{\lambda,\, i_\lambda(y)} = u_{\lambda,\, \lambda_0,\, y} \cdot t_{\lambda_0,\, i_{\lambda_0}(y)} \q \text{in} \q R \qq \text{for a unique} \q u_{\lambda,\, \lambda_0,\, y} \in R^\times.
\ee
Letting $\cY$ denote the set of the generic points of $\Spec(R \tensor_{\cO_C} k)$, for $\lambda_0 \in \Lambda$ we build the monoid
\be \lab{P-sst-def}
\tst P_{\lambda_0} \ce \p{\p{\prod_{\cY} \bN_{\ge 0} \times \prod_{\{0 \le i \le r_{\lambda_0}\} \setminus i_{\lambda_0}(\cY)} \bZ} \times \prod_{\lambda \neq \lambda_0} \p{  \prod_{0 \le i \le r_\lambda} \bZ}}\big/(\Delta_\lambda = \Delta_{\lambda_0} )_{\lambda \neq \lambda_0},
 \ee
 where the quotient means that for every $\lambda \neq \lambda_0$ we are identifying every diagonal element of $\prod_{0 \le i \le r_\lambda} \bZ$ with the corresponding diagonal element of $\prod_{\{0 \le i \le r_{\lambda_0}\} \setminus i_{\lambda_0}(\cY)} \bZ$ (interpreted to be $0$ if the indexing set is empty). The assignment (as in \S\ref{conv-chart-sm}, subscripts indicate factors in \eqref{P-sst-def})
\[
 (\bN_{\ge 0})_{y} \ni 1 \mapsto t_{\lambda_0,\, i_{\lambda_0}(y)}, \qq \bZ_{(\lambda,\, i)} \ni 1 \mapsto t_{\lambda,\, i} \ \  \text{if} \ \ i \not\in i_\lambda(\cY), \qq \bZ_{(\lambda,\, i)} \ni 1 \mapsto u_{\lambda,\, \lambda_0,\, y} \ \ \text{if} \ \ i = i_\lambda(y)
\]
determines a chart 
\[
 \tst P_{\lambda_0}  \ra R/p
 \]
for the source $\Spec(R/p)$ of a fine version of the log closed immersion \eqref{ci-chart}. In terms of this chart, the Frobenius of $R/p$ multiplies each element of $P_{\lambda_0}$ by $p$.

Due to the relation \eqref{a-def}, the images in $R/p$ of the generators of $Q$ are evidently expressible in terms of the images of the elements of $P_{\lambda_0}$ (without knowing the ``values'' of these images), so, as in the smooth case, there is a natural Frobenius-equivariant chart for a fine version of \eqref{ci-chart}:
\[
Q \ra P_{\lambda_0}
\]
that, for instance, for $\lambda \neq \lambda_0$ and $y \in \cY$, sends $1 \in (\bN_{\ge 0})_{(\lambda,\, i_\lambda(y))}$ to $(1, 1) \in (\bN_{\ge 0})_{y} \times \bZ_{(\lambda,\, i_\lambda(y))}$.

The resulting $A_{\Sigma,\, \Lambda}^\square$-algebra 
\[
A_{\Sigma,\, \Lambda}^\square \tensor_{\bZ[Q]} \bZ[P_{\lambda_0}]
\]
comes equipped with an $A_{\Sigma,\, \Lambda}^\square$-semilinear Frobenius endomorphism and is initial among the $A_{\Sigma,\, \Lambda}^\square$-algebras $B$ for which $X_{\lambda,\, i} \in B^\times$ when $ i \not\in i_\lambda(\cY)$ and that are equipped with, for each $y \in \cY$ and $\lambda\in \Lambda$, a unit $U_{\lambda,\, \lambda_0,\, y} \in B^\times$ subject to the relations 
\be \lab{nonsm-rels}
\ba\tst X_{\lambda,\, i_\lambda(y)} &= U_{\lambda,\, \lambda_0,\, y} \cdot X_{\lambda_0,\, i_{\lambda_0}(y)}, \qq U_{\lambda_0,\, \lambda_0,\, y} = 1, \qq \text{and} \\ \tst \prod_{y \in \cY} U_{\lambda,\, \lambda_0,\, y} &=\tst  \p{\prod_{\{0 \le i \le r_{\lambda_0}\}\setminus i_{\lambda_0}(\cY)  } X_{\lambda_0,\, i}} \big/ \p{\prod_{\{0 \le i \le r_\lambda\}\setminus i_{\lambda}(\cY)  } X_{\lambda,\, i}} \q \text{for}\q \lambda \in \Lambda.
 \ea
\ee
For a $\lambda_0' \in \Lambda$, we may set $U_{\lambda,\, \lambda_0',\, y} = U_{\lambda,\, \lambda_0,\, y} \cdot U_{\lambda_0',\, \lambda_0,\, y}\i$ to express the $U_{\lambda,\, \lambda_0',\, y}$ in terms of the $U_{\lambda,\, \lambda_0,\, y}$, so, up to a canonical $A_{\Sigma,\, \Lambda}^\square$-isomorphism, $A_{\Sigma,\, \Lambda}^\square \tensor_{\bZ[Q]} \bZ[P_{\lambda_0}]$ does not depend on $\lambda_0$. Moreover, 
\be \lab{R-nat-alg-2}
R \q  \text{is naturally an} \q  \text{$(A_{\Sigma,\, \Lambda}^\square\tensor_{\bZ[Q]} \bZ[P_{\lambda_0}])$-algebra} \q  \text{(with $U_{\lambda,\, \lambda_0,\, y} = u_{\lambda,\, \lambda_0,\, y}$).}
\ee
}As in the smooth case, we equip $\Spec\p{A_{\Sigma,\, \Lambda}^\square \tensor_{\bZ[Q]} \bZ[P_{\lambda_0}]}$ with the log structure determined by $P_{\lambda_0}$, so a fine version of the log closed immersion \eqref{ci-chart} factors Frobenius-equivariantly as follows:
\be \lab{ci-fact}
\xymatrix@C=18pt{
\tst \Spec(R/p) \ar@{^(->}[r]^-{j_{\lambda_0}} & \Spec\p{A_{\Sigma,\, \Lambda}^\square \tensor_{\bZ[Q]} \bZ[P_{\lambda_0}]} \ar[r]^-{q_{\lambda_0}} & \Spec(A_{\Sigma,\, \Lambda}^\square),
}
\ee
where $j_{\lambda_0}$ is an exact closed immersion and, by \cite{Kat89}*{3.5}, the projection $q_{\lambda_0}$ is log \'{e}tale. As in \S\ref{conv-chart-sm}, we have natural isomorphisms $P_{\lambda_0} \simeq P_{\lambda_0'}$ over $Q$ and the compatibility diagram \eqref{indep-0}.
\epp

We now use the charts $Q \ra P_{\lambda_0}$ to build a (non log) PD envelope whose $p$-adic completion is $D_{\Sigma,\, \Lambda}$.

\bpp[The divided power envelope of $j_{\lambda_0}$] \lab{D-j-def}
\ready{For $\lambda_0 \in \Lambda$, we let $D_{j_{\lambda_0}}$ be the divided power envelope  over $(\bZ_p, p\bZ_p)$ of the closed immersion $j_{\lambda_0}$ defined in \eqref{ci-fact-0} and \eqref{ci-fact}. The universal property of $A_\cris^0$ (see \S\ref{Acris-ring}) identifies $D_{j_{\lambda_0}}$ with the divided power envelope of the closed immersion
\[
j_{\lambda_0,\, \cris} \colon \Spec(R/p) \hra \Spec((A_{\Sigma,\, \Lambda}^\square \tensor_{A_\Inf} A_\cris^0) \tensor_{\bZ[Q]} \bZ[P_{\lambda_0}] ) \q \text{over} \q \Spec(\cO_C/p) \hra \Spec(A_\cris^0)
\]
(compare with \S\ref{log-PD}).
Since $j_{\lambda_0}$ underlies an \emph{exact} closed immersion of log schemes, we may, in addition, identify $D_{j_{\lambda_0}}$ endowed with the log structure determined by $P_{\lambda_0}$ with the log PD envelope of $j_{\lambda_0}$ over $\bZ_p$, or of $j_{\lambda_0,\, \cris}$ over $A_\cris^0$ (compare with \cite{Kat89}*{5.5.1}).} For $\lambda_0' \in \Lambda$, the vertical isomorphism in \eqref{indep-0} induces an isomorphism
\be \lab{D-lambda-pr}
D_{j_{\lambda_0}} \cong D_{j_{\lambda_0'}}.
\ee
By functoriality, $D_{j_{\lambda_0}}$ comes equipped with an $A_\cris^0$-semilinear Frobenius endomorphism, and the isomorphisms \eqref{D-lambda-pr} are Frobenius equivariant. Due to \eqref{R-nat-alg-1} and \eqref{R-nat-alg-2}, there is a map
\be \lab{D-j-map-to-R}
D_{j_{\lambda_0}} \ra R \qq \text{that lifts} \qq D_{j_{\lambda_0}} \surjects R/p;
\ee
its  formation is compatible with the isomorphisms \eqref{D-lambda-pr}.

\ready{By the universal property of $A_{\Sigma,\, \Lambda}^\square\tensor_{\bZ[Q]} \bZ[P_{\lambda_0}]$ (see \eqref{chart-rels} and \eqref{nonsm-rels}), the continuous $\Delta_{\Sigma,\, \Lambda}$-action on $A_{\Sigma,\, \Lambda}^\square$ extends to a $(p, \mu)$-adically continuous $\Delta_{\Sigma,\, \Lambda}$-action on $A_{\Sigma,\, \Lambda}^\square\tensor_{\bZ[Q]} \bZ[P_{\lambda_0}]$, so it induces an $A_\cris^0$-linear $\Delta_{\Sigma,\, \Lambda}$-action on $D_{j_{\lambda_0}}$. As an $(A_{\Sigma,\, \Lambda}^\square\tensor_{\bZ[Q]} \bZ[P_{\lambda_0}])$-algebra,  $D_{j_{\lambda_0}}$ is generated by the divided powers of the elements of the ideal of $j_{\lambda_0}$, so this action is $p$-adically continuous.}
\epp

\blem \lab{D-no-log}
For $\lambda_0 \in \Lambda$, the map $q_{\lambda_0}$ induces Frobenius- and $\Delta_{\Sigma,\, \Lambda}$-equivariant isomorphisms
\be \lab{D-no-log-eq}
D_{\Sigma,\, \Lambda,\, n} \cong D_{j_{\lambda_0}}/p^n \q \text{for} \q n > 0 \qq \text{{\upshape(}resp.,} \q D_{\Sigma,\, \Lambda} \cong \wh{D_{j_{\lambda_0}}})
\ee
that are $A_\cris$-linear and compatible with divided powers, maps to $R/p^n$ \up{resp.,~$R$\uscolon see \eqref{R-D-map} and \eqref{D-j-map-to-R}}, and the isomorphisms \eqref{D-lambda-pr}. In particular, $D_{\Sigma,\, \Lambda}$ is $p$-adically complete,
\[
D_{\Sigma,\, \Lambda}/p^n \isomto D_{\Sigma,\, \Lambda,\, n} \q \text{for} \q n > 0
\]
and the $\Delta_{\Sigma,\, \Lambda}$-action on $D_{\Sigma,\, \Lambda}$ is $p$-adically continuous.
\elem

\bpf
\ready{We may identify $D_{\Sigma,\, \Lambda,\, n}$ with the log PD envelope of
\[
\Spec(R/p) \hra \Spec(A_{\Sigma,\, \Lambda,\, \cris}^\square/p^n) \qq \text{over} \qq \Spec(\cO_C/p) \hra \Spec(A_\cris/p^n)
\]
defined using fine log structures (see \S\ref{log-PD}). On the other hand, we may identify $D_{j_{\lambda_0}}/p^n$ with the (log) divided power envelope of $j_{\lambda_0,\, \cris} \tensor_{A^0_\cris}A_\cris^0/p^n$, that is, of
\[
\Spec(R/p) \hra \Spec((A_{\Sigma,\, \Lambda,\, \cris}^\square/p^n) \tensor_{\bZ[Q]} \bZ[P_{\lambda_0}]), \qq \text{over} \qq \Spec(\cO_C/p) \hra \Spec(A_\cris/p^n)
\]
(see \S\ref{D-j-def} and \cite{SP}*{\href{http://stacks.math.columbia.edu/tag/07HB}{07HB}}). Consider a commutative square
\be\ba \lab{lift-?}
\xymatrix{
T_0 \ar@{^(->}[d] \ar[r] & \Spec((A_{\Sigma,\, \Lambda,\, \cris}^\square/p^n) \tensor_{\bZ[Q]} \bZ[P_{\lambda_0}]) \ar[d]^-{q_{\lambda_0} \tensor_{A_\Inf}A_\cris/p^n} \\ 
T \ar[r] \ar@{-->}[ru]^-{?} & \Spec(A_{\Sigma,\, \Lambda,\, \cris}^\square/p^n)
}
\ea \ee
of log schemes over $A_\cris/p^n$ in which $T_0 \hra T$ is a log PD thickening such that the log structure $\cN_T$ of $T$ (and hence also $\cN_{T_0}$ of $T_0$) is integral and quasi-coherent and the log structures of $\Spec((A_{\Sigma,\, \Lambda,\, \cris}^\square/p^n) \tensor_{\bZ[Q]} \bZ[P_{\lambda_0}])$ and $\Spec(A_{\Sigma,\, \Lambda,\, \cris}^\square/p^n)$ are determined by the charts $P_{\lambda_0}$ and $Q$, respectively (see \S\S\ref{chart-target}--\ref{conv-chart}). By \cite{Bei13a}*{1.1~Exercises~(iii)}, for any $t, t' \in \Gamma(T, \cN_T)$ and $u_0 \in \cO_{T_0}^\times$ with $t|_{T_0} = u_0 \cdot t'|_{T_0}$, there exists a unique lift $u \in \cO_{T}^\times$ of $u_0$ such that $t = ut'$.} 
Thus, by the construction of $P_{\lambda_0}$ and the universal property described by the equations \eqref{chart-rels} and \eqref{nonsm-rels}, there is a unique log morphism indicated by the dashed arrow in \eqref{lift-?} that makes the diagram commute. Consequently, $q_{\lambda_0}$ induces an isomorphism between the log PD envelopes: 
\[
D_{j_{\lambda_0}}/p^n \cong D_{\Sigma,\, \Lambda,\, n}, \qq \text{and, by letting $n$ vary, also} \qq \wh{D_{j_{\lambda_0}}} \cong D_{\Sigma,\, \Lambda}.
\]
Functoriality implies the claimed compatibilities, 
and \eqref{D-no-log-eq} implies the ``in particular'' claim.
\epf

We now use \Cref{D-no-log} to define the rings $D_{\Sigma,\, \Lambda}^{(m)}$ that are analogous to the rings $A_\cris^{(m)}(R)$ of \S\ref{AcrisR-m}.

\bpp[The rings $D^{(m)}_{\Sigma,\, \Lambda}$] \lab{D-m}
\ready{For $\lambda_0 \in \Lambda$, the divided powers 
of the elements of the ideal of $j_{\lambda_0}$ generate $D_{j_{\lambda_0}}$ as an $(A_{\Sigma,\, \Lambda}^\square \tensor_{\bZ[Q]} \bZ[P_{\lambda_0}])$-algebra. In turn, for a fixed $m \in \bZ_{\ge 1}$, the divided powers of degree at most $m$ 
generate a Frobenius-stable $(A_{\Sigma,\, \Lambda}^\square \tensor_{\bZ[Q]} \bZ[P_{\lambda_0}])$-subalgebra
\[
\tst D_{j_{\lambda_0}}^{(m)} \subset D_{j_{\lambda_0}}, \qq \text{and} \qq D_{j_{\lambda_0}} = \bigcup_{m \ge 1} D_{j_{\lambda_0}}^{(m)}.
\] 
Since $D_{j_{\lambda_0}}$ is naturally and Frobenius-semilinearly an $A_\cris^0$-algebra (see \S\ref{D-j-def}), $D_{j_{\lambda_0}}^{(m)}$ is naturally and Frobenius-semilinearly 
an algebra over the subring $A_\cris^{0,\, (m)} \subset A_\cris^0$ defined in \S\ref{Acris-m}. By \Cref{D-no-log}, 
\[
\text{the image} \q D_{\Sigma,\, \Lambda}^{0} \q \text{of}\q D_{j_{\lambda_0}} \q \text{in} \q D_{\Sigma,\, \Lambda} \overset{\eqref{D-no-log-eq}}{\cong} \wh{D_{j_{\lambda_0}}} \q \text{is Frobenius-stable and 
independent of $\lambda_0$,}
\]
and the same holds for the image $D_{\Sigma,\, \Lambda}^{0,\, (m)} \subset D_{\Sigma,\, \Lambda}^0$ of $D_{j_{\lambda_0}}^{(m)}$ in $D_{\Sigma,\, \Lambda}$. For $m \ge p$, the $p$-adic completion
\[
D_{\Sigma,\, \Lambda}^{(m)} \ce (D_{\Sigma,\, \Lambda}^{0,\, (m)})\wh{\ \ } \qq \text{is naturally an algebra over} \qq A_\cris^{(m)}
\]
and comes equipped with an $A_\cris^{(m)}$-semilinear Frobenius.

By \Cref{D-no-log}, the composition $D_{j_{\lambda_0}} \surjects D_{\Sigma,\, \Lambda}^0 \hra D_{\Sigma,\, \Lambda}$ induces an isomorphism modulo $p^n$, and hence so do both maps that comprise it. Thus, since  $D_{\Sigma,\, \Lambda}^0 = \bigcup_{m \ge p} D_{\Sigma,\, \Lambda}^{0,\, (m)}$, we have  
\be \lab{D-Dm}
\tst D_{\Sigma,\, \Lambda} \cong (D_{\Sigma,\, \Lambda}^0)\wh{\ \ } \cong  \p{\varinjlim_{m \ge p} D_{\Sigma,\, \Lambda}^{(m)}}\wh{\ \ } \qq \text{over} \qq A_\cris
\ee
compatibly with the Frobenii. In what follows, $D_{\Sigma,\, \Lambda}^0$ plays the role of the ring that underlies the hypothetical log PD envelope of the log closed immersion \eqref{ci-chart} (see \cref{no-env-foot}).}

By \Cref{D-no-log}, the $\Delta_{\Sigma,\, \Lambda}$-action on $D_{\Sigma,\, \Lambda}$ respects the subrings $D_{\Sigma,\, \Lambda}^{0,\, (m)} \subset D_{\Sigma,\, \Lambda}$. The induced $A_\cris^{(m)}$-linear $\Delta_{\Sigma,\, \Lambda}$-action on  $D_{\Sigma,\, \Lambda}^{(m)}$ is $p$-adically continuous and compatible as $m$ varies. The identifications in \eqref{D-Dm} are $\Delta_{\Sigma,\, \Lambda}$-equivariant.
\epp

\bpp[The derivations $\f{\del}{\del\log(X_\tau)}$] \lab{log-PD-der}
For brevity, let $\tau$ denote either the index ``$\sigma$'' for some $\sigma \in \Sigma$ or the index  ``$\lambda, i$'' for some $\lambda \in \Lambda$  and $i = 1, \ldots, d$. The log derivations $\f{\del}{\del\log(X_\tau)}$ defined in \eqref{derivations} with $R_{\Sigma}^\square$ or $R_\lambda^\square$ in place of $R$ give rise to the log $A_\Inf$-derivations
\be \lab{more-der}
\tst \frac{\partial}{\partial\log(X_\tau)}\colon A_{\Sigma,\, \Lambda}^\square \ra A_{\Sigma,\, \Lambda}^\square \ee
(as in \S\ref{log-dR-cx}, we do not explicate the accompanying homomorphisms from the log structure). These, in turn, induce the divided power $A_\cris$-derivations
\be \lab{log-PD-D}
\tst \frac{\partial}{\partial\log(X_\tau)} \colon D_{\Sigma,\, \Lambda} \ra  D_{\Sigma,\, \Lambda}
\ee
(compare with \Cref{D-compute} and its proof, especially, \eqref{PD-dif-id}), where a \emph{divided power $A_\cris$-derivation} $\del$ is required to satisfy $\del(x^{[m]}) = x^{[m - 1]}\del(x)$ for divided powers $x^{[m]}$ with $m \ge 1$, in addition to the $A_\cris$-linearity and the Leibniz rule.  
 
Since $q_{\lambda_0}$ is log \'{e}tale (see \S\S\ref{conv-chart-sm}--\ref{conv-chart}), the derivations \eqref{more-der} uniquely extend to log $A_\Inf$-derivations
\be \lab{LPD-0}
\tst \frac{\partial}{\partial\log(X_\tau)}\colon A_{\Sigma,\, \Lambda}^\square \tensor_{\bZ[Q]} \bZ[P_{\lambda_0}] \ra A_{\Sigma,\, \Lambda}^\square \tensor_{\bZ[Q]} \bZ[P_{\lambda_0}] \qq \text{for every} \qq \lambda_0 \in \Lambda.
\ee
These, in turn, induce divided power $A_\cris^0$-derivations (see \cite{SP}*{\href{http://stacks.math.columbia.edu/tag/07HW}{07HW}})
\be \lab{DL-der}
\tst  \frac{\partial}{\partial\log(X_\tau)} \colon D_{j_{\lambda_0}} \ra D_{j_{\lambda_0}}. 
\ee
By construction, the derivations \eqref{log-PD-D} and \eqref{DL-der} are compatible, so they induce 
$A_\cris^{(m)}$-derivations
\[
\tst \frac{\partial}{\partial\log(X_\tau)} \colon D_{\Sigma,\, \Lambda}^{(m)} \ra D_{\Sigma,\, \Lambda}^{(m)} \qq \text{for} \qq  m \ge p
\]
that are compatible as $m$ varies and recover \eqref{log-PD-D} under the identification $D_{\Sigma,\, \Lambda} \cong \p{\varinjlim D_{\Sigma,\, \Lambda}^{(m)}}\wh{\ \ }$. Consequently, we may reexpress the identification \eqref{D-Koszul} as the Frobenius-equivariant identification
\be \lab{Dm-Koszul}
\tst R\Gamma_{\log\cris}(\cO_{\fX_{\cO_C/p}/A_{\cris}}) \cong \p{\varinjlim_{m \ge p}\p{ K_{D^{(m)}_{\Sigma,\, \Lambda}}\p{\p{\f{\partial}{\partial \log (X_\tau)}}_\tau}}}\wh{\ \ },
\ee
where in degree $j$ of $K_{D^{(m)}_{\Sigma,\, \Lambda}}((\f{\partial}{\partial \log (X_\tau)})_\tau)$ the Frobenius acts as $p^j$ times the Frobenius of $D^{(m)}_{\Sigma,\, \Lambda}$.
\epp

\bpp[A functorial complex that computes $R\Gamma_{\log\cris}(\fX_{\cO_C/p}/A_\cris)$] \lab{functorial-cris}
\ready{For a fixed $R$, the formation of the rings $D_{\Sigma,\, \Lambda}$, $D_{j_{\lambda_0}}$, $D_{\Sigma,\, \Lambda}^0$, and $D_{\Sigma,\, \Lambda}^{(m)}$ and the morphisms $j_{\lambda_0}$ and $q_{\lambda_0}$ is compatible with enlarging $\Sigma$ and $\Lambda$, and the same holds for the identification \eqref{Dm-Koszul}. Thus, we may form the filtered direct limit over all the closed immersions \eqref{more-coord} for varying $\Sigma$ and $\Lambda$ to build the complex
\be \lab{source-colim}
\tst \varinjlim_{\Sigma,\, \Lambda}\p{ \p{\varinjlim_{m \ge p} \p{K_{D^{(m)}_{\Sigma,\, \Lambda}}\p{\p{\frac{\partial}{\partial\log(X_\sigma)}}_{\sigma\in\Sigma}, \p{\frac{\partial}{\partial\log(X_{\lambda,\, i})}}_{\lambda\in\Lambda,\, 1 \le i \le  d}}}}\wh{\ \ }},
\ee
where the direct limits and the $p$-adic completion are termwise. By construction, this complex comes equipped with an $A_\cris$-semilinear Frobenius endomorphism (see the end of \S\ref{log-PD-der}) and, by \eqref{Dm-Koszul}, in the derived category it is canonically and Frobenius-equivariantly identified with
\[
R\Gamma_{\log\cris}(\cO_{\fX_{\cO_C/p}/A_{\cris}}).
\]
Moreover, if $R'$ is a $p$-adically formally \'{e}tale $R$-algebra equipped with data as in \S\ref{all-possible-coordinates} for some sets $\Sigma'$ and $\Lambda'$, then the term indexed by $\Sigma$, $\Lambda$ (and by the closed immersion \eqref{more-coord}) of the direct limit \eqref{source-colim} maps to\footnote{\ready{One uses the universal properties described in \eqref{chart-rels} and \eqref{nonsm-rels} and keeps in mind the case when $R \tensor_{\cO_C} k$ is not $k$-smooth but $R' \tensor_{\cO_C} k$ is.}} the term indexed by  $\Sigma \cup \Sigma'$, $\Lambda \cup \Lambda'$ (and by a closed immersion of $\Spf(R')$) of the analogous direct limit for $R'$, compatibly with the transition maps in \eqref{source-colim} and with the Frobenius. Thus, the complex \eqref{source-colim} equipped with its Frobenius is functorial in $R$, and so is its identification with $R\Gamma_{\log\cris}(\cO_{\fX_{\cO_C/p}/A_{\cris}})$.}

Since the formation of the maps \eqref{cris-dR-map} is compatible with enlarging $\Sigma$ and $\Lambda$, and then also with replacing $R$ by $R'$, the map $R\Gamma_{\log\cris}(\cO_{\fX_{\cO_C/p}/A_{\cris}}) \ra R\Gamma_{\log\dR}(\fX/\cO_C)$ is identified with a map 
\[
\tst \varinjlim_{\Sigma,\, \Lambda}\p{ \p{\varinjlim_{m \ge p} \p{ K_{D^{(m)}_{\Sigma,\, \Lambda}}\p{\p{\frac{\partial}{\partial\log(X_\sigma)}}_{\sigma \in \Sigma}, \p{\frac{\partial}{\partial\log(X_{\lambda,\, i})}}_{\lambda\in\Lambda,\, 1 \le i \le  d}}}}\wh{\ \ }} \ra \Omega^\bullet_{\Spf(R)/\cO_C,\, \log}
\]
whose formation is compatible is compatible with replacing $R$ by $R'$.
\epp

Having constructed the functorial complexes \eqref{FAO-colim} and \eqref{source-colim}, we seek to exhibit a natural map between them and to prove that it is an isomorphism. These tasks, which will be completed in \S\ref{comp-map} and \Cref{alphaR-qiso}, are the last stepping stones to the proof of \Cref{abs-crys-comp} given in \S\ref{Acris-main-pf}. We begin with the following variants of \Cref{delta-exp} and \Cref{local-isom}.

\blem \lab{more-exp}
For $m \ge p^2$, the element $\delta_\tau \in \Delta_{\Sigma,\, \Lambda}$, where the index $\tau$ is either ``$\sigma$'' for some $\sigma \in \Sigma$ or ``$\lambda, i$'' for some $\lambda \in \Lambda$  and $i = 1, \ldots, d$ \up{see {\upshape\S\ref{more-cover}}}, acts on $D_{\Sigma,\, \Lambda}^{(m)}$ as the endomorphism
\be \lab{del-tau-sum}
\tst \sum_{n \ge 0} \f{(\log([\eps]))^n}{n!} (\f{\del}{\del \log(X_\tau)})^n,
\ee
where $\f{(\log([\eps]))^n}{n!}$ lies in $A_\cris^{(m)}$ and $p$-adically tends to $0$ \up{see {\upshape\S\ref{log-eps}}}.
\elem 

\bpf
Analogously to the proof of \Cref{delta-exp}, the series \eqref{del-tau-sum} \emph{a priori} defines an $A_\cris$-algebra endomorphism of $D_{\Sigma,\, \Lambda}$. Moreover, by \Cref{delta-exp}, the action of $\delta_\tau$ on the ring $A_{\Sigma,\, \Lambda,\, \cris}^\square$ defined in \eqref{Acris-fact} is given by \eqref{del-tau-sum}. Thus, due to the universal properties \eqref{chart-rels} and \eqref{nonsm-rels}, the same holds for the action of $\delta_\tau$ on $A_{\Sigma,\, \Lambda,\, \cris}^\square \wh{\tensor}_{\bZ[Q]} \bZ[P_{\lambda_0}] \cong ((A_{\Sigma,\, \Lambda}^\square \tensor_{A_\Inf} A_\cris^0) \tensor_{\bZ[Q]} \bZ[P_{\lambda_0}])\wh{\ \ }$ (see \eqref{LPD-0}), 
where the completion is $p$-adic. Then, by the universal property of $D_{j_{\lambda_0}}$ (see \S\ref{D-j-def}) and \eqref{D-no-log-eq}, the element $\delta_\tau$ acts on $D_{\Sigma,\, \Lambda}$, and hence also on $D_{\Sigma,\, \Lambda}^{0,\, (m)}$ and $D_{\Sigma,\, \Lambda}^{(m)}$, by the series \eqref{del-tau-sum}.
\epf

\bprop \lab{comp-map-0}
In the local setting of {\upshape\S\ref{all-possible-coordinates}}, for $m \ge p^2$, the additive morphisms
\be \lab{CM-0-eq}
\tst
 \p{D_{\Sigma,\, \Lambda}^{(m)} \xra{\frac{\partial}{\partial\log(X_\tau)}} D_{\Sigma,\, \Lambda}^{(m)}} \xra{\p{\id,\,\, \sum_{n\geq 1}\frac{(\log([\eps]))^n}{ n!}(\frac{\partial}{\partial\log (X_\tau)})^{n - 1}}}  \p{D_{\Sigma,\, \Lambda}^{(m)} \xra{\delta_\tau - 1} D_{\Sigma,\, \Lambda}^{(m)}}
\ee
of complexes in degree $0$ and $1$, where the index $\tau$ ranges over ``$\sigma$'' for $\sigma \in \Sigma$ and ``$\lambda, i$'' for $\lambda \in \Lambda$  and $i = 1, \ldots, d$, define a Frobenius-equivariant morphism \up{whose target is defined as in \eqref{eta-def}}
\[
\tst K_{D_{\Sigma,\, \Lambda}^{(m)}}\p{(\frac{\partial}{\partial\log(X_\sigma)})_{\sigma\in\Sigma}, (\frac{\partial}{\partial\log(X_{\lambda,\, i})})_{\lambda\in\Lambda,\, 1 \le i \le  d} } \ra \eta_{(\mu)}\p{K_{D_{\Sigma,\, \Lambda}^{(m)}}( (\delta_\sigma -1)_{\sigma\in\Sigma}, (\delta_{\lambda,\, i}-1)_{\lambda\in\Lambda,\, 1 \le i \le  d})},
\]
where in each degree $j$ of the source the Frobenius acts as $p^j$ times the Frobenius of $D_{\Sigma,\, \Lambda}^{(m)}$.
\eprop

\bpf
By \Cref{more-exp}, the morphism \eqref{CM-0-eq} is well defined. Moreover, the image of its degree $1$ component lies in $\mu \cdot D_{\Sigma,\, \Lambda}^{(m)}$ because, by \S\ref{log-eps},  $\frac{(\log([\eps]))^n}{\mu\cdot n!}$ lies in $A_\cris^{(m)}$ and $p$-adically tends to $0$. The rest of the claim then follows from the definitions \eqref{Koszul-def} and \eqref{eta-def}, granted that one argues the Frobenius-equivariance as in the proof of \Cref{local-isom}.
\epf

\Cref{comp-map-0} reduces the task of exhibiting a natural map from the complex \eqref{source-colim} to the complex \eqref{FAO-colim} to exhibiting a natural $\Delta_{\Sigma,\, \Lambda}$-equivariant ring morphism $D_{\Sigma,\, \Lambda}^{(m)} \ra \bA_\cris^{(m)}(R_{\Sigma,\, \Lambda,\, \infty})$. To build the latter, we will realize $\bA_\cris^{(m)}(R_{\Sigma,\, \Lambda,\, \infty})$ inside the following ring $\bA_\cris(R_{\Sigma,\, \Lambda,\, \infty})$.

\bpp[The ring $\bA_\cris(R_{\Sigma,\, \Lambda,\, \infty})$] \lab{AA-cris}
For an affinoid perfectoid $\Spa(R_\infty'[\f{1}{p}], R_\infty')$ over $\Spa(C, \cO_C)$ (such as the one with $R_\infty' = R_{\Sigma,\, \Lambda,\, \infty}$), we consider the $\bA_\Inf(R_\infty')$-subalgebra
\[
\tst \bA_\cris^0(R_\infty') \subset \bA_\Inf(R_\infty')[\f{1}{p}]
\]
generated by the elements $\f{\xi^n}{n!}$ for $n \ge 1$. Analogously to \S\ref{Acris-R}, by \cite{Tsu99}*{proof of A2.8}, we have
\[
\tst \bA_\cris^0(R_\infty') \cong (\bA_\Inf(R_\infty')[\f{T^n}{n!}]_{n \ge 1})/(T - \xi), \qq \text{so} \qq \bA_\cris^0(R_\infty') \cong \bA_\Inf(R_\infty') \tensor_{A_\Inf} A_\cris^0.
\]
Thus, analogously to \S\ref{Acris-R}, the ring $\bA_\cris^0(R_\infty')$ is identified with the divided power envelope of $(\bA_\Inf(R_\infty'), \Ker(\theta) + p
\cdot \bA_\Inf(R_\infty'))$ over $(\bZ_p, p\bZ_p)$, and
\[
\bA_\cris(R_\infty') \ce \bA_{\Inf}(R_\infty') \wh{\tensor}_{A_\Inf} A_\cris \qq \text{is identified with} \qq (\bA_\cris^0(R_\infty'))\wh{\ \ }.
\]
For an $m \in \bZ_{\ge 1}$, we let $\bA_\cris^{0,\, (m)}(R_\infty') \subset \bA_\cris^0(R_\infty')$ be the $\bA_\Inf(R_\infty')$-subalgebra generated by the elements $\f{\xi^n}{n!}$ with $n \le m$ (compare with \S\ref{Acris-m}). For a fixed $m$, the subalgebra 
\[
\tst \bA_\Inf(R_\infty')[\f{T^n}{n!}]_{m \ge n \ge 1} \subset \bA_\Inf(R_\infty')[\f{T^n}{n!}]_{n \ge 1} \subset (\bA_\Inf(R_\infty')[\f{1}{p}])[T]
\]
is described by explicit lower bounds on the ``$p$-adic valuations'' of the coefficients of $T^N$ for $N \ge 1$.
Thus, since the sequence $(p, \xi)$ is $\bA_\Inf(R_\infty')$-regular (compare with \Cref{A-Apr}), the quotient of $\bA_\Inf(R_\infty')[\f{T^n}{n!}]_{n \ge 1}$ by $\bA_\Inf(R_\infty')[\f{T^n}{n!}]_{m \ge n \ge 1}$ has no nonzero $(T - \xi)$-torsion. Consequently, 
\be \lab{Acris-0-m}
\tst \bA_\cris^{0,\, (m)}(R_\infty') \cong (\bA_\Inf(R_\infty')[\f{T^n}{n!}]_{m \ge n \ge 1})/(T - \xi), 
\ee
to the effect that
\[
\bA_\cris^{0,\, (m)}(R_\infty') \cong \bA_\Inf(R_\infty') \tensor_{A_\Inf} A_\cris^{0,\, (m)}.
\]
Thus, by letting the completion be $p$-adic if $m \ge p$ and $(p, \mu)$-adic if $m < p$, we obtain the following identification with the $A_\cris^{(m)}$-algebra $\bA_\cris^{(m)}(R_\infty')$ defined as in \eqref{more-Acrism}:
\be \lab{Acrism-comp}
(\bA_\cris^{0,\, (m)}(R_\infty'))\wh{\ \ } \cong \bA_\cris^{(m)}(R_\infty') \ce \bA_\Inf(R_\infty') \wh{\tensor}_{A_\Inf} A_\cris^{(m)}.
\ee
\epp

\bprop \lab{Acrism-inj}
For an affinoid perfectoid $\Spa(R_\infty'[\f{1}{p}], R_\infty')$ over $\Spa(C, \cO_C)$ and an $m \ge 1$, the following ring homomorphisms are injective\ucolon
\be \lab{AI-eq}
\tst \bA_\cris^{0,\, (m)}(R_\infty') \hra \bA_\cris^{(m)}(R_\infty') \hra \bA_\cris(R_\infty') \hra \bB_\dR^+(R_\infty') \ce (\bA_\Inf(R_\infty')[\f{1}{p}])\wh{\ \ }
\ee
where the completion is $\xi$-adic and the definition of the last map is explained in the proof. In particular,  the $A_\Inf$-algebras in \eqref{AI-eq} have no nonzero $\mu$-torsion.
\eprop


\bpf
The assertion about the $\mu$-torsion follows from the rest because $\mu/\xi$ is a unit in $\bB_\dR^+(R_\infty')$ (see \eqref{mu-intro}--\eqref{xi-intro}) and $\bB_\dR^+(R_\infty')$ inherits $\xi$-torsion freeness from $\bA_\Inf(R_\infty')$.

\ready{The sequence $(p, \xi)$ is $\bA_\Inf(R_\infty')$-regular and $\bA_\Inf(R_\infty')$ is $\xi$-adically separated (see \cite{SP}*{\href{http://stacks.math.columbia.edu/tag/090T}{090T}}), so the ring $\bA_\Inf(R_\infty')[\f{1}{p}]$ is also $\xi$-adically separated. Thus, we obtain the injection
\[
\tst \bA_\Inf(R_\infty')[\f{1}{p}] \hra \bB_\dR^+(R_\infty'), \qq \text{and hence also} \qq \bA_\cris^{0,\, (m)}(R_\infty') \hra \bB_\dR^+(R_\infty'),
\] 
which, in particular, allows us to assume that $m \ge p$. For varying $n \ge 0$, the $\bA_\Inf (R_\infty')$-submodules 
\[
\tst \Fil_n^0 \subset \bA_\cris^0(R_\infty') \qq \text{generated by the} \q \f{\xi^{n'}}{n'!} \q \text{for} \q n' \ge n
\]
form a decreasing filtration of $\bA_\cris^0(R_\infty')$ by ideals. By \cite{Tsu99}*{A2.9 (2)},\footnote{\ready{\lab{caveat}\emph{Loc.~cit.}~is written in a different setting, but its proof continues to work if $A$ there is replaced by our $R_\infty'$.}} each 
\be \lab{p-complete}
\bA_\cris^0(R_\infty')/\Fil_n^0 \qq \text{is $p$-torsion free and $p$-adically complete.}
\ee
Thus, the $p$-adic completions $\Fil_n \ce (\Fil_n^0)\wh{\ \ }$ form a decreasing filtration of $\bA_\cris(R_\infty')$ by ideals~with 
\be \lab{qq-fil-fil}
\bA_\cris (R_\infty') /\Fil_n \cong \bA_\cris^0 (R_\infty') /\Fil_n^0.
\ee
The $p$-torsion freeness also supplies a decreasing filtration modulo $p$:
\[
 \Fil_n^0/p\Fil_n^0 \hra \bA_\cris^0 (R_\infty')/p\bA_\cris^0 (R_\infty').
\]
The isomorphism $\bA_\cris^0(R_\infty') \cong (\bA_\Inf(R_\infty')[\f{T^n}{n!}]_{n \ge 1})/(T - \xi)$ gives the explicit description
\be \lab{BC-desc}
\tst \bA_\cris^0 (R_\infty')/p\bA_\cris^0 (R_\infty') \cong (R_\infty'^\flat/\xi^p)[Y_1, Y_2, \ldots]/(Y_1^p, Y_2^p, \ldots) \q \text{with} \q Y_j = \f{T^{p^j}}{(p^j)!}
\ee
(compare with \cite{BC09}*{9.4.1 (3)}), so the filtration $\{ \Fil_n^0/p\Fil_n^0 \}_{n \ge 0}$ is separated. Thus, since 
\[
\Fil_n^0/p\Fil_n^0 \cong \Fil_n/p\Fil_n \q \text{compatibly with} \q \bA_\cris^0 (R_\infty')/p\bA_\cris^0 (R_\infty') \cong \bA_\cris (R_\infty')/p\bA_\cris (R_\infty'),
\]
the $p$-adic separatedness of $\bA_\cris(R_\infty')$ and \eqref{p-complete} force the filtration $\{\Fil_n\}_{n \ge 0}$ to be separated,~too:
\[
\tst \bA_\cris (R_\infty') \hra \varprojlim_n\, (\bA_\cris (R_\infty') / \Fil_n) \overset{\eqref{qq-fil-fil}}{\cong}  \varprojlim_n\, (\bA_\cris^0 (R_\infty') / \Fil_n^0).
\] 
Moreover, we have the injection $\varprojlim_n\, (\bA_\cris^0 (R_\infty') / \Fil_n^0) \hra \bB_{\dR}^{+}(R_\infty')$ that results from the injections
\[
\tst \bA_\cris^0 (R_\infty') / \Fil_n^0 \hra (\bA_\cris^0 (R_\infty') / \Fil_n^0)[\f{1}{p}] \cong (\bA_\Inf(R_\infty')[\f{1}{p}])/\xi^n \cong \bB_{\dR}^{+}(R_\infty')/\xi^n.
\]
Thus, we obtain the desired natural injection $\bA_\cris(R_\infty') \hra \bB_{\dR}^{+}(R_\infty')$ of $\bA_\cris^0(R_\infty')$-algebras.

We turn to the remaining injection $\bA_\cris^{(m)}(R_\infty') \hra \bA_\cris(R_\infty')$. For $n \ge 0$, we consider the ideal
\[
\tst \Fil_n^{0,\, (m)}  \ce \bA_\cris^{0,\, (m)} (R_\infty') \bigcap \Fil_n^0   = \Ker\p{\bA_\cris^{0,\, (m)} (R_\infty') \ra (\bA_\Inf(R_\infty')[\f{1}{p}])/\xi^n} \subset \bA_\cris^{0,\, (m)} (R_\infty'),
\]
so that $\{ \Fil_n^{0,\, (m)} \}_{n \ge 0}$ forms a filtration of $\bA_\cris^{0,\, (m)}(R_\infty')$. Explicitly,  as an $\bA_\Inf(R_\infty')$-module, $\Fil_n^{0,\, (m)}$ is generated by the products $\f{\xi^{n_1}}{n_1!} \cdots \f{\xi^{n_s}}{n_s!}$ with $n_1 + \ldots + n_s \ge n$ and $0 \le n_i \le m$. By \eqref{p-complete}, each
\be \lab{Fil-m-tor}
\bA_\cris^{0,\, (m)} (R_\infty')/\Fil_n^{0,\, (m)} \qq  \text{is $p$-torsion free,}
\ee
so we again get the induced filtration modulo $p$:
\[
\Fil_n^{0,\, (m)}/p\Fil_n^{0,\, (m)} \hra \bA_\cris^{0,\, (m)}(R_\infty')/p\bA_\cris^{0,\, (m)}(R_\infty').
\]
As in the case of $\{\Fil_n^0/p\Fil_n^0\}_{n \ge 0}$, the analogous to \eqref{BC-desc} description of $\bA_\cris^{0,\, (m)}(R_\infty')/p\bA_\cris^{0,\, (m)}(R_\infty')$ supplied by the isomorphism \eqref{Acris-0-m} shows that the filtration $\{ \Fil_n^{0,\, (m)}/p\Fil_n^{0,\, (m)}\}_{n \ge 0}$ is separated.}

For each $n > 0$, there is a $j_n > 0$ such that $p^{j_n}$ kills 
\[
\bA_\cris^0(R_\infty')/(\bA_\cris^{0,\, (m)}(R_\infty') + \Fil_n^0)
\]
 (for instance, $j_n \ce \ord_p(n!)$ has this property). Consequently, $p^{j_n}$ kills the kernel of the map
\[
\tst \qq  \f{\bA_\cris^{0,\, (m)}(R_\infty')/\Fil^{0,\, (m)}_n }{p^j \cdot (\bA_\cris^{0,\, (m)}(R_\infty')/\Fil^{0,\, (m)}_n)} \ra \f{\bA_\cris^{0}(R_\infty')/\Fil_n^0}{p^j \cdot (\bA_\cris^{0}(R_\infty')/\Fil_n^0)} \qq \text{for each} \qq j > 0,
\]
so, for $j > j_n$, every element of this kernel is a multiple of $p^{j - j_n}$. The short exact sequences
\[
\tst \qq 0 \ra \f{\Fil^{0,\, (m)}_n}{p^j \cdot \Fil^{0,\, (m)}_n} \xra{\eqref{Fil-m-tor}} \f{\bA_\cris^{0,\, (m)}(R_\infty')}{p^j \cdot \bA_\cris^{0,\, (m)}(R_\infty')} \ra  \f{\bA_\cris^{0,\, (m)}(R_\infty')/\Fil^{0,\, (m)}_n }{p^j \cdot (\bA_\cris^{0,\, (m)}(R_\infty')/\Fil^{0,\, (m)}_n)} \ra 0
\]
then show that for each $n > 0$, every element of the kernel
\be \lab{ker-ker}
\tst \Ker(\bA_\cris^{(m)}(R_\infty') \ra \bA_\cris(R_\infty')) \cong \Ker\p{\varprojlim_{j > 0} (\bA_\cris^{0,\, (m)}(R_\infty')/p^j) \ra \varprojlim_{j > 0}(\bA_\cris^0(R_\infty')/p^j)}
\ee
lies in $\varprojlim_{j > 0} \p{\Fil^{0,\, (m)}_n/p^j\Fil^{0,\, (m)}_n}$. In particular, this kernel maps to $\Fil_n^{0,\, (m)}/p \subset \bA_\cris^{0,\, (m)}(R_\infty')/p$ for each $n > 0$. However, by the previous paragraph, $\bigcap_{n > 0}\, (\Fil_n^{0,\, (m)}/p) = 0$, so the kernel \eqref{ker-ker} lies in $p\cdot \bA_\cris^{(m)}(R_\infty')$. Since $\bA_\cris(R_\infty')$ has no nonzero $p$-torsion and $\bA_\cris^{(m)}(R_\infty')$ is $p$-adically separated, this implies that the map $\bA_\cris^{(m)}(R_\infty') \ra \bA_\cris(R_\infty')$ is injective, as desired.
\epf

The following lemma is the final step to building the desired map $D_{\Sigma,\, \Lambda}^{(m)} \ra \bA_\cris^{(m)}(R_{\Sigma,\, \Lambda,\, \infty})$ in \S\ref{comp-map}.

\blem \lab{non-log-env}
For $\lambda_0 \in \Lambda$, there is a $\Delta_{\Sigma,\,\Lambda}$-equivariant divided power morphism
\be \lab{NLE-eq}
D_{j_{\lambda_0}} \ra \bA_\cris^0(R_{\Sigma,\, \Lambda,\, \infty})
\ee
that is compatible with the isomorphisms $D_{j_{\lambda_0}} \cong D_{j_{\lambda_0'}}$ of \eqref{D-lambda-pr} and is Frobenius-equivariant.
\elem 

\bpf
\ready{By construction, $\bA_\Inf(R_{\Sigma,\, \Lambda,\, \infty})$ is an 
$A_{\Sigma,\, \Lambda}^\square$-algebra. Moreover, since $\bA_\Inf(R_{\Sigma,\, \Lambda,\, \infty})$ is $\xi$-adically complete with $\bA_\Inf(R_{\Sigma,\, \Lambda,\, \infty})/\xi \cong R_{\Sigma,\, \Lambda,\, \infty}$ (see \S\ref{more-Ainf}), if $t_{\lambda,\, i}$ is a unit in $R \subset R_{\Sigma,\, \Lambda,\, \infty}$, then $X_{\lambda,\, i}$ is a unit in $\bA_\Inf(R_{\Sigma,\, \Lambda,\, \infty})$ (see \eqref{AR-surj}). Thus, if $R \tensor_{\cO_C} k$ is $k$-smooth, then the equations \eqref{chart-rels} have a unique solution in $\bA_\Inf(R_{\Sigma,\, \Lambda,\, \infty})$, to the effect that, in this case, $\bA_\Inf(R_{\Sigma,\, \Lambda,\, \infty})$ is naturally and $\Delta_{\Sigma,\, \Lambda}$-equivariantly an $(A_{\Sigma,\, \Lambda}^\square \tensor_{\bZ[Q]} \bZ[P_{\lambda_0}])$-algebra, compatibly with the ``change of $\lambda_0$'' isomorphisms of \eqref{indep-0}, the maps \eqref{R-nat-alg-1} and \eqref{more-theta} to $R$ and $R_{\Sigma,\, \Lambda,\, \infty}$, and the Frobenii.

If $R \tensor_{\cO_C} k$ is not $k$-smooth, then, in the notation of \S\ref{conv-chart}, for each $m \ge 0$ and a generic point $y \in \cY$ of $\Spec(R \tensor_{\cO_C} k)$, the element $t^{1/p^m}_{\lambda_0,\, i_{\lambda_0}(y)}$ is not a zero divisor in $R_{\Sigma,\, \Lambda,\, \infty}$ and is a unit in $R_{\Sigma,\, \Lambda,\, \infty}[\f{1}{p}]$. Thus, since $R_{\Sigma,\, \Lambda,\, \infty}$ is integrally closed in $R_{\Sigma,\, \Lambda,\, \infty}[\f{1}{p}]$, we conclude from \eqref{a-def} that 
\[
t_{\lambda,\, i_\lambda(y)}^{1/p^m}/t^{1/p^m}_{\lambda_0,\, i_{\lambda_0}(y)} \in R_{\Sigma,\, \Lambda,\, \infty}^\times \qq \text{for every} \qq \lambda_0, \lambda \in \Lambda, \q m \ge 0.
\]
In other words, 
\[
t_{\lambda,\, i_\lambda(y)}^{1/p^m} = u_{\lambda,\, \lambda_0,\, y}^{(m)} \cdot  t^{1/p^m}_{\lambda_0,\, i_{\lambda_0}(y)} \qq \text{in} \qq R_{\Sigma,\, \Lambda,\, \infty} \qq \text{for a unique} \qq u_{\lambda,\, \lambda_0,\, y}^{(m)} \in R_{\Sigma,\, \Lambda,\, \infty}^\times.
\]
}By the uniqueness, $(u^{(m + 1)}_{\lambda,\, \lambda_0,\, y})^p = u^{(m)}_{\lambda,\, \lambda_0,\, y}$, so  $u^\flat_{\lambda,\, \lambda_0,\, y} \ce  (u^{(m)}_{\lambda,\, \lambda_0,\, y})_{m \ge 0} \in (R_{\Sigma,\, \Lambda,\, \infty}^\flat)^\times$ satisfies
\[
X_{\lambda,\, i_{\lambda}(y)} =  [u^\flat_{\lambda,\, \lambda_0,\, y}] \cdot X_{\lambda_0,\, i_{\lambda_0}(y)} \qq \text{in} \qq \bA_\Inf(R_{\Sigma,\, \Lambda,\, \infty})
\]
(see \S\ref{tilt} and \S\ref{Ainf}). Thus, since each $X_{\lambda,\, i}$ is a nonzero divisor in $\bA_\Inf(R_{\Sigma,\, \Lambda,\, \infty})$, the $[u^\flat_{\lambda,\, \lambda_0,\, y}]$ solve the equations \eqref{nonsm-rels} in $\bA_\Inf(R_{\Sigma,\, \Lambda,\, \infty})$, compatibly with the solution in $R \subset R_{\Sigma,\, \Lambda,\, \infty}$ of \eqref{R-nat-alg-2}. Thus, in  the nonsmooth case as well, $\bA_\Inf(R_{\Sigma,\, \Lambda,\, \infty})$ is $\Delta_{\Sigma,\, \Lambda}$-equivariantly an $(A_{\Sigma,\, \Lambda}^\square \tensor_{\bZ[Q]} \bZ[P_{\lambda_0}])$-algebra, compatibly with the change of $\lambda_0$, the maps \eqref{R-nat-alg-2} and \eqref{more-theta} to $R$ and $R_{\Sigma,\, \Lambda,\, \infty}$, and the Frobenii.

In conclusion, in all the cases we obtain a compatible with the change of $\lambda_0$ commutative square
\[\ba
\xymatrix{
A_{\Sigma,\, \Lambda}^\square \tensor_{\bZ[Q]} \bZ[P_{\lambda_0}] \ar@{->>}[r]^-{j_{\lambda_0}} \ar[d] & R \ar[d] &  & \ar@{}[d]_{\text{\normalsize{\q so also}}} &  A_{\Sigma,\, \Lambda}^\square \tensor_{\bZ[Q]} \bZ[P_{\lambda_0}] \ar@{->>}[r]^-{j_{\lambda_0}} \ar[d] & R \ar[d] \\
\bA_\Inf(R_{\Sigma,\, \Lambda,\, \infty}) \ar@{->>}[r]^-{\theta} & R_{\Sigma,\, \Lambda,\, \infty}, &&& \bA_\cris^0(R_{\Sigma,\, \Lambda,\, \infty}) \ar@{->>}[r]^-{\theta} & R_{\Sigma,\, \Lambda,\, \infty}.
}  
\ea \]
The universal property of $D_{j_{\lambda_0}}$ then supplies the desired divided power morphism \eqref{NLE-eq}.
\epf

\bpp [The comparison map] \lab{comp-map}
\ready{The $p$-adic completion of the morphism \eqref{NLE-eq} is the morphism
\be \lab{CM-1}
D_{\Sigma,\, \Lambda} \ra \bA_\cris(R_{\Sigma,\, \Lambda,\, \infty})
\ee
(see \Cref{D-no-log}), which does not depend on $\lambda_0$. By construction, it makes the diagram
\be \lab{D-cris-comp-R}
\ba
\xymatrix@C=30pt{D_{\Sigma,\, \Lambda} \ar@{->>}[r]^{\eqref{R-D-map}} \ar[d]^{\eqref{CM-1}} & R \ar@{^(->}[d] \\  
\bA_\cris(R_{\Sigma,\, \Lambda,\, \infty}) \ar@{->>}[r]^-{\theta} &R_{\Sigma,\, \Lambda,\, \infty}}
\ea\ee
commute. Its restriction to $D_{\Sigma,\, \Lambda}^{0,\, (m)}$ of \S\ref{D-m} lands in the subring $\bA_\cris^{0,\, (m)}(R_{\Sigma,\, \Lambda,\, \infty}) \overset{\ref{Acrism-inj}}{\subset} \bA_\cris(R_{\Sigma,\, \Lambda,\, \infty})$, so, by passing to $p$-adic completions and using \eqref{Acrism-comp}, we obtain the compatible morphisms
\be \lab{CM-2}
D_{\Sigma,\, \Lambda}^{(m)} \ra \bA_\cris^{(m)}(R_{\Sigma,\, \Lambda,\, \infty}) \qq \text{for} \qq m \ge p
\ee
that recover \eqref{CM-1} via the identifications \eqref{D-Dm} and \eqref{Acris-m-eq}.
By construction and \Cref{non-log-env}, the morphisms \eqref{CM-2} are $\Delta_{\Sigma,\, \Lambda}$-equivariant and Frobenius-equivariant, so they give rise to the Frobenius-equivariant morphisms of complexes
\[
\tst K_{D_{\Sigma,\, \Lambda}^{(m)}}\p{(\delta_\sigma -1)_{\sigma\in\Sigma}, (\delta_{\lambda,\, i}-1)_{\lambda\in\Lambda,\, 1 \le i \le  d}} \ra K_{\bA_\cris^{(m)}(R_{\Sigma,\, \Lambda,\, \infty})}( (\delta_\sigma -1)_{\sigma\in\Sigma}, (\delta_{\lambda,\, i}-1)_{\lambda\in\Lambda,\, 1 \le i \le  d}).
\]
After applying the functor $\eta_{(\mu)}$, these morphisms compose with the ones constructed in \Cref{comp-map-0} to give rise to the desired Frobenius-equivariant comparison map of complexes: 
\be \lab{CM-3}
\tst\p{\varinjlim_{m \ge p} \p{K_{D^{(m)}_{\Sigma,\, \Lambda}}\p{\p{\frac{\partial}{\partial\log(X_\tau)}}_\tau }}}\wh{\ \ } \ra \p{\varinjlim_{m \ge p} \p{ \eta_{(\mu)} \p{ K_{\bA_{\cris}^{(m)}(R_{\Sigma,\, \Lambda,\, \infty})}((\delta_\tau - 1)_\tau)}}}\wh{\ },
\ee
where the direct limits and the $p$-adic completions are formed termwise and, for brevity, we let $\tau$ range over the indices ``$\sigma$'' for $\sigma \in \Sigma$ and ``$(\lambda, i)$'' for $\lambda \in \Lambda$ and $i = 1, \dotsc, d$.}

The source (resp.,~the target) of the map \eqref{CM-3} is a term of the direct limit \eqref{source-colim} (resp.,~\eqref{FAO-colim}) and the formation of this map is compatible with enlarging $\Sigma$ and $\Lambda$, that is, with the transition maps of the direct limits \eqref{source-colim} and \eqref{FAO-colim}. Moreover, if $R'$ is a $p$-adically formally \'{e}tale $R$-algebra equipped with data as in \eqref{more-coord} for some sets $\Sigma'$ and $\Lambda'$, then the map \eqref{CM-3} and its analogue for $R'$ and the sets $\Sigma \cup \Sigma'$, $\Lambda \cup \Lambda'$ (and the induced closed immersion \eqref{more-coord}) are compatible with the maps between their sources (resp.,~targets) discussed in \S\ref{functorial-AO} and \S\ref{functorial-cris}.  

In conclusion, by taking the filtered direct limit of the maps \eqref{CM-3} over all the closed immersions \eqref{more-coord} for varying $\Sigma$ and $\Lambda$ (but a fixed $R$), we obtain a comparison map from the complex \eqref{source-colim} to the complex \eqref{FAO-colim}, and the formation of this map is compatible with replacing $R$ by a formally \'{e}tale $R$-algebra $R'$. The following proposition implies that this map is a quasi-isomorphism.
\epp

\bprop \lab{alphaR-qiso}
The Frobenius-equivariant comparison map \eqref{CM-3} is a quasi-isomorphism. 
\eprop

\bpf
\ready{
The proof is similar to that of \cite{BMS16}*{12.9}, and the idea is to reduce to the case of a single coordinate morphism settled in \Cref{local-isom}. More precisely, for $m \ge p$, let 
\[
\Spec(R/p) \hra \Spf(D_{\Sigma,\, \Lambda}^{(m)})
\]
be the closed immersion induced by its analogue for $D_{\Sigma,\, \Lambda}$, that is, by the first map in \eqref{Acris-fact}. For each $\lambda_0 \in \Lambda$, the ideal of $A_{\Sigma,\, \Lambda}^\square \tensor_{\bZ[Q]} \bZ[P_{\lambda_0}]$ that cuts out $R/p$ (see \eqref{ci-fact}) is finitely generated. Consequently, for each $m \ge p$, the ideal of $D_{\Sigma,\, \Lambda}^{(m)}$ that cuts out $R/p$ is finitely generated, too, and hence, due to divided powers, it is also $p$-adically topologically nilpotent. Thus, fixing a $\lambda \in \Lambda$ and, for $m \ge p$, letting $A_\cris^{(m)}(R)_{\lambda}$ denote the ring $A_\cris^{(m)}(R)$ of \S\ref{AcrisR-m} constructed using the coordinate morphism $R_\lambda^\square \ra R$, we may use the $p$-adic formal \'{e}taleness of $A_\cris^{(m)}(R_\lambda^\square) \ra A_\cris^{(m)}(R)_{\lambda}$ (see \S\ref{Ainf}) to obtain the unique indicated lifts in the commutative diagram
\[
\xymatrix@C=30pt{
\Spec(R_{\Sigma,\, \Lambda,\, \infty}/p) \ar@{^(->}[d]^{\theta} \ar[r] & \Spec(R/p) \ar@{^(->}[d] \ar[r] & \Spf(A_\cris^{(m)}(R)_\lambda) \ar[d] \\
\Spf(\bA_\cris^{(m)}(R_{\Sigma,\, \Lambda,\, \infty})) \ar@{-->}[rru] \ar[r]_-{\eqref{CM-2}} & \Spf(D_{\Sigma,\, \Lambda}^{(m)}) \ar@{-->}[ru] \ar[r] & \Spf(A_\cris^{(m)}(R_\lambda^\square))
}
\]
in which the bottom right horizontal map results from the fact that, by construction, each $D_{\Sigma,\, \Lambda}^{(m)}$ is both an $A(R_{\lambda}^\square)$-algebra and an $A_\cris^{(m)}$-algebra. By the uniqueness of such lifts, the resulting maps 
\be \lab{map-to-Dm}
A_\cris^{(m)}(R)_\lambda \ra D_{\Sigma,\, \Lambda}^{(m)} 
\ee
are compatible as $m$ varies, $\Delta_{\Sigma,\, \Lambda}$-equivariant, where $\Delta_{\Sigma,\, \Lambda}$ acts on $A_\cris^{(m)}(R)_\lambda$ through the projection $\Delta_{\Sigma,\, \Lambda} \surjects \Delta_\lambda$, and are compatible with the maps from its source and target to $\bA_\cris^{(m)}(R_{\Sigma,\, \Lambda,\, \infty})$.} By construction, the maps \eqref{map-to-Dm} are also compatible with the derivations $\f{\partial}{\partial \log(X_{\lambda,\, i})}$ for $i = 1, \dotsc, d$ discussed in \S\ref{log-dR-cx} and \S\ref{log-PD-der}. 
Consequently, the diagram
\be\ba \lab{big-comm}
\tst
\xymatrix@C=31pt{
K_{A_\cris^{(m)}(R)_{\lambda}}\p{\frac{\partial}{\partial\log(X_{\lambda,\,1})}, \dots, \frac{\partial}{\partial \log(X_{\lambda,\, d})}} \ar[d]^-{\eqref{map-to-Dm}} \ar[r]^-{\eqref{local-isom-3}\,} &  \eta_{(\mu)}\p{K_{\bA_\cris^{(m)}(R_{\lambda,\, \infty})}\p{\delta_{\lambda,\, 1} - 1, \dots, \delta_{\lambda,\,d} - 1}} \ar[d] \\
K_{D^{(m)}_{\Sigma,\, \Lambda}}\p{\p{\frac{\partial}{\partial\log(X_\tau)}}_\tau }  \ar[r]^-{\ref{comp-map-0}\text{ and }\eqref{CM-2}} & \eta_{(\mu)} \p{ K_{\bA_{\cris}^{(m)}(R_{\Sigma,\, \Lambda,\, \infty})}((\delta_\tau - 1)_\tau)}
}
\ea \ee
commutes, where the index $\tau$ ranges over ``$\sigma$'' for $\sigma \in \Sigma$ and ``$(\lambda', i)$'' for $\lambda' \in \Lambda$ and $i = 1, \dotsc, d$. By \Cref{local-isom}, for $m \ge p^2$, the top horizontal map in \eqref{big-comm} is a  quasi-isomorphism and, by \Cref{Koszul} and \Cref{more-ANB-more}, so is the right vertical map. By \Cref{RHS-rewrite} and \eqref{Dm-Koszul}, the left vertical map in \eqref{big-comm} becomes a quasi-isomorphism after applying $\varinjlim_{m \ge p}$ and forming the termwise $p$-adic completion. These operations turn the bottom horizontal map in \eqref{big-comm} into the comparison map \eqref{CM-3}, so we conclude that the latter is also a quasi-isomorphism, as desired.
\epf

\bpp[Proof of Theorem {\upshape\ref{abs-crys-comp}}] \lab{Acris-main-pf}
By \S\ref{comp-map} and \Cref{alphaR-qiso}, the functorial in $R$ complexes \eqref{FAO-colim} and \eqref{source-colim} define canonically and Frobenius-equivariantly quasi-isomorphic complexes of presheaves on a basis for the topology of $\fX_\et$. Their associated complexes of sheaves on $\fX_\et$ are then also canonically and Frobenius-equivariantly quasi-isomorphic. By \S\ref{functorial-AO} and \S\ref{functorial-cris}, these complexes of sheaves Frobenius-compatibly represent $A\Omega_\fX \wh{\tensor}^\bL_{A_\Inf} A_\cris$ and $Ru_*(\cO_{\fX_{\cO_C/p}/A_\cris})$, respectively, so that, in conclusion, \Cref{alphaR-qiso} supplies a Frobenius-equivariant isomorphism
\be \lab{main-comp-iso}
Ru_*(\cO_{\fX_{\cO_C/p}/A_\cris}) \isomto A\Omega_\fX \wh{\tensor}^\bL_{A_\Inf} A_\cris,
\ee
which gives the desired identification \eqref{abs-crys-comp-eq}.
\QED
\epp

We have obtained two ways to identify the de Rham specialization of $A\Omega_\fX$: we may either use \eqref{dR-spec-id} or combine \eqref{main-comp-iso} with the fact that the logarithmic crystalline cohomology of $\fX_{\cO_C/p}$ over $\cO_C$ is computed by $\Omega^\bullet_{\fX/\cO_C,\, \log}$. 
For use in \S\ref{lattice}, we now check that these two identifications agree.

\bprop \lab{two-dR-same}
The following diagram commutes\ucolon
\be\ba \lab{TDS-eq}
\xymatrix{
Ru_*(\cO_{\fX_{\cO_C/p}/A_\cris}) \ar[rd] \ar[rr]_-{\sim}^-{\eqref{main-comp-iso}} && A\Omega_\fX \wh{\tensor}^\bL_{A_\Inf} A_\cris \ar[ld]^-{\eqref{dR-spec-id}} \\
& \Omega^\bullet_{\fX/\cO_C,\, \log} &
}
\ea\ee
where the unlabeled map results from the identification $Ru_*(\cO_{\fX_{\cO_C/p}/A_\cris}) \wh{\tensor}^\bL_{A_\cris,\, \theta}\, \cO_C \cong \Omega^\bullet_{\fX/\cO_C,\, \log}$ of \cite{Bei13a}*{(1.8.1) and (1.11.1)} \up{compare with Remark \uref{cris-dR-bc-rem}}. 
\eprop

\bpf
\ready{The overall argument is similar to the one given in the smooth case in \cite{BMS16}*{proof of~14.1}.

The claim is local, so we assume the setup of \S\ref{all-possible-coordinates}. Then the  $d\log(X_{\sigma})$ and $d\log(X_{\lambda,\, i})$ generate the differential graded algebra $\Omega^\bullet_{D_{\Sigma,\, \Lambda}/A_\cris,\, \log,\, \PD}$ over $D_{\Sigma,\, \Lambda}$ (see \Cref{D-compute}), so, since the terms of $\Omega^\bullet_{\Spf(R)/\cO_C,\, \log}$ are $p$-torsion free and each $t_\sigma$ and $t_{\lambda,\, i}$ is a unit in $R[\f{1}{p}]$, there is at most one map of  differential graded algebras 
\be \lab{cdga-uniq}
\Omega^\bullet_{D_{\Sigma,\, \Lambda}/A_\cris,\, \log,\, \PD} \ra \Omega^\bullet_{\Spf(R)/\cO_C,\, \log} \qq \text{with} \qq  \xymatrix@C=40pt{D_{\Sigma,\, \Lambda} \ar@{->>}[r]^-{\eqref{R-D-map}}  & R} \q \text{in degree} \q  0.
\ee
By \Cref{D-compute}, the unlabeled map of \eqref{TDS-eq} is identified with this unique map, so we need to show that so is the composition in \eqref{TDS-eq}.

We recall from the proof of \Cref{dR-spec} that the right diagonal map in \eqref{TDS-eq} is defined by composing the Frobenius of $A\Omega_\fX$, the reduction modulo $\varphi(\xi)$, and the canonical identification supplied by \cite{BMS16}*{6.12} of $(L\eta_{(\varphi(\xi))}(A\Omega_\fX))/\varphi(\xi)$ with the complex\footnote{\ready{For the sake of notational simplicity, we suppress the twists inherent in the construction $H^\bullet(-)$ of \emph{loc.~cit.}}} 
$H^\bullet(A\Omega_\fX/\varphi(\xi))$ that is \emph{a posteriori} identified with $\Omega_{\Spf(R)/\cO_C,\, \log}^\bullet$ and whose differentials are \emph{a priori} given by Bockstein homomorphisms (defined in \emph{loc.~cit.} using $A\Omega_\fX/(\varphi(\xi))^2$). Letting $\tau$ range over the usual indices (see  \S\ref{log-PD-der}), we may also apply this construction to the complex $\eta_{(\mu)}\p{K_{\bA_\Inf(R_{\Sigma,\, \Lambda,\, \infty})}((\delta_\tau -1)_\tau)}$: Frobenius maps it isomorphically to $\eta_{(\varphi(\mu))}\p{K_{\bA_\Inf(R_{\Sigma,\, \Lambda,\, \infty})}((\delta_\tau -1)_\tau)}$, for which the reduction modulo $\varphi(\xi)$ map is
\[
\tst \eta_{(\varphi(\mu))}\p{K_{\bA_\Inf(R_{\Sigma,\, \Lambda,\, \infty})}((\delta_\tau -1)_\tau)} \ra H^\bullet\p{\p{\eta_{(\mu)}\p{K_{\bA_\Inf(R_{\Sigma,\, \Lambda,\, \infty})}((\delta_\tau -1)_\tau)}}/\varphi(\xi)}.
\]
Moreover, due to the isomorphism \eqref{HT-iso-psh} and \Cref{e-iso-more,ANB-c-more}, 
\[
\p{\eta_{(\mu)}\p{K_{\bA_\Inf(R_{\Sigma,\, \Lambda,\, \infty})}((\delta_\tau -1)_\tau)}}/\varphi(\xi) \isomto \eta_{(\zeta_p - 1)}\p{K_{R_{\Sigma,\, \Lambda,\, \infty}}((\delta_\tau -1)_\tau)}.
\]
The cited remarks and \Cref{Hi-OmegaX} (with \Cref{desheafify}) imply that this describes the de Rham specialization map $A\Omega_\fX \ra \Omega^\bullet_{\fX/\cO_C,\, \log}$ in terms of the complex $\eta_{(\mu)}\p{K_{\bA_\Inf(R_{\Sigma,\, \Lambda,\, \infty})}((\delta_\tau -1)_\tau)}$.

We now describe the right diagonal map of \eqref{TDS-eq} in terms of $\eta_{(\mu)} \p{ K_{\bA_{\cris}^{(m)}(R_{\Sigma,\, \Lambda,\, \infty})}((\delta_\tau - 1)_\tau)}$, which is a variable term that comprises the target of the comparison map \eqref{CM-3}. Namely, we first let $\varphi(\bA_{\cris}^{(m)}(R_{\Sigma,\, \Lambda,\, \infty}))$ be the analogue of the ring $\bA_{\cris}^{(m)}(R_{\Sigma,\, \Lambda,\, \infty})$ built using the element $\varphi(\xi)$ instead of $\xi$, so that the Frobenius gives the isomorphism $\bA_{\cris}^{(m)}(R_{\Sigma,\, \Lambda,\, \infty}) \isomto \varphi(\bA_{\cris}^{(m)}(R_{\Sigma,\, \Lambda,\, \infty}))$.\footnote{\ready{Composing with the map $\varphi(\bA_{\cris}^{(m)}(R_{\Sigma,\, \Lambda,\, \infty})) \ra \bA_{\cris}^{(m)}(R_{\Sigma,\, \Lambda,\, \infty})$ recovers the Frobenius of $\bA_{\cris}^{(m)}(R_{\Sigma,\, \Lambda,\, \infty})$. 
}}
Then the Frobenius maps the complex $\eta_{(\mu)} \p{ K_{\bA_{\cris}^{(m)}(R_{\Sigma,\, \Lambda,\, \infty})}((\delta_\tau - 1)_\tau)}$ isomorphically to the complex $\eta_{(\varphi(\mu))} \p{ K_{\varphi(\bA_{\cris}^{(m)}(R_{\Sigma,\, \Lambda,\, \infty}))}((\delta_\tau - 1)_\tau)}$, for which the reduction modulo $\varphi(\xi)$ map is
\be \lab{TDS-3}
\tst \eta_{(\varphi(\mu))} \p{ K_{\varphi(\bA_{\cris}^{(m)}(R_{\Sigma,\, \Lambda,\, \infty}))}((\delta_\tau - 1)_\tau)} \ra H^\bullet\p{\eta_{(\mu)} \p{ K_{\varphi(\bA_{\cris}^{(m)}(R_{\Sigma,\, \Lambda,\, \infty}))}((\delta_\tau - 1)_\tau)}/\varphi(\xi)}.
\ee
Via a morphism induced by $\theta \circ \varphi\i \colon \varphi(\bA_{\cris}^{(m)}(R_{\Sigma,\, \Lambda,\, \infty})) \surjects R_{\Sigma,\, \Lambda,\, \infty}$, the target of \eqref{TDS-3} maps to 
\[
H^\bullet\p{\eta_{(\zeta_p - 1)}\p{K_{R_{\Sigma,\, \Lambda,\, \infty}}((\delta_\tau -1)_\tau)}} \overset{\ref{e-iso-more}\text{ and }\ref{Hi-OmegaX}}{\cong} \Omega_{\Spf(R)/\cO_C,\,\log}^\bullet
\]
because, since each $H^i\p{\eta_{(\zeta_p - 1)}\p{K_{R_{\Sigma,\, \Lambda,\, \infty}}((\delta_\tau -1)_\tau)}}$ is $p$-torsion free, the agreement of the Bockstein differentials may be checked after inverting $p$ by using the fact that $(\bA_\Inf(R_{\Sigma,\, \Lambda,\, \infty})/\varphi(\xi)^2)[\f{1}{p}]$ is an algebra over $\varphi(\bA_{\cris}^{(m)}(R_{\Sigma,\, \Lambda,\, \infty}))$ via a map that lifts $\theta \circ \varphi\i$. 
 The resulting composition
\be \lab{TDS-5}
\tst \eta_{(\mu)} \p{ K_{\bA_{\cris}^{(m)}(R_{\Sigma,\, \Lambda,\, \infty})}((\delta_\tau - 1)_\tau)} \ra H^\bullet\p{\eta_{(\zeta_p - 1)}\p{K_{R_{\Sigma,\, \Lambda,\, \infty}}((\delta_\tau -1)_\tau)}} \cong \Omega_{\Spf(R)/\cO_C,\,\log}^\bullet
\ee
is the promised description of the right diagonal map of \eqref{TDS-eq}. In addition, by construction and \cite{BMS16}*{6.13}, 
 this composition is a morphism of differential graded algebras\footnote{\ready{The  differential graded algebra structure on the Koszul complex $K_*((\delta_\tau - 1)_\tau)$ that computes continuous group cohomology is described in \cite{BMS16}*{7.5} and its proof.}} that in degree $0$ is given by the map $\theta$ of \eqref{AAcrism-theta}. }

On the other hand, the comparison map 
\[
\tst \Omega^\bullet_{D_{\Sigma,\, \Lambda}/A_\cris,\, \log,\, \PD} \cong K_{D_{\Sigma,\, \Lambda}}((\f{\partial}{\partial\log(X_\tau)})_\tau) \xra{\eqref{CM-3}} \p{\varinjlim_{m \ge p} \p{ \eta_{(\mu)} \p{ K_{\bA_{\cris}^{(m)}(R_{\Sigma,\, \Lambda,\, \infty})}((\delta_\tau - 1)_\tau)}}}\wh{\ }
\]
would be a map of differential graded algebras if in the formula $\log([\eps])\cdot \sum_{n\geq 0} \p{\frac{(\log([\eps]))^{n}}{ (n + 1)!}(\frac{\partial}{\partial\log (X_\tau)})^{n}}$ that describes the morphism \eqref{CM-0-eq} in degree $1$ we could disregard the terms with $n \ge 1$. However, $\log([\eps])$ and $\mu$ are unit multiples of each other (see \S\ref{log-eps}) and $\theta(\f{\mu^{n}}{(n + 1)!}) = 0$ in $\cO_C$ for $n \ge 1$, so we can indeed ignore these terms if we are only interested in the composition 
\[
\tst \Omega^\bullet_{D_{\Sigma,\, \Lambda}/A_\cris,\, \log,\, \PD} \ra \p{\varinjlim_{m \ge p} \p{ \eta_{(\mu)} \p{ K_{\bA_{\cris}^{(m)}(R_{\Sigma,\, \Lambda,\, \infty})}((\delta_\tau - 1)_\tau)}}}\wh{\ } \xra{\eqref{TDS-5}} \Omega_{\Spf(R)/\cO_C,\,\log}^\bullet.
\]
that describes the composition in \eqref{TDS-eq}. In conclusion, this composition is a morphism of  differential graded algebras that, due to \eqref{D-cris-comp-R}, in degree $0$ is the map $D_{\Sigma,\, \Lambda} \surjects R$ from \eqref{R-D-map}. Therefore, as desired, it is the unique morphism \eqref{cdga-uniq}.
\epf

We conclude \S\ref{section-Acris} by using \Cref{abs-crys-comp} to analyze the crystalline specialization of $R\Gamma(\fX_\et, A\Omega_\fX)$.

\bpp[The crystalline specialization map] \lab{Wk-log}
The Witt vector functoriality gives the surjection
\[
A_\Inf \surjects W(k), \qq \text{the \emph{crystalline specialization} map of $A_\Inf$.}
\]
Since $\xi$ maps to $p$ in $W(k)$, this surjection factors through $A_\cris$ as follows: $A_\Inf \hra A_\cris \surjects W(k)$. We equip $W(k)$ with the pullback of the log structure \eqref{AL-1} on $A_\cris$. Explicitly, the resulting log structure on $W(k)$ is associated to the prelog structure $\bQ_{\ge 0} \xra{0 \,\neq\, q\, \mapsto\, 0,\, 0\, \mapsto\, 1} W(k)$. 
\epp

\bcor \lab{Acris-BC}
For  quasi-compact and quasi-separated $\fX$, we have Frobenius-equivariant identifications
\be\ba \lab{Acris-BC-eq}
R\Gamma(\fX_\et, A\Omega_\fX) \wh{\tensor}_{A_\Inf}^\bL A_\cris &\cong R\Gamma_{\log \cris}(\fX_{\cO_C/p}/A_\cris), \\
R\Gamma(\fX_\et, A\Omega_\fX) \wh{\tensor}^\bL_{A_\Inf} W(k) &\cong R\Gamma_{\log\cris}(\fX_k/W(k)).
\ea\ee
For  $\cO_C$-proper $\fX$,  we have Frobenius-equivariant identifications
\be \ba \lab{ABCA-eq}
R\Gamma(\fX_\et, A\Omega_\fX) \tensor_{A_\Inf}^\bL A_\cris &\cong R\Gamma_{\log \cris}(\fX_{\cO_C/p}/A_\cris), \\
R\Gamma(\fX_\et, A\Omega_\fX) \tensor^\bL_{A_\Inf} W(k) &\cong R\Gamma_{\log\cris}(\fX_k/W(k)),
\ea\ee
and the cohomology modules of $R\Gamma(\fX_\et, A\Omega_\fX) \tensor_{A_\Inf}^\bL A_\cris[\f{1}{p}]$ are finite free over  $A_\cris[\f{1}{p}]$.
\ecor

\bpf
\ready{By \cite{BMS16}*{4.9 (i)}, a finitely presented $A_\Inf/p^n$-module is perfect as an $A_\Inf$-module.  Consequently, any $A_\Inf/p^n$-module $M$ is a filtered direct limit of perfect $A_\Inf$-modules, so, by \cite{SP}*{\href{http://stacks.math.columbia.edu/tag/0739}{0739}}, 
\[
R\Gamma(\fX_\et, A\Omega_\fX \tensor^\bL_{A_\Inf} M) \cong R\Gamma(\fX_\et, A\Omega_\fX) \tensor^\bL_{A_\Inf} M.
\]
This applies to $M = A_\cris/p^n$, so 
the first identification in \eqref{Acris-BC-eq} follows from \Cref{abs-crys-comp}. 

For each finite subextension of $C/(W(k)[\f{1}{p}])$, we consider its ring of integers $\cO \subset \cO_C$ equipped with the (fine) log structure associated to $\cO \cap (\cO[\f{1}{p}])^\times \hra \cO$. By using \'{e}tale local semistable coordinates \eqref{sst-coord} and \Cref{log-claim,log-comp}, we employ limit arguments to find such an $\cO$ and a quasi-compact and quasi-separated, fine, log smooth log scheme $\cX$ over $\cO/p$ that descends $\fX_{\cO_C/p}$ and is of Cartier type (see \cite{Kat89}*{4.8}). 
Then the base change theorem \cite{Bei13a}*{(1.11.1)} applies\footnote{\ready{In loc.~cit.,~the map $f$ of fine log schemes is quasi-compact and \emph{separated}. One may relax this to quasi-compact and quasi-separated: once $Y$ there is affine, the iterated intersections of opens in an affine cover of $Z$ are quasi-compact and separated over $Y$, so the \v{C}ech technique (compare with \cite{SP}*{\href{http://stacks.math.columbia.edu/tag/08BN}{08BN}}) reduces to the original assumptions.}} 
and shows that 
\be \lab{Bei-id-1}
R\Gamma_{\log \cris}(\fX_{\cO_C/p}/A_\cris) \wh{\tensor}^\bL_{A_\cris} W(k) \cong R\Gamma_{\log\cris}(\fX_k/W(k)),
\ee
so the second identification in \eqref{Acris-BC-eq} follows from the first. }

If $\fX$ is $\cO_C$-proper, then, by \Cref{perfect}, the object $R\Gamma(\fX_\et, A\Omega_\fX)$ is quasi-isomorphic to a bounded complex of finite free $A_\Inf$-modules, so the identifications in \eqref{ABCA-eq} follow from those in \eqref{Acris-BC-eq}. Moreover,  then $\cX$ is $\cO$-proper and \cite{Bei13a}*{1.18, Theorem} proves that the cohomology groups of 
\[
\tst R\Gamma_{\log \cris}(\fX_{\cO_C/p}/A_\cris) \tensor_{A_\cris}^\bL A_\cris[\f{1}{p}], \qq \text{and hence also of} \qq R\Gamma(\fX_\et, A\Omega_\fX) \tensor_{A_\Inf}^\bL A_\cris[\f{1}{p}],
\]
are finite free $A_\cris[\f{1}{p}]$-modules. 
\epf

\brems \lab{log-cris-bc}
\remi
In the preceding proof, the special fiber $\cX_k$ is a descent of $\fX_k$ to a fine, log smooth log scheme of Cartier type over $k$ equipped with the log structure associated to $\bN_{\ge 0} \xra{1\,\mapsto\,0,\, 0\, \mapsto\, 1} k$ (the base change map is a ``change of log structure'' self-map of $k$ determined by the map $\bN_{\ge 0} \ra \bQ_{\ge 0}$ that sends $1$ to the valuation of a uniformizer of $\cO$). Given such a descent, the base change theorem \cite{Bei13a}*{(1.11.1)} gives the further Frobenius-equivariant identification 
\setcounter{equation}{0}
\be \lab{LCB-eq}
\qq R\Gamma_{\log\cris}(\fX_k/W(k)) \cong R\Gamma_{\log\cris}(\cX_k/W(k)),
\ee
where the $W(k)$ on the left (resp.,~right) is equipped with the log structure associated to $\bQ_{\ge 0} \xra{0 \,\neq\, q\, \mapsto\, 0,\, 0\, \mapsto\, 1} W(k)$ (resp.,~$\bN_{\ge 0} \xra{1\,\mapsto\,0,\, 0\, \mapsto\, 1} W(k)$).
Likewise, if $\fX_k$ is $k$-smooth, then \emph{loc.~cit.}~gives the Frobenius-equivariant identification 
\be \lab{sm-no-log} 
\qq R\Gamma_{\log\cris}(\fX_k/W(k)) \cong R\Gamma_{\cris}(\fX_k/W(k)).
\ee

\remi
The identification \eqref{Bei-id-1} expresses $R\Gamma_{\log\cris}(\fX_k/W(k))$ in terms of $R\Gamma_{\log \cris}(\fX_{\cO_C/p}/A_\cris)$. Further results from \cite{Bei13a} imply a converse for proper $\fX$ after extending coefficients to $B_\st^+$, see \eqref{Bei-id-2} below (when $\fX_k$ is $k$-smooth, $A_\cris[\f{1}{p}]$ in place of $B_\st^+$ suffices, see \cite{BMS16}*{13.21}).
\erems




}

%% file: BdR.tex

\section{The comparison to the $B_\dR^+$-cohomology \nopunct} \lab{section-BdR}

\revise{
\ready{The main goal of this section is to prove in \Cref{BdR-comp} that for $\cO_C$-proper $\fX$, we have
\be \lab{section-BdR-main}
R\Gamma_{\log\cris}(\fX_{\cO_C/p}/A_\cris) \tensor^\bL_{A_\cris} B_\dR^+ \cong R\Gamma_\cris(\fX_C^\ad/B_\dR^+),
\ee
where the definition of $R\Gamma_\cris(\fX_C^\ad/B_\dR^+)$, the ``crystalline cohomology of $\fX_C^\ad$ over $B_\dR^+$,'' was given in \cite{BMS16}*{\S13} (see \S\ref{BdR-et} for a brief review). This definition is purely in terms of $\fX_C^\ad$ and was engineered in \emph{op.~cit.}~to be compatible with $R\Gamma_{\log\cris}(\fX_{\cO_C/p}/A_\cris)$ when $\fX$ is smooth. Thus, for \eqref{section-BdR-main}, we only need to check that a slightly more general definition that uses the \'{e}tale topology and more general embeddings than those furnished by annuli leads to the same cohomology (see \S\S\ref{BdR-et}--\ref{BdR-more}). For this, we adapt the arguments of \emph{op.~cit.};~in fact, our $C$ is $(\ov{W(k)}[\f{1}{p}])\,\wh{\ }$ (see \S\ref{fir-setup}), so we may simplify the descent to a discretely valued base aspects of these arguments by taking advantage of a result of Huber on the local structure of \'{e}tale maps of adic spaces (see \S\ref{BdR-more}).

\bpp[The ring $B_\dR^+$] \lab{ring-BdR}
Since $\xi$ is not a zero divisor in $A_\Inf[\f{1}{p}]$ and generates $\Ker(\theta[\f{1}{p}])$, the $(\Ker(\theta[\f{1}{p}]))$-adic completion of $A_\Inf[\f{1}{p}]$ is a complete discrete valuation ring 
$B_\dR^+$ with $\xi$ as a uniformizer and $C$ as the residue field.  By \Cref{Acrism-inj}, both $A_\Inf$ and $A_\cris$ are subalgebras of $B_\dR^+$. By the ``glueing of flatness'' \cite{RG71}*{II.1.4.2.1}, 
the ring $B_\dR^+$ is flat as an $A_\Inf$-algebra. We set 
\[
B_\dR \ce \Frac(B_\dR^+).
\]
Our $A_\Inf$ is a $W(k)$-algebra (see \S\ref{Ainf-not}), so, by the infinitesimal lifting, $B_\dR^+$ is a  $\ov{W(k)}[\f{1}{p}]$-algebra.
\epp

\bpp[The $B_{\dR}^{+}$-cohomology using the \'{e}tale topology] \lab{BdR-et}
\ready{In \cite{BMS16}*{\S13}, Bhatt--Morrow--Scholze used the analytic topology of a smooth adic $C$-space $X$ to define the ``$B_\dR^+$-cohomology'' of $X$,
\[
R\Gamma_{\cris}(X/B_\dR^+) \in D^{\ge 0}(B_\dR^+).
\]
We now review their construction and show that it may also be carried out in the \'{e}tale topology.

By \cite{Hub96}*{1.6.10, 2.2.8}, the analytic (resp.,~\'{e}tale) topology of $X$ has a basis of affinoids $\Spa(A, A^\circ)$ 
each of which admits a map 
\be \lab{affd-et}
\Spa (A, A^\circ)\ra \bT_C^d \ce \Spa (C\langle T_1^{\pm 1}, \dotsc, T_d^{\pm 1} \rangle, \cO_C\langle T_1^{\pm 1}, \dotsc, T_d^{\pm 1} \rangle)  \q \text{for some} \q d \ge 0
\ee
that is a composition of a rational embedding, a finite \'{e}tale map, and a rational embedding. By localizing further, we refine the basis to consist of those $\Spa(A, A^\circ)$ as above for which there is a finite subset $\Psi \subset (A^\circ)^\times$ such that the following map is surjective: 
\be \lab{surj-onto-R}
\tst  C\In{ (X_u^{\pm 1})_{u \in \Psi} } \xra{X_u \mapsto u\,} A.
\ee
Then, by endowing each $A_\Inf/\xi^n$ with the $p$-adic topology, each $(A_\Inf/\xi^n)[\f{1}{p}]$ with the unique ring topology for which $A_\Inf/\xi^n$ is an open subring, setting
\be \lab{BdR-pow-def}
\tst B_\dR^+\In{ (X_u^{\pm 1})_{u \in \Psi} } \ce \varprojlim_{n > 0} ((B_\dR^+/\xi^n)\In{ (X_u^{\pm 1})_{u \in \Psi}}),
\ee
and composing the projection onto the $n = 1$ term with \eqref{surj-onto-R}, we obtain the surjection 
\[
\tst s\colon B_\dR^+\In{ (X_u^{\pm 1})_{u \in \Psi} } \surjects A \qq \text{and set} \qq D_{\Psi}(A) \ce \varprojlim_{n > 0} ((B_\dR^+\In{ (X_u^{\pm 1})_{u \in \Psi} })/(\Ker s)^n).
\]
By the Leibniz rule, any $(B_\dR^+\In{ (X_u^{\pm 1})_{u \in \Psi} })$-valued derivation of $B_\dR^+\In{ (X_u^{\pm 1})_{u \in \Psi} }$ induces a $D_{\Psi}(A)$-valued derivation of $D_{\Psi}(A)$. Thus, the commuting derivations $\f{\partial}{\partial\log(X_u)} \ce X_u \cdot \f{\partial}{\partial X_u}$ give rise to the Koszul complex
\[
\tst \Omega^\bullet_{D_\Psi(A)/B_\dR^+} \ce K_{D_\Psi(A)}\p{(\f{\partial}{\partial \log(X_u)})_{u \in \Psi}}
\]
that is functorial in enlarging $\Psi$. The resulting complex 
\be \lab{omg-R-def}
\tst \Omega^\bullet_{A/B_\dR^+} \ce \varinjlim_\Psi\p{\Omega^\bullet_{D_\Psi(A)/B_\dR^+}}
\ee
is contravariantly functorial in $\Spa(A, A^\circ)$. Consequently, by varying $\Spa(A, A^\circ)$, we obtain a complex of presheaves on the basis described above for the analytic (resp.,~\'{e}tale) topology of $X$. The hypercohomology of the associated complex of sheaves is, by definition, the $B_\dR^+$-cohomology of $X$:
\[
R\Gamma_{\cris}(X/B_\dR^+) \qq \text{(resp.,~its variant for the \'{e}tale topology} \qq R\Gamma_{\cris}(X_\et/B_\dR^+)).
\]
By \cite{BMS16}*{13.12~(ii), 13.13}, if $\Spa(A, A^\circ)$ is fixed and $\Psi$ is sufficiently large, then $D_\Psi(A)$ is $\xi$-torsion free and $\xi$-adically complete, $\Omega^\bullet_{D_\Psi(A)/B_\dR^+}$ maps quasi-isomorphically to $\Omega^\bullet_{A/B_\dR^+}$, and, in the derived category, we have a canonical identification 
\[
\Omega^\bullet_{D_\Psi(A)/B_\dR^+}/\xi \cong \Omega^{\bullet,\, \cont}_{A/C},\qq \text{so also} \qq \Omega^\bullet_{A/B_\dR^+}/\xi \cong \Omega^{\bullet,\, \cont}_{A/C}.
\]
Consequently, by \cite{BMS16}*{9.15} (which we also used for proving \Cref{der-comp}), our definition of $R\Gamma_{\cris}(X/B_\dR^+)$ agrees with that of \cite{BMS16}*{\S13} (where one skips the sheafification step),
\be \lab{BdR-coho-comp}
R\Gamma_{\cris}(X/B_\dR^+) \qq \text{and} \qq R\Gamma_{\cris}(X_\et/B_\dR^+) \qq \text{are derived $\xi$-adically complete,}
\ee
and their (derived) reductions modulo $\xi$ are canonically and compatibly identified with the de Rham cohomology objects $R\Gamma(X, \Omega^{\bullet,\, \cont}_{X/C})$  and $R\Gamma(X_\et, \Omega^{\bullet,\, \cont}_{X/C})$, respectively, for instance:
\be \lab{BdR-mod-xi}
R\Gamma_\cris(X/B_\dR^+) \tensor^\bL_{B_\dR^+} C \cong R\Gamma(X, \Omega^{\bullet,\, \cont}_{X/C}) \equalscolon R\Gamma_\dR(X/C).
\ee
}Thus, since, by the Hodge-to-de Rham spectral sequence and \cite{Sch13}*{9.2~(ii)}, the formation of the de Rham cohomology is insensitive to passage to the \'{e}tale topology, we have
\be \lab{Zar-et-BdR}
R\Gamma_{\cris}(X/B_\dR^+) \isomto R\Gamma_{\cris}(X_\et/B_\dR^+)
\ee
via pullback. In addition, if $X$ is proper over $C$ and there is a complete, discretely valued subfield $K \subset C$ with a perfect residue field and a proper, smooth adic space $X_0$ over $K$ equipped with an isomorphism $X \cong X_0 \wh{\tensor}_K C$, then, by \cite{BMS16}*{13.20}, there is a canonical identification
\be \lab{cris-dR-descent}
R\Gamma_\cris(X/B_\dR^+) \cong R\Gamma_\dR(X_0/K) \tensor_K B_\dR^+, \q \text{where} \q R\Gamma_\dR(X_0/K) \ce R\Gamma(X_0, \Omega_{X_0/K}^{\bullet,\, \cont}).
\ee
In this situation, by the proof of \emph{loc.~cit.}, the reduction modulo $\xi$ of the identification \eqref{cris-dR-descent}  recovers the identification \eqref{BdR-mod-xi} under the base change identification $R\Gamma_\dR(X/C) \cong R\Gamma_\dR(X_0/K) \tensor_K^\bL C$. 
\epp

\bpp[The $B_{\dR}^{+}$-cohomology using more general embeddings] \lab{BdR-more}
\ready{To relate $R\Gamma_{\cris}(X_\et/B_\dR^+)$ to the absolute crystalline cohomology of \S\ref{abs-crys-log}, we now mildly generalize the construction of the former. 

The \'{e}tale topology of $X$ has a basis of affinoids $\Spa(A, A^\circ)$ each of which admits an \'{e}tale map
\be \ba \lab{R-et-map}
\Spa(A, A^\circ) \ra \Spa (C\langle T_0, \dotsc, T_r, &T_{r + 1}^{\pm 1}, \dotsc, T_d^{\pm 1} \rangle/(T_0\cdots T_r - p^q), \\ &\cO_C\langle T_0, \dotsc, T_r, T_{r + 1}^{\pm 1}, \dotsc, T_d^{\pm 1} \rangle/(T_0\cdots T_r - p^q))
\ea \ee
for some $d \ge r \ge 0$ and $q \in \bQ_{> 0}$ (we have seen in \S\ref{BdR-et} that even the ones with $r = 0$ would suffice).  By \cite{Hub96}*{1.7.3~iii)}\footnote{\ready{Noncomplete $A$ are allowed in \emph{loc.~cit.},~so we choose $A^+ \ce \ov{W(k)}[ T_1, \dotsc, T_r, T_{r + 1}^{\pm 1}, \dotsc, T_d^{\pm 1} ]/(T_0\cdots T_r - p^q)$ and $A^\rhd \ce A^+[\f{1}{p}]$.} } and limit arguments, there is a finite subextension $W(k)[\f{1}{p}] \subset K \subset C$ with the ring of integers $\cO$ containing $p^q$ and a finite type $(\cO[T_0, \dotsc, T_r, T_{r + 1}^{\pm 1}, \dotsc, T_d^{\pm 1}]/(T_0\cdots T_r - p^q))$-algebra $A_0$ that is \'{e}tale after inverting $p$, flat over $\cO$, normal, 
and such that the morphism \eqref{R-et-map} is the $C$-base change of an \'{e}tale $\Spa(K, \cO)$-morphism 
\be\ba \lab{adic-descent}
\tst \Spa(\wh{A_0}[\f{1}{p}], \wh{A_0}) \ra \Spa (K\langle T_0, \dotsc, T_r, &T_{r + 1}^{\pm 1}, \dotsc, T_d^{\pm 1} \rangle/(T_0\cdots T_r - p^q), \\ &\cO\langle T_0, \dotsc, T_r, T_{r + 1}^{\pm 1}, \dotsc, T_d^{\pm 1} \rangle/(T_0\cdots T_r - p^q)).
\ea\ee
By also using the reduced fiber theorem \cite{SP}*{\href{http://stacks.math.columbia.edu/tag/09IL}{09IL}}, we enlarge $K$ to ensure that, in addition, 
\be\lab{red-fib-eq}
A_0 \wh{\tensor}_\cO \ov{W(k)} \cong A^\circ.
\ee
The connected components of $\Spec(A_0)$ on which $p$ is a unit do not contribute to $\wh{A_0}$, so we lose no generality by assuming that $\Spec(A_0)$ has no such components. 

By working locally on $\Spa(\wh{A_0}[\f{1}{p}], \wh{A_0})$, we refine the basis above to consist of those  $\Spa(A, A^\circ)$ for which, in addition, there are finite subsets $\Psi_0 \subset (\wh{A_0})^\times$ and $\Xi_0 \subset \wh{A_0} \cap (\wh{A_0}[\f{1}{p}])^\times$ such that the~map
\be \lab{embed-ball-0}
\tst K\langle (x_u^{\pm 1})_{u \in \Psi_0}, (x_a)_{a \in \Xi_0} \rangle \xra{x_u \mapsto u,\, x_a \mapsto a} \wh{A_0}[\f{1}{p}]
\ee
is surjective. 
Then there are finite subsets $\Psi \subset (A^\circ)^\times$ and $\Xi \subset A^\circ \cap A^\times$ such that the map
\be \lab{embed-ball}
C\langle (X_u^{\pm 1})_{u \in \Psi}, (X_a)_{a \in \Xi} \rangle \xra{X_u \mapsto u,\, X_a \mapsto a} A
\ee
is also surjective. 
Defining the ring $B_\dR^+\In{ (X_u^{\pm 1})_{u \in \Psi},\, (X_a)_{a \in \Xi} }$ analogously to \eqref{BdR-pow-def}, so that the map \eqref{embed-ball} gives rise to the surjection 
\[
s\colon B_\dR^+\In{ (X_u^{\pm 1})_{u \in \Psi},\, (X_a)_{a \in \Xi} } \surjects A,
\]
we set
\[
\tst  D_{\Psi,\, \Xi,\, n}(A) \ce (B_\dR^+\In{ (X_u^{\pm 1})_{u \in \Psi},\, (X_a)_{a \in \Xi} })/(\Ker s)^n \qq \text{and} \qq D_{\Psi,\, \Xi}(A) \ce \varprojlim_{n > 0} D_{\Psi,\, \Xi,\, n}(A).
\]
By \cite{BMS16}*{13.4~(ii)}, each $D_{\Psi,\, \Xi,\, n}(A)$ is a complete, strongly Noetherian Tate ring (in the sense of \cite{Hub93a}*{\S1}), with the image of $(A_\Inf/\xi^n)\In{(X_u^{\pm 1})_{u \in \Sigma}, (X_a)_{a \in \Xi}}$ endowed with its $p$-adic topology as a ring of definition. By construction, $D_{\Psi,\, \Xi,\, n}(A)$ is a nilpotent thickening of $D_{\Psi,\, \Xi,\, 1}(A) \cong A$.

The ring $B_\dR^+$ is a $K$-algebra (see \S\ref{ring-BdR}), so we let $(B_\dR^+/\xi^n)_0 \subset B_\dR^+/\xi^n$ for $n > 0$ be the (module-finite) $A_\Inf/\xi^n$-subalgebra generated by the image of $\cO$. The proof of \cite{BMS16}*{13.11} shows (with $R_A$ there replaced by our $\wh{A_0}[\f{1}{p}]$)\footnote{\lab{RG-arg}\ready{In fact, in our case the argument is simpler, and we sketch it here. Since $\Spec(A_0)$ has no connected components on which $p$ is a unit, by \cite{RG71}*{I.3.3.5} and \cite{SP}*{\href{http://stacks.math.columbia.edu/tag/0593}{0593}}, the ring $A_0$ is free as an $\cO$-module. Thus, the $n$-th term of the inverse limit in \eqref{BdR-inv-lim} is a $p$-adically completed direct sum of copies of $(A_\Inf/\xi^n)[\f{1}{p}]$. This makes the multiplication by $\xi^m$ map on this $n$-th term explicit and the desired claims follow by passing to the inverse limit over $n$.}} that the $B_\dR^+$-algebra 
\be \lab{BdR-inv-lim}
\tst B_\dR^+ \wh{\tensor}_K (A_0[\f{1}{p}]) \ce \varprojlim_{n > 0} \p{( (B_\dR^+/\xi^n)_0 \wh{\tensor}_\cO A_0)[\f{1}{p}]}
\ee
has no nonzero $\xi$-torsion and is $\xi$-adically complete with 
\[
\tst (B_\dR^+ \wh{\tensor}_K (A_0[\f{1}{p}]))/\xi \cong A \q \text{and, more generally,} \q (B_\dR^+ \wh{\tensor}_K (A_0[\f{1}{p}]))/\xi^n \cong ( (B_\dR^+/\xi^n)_0 \wh{\tensor}_\cO A_0)[\f{1}{p}].
\]
The argument of \cref{RG-arg} shows that the map $(B_\dR^+/\xi^{n + 1})_0 \wh{\tensor}_\cO A_0 \ra (B_\dR^+/\xi^n)_0 \wh{\tensor}_\cO A_0$  is surjective for $n > 0$ (with the kernel of square zero, as may be seen after inverting $p$), so the subring 
\be \lab{sub-surj}
\tst \varprojlim_{n > 0} \p{ (B_\dR^+/\xi^n)_0 \wh{\tensor}_\cO A_0} \subset B_\dR^+ \wh{\tensor}_K (A_0[\f{1}{p}]) \qq \text{surjects onto} \qq \cO_C \wh{\tensor}_\cO A_0 \overset{\eqref{red-fib-eq}}{\cong} A^\circ.
\ee
Moreover, we have the following analogue of \cite{BMS16}*{13.12~(ii)} whose proof will be given in \S\ref{BMS-lem-pf}:}
\epp

\begin{sublem} \lab{BMS-lem}
If $\Spa(A, A^\circ)$ is an element of the refined basis for $X_\et$ described above, $\Psi \subset (A^\circ)^\times$ \up{resp.,~$\Xi \subset A^\circ \cap A^\times$} contains the images of the $T_i$ for $r + 1 \le i \le d$ \up{resp.,~$1 \le i \le r$} under a coordinate morphism as in \eqref{R-et-map}, and $\Psi$ and $\Xi$ are large enough \up{see the proof for the precise meaning}, then
\be \lab{BMS-lem-eq}
\tst D_{\Psi,\, \Xi}(A) \cong (B_\dR^+ \wh{\tensor}_K (A_0[\f{1}{p}]))\llb (X_a - \wt{a})_{a \in (\Psi \cup \Xi) \setminus \{ T_1, \dotsc, T_d\}} \rrb
\ee
where $\wt{a} \in \varprojlim_{n > 0} \p{ (B_\dR^+/\xi^n)_0 \wh{\tensor}_\cO A_0} \subset  B_\dR^+ \wh{\tensor}_K (A_0[\f{1}{p}])$ is a fixed lift of $a$ \up{see \eqref{sub-surj}}. In particular, for large $\Psi$ and $\Xi$, the $B_\dR^+$-algebra $D_{\Psi,\, \Xi}(A)$ has no nonzero $\xi$-torsion and is $\xi$-adically complete.
\end{sublem}

Similarly to \S\ref{BdR-et}, for any $\Psi$ and $\Xi$ as in \eqref{embed-ball}, the derivations $\f{\partial}{\partial\log(X_a)} \ce X_a \cdot \f{\partial}{\partial X_a}$ with $a \in \Psi\cup \Xi$ extend to $D_{\Psi,\, \Xi}(A)$, and we may define the Koszul complex
\[
\tst \Omega^\bullet_{D_{\Psi,\, \Xi}(A)/B_\dR^+} \ce K_{D_{\Psi,\, \Xi}(A)}\p{(\f{\partial}{\partial \log(X_u)})_{u \in \Psi}, (\f{\partial}{\partial \log(X_a)})_{a \in \Xi}}
\]
that is functorial in replacing $\Psi$ and $\Xi$ by larger $\Psi'$ and $\Xi'$. Since $a \in A^\times$ for $a \in \Psi \cup \Xi$, the proof of \cite{BMS16}*{13.13} shows that for $\Psi$ and $\Xi$ to which \Cref{BMS-lem} applies,
\be \lab{mod-xi-dR}
\Omega^\bullet_{D_{\Psi,\, \Xi}(A)/B_\dR^+}/\xi \cong \Omega^{\bullet,\, \cont}_{A/C} \qq \text{in the derived category,}
\ee
compatibly with enlarging $\Psi$ and $\Xi$. In particular, due to the derived $\xi$-adic completeness supplied by \Cref{BMS-lem}, for such large enough $\Psi$ and $\Xi$, the map 
\[
\Omega^\bullet_{D_{\Psi,\, \Xi}(A)/B_\dR^+} \ra \Omega^\bullet_{D_{\Psi',\, \Xi'}(A)/B_\dR^+} \qq \text{is a quasi-isomorphism.}
\]
Thus, if the element $\Spa(A, A^\circ)$ of the refined basis above also belongs to the basis considered in \S\ref{BdR-et}, that is, if it has an \'{e}tale coordinate map as in \eqref{affd-et} and a surjection \eqref{surj-onto-R}, then 
we obtain the functorial in $\Spa(A, A^\circ)$ quasi-isomorphism with the complex $\Omega^\bullet_{A/B_\dR^+}$ of \eqref{omg-R-def}:
\be \lab{omg-comp-iso}
\tst \Omega^\bullet_{A/B_\dR^+} \isomto \varinjlim_{\Psi,\, \Xi} \p{ \Omega^\bullet_{D_{\Psi,\, \Xi}(A)/B_\dR^+}}.
\ee
Such $\Spa(A, A^\circ)$ still form a basis for $X_\et$ (see the parenthetical remark after \eqref{R-et-map}), so we  conclude that the hypercohomology of the sheafification of the complex of presheaves furnished by the target of \eqref{omg-comp-iso} is identified with $R\Gamma_\cris(X_\et/B_\dR^+)$. In conclusion, we may summarize informally:
\be \lab{informal}
\text{the complexes} \q \Omega^\bullet_{D_{\Psi,\, \Xi}(A)/B_\dR^+} \q \text{also compute the $B_\dR^+$-cohomology} \q R\Gamma_\cris(X/B_\dR^+)
\ee
and the maps \eqref{mod-xi-dR} recover the identification \eqref{BdR-mod-xi}.

\bpp[Proof of Lemma {\upshape \ref{BMS-lem}}] \lab{BMS-lem-pf}
\ready{We adapt the proof of \cite{BMS16}*{13.12 (ii)} as follows. In addition to the coordinate morphism \eqref{R-et-map} and its descent \eqref{adic-descent} used in the statement, we fix subsets $\Psi_0 \subset (\wh{A_0})^\times$ and $\Xi_0 \subset \wh{A_0} \cap (\wh{A_0}[\f{1}{p}])^\times$ such that, as in \eqref{embed-ball-0}, the map 
\[
\tst s_0\colon K \In{ (x_u^{\pm 1})_{u \in \Psi_0}, (x_{ a})_{a \in \Xi_0} } \xra{x_{ u} \mapsto u,\, x_{ a} \mapsto a} \wh{A_0}[\f{1}{p}]
\]
is surjective and $\Psi_0$ (resp.,~$\Xi_0$) contains the images of the $T_i$ for $r + 1 \le i \le d$ (resp.,~$1 \le i \le r$) under the map \eqref{adic-descent}. We require that $\Psi$ (resp.,~$\Xi$) contains the image of this $\Psi_0$ (resp.,~$\Xi_0$) in $A^\circ$---this is the meaning of ``large enough'' in the statement. We set 
\[
\tst D_{0,\, n} \ce K \In{ (x_u^{\pm 1})_{u \in \Psi_0}, (x_{ a})_{a \in \Xi_0} }/(\Ker s_0)^n \q \text{for} \q n > 0 \qq   \text{and} \qq D_0 \ce \varprojlim_{n > 0} D_{0,\, n}.
\]
The continuous map $K \ra B_\dR^+$ (see \S\ref{ring-BdR}) gives a compatible with $s_0$ and $s$ continuous map
\[
K \In{ (x_u^{\pm 1})_{u \in \Psi_0}, (x_{ a})_{a \in \Xi_0} } \xra{x_u \mapsto X_u,\, x_a \mapsto X_a} B_\dR^+ \In{ (X_u^{\pm 1})_{u \in \Psi}, (X_{ a})_{a \in \Xi} }, \q \text{so also} \q D_{0,\, n} \ra D_{\Psi,\, \Xi,\, n}(A).
\]
By the $K[T_1, \dotsc, T_r, T_{r + 1}^{\pm 1}, \dotsc, T_d^{\pm 1}]$-\'{e}taleness 
of $A_0[\f{1}{p}]$, the map
\be \lab{lift-0}
\tst A_0[\f{1}{p}] \ra \wh{A_0}[\f{1}{p}] \qq \text{lifts to a map} \qq A_0[\f{1}{p}] \ra D_0 \q \text{with} \q T_i \mapsto x_{T_i}.
\ee
By \cite{GR03}*{7.3.15}, 
for each $n > 0$, the subring $D_{0,\, n}^\circ \subset D_{0,\, n}$ of powerbounded elements is the preimage of its counterpart $(\wh{A_0}[\f{1}{p}])^\circ \subset \wh{A_0}[\f{1}{p}]$. Thus, the lift \eqref{lift-0} maps $A_0$ to $D_{0,\, n}^\circ$, so also to some ring of definition of $D_{0,\, n}$. By composing with the map $D_{0,\, n} \ra D_{\Psi,\, \Xi,\, n}(A)$, we obtain the map $A_0 \ra D_{\Psi,\, \Xi,\, n}(A)$ whose image lies in some ring of definition, so, as $n$ varies, also the map
\[ 
\tst B_\dR^+ \wh{\tensor}_K (A_0[\f{1}{p}]) \ra D_{\Psi,\, \Xi}(A) \qq \text{that is compatible with the maps to $A$.}
\]
This gives rise to the continuous map $y$ in the diagram 
\[
\xymatrix@C=40pt@R=6pt{
 & (B_\dR^+ \wh{\tensor}_K (A_0[\f{1}{p}]))\llb (X_a - \wt{a})_{a \in (\Psi \cup \Xi) \setminus \{ T_1, \dotsc, T_d\}} \rrb \ar@/_1.0pc/[dd]_-{y} \ar@{->>}@/^1.0pc/[rd]^-{\, X_a\, \mapsto\, a}  & \\
B_\dR^+ \In{ (X_u^{\pm 1})_{u \in \Psi}, (X_{ a})_{a \in \Xi} } \ar@/_0.5pc/[rd] \ar@/^1.0pc/[ru]^-{X_{T_i}\, \mapsto\, T_i} & & A \\ 
& D_{\Psi,\, \Xi}(A) \ar@/_1.0pc/[uu]_-{z} \ar@/_0.5pc/@{->>}[ru] &
}
\]
whose maps ``$X_{T_i} \mapsto T_i$'' and $z$ are defined as follows.
\begin{itemize}
\item
To define the map ``$X_{T_i} \mapsto T_i$,'' one first forms the inverse limit over $N$ of the maps 
\[
\tst \qqq (A_\Inf/\xi^n)\In{ (X_u^{\pm 1})_{u \in \Psi}, (X_{ a})_{a \in \Xi} } \xra{X_{T_i} \mapsto T_i} ((B_\dR^+/\xi^n)_0 \tensor_\cO A_0/p^N)\llb (X_a - \wt{a})_{a \in (\Psi \cup \Xi) \setminus \{ T_1, \dotsc, T_d\}}\rrb
\]
defined by using the fact that each $\wt{u}$ with $u \in \Psi$ is a unit in $(B_\dR^+/\xi^n)_0 \wh{\tensor}_\cO A_0$ (see the sentence of \eqref{sub-surj}) and the identity $X_u\i = \wt{u}\i(1 - \wt{u}\i(X_u - \wt{u}) + \wt{u}^{-2}(X_u - \wt{u})^2 - \ldots)$. Then one inverts $p$ and forms the inverse limit over $n$.

\item
The continuous map $z$ is defined by combining the top part of the diagram, the $\xi$-adic completeness of $B_\dR^+ \wh{\tensor}_K (A_0[\f{1}{p}])$ (see \S\ref{BdR-more}), and the definition of $ D_{\Psi,\, \Xi}(A)$.
\end{itemize}
By construction, the diagram commutes, so that $y \circ z = \id$. By the $K[T_1, \dotsc, T_r, T_{r + 1}^{\pm 1}, \dotsc, T_d^{\pm 1}]$-\'{e}taleness of $A_0[\f{1}{p}]$, the $K[T_1, \dotsc, T_r, T_{r + 1}^{\pm 1}, \dotsc, T_d^{\pm 1}]$-algebra endomorphism $z \circ y$ of the pro-thickening $(B_\dR^+ \wh{\tensor}_K (A_0[\f{1}{p}]))\llb (X_a - \wt{a})_{a \in (\Psi \cup \Xi) \setminus \{ T_1, \dotsc, T_d\}} \rrb$ of $A$ is the identity on $A_0[\f{1}{p}]$, so also on $(B_\dR^+ \wh{\tensor}_K (A_0[\f{1}{p}]))$. It also fixes every $X_a$, so it must be the identity. Thus, $z$ is the desired isomorphism \eqref{BMS-lem-eq}. \QED}
\epp

\bpp[The map from the absolute crystalline cohomology] \lab{map-from-abs}
\ready{Returning to the $\fX$ of \S\ref{fir-setup}, our next goal is to use the preceding discussion 
to exhibit a map 
\be \lab{abs-to-BdR}
R\Gamma_{\log\cris}(\fX_{\cO_C/p}/A_\cris)\ra R\Gamma_\cris (\fX^\ad_C /B_{\dR}^{+}) \qq \text{over} \qq A_\cris \xra{\text{\S\ref{ring-BdR}}} B_\dR^+.
\ee
For this, we use the basis of $\fX_\et$ consisting of the affine opens $\Spf (R)$ as in the ``all possible coordinates'' setting of \S\ref{all-possible-coordinates} and adopt the subsequent notation of \S\S\ref{all-possible-coordinates}--\ref{Acris-main-pf}. To relate to \S\ref{BdR-more},~we~set 
\be \lab{PX-def}
\tst A \ce R[\f{1}{p}], \q \Psi \ce \{t_\sigma\}_{\sigma \in\Sigma} \cup \bigcup_{\lambda \in \Lambda} \{ t_{\lambda,\, r_\lambda + 1}, \dotsc, t_{\lambda,\, d} \}, \q \text{and} \q\Xi \ce \bigcup_{\lambda \in \Lambda} \{ t_{\lambda,\, 1}, \dotsc, t_{\lambda,\, r_{\lambda}} \}
\ee
(so that $A^\circ \cong R$ 
and the $t_{\lambda,\, 0}$ are omitted). We may descend the \'{e}tale map \eqref{lambda-proj} for $\lambda \in \Lambda$ to the ring of integers of a finite subextension $W(k)[\f{1}{p}] \subset K \subset C$ (see \eqref{sst-approx}) and then obtain the descended coordinate map \eqref{adic-descent} on the generic fiber. In addition, by enlarging $K$ and using the closed immersion \eqref{Sigma-ci}, we may ensure that the descent has a closed immersion \eqref{embed-ball-0} (with $\Xi_0 = \emptyset$); 
we then enlarge $\Psi$ by adjoining the image in $(A^\circ)^\times$ of the resulting $\Psi_0$. Thus, the above choices of $A$, the enlarged $\Psi$, and $\Xi$ satisfy the assumptions of \S\ref{BdR-more}: specifically, $\Spa(A, A^\circ)$ is an element of the (refined) basis of $(\fX_C^\ad)_\et$ considered there, $\Psi$ (resp.,~$\Xi$) contains $t_{\lambda,\, i}$ for $r_\lambda + 1 \le i \le d$ (resp.,~$1 \le i \le r_\lambda$), and \Cref{BMS-lem} applies to (the enlarged) $\Psi$ and $\Xi$. In conclusion, with these choices, the entire \S\ref{BdR-more} applies.

By using descents and \cite{GR03}*{7.3.15} as in the proof of \Cref{BMS-lem}, we see that the elements $\f{[(p^{1/p^\infty})^{q_\lambda}]}{X_{t_{\lambda,\, 1}} \cdots X_{t_{\lambda,\, r_\lambda}}}$ of $D_{\Psi,\, \Xi,\, n}(A)$ lie in $(D_{\Psi,\, \Xi,\, n}(A))^\circ$. 
Thus, each $(D_{\Psi,\, \Xi,\, n}(A))^\circ$, so also $D_{\Psi,\, \Xi}(A)$, is naturally an algebra over the ring $A_{\Sigma,\, \Lambda}^\square$ defined in \eqref{log-ci}. In fact, since each $D_{\Psi,\, \Xi,\, n}(A)$ is a $\bQ$-algebra in which $\xi^m$ vanishes for $m \ge n$ and each $X_a$ is a unit in $D_{\Psi,\, \Xi}(A)$, the universal relations \eqref{chart-rels} and \eqref{nonsm-rels} imply that each $(D_{\Psi,\, \Xi,\, n}(A))^\circ$, so also $D_{\Psi,\, \Xi}(A)$, is naturally an algebra over
\[
(A_{\Sigma,\, \Lambda}^\square \tensor_{A_\Inf} A_\cris^0) \tensor_{\bZ[Q]} \bZ[P_{\lambda_0}] \qq \text{for} \qq \lambda_0 \in \Lambda,
\]
compatibly with the maps 
\be\lab{maps-to-R}
\xymatrix@C=85pt{(A_{\Sigma,\, \Lambda}^\square \tensor_{A_\Inf} A_\cris^0) \tensor_{\bZ[Q]} \bZ[P_{\lambda_0}] \ar@{->>}[r]^-{\text{\eqref{R-nat-alg-1} and \eqref{R-nat-alg-2}}} & R} \q \text{and} \q (D_{\Psi,\, \Xi,\, n}(A))^\circ \ra A^\circ \cong R
\ee
and the ``change of $\lambda_0$'' isomorphisms \eqref{indep-0}. The resulting algebra structure map factors through some (necessarily $p$-adically complete) ring of definition $(D_{\Psi,\, \Xi,\, n}(A))_0$: 
\[ 
(A_{\Sigma,\, \Lambda}^\square \tensor_{A_\Inf} A_\cris^0) \tensor_{\bZ[Q]} \bZ[P_{\lambda_0}] \ra (D_{\Psi,\, \Xi,\, n}(A))_0 \hra  (D_{\Psi,\, \Xi,\, n}(A))^\circ \hra D_{\Psi,\, \Xi,\, n}(A),
\]
so the map $(D_{\Psi,\, \Xi,\, n}(A))_0 \ra  R$ is surjective. In addition, by \cite{SP}*{\href{http://stacks.math.columbia.edu/tag/07GM}{07GM}}, the kernel of the map $(D_{\Psi,\, \Xi,\, n}(A))^\circ \surjects R/p$ has a unique divided power structure, so we obtain a map
\be \lab{D-D-map}
D_{j_{\lambda_0}} \ra (D_{\Psi,\, \Xi,\, n}(A))^\circ
\ee
from the divided power envelope $D_{j_{\lambda_0}}$ defined in \S\ref{D-j-def}. For a fixed $n$ and modulo the $\f{\xi^m}{m!}$ with $m \ge n$, the kernel of the map $(A_{\Sigma,\, \Lambda}^\square \tensor_{A_\Inf} A_\cris^0) \tensor_{\bZ[Q]} \bZ[P_{\lambda_0}] \surjects R/p$ is finitely generated and, due to the surjectivity of $(A_{\Sigma,\, \Lambda}^\square \tensor_{A_\Inf} A_\cris^0) \tensor_{\bZ[Q]} \bZ[P_{\lambda_0}] \surjects R$, the generating set may be arranged to consist of $p$ and a finite set of elements that vanish already in $R$. Thus, since $D_{j_{\lambda_0}}$ is generated as an $((A_{\Sigma,\, \Lambda}^\square \tensor_{A_\Inf} A_\cris^0) \tensor_{\bZ[Q]} \bZ[P_{\lambda_0}])$-algebra by the divided powers of the elements of this kernel, by enlarging $(D_{\Psi,\, \Xi,\, n}(A))_0$ we may factor the map \eqref{D-D-map} as follows:
\[
D_{j_{\lambda_0}} \ra (D_{\Psi,\, \Xi,\, n}(A))_0 \hra (D_{\Psi,\, \Xi,\, n}(A))^\circ \hra D_{\Psi,\, \Xi,\, n}(A).
\]
Consequently, we obtain independent of $\lambda_0$ and compatible as $n$ varies continuous maps
\be \lab{D-D-no-n}
\wh{D_{j_{\lambda_0}}} \overset{\eqref{D-no-log-eq}}{\cong} D_{\Sigma,\, \Lambda} \ra (D_{\Psi,\, \Xi,\, n}(A))^\circ \hra D_{\Psi,\, \Xi,\, n}(A), \qq \text{so also} \qq D_{\Sigma,\, \Lambda} \ra D_{\Psi,\, \Xi}(A).
\ee
Via the last map, the $D_{j_{\lambda_0}}$-valued derivations $\f{\partial}{\partial \log(X_\sigma)}$ for $\sigma \in \Sigma$ and $\f{\partial}{\partial \log(X_{\lambda,\, i})}$ for $\lambda \in \Lambda$ and $1 \le i \le d$ of $D_{j_{\lambda_0}}$ are compatible with the corresponding $D_{\Psi,\, \Xi}(A)$-valued derivations of $D_{\Psi,\, \Xi}(A)$ (see \eqref{PX-def}, \S\ref{log-PD-der}, and \S\ref{BdR-more}). 
Thus, due to the density of $D_{j_{\lambda_0}}$ in $D_{\Sigma,\, \Lambda}$, the same compatibility holds for the map $D_{\Sigma,\, \Lambda} \ra D_{\Psi,\, \Xi}(A)$, 
to the effect that we obtain a map of complexes 
\[
\tst K_{D_{\Sigma,\, \Lambda}}\p{\p{\frac{\partial}{\partial\log(X_\sigma)}}_{\sigma\in\Sigma}, \p{\frac{\partial}{\partial\log(X_{\lambda,\, i})}}_{\lambda\in\Lambda,\, 1 \le i \le  d}} \ra  K_{D_{\Psi,\, \Xi}(A)}\p{(\f{\partial}{\partial \log(X_a)})_{a \in \Psi \cup \Xi}}.
\]
Its formation commutes with enlarging $\Sigma$ and $\Lambda$ (and, respectively, $\Psi$ and $\Xi$), so we obtain the~map
\be \lab{cris-BdR-cx}
\tst \varinjlim_{\Sigma,\, \Lambda} \p{K_{D_{\Sigma,\, \Lambda}}\p{\p{\frac{\partial}{\partial\log(X_\sigma)}}_{\sigma\in\Sigma}, \p{\frac{\partial}{\partial\log(X_{\lambda,\, i})}}_{\lambda\in\Lambda,\, 1 \le i \le  d}}} \ra
 \varinjlim_{\Psi,\, \Xi} \p{ \Omega^\bullet_{D_{\Psi,\, \Xi}(R[\f{1}{p}])/B_\dR^+}}.
\ee
The formation of this map is compatible with replacing $R$ by a $p$-adically formally \'{e}tale $R$-algebra $R'$ equipped with data as in \S\ref{all-possible-coordinates}. The resulting map of complexes of presheaves gives rise to the map of complexes of sheaves on $\fX_\et$ from the complex whose $R\Gamma(\fX_\et, -)$ is identified with $R\Gamma_{\log\cris}(\fX_{\cO_C/p}/A_\cris)$ (see \eqref{D-Koszul} and \S\ref{functorial-cris}) to the pushforward of the complex whose $R\Gamma((\fX_C^\ad)_\et, -)$ is identified with $R\Gamma_\cris (\fX^\ad_C /B_{\dR}^{+})$ (see \S\ref{BdR-more} and \eqref{informal}). Thus, by applying $R\Gamma(\fX_\et, -)$, we obtain the desired map \eqref{abs-to-BdR}:
\[
R\Gamma_{\log\cris}(\fX_{\cO_C/p}/A_\cris)\ra R\Gamma_\cris (\fX^\ad_C /B_{\dR}^{+}).
\]
In addition, by its construction and \Cref{D-no-log}, the map $D_{\Sigma,\, \Lambda} \ra D_{\Psi,\, \Xi}(R[\f{1}{p}])$ of \eqref{D-D-no-n} is compatible with the maps to $R[\f{1}{p}]$ (see \eqref{maps-to-R}).} Thus, \cite{BMS16}*{13.13} 
used to obtain \eqref{mod-xi-dR} implies 
that the map \eqref{cris-BdR-cx} is compatible with the maps in the derived category to $\Omega^{\bullet,\, \cont}_{R[\f{1}{p}]/C}$ described in the last display of \S\ref{functorial-cris} and \eqref{mod-xi-dR}. In conclusion, the map \eqref{abs-to-BdR} fits into the commutative square:
\be\ba \lab{cris-BdR-comp}
\xymatrix{
R\Gamma_{\log\cris}(\fX_{\cO_C/p}/A_\cris) \ar[d]^-{\eqref{cris-dR-map}} \ar[r]^-{\eqref{abs-to-BdR}} &  R\Gamma_\cris (\fX^\ad_C /B_{\dR}^{+}) \ar[d]^-{\eqref{BdR-mod-xi}} \\
R\Gamma_{\log\dR}(\fX/\cO_C)  \ar[r] &     R\Gamma_\dR(\fX_C^\ad/C).
}
\ea\ee
\epp

Having constructed the map \eqref{abs-to-BdR}, we are ready for the following extension of \cite{BMS16}*{13.23}.

\bthm \lab{BdR-comp}
If $\fX$ is $\cO_C$-proper, then the map \eqref{abs-to-BdR} induces the identification
\be \lab{BdR-comp-eq}
R\Gamma_{\log\cris}(\fX_{\cO_C/p}/A_\cris) \tensor^\bL_{A_\cris} B_\dR^+ \cong R\Gamma_\cris (\fX^\ad_C /B_{\dR}^{+})
\ee
and the cohomology modules of $R\Gamma_\cris(\fX_C^\ad/B_\dR^+)$ are finite free over $B_\dR^+$. In particular, then
\be \lab{BdR-comp-cor-eq}
R\Gamma(\fX_\et, A\Omega_\fX) \tensor_{A_\Inf}^\bL B_\dR^+ \cong R\Gamma_\cris (\fX^\ad_C /B_{\dR}^{+}),
\ee
compatibly with the identifications modulo $\xi$ with $R\Gamma_\dR(\fX_C^\ad/C)$ given by \eqref{G-dR-spec-eq} and \eqref{BdR-mod-xi}.
\ethm

\bpf
By \Cref{perfect,Acris-BC}, the object $R\Gamma_{\log\cris}(\fX_{\cO_C/p}/A_\cris)$ of $D(A_\cris)$ is perfect and its cohomology modules become finite free after inverting $p$. Therefore, due to \eqref{cris-dR-bc} and the derived $p$-adic completeness, we have the identification
\[
\xymatrix@C=36pt{ R\Gamma_{\log\cris}(\fX_{\cO_C/p}/A_\cris) \tensor^\bL_{A_\cris} \cO_C  \ar[r]_-{\sim}^-{\eqref{cris-dR-map}\,} & R\Gamma_{\log\dR}(\fX/\cO_C)}.
\]
Consequently, both sides of \eqref{BdR-comp-eq} are derived $\xi$-adically complete (see \eqref{BdR-coho-comp}) and, due to \eqref{BdR-mod-xi} and the commutativity of the diagram \eqref{cris-BdR-comp}, the map \eqref{abs-to-BdR} identifies their reductions modulo $\xi$. In conclusion, \eqref{abs-to-BdR} induces the desired identification \eqref{BdR-comp-eq} and the $B_\dR^+$-freeness claim follows from the first sentence of the proof. The combination of \eqref{ABCA-eq} and \eqref{BdR-comp-eq} gives \eqref{BdR-comp-cor-eq} and the asserted compatibility follows from \Cref{two-dR-same} and the commutativity of \eqref{cris-BdR-comp}.
\epf

\bpp[The $B_\dR^+$-cohomology and the \'{e}tale cohomology]
\ready{For any proper, smooth adic space $X$ over $C$, in \cite{BMS16}*{13.1} Bhatt--Morrow--Scholze constructed the functorial in $X$ identification
\be \lab{BdR-et-id}
R\Gamma_\cris(X/B_\dR^+) \tensor_{B_\dR^+} B_\dR \cong R\Gamma_\et(X, \bZ_p) \tensor_{\bZ_p} B_\dR.
\ee
 Due to the identification \eqref{cris-dR-descent}, when $X \cong X_0 \wh{\tensor}_K C$ for a proper, smooth adic space $X_0$ defined over a complete, discretely valued subfield $K \subset C$ that has a perfect residue field, the inverse of \eqref{BdR-et-id} supplies the functorial in $X_0$ de Rham comparison isomorphism
\be \lab{dR-comp-iso}
R\Gamma_\et(X_0 \wh{\tensor}_K C, \bZ_p) \tensor_{\bZ_p} B_\dR \cong R\Gamma_\dR(X_0/K) \tensor_{K} B_\dR.
\ee
If $C \cong \wh{\ov{K}}$, then, by transport of structure, the identification \eqref{dR-comp-iso} is $\Gal(\ov{K}/K)$-equivariant (with $\Gal(\ov{K}/K)$ acting trivially on $R\Gamma_\dR(X_0/K)$) and, by \emph{loc.~cit.},~it recovers the isomorphism constructed in \cite{Sch13}*{8.4}. In particular, in this case, \eqref{dR-comp-iso} is compatible with filtrations, where $B_\dR$ is filtered by its discrete valuation and $R\Gamma_\dR(X_0/K)$ (resp.,~$R\Gamma_\et(X_0 \wh{\tensor}_K C, \bZ_p)$) is equipped with the the Hodge (resp.,~trivial) filtration.}
\epp

\ready{For proper $\fX$, we now have two ways to identify $R\Gamma(\fX_\et, A\Omega_\fX) \tensor^\bL_{A_\Inf} B_\dR$ with $R\Gamma_\et(\fX_C^\ad, \bZ_p) \tensor_{\bZ_p}^\bL B_\dR$: we can either base change \eqref{RG-et-id-eq} to $B_\dR$ or combine \eqref{BdR-comp-cor-eq} and \eqref{BdR-et-id}. We now prove that the two ways give the same identification; this will be important in the proof of \Cref{lat-dR-id}.}

\bprop \lab{BdR-et-comp}
If $\fX$ is $\cO_C$-proper, then the map $R\Gamma_{\cris} (\fX^\ad_C / B_{\dR}^{+}) \ra R\Gamma (\fX^\ad_C, \bZ_p)\otimes_{\bZ_p} B_{\dR}^{+}$ of \cite{BMS16}*{proof of 13.1} that underlies the identification \eqref{BdR-et-id} for $X = \fX_C^\ad$ makes the diagram
\[
\xymatrix{
R\Gamma_{\log\cris}(\fX_{\cO_C/p}/ A_\cris) \ar[d]_-{\wr}^-{\eqref{main-comp-iso}} \ar[rrr]^-{\eqref{abs-to-BdR}} &&&  R\Gamma_{\cris} (\fX^\ad_C / B_{\dR}^{+})  \ar[d] \\
R\Gamma (\fX_\et, A\Omega_\fX) \otimes^\bL_{A_\Inf} A_\cris \ar[rr]^-{\text{\cite{BMS16}*{6.10}\,\,}} && R\Gamma_\et(\fX_C^\ad, \bA_{\Inf,\, \fX_C^\ad}) \tensor^\bL_{A_\Inf} B_\dR^+  &R\Gamma (\fX^\ad_C, \bZ_p)\otimes^\bL_{\bZ_p} \ar[l]^-{\sim}_-{\,\,\eqref{RG-under-map}} B_{\dR}^{+}. 
}
\]
commute\uscolon in particular, the identification of $R\Gamma(\fX_\et, A\Omega_\fX) \tensor^\bL_{A_\Inf} B_\dR$ with $R\Gamma_\et(\fX_C^\ad, \bZ_p) \tensor_{\bZ_p}^\bL B_\dR$ that results from \eqref{RG-et-id-eq} \up{and is encoded by the bottom part of the above diagram} agrees with the identification that results from \eqref{BdR-comp-cor-eq} and \eqref{BdR-et-id} \up{and is encoded by the top part of the diagram}.
\eprop

\bpf
\ready{Since $\varphi^{-1}(\mu)$ lies in $W(\fm^\flat)$ and is a unit in $B_\dR^+$, the discussion after \Cref{RG-et-id} implies that the map labeled ``\eqref{RG-under-map}'' in the diagram is an isomorphism. In particular, due to \cite{Sch13}*{5.1}, the object $R\Gamma_\et(\fX_C^\ad, \bA_{\Inf,\, \fX_C^\ad}) \tensor^\bL_{A_\Inf} B_\dR^+$ of $D(B_\dR^+)$ is perfect. We will now review the definition given in \cite{BMS16}*{proof of 13.1} of the composition $f$ of the right vertical map with this map ``\eqref{RG-under-map}.'' 

Let $\Spa(A, A^\circ)$ be an element of the basis for the analytic topology of $\fX_C^\ad$ discussed in \S\ref{BdR-et}. For a large enough set $\Psi$ as in \S\ref{BdR-et}, we consider the surjection $C\In{ (X_u^{\pm 1})_{u \in \Psi} } \xra{X_u \mapsto u\,} A$ from \eqref{surj-onto-R}, as well as the perfectoid $(\prod_\Psi \bZ_p(1))$-cover $C\In{ (X_u^{\pm 1/p^\infty})_{u \in \Psi} }$ of $C\langle (X_u^{\pm 1})_{u \in \Psi} \rangle$. Granted that $\Psi$ contains the images of the $T_i$ under some \'{e}tale coordinate map \eqref{affd-et}, the base change of this cover to $\Spa(A, A^\circ)$ is a perfectoid $(\prod_\Psi \bZ_p(1))$-cover 
\be \lab{Psi-cover}
\Spa(A_{\Psi,\, \infty}, A_{\Psi,\, \infty}^+) \ra \Spa(A, A^\circ).
\ee
Each $u \in \Psi$ has a canonical system $u^{1/p^\infty}$ of $p$-power roots in $A_{\Psi,\, \infty}^+$, which gives the unit $[u^{1/p^\infty}]$ in the $B_\dR^+$-algebra $\bB_\dR^+(A_{\Psi,\, \infty}^+)$ (see \Cref{Acrism-inj}). Since $\bB_\dR^+(A_{\Psi,\, \infty}^+)$ may be viewed as a pro-(infinitesimal thickening) of $A_{\Psi,\, \infty}$, the map $X_u \mapsto [u^{1/p^\infty}]$ extends to a continuous $B_\dR^+$-morphism 
\be \lab{DPsi-BdR}
D_\Psi(A) \ra \bB_\dR^+(A_{\Psi,\, \infty}^+) \qq \text{over} \qq A \ra A_{\Psi,\, \infty}.
\ee
By construction, for each $u\in \Psi$, this morphism intertwines $\exp\p{\log([\eps]) \cdot \f{\partial}{ \partial \log(X_u)}}$ defined by the formula \eqref{exp-def} and viewed as a ring endomorphism of $D_\Psi(A)$ 
with the action of the generator $[\eps]$ of the $u$-th copy of $\bZ_p(1)$ on $\bB_\dR^+(A_{\Psi,\, \infty}^+)$. In particular, letting $\gamma_u$ denote this generator, we may use the same formula as in \eqref{local-isom-1} to define the morphism of complexes
\be \lab{131-map}
\tst \Omega^\bullet_{D_\Psi(A)/B_\dR^+} = K_{D_\Psi(A)}\p{(\f{\partial}{\partial \log(X_u)})_{u \in \Psi}} \ra K_{\bB_\dR^+(A_{\Psi,\, \infty}^+)}((\gamma_u - 1)_{u \in \Psi}),
\ee
whose formation is functorial in $\Psi$ and, after passing to the direct limit over all $\Psi$, also in $\Spa(A, A^\circ)$. 
The almost purity theorem identifies the  cohomology of the sheaf of complexes determined by the target of \eqref{131-map} with $R\Gamma_\et(\fX_C^\ad, \bA_{\Inf,\, \fX_C^\ad}) \tensor^\bL_{A_\Inf} B_\dR^+$ (see \emph{loc.~cit.}). 
The cohomology of the sheaf of complexes determined by the source of \eqref{131-map} is, by definition, $R\Gamma_{\cris} (\fX^\ad_C / B_{\dR}^{+})$ (see \S\ref{BdR-et}). Therefore, by passing to the direct limit over all $\Psi$, sheafifying, and forming cohomology, the maps \eqref{131-map} produce the aforementioned composition $f$ defined in \emph{loc.~cit.}

The same construction gives the morphisms \eqref{131-map} for the objects $\Spa(A, A^\circ)$ of the basis of the \'{e}tale topology  of $\fX_C^\ad$ considered in \S\ref{BdR-et}. Due to \eqref{Zar-et-BdR}, this leads to the same map $f$. In addition, we may generalize the construction of the morphisms \eqref{131-map} further by using the basis for the \'{e}tale topology of $\fX_C^\ad$ considered in \S\ref{BdR-more}: the cover \eqref{Psi-cover} gets replaced by the cover
\[
\Spa(A_{\Psi,\, \Xi,\, \infty}, A_{\Psi,\, \Xi,\, \infty}^+) \ra \Spa(A, A^\circ)
\]
that is the base change of the perfectoid $(\prod_\Psi \bZ_p(1) \times \prod_\Xi \bZ_p(1))$-cover $C\langle (X_u^{\pm 1/p^\infty})_{u \in \Psi}, (X_a^{1/p^\infty})_{a \in \Xi} \rangle$ of $C\langle (X_u^{\pm 1})_{u \in \Psi}, (X_a)_{a \in \Xi} \rangle$ for large enough $\Psi \subset (A^\circ)^\times$ and $\Xi \subset A^\circ \cap A^\times$, and the rest is (mildly) modified accordingly. Due to \eqref{informal}, this variant of the construction gives the same map $f$.

In conclusion, since the construction of $f$ may be carried out in the setting of \S\ref{BdR-more} and follows the same pattern as the construction of the map \eqref{main-comp-iso}, namely, is based on the map as in \eqref{local-isom-1}, all we need to check is that, in the notation of \S\ref{map-from-abs}, the following diagram commutes:
\be \lab{D-D-com-d} \ba
\xymatrix@C=40pt{
D_{\Sigma,\, \Lambda} \ar[d]_{\eqref{CM-1}} \ar[r]^-{\eqref{D-D-no-n}} & D_{\Psi,\, \Xi}(A) \ar[d]^-{\eqref{DPsi-BdR}} \\
\bA_\cris(R_{\Sigma,\, \Lambda,\, \infty}) \ar[r]^-{\eqref{AI-eq}} & \bB_\dR^+(A_{\Psi,\, \Xi,\, \infty}^+).
}
\ea\ee
For this desired commutativity, we may first replace $\bB_\dR^+(A_{\Psi,\, \Xi,\, \infty}^+)$ by $\bB_\dR^+(A_{\Psi,\, \Xi,\, \infty}^+)/\xi^n$ for a variable $n > 0$, then replace $D_{\Sigma,\, \Lambda}$ by $D_{j_{\lambda_0}}$ for some $\lambda_0 \in \Lambda$, and, finally, since $\bB_\dR^+(A_{\Psi,\, \Xi,\, \infty}^+)/\xi^n$ is a $\bQ$-algebra and $D_{j_{\lambda_0}}$ is generated by divided powers, replace $D_{j_{\lambda_0}}$ by $(A_{\Sigma,\, \Lambda}^\square \tensor_{A_\Inf} A_\cris^0) \tensor_{\bZ[Q]} \bZ[P_{\lambda_0}]$.} However, each $X_{\tau}$ of \eqref{A-box-def} with either $\tau = \sigma$ for $\sigma \in \Sigma$ or $\tau = (\lambda, i)$ for $\lambda \in \Lambda$ and $1 \le i \le d$  maps to the unit $[X_\tau^{1/p^\infty}] \in (\bB_\dR^+(A_{\Psi,\, \Xi,\, \infty}^+))^\times$ 
under either of the two maps from $(A_{\Sigma,\, \Lambda}^\square \tensor_{A_\Inf} A_\cris^0) \tensor_{\bZ[Q]} \bZ[P_{\lambda_0}]$ to $\bB_\dR^+(A_{\Psi,\, \Xi,\, \infty}^+)/\xi^n$ supplied by the diagram \eqref{D-D-com-d}, so these two maps indeed agree, as desired. 
\epf
}}

%% file: RG.tex

\section{The $A_\Inf$-cohomology modules $H^i_{A_\Inf}(\fX)$ and their specializations \nopunct} \lab{RG}

\revise{
\ready{In this section, we define and analyze the $A_\Inf$-cohomology groups $H^i_{A_\Inf}(\fX)$ of an $\cO_C$-proper $\fX$. We show that each $H^i_{A_\Inf}(\fX)$ is a Breuil--Kisin--Fargues module (see \Cref{Hinf-BKF}) and deduce that, loosely speaking, the $p$-adic \'{e}tale cohomology of $\fX_C^\ad$ has at most the amount of torsion contained in the logarithmic crystalline cohomology of $\fX_k$ or the  logarithmic de Rham cohomology of $\fX$ (see \Cref{et-cris-tors,et-dR-tors}). Most of these results are variants of their analogues established in the smooth case in \cite{BMS16}. Their proofs, granted inputs from \S\ref{etale-section} and \S\S\ref{AOmegaX}--\ref{section-Acris},  are generally similar to those of \emph{op.~cit.} and in large part rely on commutative algebra over $A_\Inf$.

\bpp[Properness of $\fX$]
Throughout \S\ref{RG}, we assume that $\fX$ is proper and $\fX_k$ is purely $d$-dimensional.
\epp

\bpp[The $A_\Inf$-cohomology $R\Gamma_{A_\Inf}(\fX)$] \lab{Ainf-coho}
\ready{We use the object $A\Omega_\fX \in D^{\ge 0}(\fX_\et, A_\Inf)$ of \S\ref{AOX-def} to set
\[
R\Gamma_{A_\Inf}(\fX) \ce R\Gamma(\fX_\et, A\Omega_\fX) \in D^{\ge 0}(A_\Inf) \q \text{and} \q H^i_{A_\Inf}(\fX) \ce H^i(R\Gamma(\fX_\et, A\Omega_\fX)) \q \text{for $i \in \bZ$.}
\]
Since $L\eta$ commutes with pullback along a flat morphism of ringed topoi (see \cite{BMS16}*{6.14}), the object $R\Gamma_{A_\Inf}(\fX)$ is contravariantly functorial in $\fX$: an $\cO_C$-morphism $\fX' \ra \fX$ induces a morphism 
\[
R\Gamma_{A_\Inf}(\fX) \ra R\Gamma_{A_\Inf}(\fX') \q \text{in} \q D^{\ge 0}(A_\Inf), \qq \text{so also} \qq  H^i_{A_\Inf}(\fX) \ra H^i_{A_\Inf}(\fX') \q \text{for} \q i \in \bZ.
\]
\Cref{perfect} ensures that $R\Gamma_{A_\Inf}(\fX)$ is perfect, that is, isomorphic to a bounded complex of finite free $A_\Inf$-modules. Moreover, by \eqref{RG-et-id}, \eqref{G-dR-spec-eq}, and \eqref{ABCA-eq}, we have the following identifications:
\be \lab{spec-isos} \ba
\tst R\Gamma_{A_\Inf}(\fX) \tensor_{A_\Inf}^\bL A_\Inf[\f{1}{\mu}]  &\cong \tst R\Gamma_\et(\fX_C^\ad, \bZ_p) \tensor_{\bZ_p}^\bL A_\Inf[\f{1}{\mu}]; \\
 R\Gamma_{A_\Inf}(\fX) \tensor^\bL_{A_\Inf,\, \theta} \cO_C &\cong R\Gamma_{\log\dR}(\fX/\cO_C); \\
 R\Gamma_{A_\Inf}(\fX) \tensor^\bL_{A_\Inf} W(k) &\cong R\Gamma_{\log\cris}(\fX_k/W(k)).
\ea\ee
If $\fX_k$ is $k$-smooth, then we may drop ``$\log$'' from the subscripts (compare with \eqref{sm-no-log}).}

The Frobenius morphism \eqref{Frob-mor} gives rise to the Frobenius morphism
\[
\tst R\Gamma_{A_\Inf}(\fX) \tensor_{A_\Inf,\, \varphi} A_\Inf \ra R\Gamma_{A_\Inf}(\fX) \qq \text{in} \qq D^{\ge 0}(A_\Inf)
\]
that becomes an isomorphism after inverting $\varphi(\xi)$ (see \eqref{Frob-iso}). 
Consequently the cohomology modules $H^i_{A_\Inf}(\fX)$ come equipped with the $A_\Inf$-module morphism
\[
\tst \varphi\colon H^i_{A_\Inf}(\fX) \tensor_{A_\Inf,\, \varphi} A_\Inf \ra H^i_{A_\Inf}(\fX) 
\]
that becomes an isomorphism after inverting $\varphi(\xi)$. We will prove in \Cref{Hinf-BKF} that these morphisms make each $H^i_{A_\Inf}(\fX)$ a Breuil--Kisin--Fargues module in the following sense of \cite{BMS16}*{4.22}.
\epp

\bpp[Breuil--Kisin--Fargues modules] \lab{BKF-mod}
A \emph{Breuil--Kisin--Fargues module} is a finitely presented $A_\Inf$-module $M$ equipped with an $A_\Inf[\f{1}{\varphi(\xi)}]$-module isomorphism
\[
\tst \varphi_M \colon (M \tensor_{A_\Inf,\, \varphi} A_\Inf)[\f{1}{\varphi(\xi)}] \isomto M[\f{1}{\varphi(\xi)}]
\]
such that $M[\f{1}{p}]$ is $A_\Inf[\f{1}{p}]$-free. By \cite{BMS16}*{4.9~(i)}, any such $M$ is perfect as an $A_\Inf$-module, that is, $M$ has a finite resolution by finite free $A_\Inf$-modules. A morphism of Breuil--Kisin--Fargues modules is an $A_\Inf$-module morphism that commutes with the isomorphisms~$\varphi_M$. 
\epp

\bthm \lab{Hinf-BKF}
Each $(H^i_{A_\Inf}(\fX), \varphi)$ is a Breuil--Kisin--Fargues module and vanishes unless $i \in [0, 2d]$. In particular, each $H^i_{A_\Inf}(\fX)$ is perfect as an $A_\Inf$-module and each $(H^i_{A_\Inf}(\fX))[\f{1}{p}]$ is $A_\Inf[\f{1}{p}]$-free. 
\ethm

\bpf
Due to the relation with $R\Gamma_\et(\fX_C^\ad, \bZ_p)$, each $(H^i_{A_\Inf}(\fX))[\f{1}{p\mu}]$ is a free $A_\Inf[\f{1}{p\mu}]$-module. Moreover, by \Cref{Acris-BC}, the cohomology modules of $R\Gamma_{A_\Inf}(\fX) \tensor^\bL_{A_\Inf} A_\cris[\f{1}{p}]$ are free over $A_\cris[\f{1}{p}]$. Therefore, \cite{BMS16}*{4.20} applies and proves that each $H^i_{A_\Inf}(\fX)$ is a finitely presented $A_\Inf$-module that becomes free after inverting $p$, so $(H^i_{A_\Inf}(\fX), \varphi)$ is a Breuil--Kisin--Fargues module.

Since $R\Gamma_{A_\Inf}(\fX)$ is perfect, its top degree cohomology is finitely presented and of formation compatible with base change. Thus, by the de Rham specialization of \eqref{spec-isos} and the Nakayama lemma, $H^i_{A_\Inf}(\fX) = 0 $ for $i > 2d$. The same holds for $i < 0$ because $R\Gamma_{A_\Inf}(\fX) \in D^{\ge 0}(A_\Inf)$.
\epf

 For completeness sake, we mention the following corollary, which may also be proved more directly.

\bcor \lab{same-rk}
For each $i \in \bZ$, the rank of the finitely presented $\bZ_p$-module $H^i_\et(\fX_C^\ad, \bZ_p)$ is equal to the rank of the finitely presented $W(k)$-module $H^i_{\log\cris}(\fX_k/W(k))$, and is also equal to the rank of the finitely presented $\cO_C$-module $H^i_{\log \dR}(\fX/\cO_C) \ce R^i\Gamma(\fX_\et, \Omega^\bullet_{\fX/\cO_C,\, \log})$ \up{see also \eqref{str-thm-iso} below}.
\ecor

\bpf
The finite presentation assertions follow, for instance, from the perfectness of $R\Gamma_{A_\Inf}(\fX)$, the comparisons \eqref{spec-isos}, and the coherence of the ring $\cO_C$. Due to \Cref{Hinf-BKF} and the comparisons \eqref{spec-isos}, all the ranks in question are equal to the rank of the free $A_\Inf[\f{1}{p}]$-module $(H^i_{A_\Inf}(\fX))[\f{1}{p}]$.
\epf

\bpp[Base change for individual $H^i_{A_\Inf}(\fX)$] \lab{indiv-BC}
Since $A_\Inf[\f{1}{\mu}]$ is $A_\Inf$-flat, \eqref{spec-isos} implies that
\be \lab{et-iso-Hi}
\tst (H^i_{A_\Inf}(\fX))[\f{1}{\mu}] \cong H^i_\et(\fX_C^\ad, \bZ_p) \tensor_{\bZ_p} A_\Inf[\f{1}{\mu}] \q \text{for each} \q i \in \bZ.
\ee
In particular, since $\mu$ is a unit in $W(C^\flat)$ and $W(C^\flat)$ is $A_\Inf$-flat (the localization of $A_\Inf$ at the prime ideal $(p)$ is a discrete valuation ring whose completion is $W(C^\flat)$)
\be \lab{WCflat-iso}
H^i_{A_\Inf}(\fX) \tensor_{A_\Inf} W(C^\flat) \cong H^i_\et(\fX_C^\ad, \bZ_p) \tensor_{\bZ_p} W(C^\flat).
\ee
A similar de Rham comparison consists of exact sequences that result from \eqref{spec-isos} and \cite{SP}*{\href{http://stacks.math.columbia.edu/tag/0662}{0662}}: 
\be \lab{dR-bc-SES}
0 \ra H^i_{A_\Inf}(\fX) \tensor_{A_\Inf,\, \theta} \cO_C \ra H^i_{\log \dR}(\fX/\cO_C) \ra (H^{i + 1}_{A_\Inf}(\fX))[\xi] \ra 0 \q \text{for each} \q i \in \bZ.
\ee
Similarly, by \Cref{Hinf-BKF} and \cite{BMS16}*{4.9}, we have a Frobenius-equivariant exact sequence 
\be \lab{cris-bc-SES}
0 \ra H^i_{A_\Inf}(\fX) \tensor_{A_\Inf} W(k) \ra H^i_{\log \cris}(\fX_k/W(k)) \ra \Tor^1_{A_\Inf}(H^{i + 1}_{A_\Inf}(\fX), W(k)) \ra 0
\ee
for each $i \in \bZ$. In particular, we have the top degree base changes
\[
H^{2d}_{A_\Inf}(\fX) \tensor_{A_\Inf,\, \theta} \cO_C \cong H^{2d}_{\log \dR}(\fX/\cO_C) \qq \text{and} \qq H^{2d}_{A_\Inf}(\fX) \tensor_{A_\Inf} W(k) \cong H^{2d}_{\log \cris}(\fX_k/W(k)).
\]
Due to \Cref{Hinf-BKF}, the injections in the sequences \eqref{dR-bc-SES}--\eqref{cris-bc-SES} become isomorphisms after inverting $p$. The same holds without inverting $p$ in the case when $H^{i + 1}_{A_\Inf}(\fX)$ is $A_\Inf$-free. For such freeness, we have the following consequence of \Cref{Hinf-BKF} and \cite{BMS16}*{\S4.2}.
\epp

\bprop \lab{free-free}
For each $i \in \bZ$, the $\cO_C$-module $H^i_{\log \dR}(\fX/\cO_C)$ is $p$-torsion free \up{equivalently, free} if and only if the $W(k)$-module $H^i_{\log \cris}(\fX_k/W(k))$ is $p$-torsion free \up{equivalently, free}, in which case $H^i_{A_\Inf}(\fX)$ is free as an $A_\Inf$-module and $H^i_\et(\fX_C^\ad, \bZ_p)$ is free as a $\bZ_p$-module.
\eprop

\bpf
Due to \Cref{Hinf-BKF}, we may apply \cite{BMS16}*{4.18} and combine it with \eqref{spec-isos} to conclude that $H^i_{\log \dR}(\fX/\cO_C)$ is $p$-torsion free if and only if so is $H^i_{\log \cris}(\fX_k/W(k))$. When these conditions hold, the freeness of $H^i_{A_\Inf}(\fX)$ and $H^i_\et(\fX_C^\ad, \bZ_p)$ follows from \cite{BMS16}*{4.17} and \eqref{et-iso-Hi}. 
\epf



\brem
As was observed by Jesse Silliman and Ravi Fernando during the Arizona Winter School 2017, the first assertion of \Cref{free-free} may be strengthened as follows: for each $i \in \bZ$, 
\be \lab{same-dim-k}
\dim_k \p{H^i_{\log \dR}(\fX/\cO_C)_\tors \tensor_{\cO_C} k} = \dim_k \p{H^i_{\log \cris}(\fX_k/W(k))_\tors \tensor_{W(k)} k}, 
\ee
that is, $H^i_{\log \dR}(\fX/\cO_C)$ and $H^i_{\log \cris}(\fX_k/W(k))$ have the same number of cyclic summands (in the sense of \eqref{str-thm-iso} below). Indeed, by \Cref{same-rk}, the ranks of $H^i_{\log \dR}(\fX/\cO_C)$ and $H^i_{\log \cris}(\fX_k/W(k))$ agree and, by \cite{Bei13a}*{(1.8.1)}, so do the $k$-fibers of $R\Gamma_{\log\dR}(\fX/\cO_C)$ and $R\Gamma_{\log\cris}(\fX_k/W(k))$, so the claim follows by descending induction on $i$ from the following exact sequences  supplied by \cite{SP}*{\href{http://stacks.math.columbia.edu/tag/0662}{0662}}:
\[\ba
0 \ra H^i_{\log\dR}(\fX/\cO_C) \tensor_{\cO_C} k \ra &H^i(R\Gamma_{\log\dR}(\fX/\cO_C) \tensor^\bL_{\cO_C} k) \ra \Tor_1^{\cO_C}(H^{i + 1}_{\log\dR}(\fX/\cO_C), k) \ra 0, \\
0 \ra H^i_{\log\cris}(\fX_k/W(k)) \tensor_{W(k)} k \ra &H^i(R\Gamma_{\log\cris}(\fX_k/W(k)) \tensor^\bL_{W(k)} k) \ra H^{i + 1}_{\log\cris}(\fX_k/W(k))[p] \ra 0.
\ea \]
\erem

The following variant of \cite{BMS16}*{14.5~(ii)} strengthens the relationship between the freeness of $H^i_\et(\fX_C^\ad, \bZ_p)$ and that of $H^i_{\log \cris}(\fX_k/W(k))$ supplied by \Cref{free-free}.

\bthm \lab{et-cris-tors}
For every $i \in \bZ$ and $n \in \bZ_{\ge 0}$, we have 
\be \lab{ECT-eq}\ba
\length_{\bZ_p}((H^i_\et(\fX_C^\ad, \bZ_p)_\tors)/p^n) &\le \length_{W(k)}((H^i_{\log \cris}(\fX_k/W(k))_\tors)/p^n), \\
\length_{\bZ_p}(H^i_\et(\fX_C^\ad, \bZ/p^n\bZ)) &\le \length_{W(k)}(H^i_{\log \cris}(\fX_k/W_n(k))).
\ea\ee
\ethm

\bpf
The proof of the first inequality is analogous to the proof of \emph{loc.~cit.} Namely, by \Cref{same-rk}, we may drop the subscripts ``$\tors$'' and, by \Cref{Hinf-BKF}, \eqref{WCflat-iso}, and \cite{BMS16}*{4.15~(ii)}, we have
\be \lab{ECT-1}
\length_{\bZ_p}(H^i_\et(\fX_C^\ad, \bZ_p)/p^n) \le \length_{W(k)}((H^i_{A_\Inf}(\fX) \tensor_{A_\Inf} W(k))/p^n).
\ee
Since $\length_{W(k)}(Q/p^n) = \length_{W(k)}(\Tor_1^{W(k)}(Q, W(k)/p^n))$ for every $W(k)$-module $Q$ that is finite and torsion, the short exact sequence \eqref{cris-bc-SES} yields the inequality
\[
\length_{W(k)}((H^i_{A_\Inf}(\fX) \tensor_{A_\Inf} W(k))/p^n) \le \length_{W(k)}(H^i_{\log\cris}(\fX_k/W(k))/p^n),
\]
and the first inequality in \eqref{ECT-eq} follows. Due to the short exact sequences 
\[ \ba
0 \ra H^i_\et(\fX_C^\ad, \bZ_p)/p^n \ra &H^i_\et(\fX_C^\ad, \bZ/p^n\bZ) \ra (H^{i + 1}_\et(\fX_C^\ad, \bZ_p))[p^n] \ra 0,\\
0 \ra H^i_{\log\cris}(\fX_k/W(k))/p^n \ra &H^i_{\log\cris}(\fX_k/W_n(k)) \ra (H^{i + 1}_{\log\cris}(\fX_k/W(k)))[p^n] \ra 0
\ea
\]
that result from  \cite{SP}*{\href{http://stacks.math.columbia.edu/tag/0662}{0662}}, 
the second inequality in \eqref{ECT-eq} follows from the first.
\epf

The de Rham analogue of \Cref{et-cris-tors} is \Cref{et-dR-tors} below and uses the following formalism.

\bpp[The normalized length] \lab{norm-length}
Let $\fo$ be a valuation ring of rank $1$ and mixed characteristic $(0, p)$. We normalize its valuation $\val_\fo$ by requiring that $\val_\fo(p) = 1$. By the structure theorem \cite{SP}*{\href{http://stacks.math.columbia.edu/tag/0ASP}{0ASP}} (see also \cite{GR03}*{6.1.14}), every finitely presented $\fo$-module $M$ is of the form
\be \lab{str-thm-iso}
\tst M \cong \bigoplus_{i = 1}^n \fo/(a_i) \qq \text{with} \qq a_i \in \fo. 
\ee
If $M$ is, in addition, torsion, to the effect that the $a_i$ are nonzero, then we set 
\[
\tst \val_\fo(M) \ce \sum_{i = 1}^n \val(a_i).
\]
More intrinsically, $\val_\fo(M)$ is the valuation of any generator of the $0\th$ Fitting ideal $\Fitt_0(M) \subset \fo$ of $M$, so it depends only on $M$. If $\fo$ is a discrete valuation ring for which $p$ is a uniformizer, then $\val_\fo(M) = \length_\fo(M)$. In general, $\val_\fo$ has the advantage of being invariant under the extension of scalars to a larger $\fo$. Any short exact sequence 
\[
0 \ra M_1 \ra M_2 \ra M_3 \ra 0
\]
of finitely presented, torsion $\fo$-modules gives rise to the equality $\Fitt_0(M_2) = \Fitt_0(M_1)\Fitt_0(M_3)$ (see \cite{GR03}*{6.3.1 and 6.3.5 (i)}), so the assignment $\val_\fo(-)$ satisfies
\be \lab{norm-length-add}
\val_\fo(M_2) = \val_\fo(M_1) + \val_\fo(M_3).
\ee
\epp

The following lemma is the de Rham version of \cite{BMS16}*{4.14}, which gave the inequality \eqref{ECT-1}.

\blem \lab{spec-ineq}
For a finitely presented $W_n(\cO_C^\flat)$-module $M$ for some $n \ge 1$, we have
\be \lab{SS-eq}
\val_{W(C^\flat)}(M \tensor_{A_\Inf} W(C^\flat)) = \val_{\cO_C}(M/\xi M) - \val_{\cO_C}(M [\xi]).
\ee
\elem

\bpf
Since the ring $W_n(\cO_C^\flat)$ is coherent (see \cite{BMS16}*{3.24}), the $W_n(\cO_C^\flat)$-module $M[\xi]$ is finitely presented. Moreover, due to \eqref{norm-length-add}, the flatness of $A_\Inf \ra W(C^\flat)$ (see \S\ref{indiv-BC}), and the snake lemma, both sides of \eqref{SS-eq} are additive in short exact sequences. Therefore, we may assume that $n = 1$ and, due to the structure theorem \cite{SP}*{\href{http://stacks.math.columbia.edu/tag/0ASP}{0ASP}}, that $M = \cO_C^\flat/(x)$ for some $x \in \cO_C^\flat$. 

If $x = 0$, then both sides of \eqref{SS-eq} are equal to $1$. If $x \neq 0$, then the left side vanishes, and so does the right side because $M[\xi] \cong \Tor^1_{\cO_C^\flat}(M, \cO_C/p)$ and the following sequence is exact:
\[
0 \ra \Tor^1_{\cO_C^\flat}(\cO_C^\flat/(x), \cO_C/p) \ra \cO_C/p \xra{\theta([x])} \cO_C/p \ra M/\xi M \ra 0. \qedhere
\]
\epf

\bthm \lab{et-dR-tors}
For every $i \in \bZ$ and $n \in \bZ_{\ge 0}$, we have \up{recall from \uS\uref{norm-length} that $\val_{\bZ_p} = \length_{\bZ_p}$}
\be \lab{EDT-eq}\ba
\val_{\bZ_p}((H^i_\et(\fX_C^\ad, \bZ_p)_\tors)/p^n) &\le \val_{\cO_C}((H^i_{\log \dR}(\fX/\cO_C)_\tors)/p^n), \\
\val_{\bZ_p}(H^i_\et(\fX_C^\ad, \bZ/p^n\bZ)) &\le \val_{\cO_C}(R^i\Gamma(\fX_{\cO_C/p^n,\, \et}, \Omega^\bullet_{\fX_{\cO_C/p^n}/(\cO_C/p^n),\,\log})).
\ea\ee
\ethm

\bpf
The proof is analogous to that of \Cref{et-cris-tors}. Namely, by \Cref{same-rk}, we may drop the subscripts ``$\tors$'' and, by \Cref{Hinf-BKF}, \eqref{WCflat-iso}, and \Cref{spec-ineq}, we have
\[
\val_{\bZ_p}(H^i_\et(\fX_C^\ad, \bZ_p)/p^n) \le \val_{\cO_C}(H^i_{A_\Inf}(\fX)/(p^n, \xi)).
\]
The presentation \eqref{str-thm-iso} implies that $\val_{\cO_C}(Q/p^n) = \val_{\cO_C}(\Tor_1^{\cO_C}(Q, \cO_C/p^n))$ for every finitely presented, torsion $\cO_C$-module $Q$, so the short exact sequence \eqref{dR-bc-SES} yields the inequality
\[
\val_{\cO_C}(H^i_{A_\Inf}(\fX)/(p^n, \xi)) \le \val_{\cO_C}(H^i_{\log \dR}(\fX/\cO_C)/p^n).
\]
This proves the first inequality in \eqref{EDT-eq} and, analogously to the proof of \Cref{et-cris-tors}, the second one follows from the first.
\epf

The results above, specifically, \eqref{same-dim-k} and \Cref{et-cris-tors,et-dR-tors} prompt the following question.

\bq
Are there examples of $\cO_C$-proper $\fX$ satisfying the assumptions of \uS\uref{fir-setup} for which
\[
\val_{W(k)} (H^i_{\log\cris}(\fX/W(k))_\tors) \neq \val_{\cO_C}(H^i_{\log\dR}(\fX/\cO_C)_\tors)?
\]
\eq

}}


%% file: lattice.tex
\section{A functorial lattice inside the de Rham cohomology \nopunct}  \lab{lattice}

\revise{
\ready{To a proper, smooth scheme $X$ over a complete, discretely valued extension $K$ of $\bQ_p$ with a perfect residue field, in \Cref{lat-main-eg} we functorially associate an $\cO_K$-lattice 
\[
L^i_\dR(X) \subset H^i_\dR(X/K) \qq \text{for every} \qq i \in \bZ.
\]
In fact, $L^i_\dR(X)$ functorially depends only on $H^i_\et(X_{\ov{K}}, \bZ_p)$ and its construction, which relies on the theory of Breuil--Kisin--Fargues modules, proceeds along familiar lines of integral $p$-adic Hodge theory, compare, for instance, with \cite{Liu17}*{\S4}. The work of the preceding sections allows us to interpret $L^i_\dR(X)$ geometrically: we show in \Cref{lat-dR-id} that if $X$ has a proper, flat, semistable $\cO_K$-model $\cX$ for which $H^i_{\log\dR}(\cX/\cO_K)$ and $H^{i + 1}_{\log\dR}(\cX/\cO_K)$ are $\cO_K$-free, then 
\[
L^i_\dR(X) = H^i_{\log\dR}(\cX/\cO_K) \qq \text{inside} \qq H^i_\dR(X/K).
\]
We do not know whether the same holds ``modulo torsion'' if one drops the $\cO_K$-freeness assumption.

\bpp[The base field $K$] \lab{K-intro}
\ready{Throughout \S\ref{lattice}, we assume that $C \cong \wh{\ov{K}}$ for a fixed complete, discretely valued field $K$ that is of mixed characteristic $(0, p)$ and has a perfect residue field $k_0$. 
We set
\[
G \ce \Gal(\ov{K}/K),
\]
so that $G$ acts continuously on $C$, and hence also on $A_\Inf$. The continuous maps $\varphi$ and $\theta$ are $G$-equivariant, and the ideals $(\xi)$, $(\varphi(\xi))$, and $(\mu)$ of $A_\Inf$ are $G$-stable (see \S\ref{Ainf-not}). }

If $\cX$ is a $p$-adic formal $\cO_K$-scheme for which $\fX \ce \cX \wh{\tensor}_{\cO_K} \cO_C$ satisfies the assumptions of \S\ref{fir-setup}, then, by the functoriality of $R\Gamma_{A_\Inf}(\fX)$ (see \S\ref{Ainf-coho}), $G$ acts $A_\Inf$-semilinearly on each $H^i_{A_\Inf}(\fX)$.
\epp

\bpp[The Fargues equivalence] \lab{Fargues-eq}
\ready{By \cite{BMS16}*{4.26}, for any Breuil--Kisin--Fargues module $(M, \varphi_M)$ (see \S\ref{BKF-mod}), its \emph{\'{e}tale realization}, namely, 
\[
M_\et \ce (M \tensor_{A_\Inf} W(C^\flat))^{\varphi_M \tensor\, \varphi\, =\,  1},
\]
is a finitely generated $\bZ_p$-module that comes equipped with an identification
\[
\tst M \tensor_{A_\Inf} W(C^\flat) \cong M_\et \tensor_{\bZ_p} W(C^\flat) \qq \text{under which} \qq M \tensor_{A_\Inf} A_\Inf[\f{1}{\mu}] \cong M_\et \tensor_{\bZ_p} A_\Inf[\f{1}{\mu}]. 
\]
Thus, $M_\et$ is $\bZ_p$-free if $M$ is $A_\Inf$-free and, for any $(M, \varphi_M)$, we have $M \tensor_{A_\Inf} B_\dR \cong M_\et \tensor_{\bZ_p} B_\dR$,
so that $M_\et$ comes equipped with a $B_\dR^+$-sublattice (recall that $M[\f{1}{p}]$ is $A_\Inf[\f{1}{p}]$-free, see \S\ref{BKF-mod})
\[
M \tensor_{A_\Inf} B_\dR^+ \subset M_\et \tensor_{\bZ_p} B_\dR.
\]
}By a theorem of Fargues \cite{BMS16}*{4.28}, 
the category of $A_\Inf$-free Breuil--Kisin--Fargues modules $(M, \varphi_M)$ is equivalent to that of pairs $(T, \Xi)$ consisting of a finite free $\bZ_p$-module $T$ and a $B_\dR^+$-lattice $\Xi \subset T \tensor_{\bZ_p} B_\dR$ via the functor
\[
(M, \varphi_M) \mapsto (M_\et, M \tensor_{A_\Inf} B_\dR^+).
\]
\epp

\bpp[Breuil--Kisin--Fargues $G$-modules] \lab{BKF-G}
Due to the origin of our $C$ (see \S\ref{K-intro}), we may consider \emph{Breuil--Kisin--Fargues $G$-modules}, that is, Breuil--Kisin--Fargues modules $(M, \varphi_M)$ equipped with an
$A_\Inf$-semilinear $G$-action on $M$ for which $\varphi_M$ is $G$-equivariant. A morphism of Breuil--Kisin--Fargues $G$-modules is a $G$-equivariant $A_\Inf$-module morphism that commutes with the isomorphisms~$\varphi_M$. 

For instance, if an $\cX$ as in \S\ref{K-intro} is proper, then each $H^i_{A_\Inf}(\fX)$ is a Breuil--Kisin--Fargues $G$-module (see \Cref{Hinf-BKF}). The \'{e}tale realization $M_\et$ of a Breuil--Kisin--Fargues $G$-module $(M, \varphi_M)$ carries the induced $\bZ_p$-linear~$G$-action.
\epp

\bprop \lab{dR-eq-cat}
The category of $A_\Inf$-free Breuil--Kisin--Fargues $G$-modules $(M, \varphi_M)$ is equivalent to that of pairs $(T, \Xi)$ consisting of a finite free $\bZ_p$-module $T$ equipped with a 
$G$-action and a $G$-stable $B_\dR^+$-lattice $\Xi \subset T \tensor_{\bZ_p} B_\dR$ via the functor
\[
(M, \varphi_M) \mapsto (M_\et, M \tensor_{A_\Inf} B_\dR^+).
\]
\eprop

\bpf
The claim follows from the Fargues equivalence reviewed in \S\ref{Fargues-eq}.
\epf

\bpp[An \'{e}tale lattice determines a de Rham lattice] \lab{et-lat-dR-lat}
Let $T$ be a finite free $\bZ_p$-module endowed with a continuous $G$-action for which the $G$-representation $T[\f{1}{p}]$ is de Rham, so that there is a $G$-equivariant identification
\[
T \tensor_{\bZ_p} B_\dR \cong D_\dR(T) \tensor_K B_\dR, \qq \text{where} \qq D_\dR(T) \ce (T \tensor_{\bZ_p} B_\dR)^G.
\]
For such $T$, the $B_\dR^+$-lattice $D_\dR(T) \tensor_K B_\dR^+$ is evidently $G$-stable in $T \tensor_{\bZ_p} B_\dR$. Thus, by \Cref{dR-eq-cat}, the pair $(T, D_\dR(T) \tensor_K B_\dR^+)$, so $T$, determines an $A_\Inf$-free Breuil--Kisin--Fargues $G$-module 
\[
(M(T), \varphi_{M(T)})
\]
that depends functorially on $T$. 
The de Rham realization 
\[
M(T)_\dR \ce M(T) \tensor_{A_\Inf,\, \theta} \cO_C \qq \text{of} \qq (M(T), \varphi_{M(T)})
\]
is an $\cO_C$-lattice in 
\[
(M(T) \tensor_{A_\Inf} B_\dR^+)/\xi \cong (D_\dR(T) \tensor_K B_\dR^+)/\xi \cong D_\dR(T) \tensor_K C.
\]
Therefore, functorially in $T$, we obtain the $\cO_K$-lattice
\[
\tst (M(T)_\dR)^G \qq \text{inside the $K$-vector space} \qq D_\dR(T).
\]
\epp

\beg \lab{lat-main-eg}
\ready{We fix a $K$-scheme $X$ (or even a $K$-rigid space, which we view as an adic space, see \cite{Hub96}*{1.1.11 (d)}) that is proper and smooth, and~set 
\[
L^i_\et(X) \ce H^i_\et(X_{\ov{K}}, \bZ_p)/H^i_\et(X_{\ov{K}}, \bZ_p)_\tors \cong H^i_\et(X_C, \bZ_p)/H^i_\et(X_C, \bZ_p)_\tors \qq \text{for} \qq i \ge 0.
\]
As is well known and follows from \eqref{dR-comp-iso}, the $G$-representation $(L^i_\et(X))[\f{1}{p}]$ is de Rham and 
\be \lab{D-dR-coh}
D_\dR(L^i_\et(X)) \cong  (L^i_\et(X) \tensor_{\bZ_p} B_\dR)^G \overset{\eqref{dR-comp-iso}}{\cong} (H^i_\dR(X/K) \tensor_K B_\dR)^G \cong H^i_\dR(X/K)
\ee
functorially in $X$. Thus, using the discussion of \S\ref{et-lat-dR-lat}, we obtain the $\cO_K$-lattice  
\[
L^i_\dR(X) \ce (M(L^i_\et(X))_\dR)^G \subset H^i_\dR(X/K)
\]
that is functorial in $X$ (even in $L^i_\et(X)$). Its definition implies that for a finite Galois extension~$K'/K$,
\[
L^i_\dR(X) =  (L^i_\dR(X_{K'}))^{\Gal(K'/K)} \qq \text{inside} \qq H^i_\dR(X/K) = (H^i_\dR(X_{K'}/K'))^{\Gal(K'/K)}.
\]
}If $X$ extends to a proper, flat, semistable $\cO_K$-scheme $\cX$ such that $H^i_{\log\dR}(\cX/\cO_K)$ and $H^{i + 1}_{\log\dR}(\cX/\cO_K)$ are $\cO_K$-free (where $\cX$ is endowed with the log structure $\cO_{\cX,\, \et} \cap (\cO_{\cX,\,\et}[\f{1}{p}])^\times$), then, by the following \Cref{lat-dR-id} (and GAGA techniques, similarly to \Cref{comp-sch}), 
\[
L^i_\dR(X) = H^i_{\log\dR}(\cX/\cO_K) \qq \text{inside} \qq H^i_\dR(X/K);
\]
in particular, if $\cX'$ is another such $\cO_K$-model of $X$, then 
\be \lab{log-dR-X-Xpr}
H^i_{\log\dR}(\cX/\cO_K) = H^i_{\log\dR}(\cX'/\cO_K) \qq \text{inside}\qq H^i_\dR(X/K).
\ee
\eeg

\bthm \lab{lat-dR-id}
Let $\cX$ be a proper, flat $p$-adic formal $\cO_K$-scheme endowed with the log structure $\cO_{\cX,\,\et} \cap (\cO_{\cX,\,\et}[\f{1}{p}])^\times$ such that $\cX$ has an \'{e}tale cover by affines $\cU$ each of which has an \'{e}tale morphism 
\be  \lab{lat-dR-id-sst} 
\cU \ra \Spf (\cO_K\{t_0, \dotsc, t_r, t_{r + 1}, \dotsc,  t_d\}/(t_0\cdots t_r - \pi)) \q \text{for some nonunit} \q \pi \in \cO_K \setminus \{ 0 \}
\ee
\up{where $d$, $r$, and $\pi$ depend on $\cU$}. If $H^i_{\log\dR}(\cX/\cO_K)$ and $H^{i + 1}_{\log\dR}(\cX/\cO_K)$ are $\cO_K$-free,~then 
\be \lab{lat-dR-id-eq}
L^i_\dR(\cX_K^\ad) = H^i_{\log\dR}(\cX/\cO_K) \qq \text{inside} \qq H^i_\dR(\cX_K^\ad/K);
\ee
in fact, then, setting $\fX \ce \cX \wh{\tensor}_{\cO_K} \cO_C$, we have the identification 
\be \lab{et-det-BKF}
M(L^i_\et(\cX_K^\ad)) \cong H^i_{A_\Inf}(\fX)
\ee
of Breuil--Kisin--Fargues $G$-modules.
\ethm

\bpf
\ready{By working locally on $\cU$, in the target of \eqref{lat-dR-id-sst} we may replace each $t_i$ by $t_i^{\pm 1}$ for $r + 1 \le i \le d$, so $\fX$ satisfies the assumptions of \S\ref{fir-setup}. Moreover, by the Grothendieck comparison theorem 
and flat base change (compare with \Cref{comp-sch}), for $j = i$ and $j = i + 1$, we have
\be \lab{lat-dR-id-1}
H^j_{\log\dR}(\fX/\cO_C) \cong H^j_{\log\dR}(\cX/\cO_K) \tensor_{\cO_K} \cO_C, \q \text{so} \q H^j_{\log\dR}(\cX/\cO_K) \cong (H^j_{\log\dR}(\fX/\cO_C))^G.
\ee
Thus, by \Cref{free-free}, the Breuil--Kisin--Fargues $G$-modules $H^i_{A_\Inf}(\fX)$ and $H^{i + 1}_{A_\Inf}(\fX)$ (see \S\ref{BKF-G}) are $A_\Inf$-free. By \eqref{WCflat-iso}, we have the $G$-equivariant identification of the \'{e}tale realization: 
\[
(H^i_{A_\Inf}(\fX))_\et \cong H^i_\et(\cX_C^\ad, \bZ_p),
\]
which, consequently, is $\bZ_p$-free. By \Cref{BdR-et-comp}, the $B_\dR$-base change of this identification agrees with the identification $H^i_{A_\Inf}(\fX) \tensor_{A_\Inf} B_\dR \cong H^i_\et(\cX_C^\ad, \bZ_p) \tensor_{\bZ_p} B_\dR$ that results by combining
\[
H^i_{A_\Inf}(\fX) \tensor_{A_\Inf} B_\dR^+ \overset{\eqref{BdR-comp-cor-eq}}{\cong} H^i_\cris(\cX_C^\ad/B_\dR^+) \overset{\eqref{cris-dR-descent}}{\cong} H^i_\dR(\cX_K^\ad/K) \tensor_K B_\dR^+ 
\]
and
\[
H^i_\dR(\cX_K^\ad/K) \tensor_K B_\dR \overset{\eqref{dR-comp-iso}}{\cong} H^i_\et(\cX_C^\ad, \bZ_p) \tensor_{\bZ_p} B_\dR.
\]
This compatibility and \S\ref{et-lat-dR-lat} (see also \eqref{D-dR-coh}) supply the desired $G$-equivariant identification \eqref{et-det-BKF}:
\[
M(L^i_\et(\cX_K^\ad)) \cong H^i_{A_\Inf}(\fX).
\]
Under this identification, by \Cref{BdR-comp} and the sentence after \eqref{cris-dR-descent}, the identifications
\[
M(L^i_\et(\cX_K^\ad)) \tensor_{A_\Inf,\, \theta} C \overset{\eqref{D-dR-coh}}{\cong} H^i_\dR(\cX_C^\ad/C) \qq \text{and} \qq H^i_{A_\Inf}(\fX) \tensor_{A_\Inf,\, \theta} C \overset{\eqref{G-dR-spec-eq}}{\cong} H^i_\dR(\cX_C^\ad/C)
\]
agree.} Thus, \eqref{dR-bc-SES} implies the following equality inside $H^i_\dR(\cX_C^\ad/C)$:
\[
M(L^i_\et(\cX_K^\ad))_\dR = M(L^i_\et(\cX_K^\ad)) \tensor_{A_\Inf,\, \theta} \cO_C = H^i_{A_\Inf}(\fX) \tensor_{A_\Inf,\, \theta} \cO_C = H^i_{\log\dR}(\fX/\cO_C),
\]
which, together with the second identification in \eqref{lat-dR-id-1}, gives the desired \eqref{lat-dR-id-eq}.
\epf

\brem
In the proof above, we have seen that both $H^i_{A_\Inf}(\fX)$ and $H^{i + 1}_{A_\Inf}(\fX)$ are $A_\Inf$-free, so, by \eqref{cris-bc-SES}, we have the $G$-equivariant and Frobenius-equivariant identifications
\[
H^i_{A_\Inf}(\fX) \tensor_{A_\Inf} W(k) \cong H^i_{\log\cris}(\fX_k/W(k)) \overset{\eqref{LCB-eq}}{\cong} H^i_{\log\cris}(\cX_k/W(k)),
\]
and hence also the Frobenius-equivariant identification
\be \lab{BKF-det-cris}
(H^i_{A_\Inf}(\fX) \tensor_{A_\Inf} W(k))^G \cong H^i_{\log\cris}(\cX_{k_0}/W(k_0)).
\ee
In particular, \eqref{et-det-BKF} and \eqref{BKF-det-cris} show that, under the assumptions of \Cref{lat-dR-id}, the integral $p$-adic \'{e}tale cohomology $H^i_\et(\cX_C^\ad, \bZ_p)$ endowed with its Galois action functorially determines the integral logarithmic crystalline cohomology $H^i_{\log\cris}(\cX_{k_0}/W(k_0))$ endowed with its Frobenius.
\erem
}
}

%% file: semistable.tex

\section{The semistable comparison isomorphism} \lab{semistable}

\revise{
\ready{Our final goal is to deduce the semistable comparison isomorphism for suitable proper, ``semistable'' formal schemes (see \Cref{sst-comp}). This extends \cite{BMS16}*{1.1~(i)}, which treated the good reduction case (see also \cite{TT15}*{1.2} for a result ``with coefficients'' over an absolutely unramified base), and is similar to the semistable comparison established by Colmez--Niziol \cite{CN17}*{5.26}. More precisely, \emph{loc.~cit.}~also includes cases in which the log structures are not ``vertical.''

\bpp[The ring $B_\st$]
\ready{We consider the log PD thickenings $A_\cris/p^n$ of $\cO_C/p$ (see \S\ref{Acris-log}) and set
\[
\tst J_n \ce \Ker(A_\cris/p^n \surjects \cO_C/p) \qq \text{and} \qq J \ce \varprojlim_{n\ge 1} J_n \cong \Ker(A_\cris \surjects \cO_C/p).
\]
The element $p \in \cO_C \setminus \{ 0 \}$ belongs to the log structure of $\cO_C/p$ (see \S\ref{log-str-def}~\ref{OC-log-str}), so its preimage in the log structure of $A_\cris/p^n$ is a $(1 + J_n, \times)$-torsor, which is necessarily trivial\footnote{\ready{Quasi-coherent cohomology of affine schemes vanishes, so, for a finitely generated, and hence nilpotent, ideal $J' \subset J_n$ of $A_\cris/p^n$, the \'{e}tale sheaf on $\Spec(A_\cris/p^n)$ associated to $(1 + J', \times)$ has no nontrivial torsors. The filtered direct limit of these sheaves is the analogous sheaf associated to $(1 + J_n, \times)$, so it, too, has no nontrivial torsors.}} (compare with \cite{Bei13a}*{\S1.1, Exercises, (iii)}). Consequently, as $n$ varies, these torsors comprise a trivial $(1 + J, \times)$-torsor $\tau_0$, whose base change along the logarithm map $(1 + J, \times) \ra (J, +) \subset (A_\cris, +)$ furnished by the divided power structure on $J$ is a trivial $(A_\cris, +)$-torsor $\tau$, the \emph{Fontaine--Hyodo--Kato torsor} (compare with \cite{Bei13a}*{\S1.15, p.~23}). The functor which to an $A_\cris$-algebra $A$ assigns the underlying set of the $(A, +)$-torsor $\tau \times_{(A_\cris, +)} (A, +)$ is represented by the $A_\cris$-algebra $A_\st$, so $A_\st$ is the initial $A_\cris$-algebra over which the Fontaine--Hyodo--Kato torsor is canonically trivialized. 

 We may noncanonically trivialize $\tau_0$ (for instance, $[p^{1/p^\infty}]$ is a trivialization, see \eqref{AL-1}) to obtain an isomorphism $A_\st \simeq A_\cris[T]$, which, upon adjusting the trivialization by an $a \in 1 + J$, gets postcomposed with the $A_\cris$-automorphism $A_\cris[T] \xra{T \mapsto T + \log(a)} A_\cris[T]$. 
 The $A_\cris$-derivation $-\f{d}{dT}$ respects these automorphisms, so it induces a canonical $A_\cris$-derivation, the \emph{monodromy operator},
\[
N \colon A_\st \ra A_\st \qq \text{for which} \qq (A_\st)^{N = 0} = A_\cris
\]
(our $N$ agrees with that of \emph{op.~cit.}, see \cite{Bei13a}*{\S1.15, Remarks (i)}; compare also with \cite{Tsu99}*{4.1.1}). 

By \cite{Bei13a}*{(1.15.2)}, the Frobenius pullback of $\tau_0$ is isomorphic to the $p$-fold self-product of $\tau_0$, and hence likewise for the base change $\tau$ of $\tau_0$ to $(A_\cris, +)$. Thus, we have an 
 $A_\cris$-semilinear Frobenius
\be \lab{Ast-Frob}
\varphi \colon A_\st \ra A_\st
\ee
that in terms of an isomorphism $A_\st \simeq A_\cris[T]$ obtained by trivializing $\tau_0$ is described by $T \mapsto pT$. 
The interaction of $\varphi$ and $N$ is described by the formula $N \varphi = p \varphi N$. 

}Since $\mu$ and $\log([\eps])$ are unit multiples of each other in $A_\cris$ (see \S\ref{log-eps}) and $\varphi(\log([\eps])) = p\log([\eps])$, the Frobenius \eqref{Ast-Frob} and, evidently, also the derivation $N$ induce their counterparts on
\[
\tst B_\st^+ \ce A_\st[\f{1}{p}] \qq \text{and} \qq B_\st \ce A_\st[\f{1}{p\mu}].
\]
The relation $N\varphi = p\varphi N$ continues to hold for $B_\st^+$ and $B_\st$. As is explained in \cite{Bei13a}*{\S1.17}, the $A_\cris$-algebras $B_\st^+$ and $B_\st$ reviewed above agree with the ones constructed in \cite{Fon94}*{\S3}.
\epp

 For us, the significance of the period ring $B_\st^+$ lies in the following comparison between the logarithmic crystalline cohomology of $\fX_k$ over $W(k)$ and of $\fX_{\cO_C/p}$ over $A_\cris$ (compare with \cite{BMS16}*{13.21}).

\bprop \lab{Y-and-X}
If $\fX$ is $\cO_C$-proper, then 
\be \lab{cris-cris-rel}
R\Gamma_{\log\cris}(\fX_{k} /W(k)) \otimes^\bL_{W(k)} B_\st^+ \cong
R\Gamma_{\log\cris}(\fX_{\cO_C/p}/A_\cris) \otimes^\bL_{A_\cris} B_\st^+,
\ee
where the log structures are those of \uS\uref{log-str-def}~\ref{X-log-str}, \uS\uref{Acris-log}, and \uS\uref{Wk-log}. In particular, if $\fX$ is $\cO_C$-proper and $\cY$ is a descent of $\fX_{\cO_C/p}$ to a proper, log smooth, fine log $\cO/p$-scheme of Cartier type for some discrete valuation subring $\cO \subset \cO_C$ with a perfect residue field $k_0$ and $C \cong \p{\ov{\cO[\f{1}{p}]}}\, \wh{}\,$ \up{where $\cO/p$ is equipped with the log structure associated to the chart $\cO \setminus \{ 0 \} \ra \cO/p$}, then we have the following identification that is compatible with the actions of $\varphi$ and $N$ \up{described in the proof}\ucolon
\be \lab{Bei-id-2}
R\Gamma_{\log\cris}(\cY_{k_0} /W(k_0)) \otimes^\bL_{W(k_0)} B_\st^+ \cong
R\Gamma_{\log\cris}(\fX_{\cO_C/p}/A_\cris) \otimes^\bL_{A_\cris} B_\st^+,
\ee
where $W(k_0)$ is endowed with the log structure associated to $\bN_{\ge 0} \xra{1\,\mapsto\,0,\, 0\, \mapsto\, 1} W(k_0)$.
\eprop

\bpf
A descent $\cY$ exists (see the proof of Corollary \uref{Acris-BC}), so \eqref{Bei-id-2} follows from  \cite{Bei13a}*{(1.16.2) and (1.18.5)} and, due to \eqref{LCB-eq}, it implies \eqref{cris-cris-rel}. On the left side of \eqref{Bei-id-2}, the operator $N$ combines the monodromy of $R\Gamma_{\log\cris}(\cX_{k_0} /W(k_0))$ and $B_\st^+$, so is ``$N \tensor 1 + 1 \tensor N$''; on the right, $N$ is the monodromy of $B_\st^+$. On both sides of \eqref{Bei-id-2}, the Frobenius $\varphi$ acts on both factors.
\epf

\brem
One may eliminate the dependence of \eqref{Bei-id-2} on the choice of $\cY$ by forming a direct limit over all the possible $\cY$, see \cite{Bei13a}*{\S1.18, Remarks (i)}.
\erem

\bpp[The base field $K$] \lab{base-K}
For the rest of \S\ref{semistable}, we assume that $C = \wh{\ov{K}}$ for a fixed complete, discretely valued subfield $K \subset C$ with a perfect residue field $k_0$, set $G \ce \Gal(\ov{K}/K)$, and endow $\cO_K$ (resp.,~$\cO_K/p$) with the log structure associated to the chart $\cO_K \setminus \{0\} \hra \cO_K$ (resp.,~its base change). 
By functoriality, $G$ acts on $A_\cris$, $A_\st$, $B_\st^+$, and, since the ideal $(\mu)$ does not depend on the choice of $\eps$ (see \S\ref{Ainf-not}), also on $B_\st$. These $G$-actions commute with $\varphi$ and $N$. When $\cO$ of \Cref{Y-and-X} is our $\cO_K$, the identification \eqref{Bei-id-2} is $G$-equivariant granted that $G$ acts on both sides by functoriality. 
\epp

\bthm \lab{sst-comp}
Let $\cX$ be a proper $p$-adic formal $\cO_K$-scheme that has an \'{e}tale cover by affines $\cU$ each of which has an \'{e}tale $\cO_K$-morphism 
\[
\cU \ra \Spf (\cO_K\{t_0, \dotsc, t_r, t_{r + 1}, \dotsc,  t_d\}/(t_0\cdots t_r - \pi)) \q \text{for some nonunit} \q \pi \in \cO_K \setminus \{ 0 \}
\]
\up{where $d$, $r$, and $\pi$ depend on $\cU$} and endow $\cX$ with the log structure $\cO_{\cX,\,\et} \cap (\cO_{\cX,\,\et}[\f{1}{p}])^\times$. There is the following natural, $G$-equivariant isomorphism that is compatible with the actions of $\varphi$ and $N$\ucolon
\be \lab{sst-comp-eq}
R\Gamma_\et (\cX_{C}^\ad, \bZ_p)\otimes^\bL_{\bZ_p} B_\st \cong
R\Gamma_{\log\cris}(\cX_{k_0} /W(k_0)) \otimes^\bL_{W(k_0)} B_\st,
\ee
where $W(k_0)$ is endowed with the log structure associated to $\bN_{\ge 0} \xra{1\,\mapsto\,0,\, 0\, \mapsto\, 1} W(k_0)$. In particular, 
\[
\text{the $G$-representation} \q H^i_\et(\cX_C^\ad, \bQ_p) \q \text{is semistable for every} \q i \in \bZ.
\]
\ethm

\bpf
We set $\fX \ce \cX \wh{\tensor}_{\cO_K} \cO_C$, so that $\fX$ meets the requirements of \uS\uref{fir-setup}. By \Cref{log-claim,log-comp} and \cite{Kat89}*{4.8}, the base change $\cX_{\cO_K/p}$ is fine, log smooth, and of Cartier type over $\cO_K/p$, so \Cref{Y-and-X} applies to it and gives the $G$-equivariant (see \S\ref{base-K}) identification
\[
R\Gamma_{\log\cris}(\cX_{k_0} /W(k_0)) \otimes^\bL_{W(k_0)} B_\st^+ \cong
R\Gamma_{\log\cris}(\fX_{\cO_C/p}/A_\cris) \otimes^\bL_{A_\cris} B_\st^+ \overset{\eqref{ABCA-eq}}{\cong}
R\Gamma(\fX_\et, A\Omega_\fX) \tensor^\bL_{A_\Inf} B_\st^+
\]
that is compatible with $\varphi$ and $N$. In addition, by \eqref{RG-et-id-eq}, we have the $G$-equivariant identification 
\[
R\Gamma (\fX_\et, A\Omega_\fX) \otimes^\bL_{A_\Inf} B_\st \cong
R\Gamma (\fX_C^\ad, \bZ_p)\otimes^\bL_{\bZ_p} B_\st \cong R\Gamma (\cX_{C}^\ad, \bZ_p)\otimes^\bL_{\bZ_p} B_\st
\]
that is trivially compatible with $N$ and is compatible with $\varphi$ by the discussion after \Cref{RG-et-id}. The desired \eqref{sst-comp-eq} follows by combining the displayed identifications.
\epf

\brem
The isomorphism \eqref{sst-comp-eq} is compatible with filtrations in the following sense: by \cite{Fon94}*{\S4.2}, there is a (noncanonical) $A_\cris$-algebra homomorphism $B_\st \ra B_\dR$ and, by the proof above and \Cref{BdR-et-comp}, the $B_\dR$-base change of the isomorphism \eqref{sst-comp-eq} is identified with the de Rham comparison isomorphism \eqref{dR-comp-iso} (with $X_0 = \cX_K^\ad$) that is compatible with filtrations.
\erem

}}